\documentclass[12pt]{article}
\usepackage{amsthm,epsfig,a4}
\usepackage[ansinew]{inputenc}
\usepackage{amsmath,amssymb}
\usepackage{amsfonts}
\usepackage{eufrak}
\usepackage[UglyObsolete,tight,heads=LaTeX]{diagrams}
\def\q{\nolinebreak \hfill $\Box$} 
\def\N{\\[5mm] \indent} 
\def\D{\displaystyle}
\def\:{\colon\thinspace}
\theoremstyle{plain}
\newtheorem{Lemma}{Lemma}[section]
\newtheorem{Prop}[Lemma]{Proposition}
\newtheorem{Thm}[Lemma]{Theorem}
\newtheorem{Cor}[Lemma]{Corollary}
\newtheorem{SLemma}{Lemma}[subsection]
\newtheorem{SProp}[SLemma]{Proposition}
\newtheorem{SThm}[SLemma]{Theorem}
\newtheorem{SCor}[SLemma]{Corollary}
\theoremstyle{definition}
\newtheorem{Def}[Lemma]{Definition}
\newtheorem{SDef}[SLemma]{Definition}
\newtheorem*{Ex}{Example}
\newtheorem*{Rem}{Remark}
\newtheorem*{Ac}{Acknowledgement}
\allowdisplaybreaks
\begin{document}
\title{On Low-Dimensional Solvmanifolds}
\author{\textsc{Christoph Bock}
}
\date{}
\maketitle
%
%
%
\section{Introduction}
In this note we want to study compact homogeneous spaces $G/\Gamma$,
where $G$ is a connected and simply-connected Lie group and $\Gamma$
a discrete subgroup in $G$. It is well known that the existence of
such a $\Gamma$ implies the unimodularity of the Lie group $G$.
Recall that a Lie group $G$ is called
\emph{unimodular}\index{Unimodular Lie Group} if for all $X \in
\mathfrak{g}$ holds $\mathrm{tr} \, \mathrm{ad} \, X = 0$, where
$\mathfrak{g}$ denotes the Lie algebra of $G$.

If we further demand $G/\Gamma$ to be symplectic (when $G$ is
even-dimensional), a result of Chu \cite{Chu} shows that $G$ has to
be solvable.

Therefore, we regard compact quotients of connected and
simply-connected solvable Lie groups by discrete subgroups, so
called \emph{solvmanifolds}.

First, we recall the definition of nilpotent and solvable groups
resp.\ Lie algebras.
\begin{itemize}
\item[(i)] Let $G$ be group and denote its neutral element by $e$.
We define the \emph{derived series} $(D^{(k)}G)_{k \in \mathbb{N}}$,
\index{Derived Series} \emph{descending series} $(G^{(k)})_{k \in
\mathbb{N}}$\index{Descending Series} and \emph{ascending series}
$(G_{(k)})_{k \in \mathbb{N}}$\index{Ascending Series} of subgroups
in $G$ inductively as follows:
\begin{eqnarray*}
& D^{(0)}G := G^{(0)} := G, & \\ & D^{(k)}G :=
[D^{(k-1)}G,D^{(k-1)}G],~~ G^{(k)} := [G,G^{(k-1)}], & \\ & G_{(0)}
= \{e\},~~ G_{(k)} := \{g \in G \,|\, [g,G] \subset G_{(k-1)}\}. &
\end{eqnarray*}

$G$ is called \emph{nilpotent}\index{Nilpotent!Group} if there
exists $k_0 \in \mathbb{N}$ such that $G^{(k_0)} = \{e\}$.

$G$ is called \emph{solvable}\index{Solvable!Group} if there exists
$k_0 \in \mathbb{N}$ such that $D^{(k_0)}G = \{e\}$.
\item[(ii)] Given a Lie algebra $\mathfrak{g}$, one defines the
\emph{derived}, \emph{descending} and \emph{ascending series} of
subalgebras in $\mathfrak{g}$ via
\begin{eqnarray*}
& D^{(0)}\mathfrak{g} := \mathfrak{g}^{(0)} := \mathfrak{g}, & \\ &
D^{(k)}\mathfrak{g} :=
[D^{(k-1)}\mathfrak{g},D^{(k-1)}\mathfrak{g}],~~ \mathfrak{g}^{(k)}
:= [\mathfrak{g},\mathfrak{g}^{(k-1)}], &
\\ & \mathfrak{g}_{(0)} = \{0\},~~ \mathfrak{g}_{(k)} := \{X \in \mathfrak{g}
\,|\, [X,\mathfrak{g}] \subset \mathfrak{g}_{(k-1)}\} &
\end{eqnarray*}
and calls $\mathfrak{g}$ \emph{nilpotent}\index{Nilpotent!Lie
Algebra} resp.\ \emph{solvable}\index{Solvable!Lie Algebra} if its
ascending resp.\ derived series becomes trivial for $k_0$ large
enough.
\end{itemize}

We collect some properties in the following proposition. The parts
which are not obvious can be found in \cite[Section 3.18]{V}.

\setcounter{SLemma}{0}
\begin{SProp} $\,$
\begin{itemize}
\item[(i)] The subgroups arising in the derived, descending and
ascending series of a group are normal. Moreover, they are closed
and simply-connected Lie subgroups in the case of a connected and
simply-connected Lie group.
\item[(ii)] The subalgebras arising in the derived, descending and
ascending series of a Lie algebra are ideals.
\item[(iii)] A Lie group is nilpotent resp.\ solvable if and only if its Lie
algebra is nilpotent resp.\ solvable. \q
\end{itemize}
\end{SProp}

We shall examine whether certain solvmanifolds are formal,
symplectic, K\"ahler or Lefschetz. Therefore, we give a brief
introduction of these notions.
\subsection{Formality}
A \emph{differential graded algebra (DGA)}\index{Differential Graded
Algebra (DGA)} is a graded $\mathbb{R}$-algebra $A = \bigoplus_{i
\in \mathbb{N}}~ A^i$ together with an $\mathbb{R}$-linear map $d \:
A \to A$ such that $d(A^i) \subset A^{i+1}$ and the following
conditions are satisfied:
\begin{itemize}
\item[(i)] The $\mathbb{R}$-algebra structure of $A$ is given by an
inclusion $\mathbb{R} \hookrightarrow A^0$.

\item[(ii)] The multiplication is graded commutative, i.e.\ for $a \in
A^i$ and $b \in A^j$ one has $a \cdot b = (-1)^{i \cdot j} b \cdot a
\in A^{i+j}$.

\item[(iii)] The Leibniz rule holds: $\forall_{a \in A^i} \forall_{b
\in A} ~ d(a \cdot b) = d(a) \cdot b + (-1)^i a \cdot d(b)$

\item[(iv)] The map $d$ is a differential, i.e.\ $d^2 = 0$.
\end{itemize}
Further, we define $|a| := i$ for $a \in A^i$.

The \emph{$i$-th cohomology of a DGA} $(A,d)$ is the algebra
$$ H^i(A,d) := \frac{\ker (d \: A^i \to A^{i+1})}{\mathrm{im}
(d \: A^{i-1} \to A^i)}.$$

If $(B,d_B)$ is another DGA, then a $\mathbb{R}$-linear map $f \: A
\to B$ is called \emph{morphism} if $f(A^i) \subset B^i$, $f$ is
multiplicative, and $d_B \circ f = f \circ d_A$. Obviously, any such
$f$ induces a homomorphism $f^* \: H^*(A,d_A) \to H^*(B,d_B)$. A
morphism of differential graded algebras inducing an isomorphism on
cohomology is called
\emph{quasi-isomorphism}\index{Quasi-Isomorphism}.
\begin{SDef}
A DGA $(\mathcal{M},d)$ is said to be \emph{minimal} if
\begin{itemize}
\item[(i)] there is a graded vector space $V = \big( \bigoplus_{i \in
\mathbb{N}_+} V^i \big) = \mathrm{Span} \, \{a_k ~ | ~ k \in I\}$
with homogeneous elements $a_k$, which we call the generators,

\item[(ii)] $\mathcal{M} = \bigwedge V$,

\item[(iii)]  the index set $I$ is well ordered, such that $k < l ~
\Rightarrow |a_k| \le |a_l|$ and the expression for $da_k$ contains
only generators $a_l$ with $l < k$.
\end{itemize}
\end{SDef}
We shall say that $(\mathcal{M},d)$ is a \emph{minimal model for a
differential graded algebra}\index{Minimal Model of a! Differential
Graded Algebra} $(A,d_A)$ if $(\mathcal{M},d)$ is minimal and there
is a quasi-isomorphism of DGAs $\rho \: (\mathcal{M},d) \to
(A,d_A)$, i.e.\ it induces an isomorphism $\rho^* \:
H^*(\mathcal{M},d) \to H^*(A,d_A)$  on cohomology.

The importance of minimal models is reflected by the following
theorem, which is taken from Sullivan's work \cite[Section
5]{Sullivan}.
\begin{SThm}\label{Konstruktion d m M}
 A differential graded algebra $(A,d_A)$ with $H^0(A,d_A) =
\mathbb{R}$ possesses a minimal model. It is unique up to
isomorphism of differential graded algebras.
\end{SThm}

We quote the existence-part of Sullivan's proof, which gives an
explicit construction of the minimal model. Whenever we are going to
construct such a model for a given algebra in this article, we will
do it as we do it in this proof.

\textit{Proof of the existence.} We need the following algebraic
operations to ``add'' resp.\ ``kill'' cohomology.

Let $(\mathcal{M},d)$ be a DGA. We ``add'' cohomology by choosing a
new generator $x$ and setting
$$\widetilde{\mathcal{M}} := \mathcal{M} \otimes \bigwedge(x),~~
\tilde{d}|_{\mathcal{M}} = d, ~~ \tilde{d}(x) = 0,$$ and ``kill'' a
cohomology class $[z] \in H^k(\mathcal{M},d)$ by choosing a new
generator $y$ of degree $k-1$ and setting
$$\widetilde{\mathcal{M}} := \mathcal{M} \otimes \bigwedge(y),~~
\tilde{d}|_{\mathcal{M}} = d, ~~ \tilde{d}(y) = z.$$ Note that $z$
is a polynomial in the generators of $\mathcal{M}$.

Now, let $(A,d_A)$ a DGA with $H^0(A,d_A) = \mathbb{R}$. We set
$\mathcal{M}_0 := \mathbb{R}$,  $d_0 := 0$ and $\rho_0(x) = x$.

Suppose now $\rho_k \: (\mathcal{M}_k,d_k) \to (A,d_A)$ has been
constructed so that $\rho_k$ induces isomorphisms on cohomology in
degrees $\le k$ and a monomorphism in degree $(k+ \nolinebreak 1)$.

``Add'' cohomology in degree $(k+1)$ to get morphism of differential
graded algebras $\rho_{(k+1),0} \:
(\mathcal{M}_{(k+1),0},d_{(k+1),0}) \to (A,d_A)$ which induces an
isomorphism $\rho_{(k+1),0}^*$ on cohomology in degrees $\le (k+1)$.
Now, we want to make the induced map $\rho_{(k+1),0}^*$ injective on
cohomology in degree $(k+2)$ .

We ``kill'' the kernel on cohomology in degree $(k+2)$ (by
non-closed generators of degree (k+1)) and define $\rho_{(k+1),1} \:
(\mathcal{M}_{(k+1),1},d_{(k+1),1}) \to (A,d_A)$ accordingly. If
there are generators of degree one in
$(\mathcal{M}_{(k+1),0},d_{(k+1),0})$ it is possible that this
killing process generates new kernel on cohomology in degree
$(k+2)$. Therefore, we may have to ``kill'' the kernel in degree
$(k+2)$ repeatedly.

We end up with a morphism $\rho_{(k+1),\infty} \:
(\mathcal{M}_{(k+1),\infty},d_{(k+1),\infty}) \to (A,d_A)$ which
induces isomorphisms on cohomology in degrees $\le (k+1)$ and a
monomorphism in degree $(k+2)$. Set $\rho_{k+1} :=
\rho_{(k+1),\infty}$ and $(\mathcal{M}_{k+1},d_{k+1}) :=
(\mathcal{M}_{(k+1),\infty},d_{(k+1),\infty})$.

Inductively we get the minimal model $\rho \: (\mathcal{M},d) \to
(A,d_A)$. \q
\N A \emph{minimal model} $(\mathcal{M}_M,d)$ \emph{of a connected
smooth manifold}\index{Minimal Model of a!Manifold} $M$ is a minimal
model for the de Rahm complex $(\Omega(M),d)$ of differential forms
on $M$. The last theorem implies that every connected smooth
manifold possesses a minimal model which is unique up to isomorphism
of differential graded algebras.
\N Now, we are able to introduce the notion of formality. Endowed
with the trivial differential, the cohomology of a minimal DGA is a
DGA itself, and therefore it also possesses a minimal model. In
general, this two minimal models need not to be isomorphic.

A minimal differential graded algebra $(\mathcal{M},d)$ is called
\emph{formal}\index{Formality of a!Differential Graded Algebra} if
there is a morphism of differential graded algebras
$$\psi \: (\mathcal{M},d) \to (H^*(\mathcal{M},d),d_H=
 0)$$ that induces the identity on cohomology.

This means that $(\mathcal{M},d)$ and $(H^*(\mathcal{M},d),d_H=
 0)$ share their minimal model. The following theorem gives an
equivalent characterisation.

\begin{SThm}[{\cite[Theorem 1.3.1]{TO}}]  \label{form}
A minimal model $(\mathcal{M},d)$ is formal if and only if we can
write $\mathcal{M} = \bigwedge V$ and the space $V$ decomposes as a
direct sum $V = C \oplus N$ with $d(C) = 0$, $d$ is injective on
$N$, and such that every closed element in the ideal $I(N)$
generated by $N$ in $\bigwedge V$ is exact. \q
\end{SThm}

This allows us to give a weaker version of the notion of formality.

\begin{SDef}  \label{s-form}
A minimal model $(\mathcal{M},d)$ is called
$\emph{s\mbox{-formal}}$, $s \in \mathbb{N}$, if we can write
$\mathcal{M} = \bigwedge V$ and for each $i \le s$ the space $V^i$
generated by generators of degree $i$ decomposes as a direct sum
$V^i = C^i \oplus N^i$ with $d(C^i) = 0$, $d$ is injective on $N^i$
and such that every closed element in the ideal $I(\bigoplus_{i \le
s} N^i)$ generated by $\bigoplus_{i \le s} N^i$ in $\bigwedge \big(
\bigoplus_{i \le s} V^i \big)$ is exact in $\bigwedge V$.
\end{SDef}

Obviously, formality implies $s$-formality for every $s$.

The following theorem is an immediate consequence of the last
definition.

\begin{SThm} \label{nichtformal}
Let $(\mathcal{M},d)$ be a minimal model, where 
$\mathcal{M} = \bigwedge V$, $V = C \oplus N$ with $d(C) = 0$ and 
$d$ is injective on $N$.

Assume that there exist $r,s \in \mathbb{N}_+$, $n \in N^r$ and 
$x \in \bigwedge \big( \bigoplus_{i \le s} V^i \big)$ such that holds
$$\forall_{c \in C^r}~ (n+c)\,x \mbox{~is closed and not exact.}$$
Then $(\mathcal{M},d)$ is not $\max \{r,s\}$-formal. \q
\end{SThm}

A connected smooth manifold is called \emph{formal}\index{Formality
of a!Manifold} (resp.\ \emph{$s$-formal}) if its minimal model is
formal (resp.\ $s$-formal).
\N The next theorem shows the reason for defining $s$-formality: in
certain cases $s$-formality is sufficient for a manifold to be
formal.

\begin{SThm}[{\cite[Theorem 3.1]{FMDon}}]  \label{formal = n-1 formal}
Let $M$ be a connected and orientable compact smooth manifold of
dimension $2n$ or $(2n-1)$.

Then $M$ is formal if and only if it is $(n-1)$-formal. \q
\end{SThm}

\begin{Ex}[{\cite[Corollary 3.3]{FMDon}}] $\,$ \label{Bsp s-formal}
\begin{itemize}
\item[(i)] Every connected and simply-connected compact smooth
manifold is $2$-formal.

\item[(ii)] Every connected and simply-connected compact smooth
manifold of dimension seven or eight is formal if and only if it is
$3$-formal. \q
\end{itemize}
\end{Ex}

\begin{SProp}[{\cite[Lemma 2.11]{FMDon}}]  \label{prod}
Let $M_1,M_2$ be connected smooth manifolds. They are both formal
(resp.\ $s$-formal) if and only if $M_1 \times M_2$ is formal
(resp.\ $s$-formal). \q
\end{SProp}
An important tool for detecting non-formality is the concept of
Massey products: As we shall see below, the triviality of the Massey
products is necessary for formality.

Let $(A,d)$ be a differential graded algebra and $a_i \in
H^{p_i}(A)$, $p_i > 0,~ 1\le i \le 3$, satisfying $a_j \cdot a_{j+1}
= 0$ for $j = 1,2$. Take elements $\alpha_i$ of $A$ with $a_i =
[\alpha_i]$ and write $\alpha_j \cdot \alpha_{j+1} = d \xi_{j,j+1}$
for $j = 1,2$. The \emph{(triple-)Massey product}\index{Massey
Product} $\langle a_1,a_2,a_3 \rangle$ of the classes $a_i$ is
defined as
$$ [\alpha_1 \cdot \xi_{2,3} +
(-1)^{p_1 +1} \xi_{1,2} \cdot \alpha_3] \in \frac{H^{p_1 + p_2 + p_3
- 1}(A)}{a_1 \cdot H^{p_2 + p_3 -1}(A) + H^{p_1 + p_2 -1}(A) \cdot
a_3 }.$$

\begin{Rem}
The definition of the triple-Massey product as an element of a
quotient space is well defined, see e.g.\ \cite[Section 1.6]{TO}.
\end{Rem}

The next lemma shows the relation between formality and Massey
products.

\begin{SLemma}[{\cite[Theorem 1.6.5]{TO}}]
For any formal minimal differential graded algebra all Massey
products vanish. \q
\end{SLemma}

\begin{SCor}
If the de Rahm complex $(\Omega(M),d)$ of a smooth manifold $M$
possesses a non-vanishing Massey product, then $M$ is not formal. \q
\end{SCor}

Fern\'andez and Mu\~noz considered in \cite{FMGeo} the geography of
non-formal compact manifolds. This means they examined whether there
are non-formal compact manifolds of a given dimension with a given
first Betti number. They obtained the following theorem:
\begin{SThm}  \label{nicht sympl}
Given $m \in \mathbb{N}_+$ and $b \in \mathbb{N}$, there are compact
oriented $m$-dimensional smooth manifolds with $b_1 = b$ which are
non-formal if and only if one of the following conditions holds:
\begin{itemize}
\item[(i)] $m \ge 3 \mbox{ and } b \ge 2,$

\item[(ii)] $m \ge 5 \mbox{ and } b = 1$,

\item[(iii)] $m \ge 7 \mbox{ and } b = 0.$ \q
\end{itemize}
\end{SThm}
\subsection{Symplectic, K\"ahler and Lefschetz manifolds}\label{fskl}
The main examples of formal spaces are K\"ahler manifolds. By
definition, a K\"ahler manifold possesses a Riemannian, a symplectic
and a complex structure that are compatible in a sense we are going
to explain now.

Recall that a \emph{symplectic manifold}\index{Symplectic Manifold}
is a pair $(M,\omega)$, where $M$ is a $(2n)$-\-di\-men\-sio\-nal
smooth manifold and $\omega \in \Omega^2(M)$ is a closed $2$-form on
$M$ such that $\omega$ is non-degenerate, i.e.\ $w^n_p \ne 0$ for
all $p \in M$.
\begin{SDef} $\,$
\begin{itemize}
\item[(i)] An \emph{almost complex structure} on an even-dimensional
smooth manifold $M$ is a complex structure $J$ on the tangent bundle
$TM$.

\item[(ii)] Let $M$, $J$ be as in (i) and $\omega \in \Omega^2(M)$ a
non-degenerate $2$-form on M. The $2$-form $\omega$ is called
\emph{compatible} with $J$ if the bilinear form $\langle \ldots ,
\ldots \rangle$ given by
$$\forall_{p \in M}~ \forall_{v,w \in T_p M}~~ \langle v , w \rangle
= \omega(v , Jw)$$ defines a Riemannian metric on $M$.

\item[(iii)] An almost complex structure $J$ on $M$ as in (i)
is called \emph{integrable} if there exists an atlas $\mathcal{A}_M$
on $M$ such that
$$ \forall_{u \in \mathcal{A}_M} ~ \forall_{p \in \mathrm{Domain}(u)}~~
d_p u \circ J_p = J_0 \circ d_p u \: T_pM \to \mathbb{R}^{2n},$$
where
$$ J_0 = \left( \begin{array}{cc} 0 & - Id \\ Id & 0 \end{array}
\right) .$$ $J$ is called \emph{complex structure} for $M$.

\item[(iv)] A \emph{K\"ahler manifold}\index{K\"ahler Manifold} is a
symplectic manifold $(M,\omega)$ with a complex structure $J$ on M
such that $\omega$ is compatible with $J$.
\end{itemize}
\end{SDef}

If one wants to show that a given almost complex structure is not
integrable, it may be hard to disprove the condition (iii) of the
last definition. But in \cite{NN}, Newlander and Nirenberg proved
their famous result that an almost complex structure $J$ on a smooth
manifold $M$ is integrable if and only if $N_J \equiv 0$, where the
\emph{Nijenhuis tensor}\index{Nijenhuis Tensor} $N_J$ is defined as
$$N_J(X,Y) = [JX,JY] - J[JX,Y] - J[X,JY] - [X,Y]$$ for all vector
fields $X,Y$ on $M$.
\N For a time, it was not clear whether every symplectic manifold
was not in fact K\"ahlerian. Meanwhile, many examples of
non-K\"ahlerian symplectic manifolds are known. The first such was
given by Thurston in $1976$ -- the so-called Kodaira-Thurston
manifold, see \cite{Thurston}.

The difficulty to prove non-existence of any K\"ahler structure is
obvious. Nowadays, two easily verifiable necessary conditions for
K\"ahler manifolds are known. First, we have the main theorem from
the work \cite{DGMS} of Deligne, Griffiths, Morgan and Sullivan.

\begin{SThm}[{\cite[p.\ 270]{DGMS}}]
Compact K\"ahler manifolds are formal. \q
\end{SThm}

In order to prove that his manifold is not K\"ahlerian, Thurston
used another method. His manifold has first Betti number equal to
three and the Hodge decomposition for K\"ahler manifolds implies
that its odd degree Betti numbers have to be even, see e.g.\
\cite[pp.\ 116 and 117]{GH}. This is even satisfied for every Hard
Lefschetz manifold.

We say that a symplectic manifold $(M^{2n},\omega)$ is
\emph{Lefschetz}\index{Lefschetz Manifold} if the homomorphism
$$ \begin{array}{cccc}
L^k \: & H^{n-k}(M,\mathbb{R}) & \longrightarrow & H^{n+k}(M,\mathbb{R}) \\
& [\alpha] & \longmapsto & [\alpha \wedge \omega^k]
\end{array}$$
is surjective for $k = n-1$. If $L^k$ is surjective for $k \in
\{0,\ldots,n-1\}$, then $(M,\omega)$ is called \emph{Hard
Lefschetz}\index{Hard Lefschetz Manifold}.

Note that for compact $M$ the surjectivity of $L^k$ implies its
injectivity.

Obviously, the Lefschetz property depends on the choice of the
symplectic form. It may be possible that a smooth manifold $M$
possesses two symplectic forms $\omega_1,\omega_2$ such that
$(M,\omega_1)$ is Lefschetz and $(M,\omega_2)$ not.
But as mentioned above, the existence of such an $\omega_1$ has the
following consequence that is purely topological.

\begin{SThm}  \label{b2i+1 Lefschetz}
The odd degree Betti numbers of a Hard Lefschetz manifold are even.
\end{SThm}

\textit{Proof.} Let $(M^{2n},\omega)$ be a symplectic manifold
satisfying the Lefschetz property. We us the same idea as in
\cite[p.\ 123]{GH}. For each $i \in \{0,\ldots,n-1\}$ one has a
non-degenerated skew-symmetric bilinear form
$$\begin{array}{ccc}
H^{2i+1}(M,\mathbb{R}) \times H^{2i+1}(M,\mathbb{R}) &
\longrightarrow & \mathbb{R}, \\ ([\alpha],[\beta]) & \longmapsto &
[\alpha \wedge \beta \wedge \omega^{n-2i-1}]
\end{array}$$
i.e.\ $H^{2i+1}(M,\mathbb{R})$ must be even-dimensional. \q
\N Obviously, this also proves the next corollary.

\begin{SCor}\label{b1 Lefschetz}
The first Betti number of a Lefschetz manifold is even. \q
\end{SCor}

Finally, the following shows that the statement of the last theorem
holds for K\"ahler manifolds:

\begin{SThm}[{\cite[p.\ 122]{GH}}]
Compact K\"ahler manifolds are Hard Lefschetz. \q
\end{SThm}

\section{Nilmanifolds} \label{Nilmannigfaltigkeiten}
We give a brief review of known results about a special kind of
solvmanifolds, namely nilmanifolds. For the study non-formal
symplectic manifolds, nilmanifolds form one of the best classes. On
the one hand, the non-toral nilmanifolds introduce a geometrical
complexity, while on the other hand their homotopy theory is still
amenable to study. In particular, their minimal models are very easy
to calculate and we shall see that each non-toral nilmanifold is
non-formal.

A \emph{nilmanifold}\index{Nilmanifold} is a compact homogeneous
space $G / \Gamma$, where $G$ is a connected and  simply-connected
nilpotent Lie group and $\Gamma$ a \emph{lattice}\index{Lattice} in
$G$, i.e.\ a discrete co-compact subgroup.

\begin{Ex}
Every lattice in the abelian Lie group $\mathbb{R}^n$ is isomorphic
to $\mathbb{Z}^n$. The corresponding nilmanifold is the
$n$-dimensional torus. \q
\end{Ex}

In contrast to arbitrary solvable Lie groups, there is an easy
criterion for nilpotent ones which enables one to decide whether
there is a lattice or not.

Recall that the exponential map $\exp \: \mathfrak{g} \to G$ of a
connected and simply-connected nilpotent Lie group is a
diffeomorphism. We denote its inverse by $\log \: G \to
\mathfrak{g}$.

\begin{Thm}[{\cite[Theorem 2.12]{Rag}}]  \label{rat Strukturkonst}
A simply-connected nilpotent Lie group G admits a lattice if and
only if there exists a basis $\{X_1, \ldots, X_n\}$ of the Lie
algebra $\mathfrak{g}$ of $G$ such that the structure constants
\index{Structure Constants} $C_{ij}^k$ arising in the brackets
$$ [X_i,X_j] = \sum_k C_{ij}^k \, X_k $$
are rational numbers.

More precisely we have:
\begin{itemize}
\item[(i)] Let $\mathfrak{g}$ have a basis with
respect to which the structure constants are rational. Let
$\mathfrak{g}_{\mathbb{Q}}$ be the vector space over $\mathbb{Q}$
spanned by this basis.

Then, if $\mathcal{L}$ is any lattice of maximal rank in
$\mathfrak{g}$ contained in $\mathfrak{g}_{\mathbb{Q}}$, the group
generated by $\exp(\mathcal{L})$ is a lattice in $G$.
\item[(ii)] If $\Gamma$ is a lattice in $G$, then the
$\mathbb{Z}$-span of $\log(\Gamma)$ is a lattice $\mathcal{L}$ of
maximal rank in the vector space $\mathfrak{g}$ such that the
structure constants of $\mathfrak{g}$ with respect to any basis
contained in $\mathcal{L}$ belong to $\mathbb{Q}$. \q
\end{itemize}
\end{Thm}

For a given lattice $\Gamma$ in a connected and simply-connected
nilpotent Lie group $G$, the subset $\log(\Gamma)$ need not to be an
additive subgroup of the Lie algebra $\mathfrak{g}$.

\begin{Ex}
Consider the nilpotent Lie group $G := \{ \left( \begin{array}{ccc} 1 & x & z \\
0 & 1 & y \\ 0 & 0 & 1 \end{array} \right) \,|\, x,y,z \in
\mathbb{R}
\}$. Its Lie algebra is $\mathfrak{g} := \{ \left( \begin{array}{ccc} 0 & x & z \\
0 & 0 & y \\ 0 & 0 & 0 \end{array} \right) \,|\, x,y,z \in
\mathbb{R} \}$, and the logarithm is given by
$$\log( \left( \begin{array}{ccc} 0 & x & z \\ 0
& 0 & y \\ 0 & 0 & 0 \end{array} \right) ) = \left(
\begin{array}{ccc} 1 & x & z - xy \\ 0 & 1 & y \\ 0 & 0 & 1 \end{array}
\right).$$ The set of integer matrices contained in $G$ forms a
lattice $\Gamma$ in $G$ and $$\log(\Gamma) = \{ \left( \begin{array}{ccc} 0 & a & c \\
0 & 0 & b \\ 0 & 0 & 0 \end{array} \right) \,|\, a,b \in \mathbb{Z},
\, (ab \equiv 0 (2) \Rightarrow c \in \mathbb{Z}), \, (ab \equiv 1
(2) \Rightarrow c \in \frac{1}{2}\mathbb{Z})  \}$$ is not a subgroup
of $\mathfrak{g}$.
\end{Ex}

If $\Gamma$ is a lattice such that $\log(\Gamma)$ is a subgroup of
the Lie algebra, we call $\Gamma$ a \emph{lattice
subgroup}\index{Lattice Subgroup}.

Note that in the context of general Lie groups the name ``lattice
subgroup'' has a different meaning, namely that $G/\Gamma$ has a
finite invariant measure. For nilpotent groups and discrete
$\Gamma$, the latter is the same as to require that $\Gamma$ is a
lattice.

\begin{Thm}[{\cite[Theorem 5.4.2]{CG}}]\label{lattice subgroup}
Let $\Gamma$ be a lattice in a connected and simply-connected
nilpotent Lie group.
\begin{itemize}
\item[(i)] $\Gamma$ contains a lattice subgroup of finite index.
\item[(ii)] $\Gamma$ is contained as a subgroup of finite index in a
lattice subgroup. \q
\end{itemize}
\end{Thm}

For later uses, we quote the following two results.

\begin{Prop}[{\cite[Lemma 5.1.4 (a)]{CG}}]\label{lokal
kompakt}  Let $G$ be a locally compact group, $H$ a closed normal
subgroup and $\Gamma$ a discrete subgroup of $G$. Moreover, denote
by $\pi \: G \to G/H$ the natural map.

If $\Gamma \cap H$ is a lattice in $H$, and $\Gamma$ is a lattice in
$G$, then $\pi(\Gamma)$ is a lattice in $G/H$ and $\Gamma H = H
\Gamma$ is a closed subgroup of $G$. \q
\end{Prop}

\begin{Thm}[{\cite[p.\ 208]{CG}}]  \label{Zentrum}
Let $G$ be a connected and simply-connected nilpotent Lie group with
lattice $\Gamma$ and $k \in \mathbb{N}$.

Then $\Gamma \cap D^{(k)}G$, $\Gamma \cap G^{(k)}$ resp.\ $\Gamma
\cap G_{(k)}$ are lattices in $D^{(k)}G$, $G^{(k)}$ resp.\
$G_{(k)}$. Note, $G_{(1)}$ is the center $Z(G)$ of $G$. \q
\end{Thm}

We have seen that it is easy to decide if there is a lattice in a
given connected and simply-connected nilpotent Lie group, i.e.\ if
it induces a nilmanifold. Moreover, nilmanifolds have very nice
properties which will be described now. Below, we shall see that
these properties are not satisfied for general solvmanifolds.

Note that we can associate a DGA to each Lie algebra $\mathfrak{g}$
as follows:

Let $\{X_1, \ldots, X_n\}$ be a basis of $\mathfrak{g}$ and denote
by $\{x_1, \ldots, x_n\}$ the dual basis of $\mathfrak{g}^*$. The
\emph{Chevalley-Eilenberg complex}\index{Chevalley-Eilenberg
Complex} of $\mathfrak{g}$ is the differential graded algebra
$(\bigwedge \mathfrak{g}^* , \delta)$ with $\delta$ given by $$
\delta(x_k) = - \sum_{i<j} C_{ij}^k \, x_i \wedge x_j,$$ where
$C_{ij}^k$ are the structure constants of $\{X_1, \ldots, X_n\}$.

\begin{Thm}[\cite{N}, {\cite[Theorem 2.1.3]{TO}}]
\label{minModNilmgf} Let $G / \Gamma$ be a nilmanifold and denote by
$\Omega_{l.i.}(G)$ the vector space of left-invariant differential
forms on $G$.

Then the natural inclusion $\Omega_{l.i.}(G) \to \Omega(G / \Gamma)$
induces an isomorphism on cohomology.

Moreover, the minimal model of $G / \Gamma$ is isomorphic to the
Chevalley-Eilenberg complex of the Lie algebra of $G$. \q
\end{Thm}

\begin{Cor}
Any nilmanifold satisfies $b_1 \ge 2$.
\end{Cor}

\textit{Proof.} Let $\mathfrak{g}$ be a nilpotent Lie algebra. By
\cite[Theorem 7.4.1]{Weibel} we have $H^1(\bigwedge
\mathfrak{g}^*,\delta) \cong \mathfrak{g} /
[\mathfrak{g},\mathfrak{g}]$. By \cite{Dix} any nilpotent Lie
algebra $\mathfrak{g}$ satisfies the inequality $\dim \mathfrak{g} /
[\mathfrak{g},\mathfrak{g}] \ge 2$ which then implies
$b_1(\mathfrak{g}) \ge 2$. Hence the claim follows from the
preceding theorem. \q
\N We now quote some results that show that it is easy to decide
whether a nilmanifold is formal, K\"ahlerian or Hard Lefschetz.

\begin{Thm}[{\cite[Theorem 1]{HasNil}}]  \label{formale Nilmgf}
A nilmanifold is formal if and only if it is a torus. \q
\end{Thm}

\begin{Thm}[{\cite[Theorem 2.2.2]{TO}}]
If a nilmanifold is K\"ahlerian, then it is a torus. \q
\end{Thm}

This theorem follows from Theorem \ref{formale Nilmgf}. Another
proof was given by Benson and Gordon in \cite{BG88}. In fact they
proved the following:

\begin{Thm}[{\cite[pp.\ 514 et seq.]{BG88}}] \label{NilLefschetz}
A symplectic non-toral nilmanifold is not Lefschetz. \q
\end{Thm}

\begin{Cor} \label{NilLefschetzgdw}
A symplectic nilmanifold is Hard Lefschetz if and only if it is a
torus, independent of the special choice of the symplectic form.\q
\end{Cor}

\section{Solvmanifolds in general}
A \emph{solvmanifold}\index{Solvmanifold} is a compact homogeneous
space $G / \Gamma$, where $G$ is a connected and simply-connected
solvable Lie group and $\Gamma$ a \emph{lattice}\index{Lattice} in
$G$, i.e.\ a discrete co-compact subgroup.

\begin{Rem}
It is important to note that there is a more general notion of
solvmanifold, namely a compact quotient of a connected and
simply-connected solvable Lie group by a (possibly non-discrete)
closed Lie subgroup (see \cite{Aus}), but we are only considering solvmanifolds
as in the last definition. Sometimes, such are called
special solvmanifolds in the literature.

By \cite[Theorem 2.3.11]{TO}, a solvmanifold in our sense is necessary parallelisable.
E.g.\ the Klein bottle (which can be written as compact homogeneous space of a
three-dimensional connected and simply-connected solvable Lie group)
is not a solvmanifold covered by our definition.
\end{Rem}

Obviously, every nilmanifold is also a solvmanifold. But most
solvmanifolds are not diffeomorphic to nilmanifolds: Every connected
and simply connected solvable Lie group is diffeomorphic to
$\mathbb{R}^m$ (see e.g.\ \cite{V}), hence solvmanifolds are
aspherical and their fundamental group is isomorphic to the
considered lattice. Each lattice in a nilpotent Lie group must be
nilpotent. But in general, lattices in solvable Lie group are not
nilpotent and therefore the corresponding solvmanifolds are not
nilmanifolds.

The fundamental group plays an important role in the study of
solvmanifolds.

\begin{Thm}[{\cite[Theorem 3.6]{Rag}}]
Let $G_i / \Gamma_i$ be solvmanifolds for $i \in \{1,2\}$ and
$\varphi \: \Gamma_1 \to \Gamma_2$ an isomorphism.

Then there exists a diffeomorphism $\Phi \: G_1 \to G_2$ such that
\begin{itemize}
\item[(i)] $\Phi|_{\Gamma_1} = \varphi ,$
\item[(ii)] $\forall_{\gamma \in \Gamma_1} \forall_{p \in G_1}~ \Phi(p
\gamma) = \Phi(p) \varphi(\gamma).$ \q
\end{itemize}
\end{Thm}

\begin{Cor}  \label{Wang eind}
Two solvmanifolds with isomorphic fundamental groups are
diffeomorphic. \q
\end{Cor}

The study of solvmanifolds meets with noticeably greater obstacles
than the study of nilmanifolds. Even the construction of
solvmanifolds is considerably more difficult than is the case for
nilmanifolds. The reason is that there is no simple criterion for
the existence of a lattice in a connected and simply-connected
solvable Lie group.

We shall quote some necessary criteria.

\begin{Prop}[{\cite[Lemma 6.2]{Mil}}]\label{unimod}
If a connected and simply-connected solvable Lie group admits a
lattice then it is unimodular. \q
\end{Prop}

\begin{Thm}[{\cite{Mos},\cite[Theorem 3.1.2]{TO}}]\label{Mostow}
Let $G / \Gamma$ be a solvmanifold that is not a nilmanifold and
denote by $N$ the nilradical of $G$.

Then $\Gamma_N := \Gamma \cap N$ is a lattice in N, $\Gamma N = N
\Gamma$ is a closed subgroup in $G$ and $G / (N \Gamma)$ is a torus.
Therefore, $G / \Gamma$ can be naturally fibred over a non-trivial
torus with a nilmanifold as fiber:
$$ N / \Gamma_N = (N \Gamma) / \Gamma \longrightarrow G / \Gamma
\longrightarrow G / (N \Gamma) = T^k$$ This bundle will be called
the \emph{Mostow bundle}\index{Mostow Bundle}. \q
\end{Thm}

\begin{Rem}\label{BemStGr} $\,$
\begin{itemize}
\item[(i)] The structure group action of the Mostow bundle is given by
left translations
$$ N \Gamma / \Gamma_0 \times N \Gamma / \Gamma \longrightarrow N \Gamma /
\Gamma,$$ where $\Gamma_0$ is the largest normal subgroup of
$\Gamma$ which is normal in $N \Gamma$. (A proof of the topological
version of this fact can be found in \cite[Theorem I.8.15]{Steen}.
The proof for the smooth category is analogous.)
\item[(ii)] A non-toral nilmanifold $G / \Gamma$ fibers over a non-trivial torus
with fibre a nilmanifold, too, since $\Gamma \cap [G,G]$ resp.\
$\mathrm{im} \big( \Gamma \to G / [G,G] \big)$ are lattices in
$[G,G]$ resp.\ $G / [G,G]$, see above.
\end{itemize}
\end{Rem}

In view of Theorem \ref{Mostow}, we are interested in properties of
the nilradical of a solvable Lie group. The following proposition
was first proved in \cite{Mub63}.
\begin{Prop}  \label{dim Nil}
Let $G$ be a solvable Lie group and $N$ its nilradical.

Then $\dim N \ge \frac{1}{2} \dim G$.
\end{Prop}

\textit{Proof.} Denote by $\mathfrak{n} \subset \mathfrak{g}$ the
Lie algebras of $N \subset G$ and by $\mathfrak{n}_{\mathbb{C}}
\subset \mathfrak{g}_{\mathbb{C}}$ their complexifications. Note
that $\mathfrak{g}_{\mathbb{C}}$ is solvable with nilradical
$\mathfrak{n}_{\mathbb{C}}$, so from \cite[Corollary 3.8.4]{V} it
follows that $\mathfrak{n}_{\mathbb{C}} = \{ X \in
\mathfrak{g}_{\mathbb{C}} \,|\,
\mathrm{ad}|_{[\mathfrak{g}_{\mathbb{C}},\mathfrak{g}_{\mathbb{C}}]}
\, \mbox{nilpotent} \}$. Therefore, since $\mathrm{ad} \:
\mathfrak{g}_{\mathbb{C}} \to
\mathrm{Aut}([\mathfrak{g}_{\mathbb{C}} ,
\mathfrak{g}_{\mathbb{C}}])$ is a representation of
$\mathfrak{g}_{\mathbb{C}}$ in $[\mathfrak{g}_{\mathbb{C}} ,
\mathfrak{g}_{\mathbb{C}}]$, by Lie's theorem (e.g.\ \cite[Theorem
3.7.3]{V}) there exist $\lambda_1, \ldots, \lambda_k \in
\mathfrak{g}_{\mathbb{C}}^*$ such that $\mathfrak{n}_{\mathbb{C}} =
\bigcap_{i=1}^k \ker \lambda_i$, where $k := \dim_{\mathbb{C}}
[\mathfrak{g}_{\mathbb{C}} , \mathfrak{g}_{\mathbb{C}}]$.

A straightforward calculation shows $\dim_{\mathbb{C}} \bigcap_{i=1}^k
\ker \lambda_i \ge \dim_{\mathbb{C}} \mathfrak{g}_{\mathbb{C}} - k$,
so we have proven: $\dim_{\mathbb{C}} \mathfrak{n}_{\mathbb{C}} \ge
\dim_{\mathbb{C}} \mathfrak{g}_{\mathbb{C}} - \dim_{\mathbb{C}}
[\mathfrak{g}_{\mathbb{C}} , \mathfrak{g}_{\mathbb{C}}]$.

Because $\mathfrak{g}_{\mathbb{C}}$ is solvable, we get by
\cite[Corollary 3.8.4]{V} that $[\mathfrak{g}_{\mathbb{C}} ,
\mathfrak{g}_{\mathbb{C}}] \subset \mathfrak{n}_{\mathbb{C}}$ and
hence $\dim_{\mathbb{C}} \mathfrak{n}_{\mathbb{C}} \ge
\dim_{\mathbb{C}} \mathfrak{g}_{\mathbb{C}} - \dim_{\mathbb{C}}
\mathfrak{n}_{\mathbb{C}}$, i.e.\
$$ 2 \dim_{\mathbb{C}} \mathfrak{n}_{\mathbb{C}} \ge
\dim_{\mathbb{C}} \mathfrak{g}_{\mathbb{C}}.$$ The proposition now
follows from $\dim_{\mathbb{R}} \mathfrak{g} = \dim_{\mathbb{C}}
\mathfrak{g}_{\mathbb{C}}$ and $\dim_{\mathbb{R}} \mathfrak{n} =
\dim_{\mathbb{C}} \mathfrak{n}_{\mathbb{C}}$. \q
\N In some cases, we will be able to apply the next theorem to the
situation of Theorem \ref{Mostow}. It then gives a sufficient
condition for the Mostow bundle to be a principal bundle.

\begin{Thm}  \label{SolvHauptfaserbdl}
Let $G$ be a connected and simply-connected solvable Lie group and
$\Gamma$ a lattice in $G$. Suppose that $\{e\} \ne H \varsubsetneqq
G$ is a closed normal abelian Lie subgroup of $G$ with $H \subset
N(\Gamma)$, the normalizer of $\Gamma$. (For example the latter is
satisfied if $H$ is central.) Assume further that $\Gamma_H :=
\Gamma \cap H$ is a lattice in $H$.

Then $H / \Gamma_H = H \Gamma / \Gamma$ is a torus and
\begin{eqnarray} \label{H Hauptfaserbdl}
H / \Gamma_H \longrightarrow G / \Gamma \longrightarrow G / H \Gamma
\end{eqnarray}
is a principal torus bundle over a solvmanifold.
\end{Thm}

\textit{Proof.} By assumption, $H$ is a closed normal abelian
subgroup of $G$ and $\Gamma_H$ is a lattice in $H$. We have for $h_1
\gamma_1, h_2 \gamma_2 \in H \Gamma$ with $h_i \in H$, $\gamma_i \in
\Gamma$ the equivalence
$$ (h_1 \gamma_1)^{-1} (h_2 \gamma_2) \in \Gamma \Leftrightarrow
h_1^{-1} h_2 \in \Gamma_H,$$ i.e.\ $H / \Gamma_H = H \Gamma /
\Gamma$. Therefore, Proposition \ref{lokal kompakt} implies that
(\ref{H Hauptfaserbdl}) is a fibre bundle whose fibre is clearly a
torus and its base a solvmanifold. The structure group action is
given by the left translations
$$ H \Gamma / \Gamma_0 \times H \Gamma / \Gamma \longrightarrow
H \Gamma / \Gamma,$$ where $\Gamma_0$ is the largest normal subgroup
of $\Gamma$ which is normal in $H \Gamma$. (This can be seen
analogous as in Remark (i) on page \pageref{BemStGr}.) Since $H$ is
contained in $N(\Gamma) = \{ g \in G \,|\, g \Gamma g^{-1} = \Gamma
\}$, we have for each $h \in H$ and $\gamma, \gamma_0 \in \Gamma$
$$ (h \gamma) \gamma_0 (h \gamma)^{-1} = h \gamma \gamma_0 \gamma^{-1}
h^{-1} \in  h \Gamma h^{-1} = \Gamma, $$ i.e.\ $\Gamma = \Gamma_0$
and the theorem follows. \q
\N We have seen that the Chevalley-Eilenberg complex associated to a
nilmanifold is its minimal model. In this respect, arbitrary
solvmanifolds differ in an essential way from nilmanifolds. However,
in the special case of a solvmanifold which is the quotient of
completely solvable Lie group, one has an access to the minimal
model.

\begin{Def}
Let $G$ be a Lie group with Lie algebra $\mathfrak{g}$.
\begin{itemize}
\item[(i)] $G$ and $\mathfrak{g}$ are
called \emph{completely solvable}\index{Completely Solvable!Lie
Group}\index{Completely Solvable!Lie Algebra} if the linear map
$\mathrm{ad}\, X \: \mathfrak{g} \to \mathfrak{g}$ has only real
roots\footnote{By a root of a linear map, we mean a (possibly
non-real) root of the characteristic polynomial.} for all $X \in
\mathfrak{g}$.
\item[(ii)] If $G$ is simply-connected
and $\exp \: \mathfrak{g} \to G$ is a diffeomorphism, then $G$ is
called \emph{exponential}\index{Exponential Lie Group} .
\end{itemize}
\end{Def}

\begin{Rem}
In the literature a connected and simply-connected Lie group is
sometimes called exponential if the exponential map is surjective.
This is weaker than our definition.
\end{Rem}

A nilpotent Lie group or algebra is completely solvable, and it is
easy to see that completely solvable Lie groups or algebras are
solvable. Moreover, the two propositions below show that
simply-connected completely solvable Lie groups are exponential, and
exponential Lie groups are solvable. Note that the second
proposition is simply a reformulation of results of Sait\^o and
Dixmier.

\begin{Prop}[{\cite[Theorem 2.6.3]{OV}}]
Any exponential Lie group is solvable. \q
\end{Prop}

\begin{Prop} \label{Typ (I)}
A connected and simply-connected solvable Lie group $G$ with Lie
algebra $\mathfrak{g}$ is exponential if and only if the linear map
$\mathrm{ad}\, X \: \mathfrak{g} \to \mathfrak{g}$ has no purely
imaginary roots for all $X \in \mathfrak{g}$. \q
\end{Prop}

\textit{Proof.} Let $G$ be a solvable Lie group. By
\cite[Th\'eor\`eme II.1]{Saito}, $\mathrm{ad}\,X$ has no purely
imaginary roots for all $X \in \mathfrak{g}$ if and only if the
exponential map is surjective. If this is the case,
\cite[Th\'eor\`eme I.1]{Saito} implies that the exponential map is
even bijective. For solvable Lie groups, the statement ``$(1°)
\Leftrightarrow (2°)$'' of \cite[Th\'eor\`eme 3]{Dix2} says that
this is equivalent to the exponential map being a diffeomorphism. \q
\N Let a lattice in a connected and simply-connected solvable Lie
group be given. Then Theorem \ref{Mostow} stated that its
intersection with the nilradical is a lattice in the nilradical. In
the case of completely solvable Lie groups, we have an analogue for
the commutator.

\begin{Prop}[{\cite[Proposition 1]{Gorb}}]
Let $G$ be a connected and simply-connected completely solvable Lie
group and $\Gamma$ a lattice in $G$.

Then $[\Gamma,\Gamma]$ is a lattice in $[G,G]$. In particular,
$\Gamma \cap [G,G]$ is a lattice in $[G,G]$. \q
\end{Prop}

We formulate the result that enables us to compute the minimal model
of solvmanifolds which are built by dividing a lattice out of a
completely solvable group. The main part of the next theorem is due
to Hattori \cite{Hat}.

\begin{Thm}  \label{min chev eil}
Let $G / \Gamma$ be a solvmanifold. Denote by $(\bigwedge
\mathfrak{g}^* , \delta)$ the Chevalley-Eilenberg complex of $G$ and
recall that $\mathfrak{g}^*$ is the set of left-invariant
differential $1$-forms on $G$. Then the following holds:
\begin{itemize}
\item[(i)] The natural inclusion $(\bigwedge
\mathfrak{g}^* , \delta) \to (\Omega(G / \Gamma) , d)$ induces an
injection on cohomology.
\item[(ii)] If $G$ is completely solvable, then the inclusion
in (i) is a quasi-isomorphism, i.e.\ it induces an isomorphism on
cohomology. Therefore, the minimal model $\mathcal{M}_{G / \Gamma}$
is isomorphic to the minimal model of the Chevalley-Eilenberg
complex.
\item[(iii)] If $\mathrm{Ad}\,(\Gamma) \,$ and $\, \mathrm{Ad}\,(G)$
have the same Zariski closures\footnote{A basis for the Zariski
topology on $\mathrm{GL}(m,\mathbb{C})$ is given by the sets $$U_p
:= \mathrm{GL}(m,\mathbb{C}) \setminus p^{-1}(\{0\}),$$ where $p \:
\mathrm{GL}(m,\mathbb{C}) \cong \mathbb{C}^{(m^2)} \to \mathbb{C}$
ranges over polynomials.}, then the inclusion in (i) is a
quasi-isomorphism. \q
\end{itemize}
\end{Thm}

\textit{Proof.} (i) is \cite[Theorem 3.2.10]{TO} and (iii) is
\cite[Corollary 7.29]{Rag}.

ad (ii): Denote the mentioned inclusion by $i \: (\bigwedge
\mathfrak{g}^* , \delta) \to (\Omega(G / \Gamma) , d)$. By Hattori's
Theorem (see \cite[p.\ 77]{TO}), this is a quasi-isomorphism. It
remains to show that the minimal model $\rho \: (\mathcal{M}_{CE} ,
\delta_{CE}) \to (\bigwedge \mathfrak{g}^* , \delta)$ of $(\bigwedge
\mathfrak{g}^* , \delta)$ is the minimal model of $(\Omega(G /
\Gamma) , d)$. Since the minimal model is unique and the map $i
\circ \rho \: (\mathcal{M}_{CE} , \delta_{CE}) \to (\Omega(G /
\Gamma) , d)$ is a quasi-isomorphism, this is obvious. \q
\N There are examples where the inclusion in (i) in the last theorem
is not a quasi-isomorphism: Consider the Lie group $G$ which is
$\mathbb{R}^3$ as a manifold and whose Lie group structure is given
by
$$(s,a,b) \cdot (t,x,y) = (s+t,\cos(2 \pi t) \, a - \sin(2 \pi t) \, b +
x,\sin(2 \pi t) \, a + \cos(2 \pi t) \, b + y ).$$ $G$ is not
completely solvable and one calculates for its Lie algebra
$b_1(\mathfrak{g}) =1$. $G$ contains the abelian lattice $\Gamma :=
\mathbb{Z}^3$ and $G/\Gamma$ is the $3$-torus which has $b_1 = 3$.
\N We have seen in the last section that the first Betti number of a
nilmanifold is greater than or equal to two. For arbitrary
solvmanifolds this is no longer true. Below, we shall see various
examples of solvmanifolds with $b_1=1$. The following corollary
shows that $b_1 = 0$ cannot arise.

\begin{Cor}
Any solvmanifold satisfies $b_1 \ge 1$.
\end{Cor}

\textit{Proof.} Let $\mathfrak{g}$ be a solvable Lie algebra. As in
the nilpotent case we have $b_1(\bigwedge \mathfrak{g}^*,\delta) =
\dim \mathfrak{g} / [\mathfrak{g},\mathfrak{g}]$, and $\dim
\mathfrak{g} / [\mathfrak{g},\mathfrak{g}] \ge 1$ by solvability.
The claim now follows from Theorem \ref{min chev eil} (i). \q
\N To end this section, we shortly discuss the existence problem for
K\"ahler structures on solvmanifolds. The only K\"ahlerian
nilmanifolds are tori, but in the general context we have the
hyperelliptic surfaces, which are non-toral K\"ahlerian
solvmanifolds, see Section \ref{KapSolvmgfen in Dim4} below.
(\cite[Theorem 3.4.1]{TO} states the only K\"ahlerian solvmanifolds
in dimension four are tori. This is not correct, as first noted by
Hasegawa in \cite{HasK}.) Benson and Gordon \cite{BG90} conjectured
in $1990$ that the existence of a K\"ahler structure on a
solvmanifold $G/\Gamma$ with $G$ completely solvable forces
$G/\Gamma$ to be toral and this is true. In fact, Hasegawa proved in
the first half of this decade the following:

\begin{Thm}[\cite{HasK}] \label{Kaehler vollst aufl}
A solvmanifold $G / \Gamma$ is K\"ahlerian if and only if it is a
finite quotient of a complex torus which has a structure of a
complex torus bundle over a complex torus.

If $G$ is completely solvable, then $G / \Gamma$ is K\"ahlerian if
and only if it is a complex torus. \q
\end{Thm}

In later sections we shall see that neither the Hard Lefschetz
property nor formality is sufficient for an even-dimensional
solvmanifold to be K\"ahlerian.

\section{Semidirect products}
In later sections we shall try to examine low-dimensional
solvmanifolds. Concerning this, a first step is to use the known
classification of the (connected and simply-connected)
low-dimensional solvable Lie groups. Most of them have the structure
of semidirect products. In order to define this notion, we recall
the construction of the Lie group structure of the group of Lie
group automorphisms of a simply-connected Lie group in the following
theorem. It collects results that can be found in \cite[pp.\ 117 et
seq.]{War}.

\begin{Thm}\label{Aut} $\,$
\begin{itemize}
\item[(i)]
Let $\mathfrak{h} = \big( |\mathfrak{h}| = \mathbb{R}^h ,
[\ldots,\ldots] \big)$ be an $h$-dimensional Lie algebra. Then the
set $\mathrm{A}(\mathfrak{h})$ of Lie algebra isomorphisms of
$\mathfrak{h}$ is a closed Lie subgroup of the automorphism group
$\mathrm{Aut}(|\mathfrak{h}|)$ of the $h$-dimensional vector space
$|\mathfrak{h}|$. The Lie algebra of $\mathrm{A}(\mathfrak{h})$ is
$$\mathfrak{d}(\mathfrak{h}) = \{ \varphi \in
\mathrm{End}(|\mathfrak{h}|) \,|\, \varphi \mbox{ derivation with
respect to } [\ldots,\ldots] \}.$$
\item[(ii)]
Let $H$ be a connected and simply-connected Lie group with neutral
element $e$ and Lie algebra $\mathfrak{h}$. The Lie group structure
of $\mathrm{A}(H)$, the group of Lie group automorphisms of $H$, is
given by the following group isomorphism:
$$ \mathrm{A}(H) \longrightarrow \mathrm{A}(\mathfrak{h}) ~,~ f \longmapsto d_e f$$

Moreover, if $H$ is exponential, its inverse is the map $$
\mathrm{A}(\mathfrak{h}) \longrightarrow \mathrm{A}(H) ~,~ \varphi
\longmapsto \exp^{H} \circ \varphi \circ \log^{H}.$$ \q
\end{itemize}
\end{Thm}

For given (Lie) groups $G,H$ and a (smooth) action $\mu \: G \times
H \to H$ by (Lie) group automorphisms, one defines the
\emph{semidirect product}\index{Semidirect Product!of (Lie) Groups}
\emph{of} $G$ \emph{and} $H$ \emph{via} $\mu$ as the (Lie) group $G
\ltimes_{\mu} H$ with underlying set (manifold) $G \times H$ and
group structure defined as follows:
$$ \forall_{(g_1,h_1),(g_2,h_2) \in G \times H}~ (g_1,h_1)
(g_2,h_2) = \big( g_1 g_2 , \mu(g_2^{-1} , h_1) h_2 \big)$$
Note that for $(g,h) \in G \ltimes_{\mu} H$ we have $(g,h)^{-1} =
\big( g^{-1}, \mu(g,h^{-1}) \big)$.

If the action $\mu$ is trivial, i.e.\ $\forall_{g \in G,\, h \in H}~
\mu(g,h) =h$, one obtains the ordinary direct product. In the case
of Lie groups $G$ and $H$, the exponential map $\exp^{G \times H}$
is known to be the direct product of $\exp^G$ and $\exp^H$. If the
action is not trivial, the situation becomes a little more
complicated:

\begin{Thm}\label{expsemidirekt}
Let $G,H$ be connected Lie groups and $\mu \: G \times H \to H$ a
smooth action by Lie group automorphisms. Denote the Lie algebras of
$G$ and $H$ by $\mathfrak{g}$ and $\mathfrak{h}$ and let $\phi :=
(d_{e_G} \mu_1) \: \mathfrak{g} \to \mathfrak{d}(\mathfrak{h}),$
where $\mu_1 \: G \to \mathrm{A}(\mathfrak{h})$ is given by
$\mu_1(g) = d_{e_H}\mu(g,\ldots) = \mathrm{Ad}^{G\ltimes_{\mu}H}_g$.
\begin{itemize}
\item[(i)] The Lie algebra of $G \ltimes_{\mu} H$ is $\mathfrak{g}
\ltimes_{\phi} \mathfrak{h}$. This Lie algebra is called semidirect
product \index{Semidirect Product!of Lie Algebras} of $\mathfrak{g}$
and $\mathfrak{h}$ via $\phi$. Its underlying vector space is
$\mathfrak{g} \times \mathfrak{h}$ and the bracket for
$(X_1,Y_1),(X_2,Y_2) \in \mathfrak{g} \times \mathfrak{h}$ is given
by
$$[(X_1,Y_1),(X_2,Y_2)] = \big(
[X_1,X_2]_{\mathfrak{g}} , [Y_1,Y_2]_{\mathfrak{h}} + \phi(X_1)(Y_2)
- \phi(X_2)(Y_1) \big).$$ In the sequel we shall identify $X \equiv
(X,0)$ and $Y \equiv (0,Y)$.
\item[(ii)] For $(X,Y) \in \mathfrak{g} \ltimes_{\phi} \mathfrak{h}$
one has $ \exp^{G\ltimes_{\mu}H}((X,Y)) = (\exp^G(X), \gamma(1)),$
where $\gamma \: \mathbb{R} \to H$ is the solution of $$
\dot{\gamma}(t) = (d_{e_H} R_{\gamma(t)}) \big(
\exp^{A(\mathfrak{h})}(-t \, \mathrm{ad}(X)|_{\mathfrak{h}})(Y)
\big),~~ \gamma(0) = e_H.$$ Here $R_a$ denotes the right translation
by an element $a \in H$.
\end{itemize}
\end{Thm}

\textit{Proof.} The proof of (i) can be found in \cite{V}. We give a
proof of (ii). Given a Lie group homomorphism $f$ between Lie
groups, we denote its differential at the neutral element by $f_*$.

For $(g_0,h_0),(g,h) \in G \ltimes_{\mu} H$ we have $R_{(g_0,h_0)}
(g,h) = (R_{g_0}(g) , R_{h_0}( \mu(g_0^{-1},h))$, and this yields
for $(X,Y) \in \mathfrak{g} \ltimes_{\phi} \mathfrak{h}$
$$(R_{(g_0,h_0)})_* \big( (X,Y) \big) = \big((R_{g_0})_*(X) , (R_{h_0})_*
\big( \mu_1(g_0^{-1}) (Y) \big) \big).$$ Since
$\big(\gamma_1(t),\gamma_2(t)\big) := \exp^{G\ltimes_{\mu}H}\big(
t\, (X,Y) \big)$ is the integral curve through the identity of both
the right- and left-invariant vector fields associated to $(X,Y)$,
the last equation implies that $\big(\gamma_1(t),\gamma_2(t)\big)$
is the solution of the following differential equations:
\begin{eqnarray}
\label{Komp1} \gamma_1(0) = e_G,  && \dot{\gamma_1} (t) = (R_{\gamma_1(t)})_* (X), \\
\label{Komp2} \gamma_2(0) = e_H, && \dot{\gamma_2} (t) =
(R_{\gamma_2(t)})_* (\mu_1(\gamma_1(-t))(Y)).
\end{eqnarray}
$\gamma_1(t) = \exp^G(t \, X)$ is the solution of (\ref{Komp1}), and
this implies $$\mu_1(\gamma_1(-t)) = \mathrm{Ad}^{G \ltimes_{\mu}
H}_{\gamma_1(-t)}|_{\mathfrak{h}} = \exp^{A(\mathfrak{h})}(-t \,
\mathrm{ad}(X)|_{\mathfrak{h}}),$$ i.e.\ (\ref{Komp2}) is equivalent
to $\gamma_2(0) = e_H,~ \dot{\gamma_2} (t) = (R_{\gamma_2(t)})_*
(\exp^{A(\mathfrak{h})}(-t \, \mathrm{ad}(X)|_{\mathfrak{h}})(Y))$.
So the theorem is proven. \q
\N A connected and simply-connected solvable Lie group $G$ with
nilradical $N$ is called \emph{almost nilpotent}\index{Almost
Nilpotent Lie Group} if it can be written as $G = \mathbb{R}
\ltimes_{\mu} N$. Moreover, if $N$ is abelian, i.e.\ $N =
\mathbb{R}^n$, then $G$ is called \emph{almost
abelian}.\index{Almost Abelian Lie Group}

Let $G = \mathbb{R} \ltimes_{\mu} N$ be an almost nilpotent Lie
group. Since $N$ has codimension one in $G$, we can consider $\mu$
as a one-parameter group $\mathbb{R} \to \mathrm{A}(N)$. By Theorem
\ref{Aut}, there exists $\varphi \in \mathfrak{d}(\mathfrak{n})$
with
$$\forall_{t \in \mathbb{R}}~ \mu(t) = \exp^{N} \circ
\exp^{Aut(|\mathfrak{n}|)} (t \varphi) \circ \log^{N}.$$ Choosing a
basis of $|\mathfrak{n}|$, we can identify
$\mathrm{Aut}(|\mathfrak{n}|)$ with a subset of
$\mathfrak{gl}(n,\mathbb{R})$ and get
$$\forall_{t \in \mathbb{R}} ~ d_e \big(\mu(t)\big) \in
\exp^{GL(n,\mathbb{R})} \big( \mathfrak{gl}(n,\mathbb{R}) \big) .$$
Note, if $N$ is abelian, the exponential map $\exp^N \: \mathfrak{n}
\to N$ is the identity. These considerations make it interesting to
examine the image of $\exp^{GL(n,\mathbb{R})}$.

\begin{Thm}[{\cite[Theorem 6]{Nono}}]\label{Bild exp}
$M$ is an element of
$\exp^{GL(n,\mathbb{R})}(\mathfrak{gl}(n,\mathbb{R}))$ if and only
if the real Jordan form of $M$ contains in the form of pairs the
blocks belonging to real negative eigenvalues $\lambda_i^-$,
whenever there exist real negative eigenvalues $\lambda_i^-$ of $M$.
I.e.\ the block belonging to such a $\lambda_i^-$ is of the
following form $$ \bigoplus_{j=1}^{n_i} \left(
\begin{array}{cc} J_{n_{ij}} & 0 \\ 0 & J_{n_{ij}}
\end{array} \right)$$
with $$ J_{n_{ij}} = \left(
\begin{array}{cccc} \lambda_i^- & 1  & & 0 \\  & \lambda_i^- & \ddots & \\
& & \ddots & 1 \\ 0 & & &  \lambda_i^-
\end{array} \right) \in \mathrm{M}(n_{ij},n_{ij};\mathbb{R}).$$ \q
\end{Thm}

We are now going to derive some facts that follow from the existence
of a lattice in an almost nilpotent Lie group.

\begin{Thm}[\cite{Tralle}]\label{Gitter vollst aufl}
Let $G = \mathbb{R} \ltimes_{\mu} N$ be an almost nilpotent and
completely solvable Lie group containing a lattice $\Gamma$.

Then there is a one-parameter group $\nu \: \mathbb{R} \to
\mathrm{A}(N)$ such that $\nu (k)$ preserves the lattice $\Gamma_N
:= \Gamma \cap N$ for all $k \in \mathbb{Z}$. $\Gamma$ is isomorphic
to $\mathbb{Z} \ltimes_{\nu} \Gamma_N$ and $G / \Gamma$ is
diffeomorphic to $\big( \mathbb{R} \ltimes_{\nu} N \big) / \big(
\mathbb{Z} \ltimes_{\nu} \Gamma_N$ \big).

Moreover, there are $t_1 \in \mathbb{R} \setminus \{0\}$ and an
inner automorphism $I_{n_1} \in \mathrm{A}(N)$ such that $\nu(1) =
\mu(t_1) \circ I_{n_1}$.
\end{Thm}

\textit{Proof.} We know that $\Gamma_N$ is a lattice in $N$ and
$\mathrm{im}(\Gamma \to G/N) \cong \Gamma / \Gamma_N$ is a lattice
in $G / N \cong \mathbb{R}$. Therefore, $\Gamma / \Gamma_N \cong
\mathbb{Z}$ is free, and the following exact sequence is split:
$$ \{1\} \longrightarrow \Gamma_N \longrightarrow \Gamma
\longrightarrow \mathbb{Z} \longrightarrow \{0\},$$ i.e.\ there is a
group-theoretic section $s \: \mathbb{Z} \to \Gamma$.
\cite[Th\'eor\`eme II.5]{Saito} states that a group homomorphism
from a lattice of completely solvable Lie group into another
completely solvable Lie group uniquely extends to a Lie group
homomorphism of the Lie groups. Hence, $s$ extends uniquely to a
one-parameter group $s \: \mathbb{R} \to G$. Therefore,
$$\nu \: \mathbb{R} \longrightarrow \mathrm{A}(N), ~~ \nu (t) (n) = s(t) \cdot n
\cdot s(t)^{-1}$$ is a one-parameter group with $\forall_{k \in
\mathbb{Z}}~ \nu (k) (\Gamma_N) = \Gamma_N$, the lattice $\Gamma$ is
isomorphic to $\mathbb{Z} \ltimes_{\nu} \Gamma_N$ and $G / \Gamma$
is diffeomorphic to $\big( \mathbb{R} \ltimes_{\nu} N \big) / \big(
\mathbb{Z} \ltimes_{\nu} \Gamma_N \big)$.

Let $\gamma_1 := s(1) \in (\Gamma \setminus \Gamma_N )\subset
\mathbb{R} \ltimes_{\mu} N$. There are unique $t_1 \in \mathbb{R}
\setminus \{ 0 \}$, $n_1 \in N$ with $\gamma_1 = t_1 \cdot n_1$,
where we identify $t_1 \equiv (t_1,e_N) \in G$ and $n_1 \equiv
(0,n_1) \in G$. Since $G = \mathbb{R} \ltimes_{\nu} N$ and $G =
\mathbb{R} \ltimes_{\mu} N$ with the same normal subgroup $N$ of
$G$, one has for all $n \in N$
$$ \nu(1) (n) = \gamma_1 \cdot n \cdot \gamma_1^{-1}
= t_1 \cdot n_1 \cdot n \cdot n_1^{-1} \cdot t_1^{-1} = \mu(t_1)
(n_1 \cdot n \cdot n_1^{-1}) = \mu(t_1) (I_{n_1}(n)),$$ from where
the theorem follows. \q

\begin{Cor}  \label{Gitter in fastabelsch}
Let $G = \mathbb{R} \ltimes_{\mu} N$ be an almost nilpotent (not
necessary completely solvable) Lie group containing a lattice
$\Gamma$. Again, denote by $\Gamma_N := \Gamma \cap N$ the induced
lattice in the nilradical of $G$.

Then there exist $t_1 \in \mathbb{R} \setminus \{0\}$, a group
homomorphism $\nu \: \mathbb{Z} \to \mathrm{Aut}(\Gamma_N)$, and an
inner automorphism $I_{n_1}$ of $N$ such that $\Gamma \cong
\mathbb{Z} \ltimes_{\nu} \Gamma_N$ and $\nu(1) = \mu(t_1) \circ
I_{n_1}$.

If $G$ is almost abelian, then a basis transformation yields $\Gamma
\cong t_1 \mathbb{Z} \ltimes_{\mu|_{\mathbb{Z}^n}} \mathbb{Z}^n$.
\end{Cor}

\textit{Proof.} We argue as in the last proof. But we do not use
\cite[Th\'eor\`eme 5]{Saito} and get only a group homomorphism $\nu
\: \mathbb{Z} \to \mathrm{Aut}(\Gamma_N)$ (defined on $\mathbb{Z}$
instead of $\mathbb{R}$). For general $N$, the
calculation at the end of the proof implies the claim.

Since an abelian group has only one inner automorphism, in the
almost abelian case this yields $\nu(1) = \mu(t_1)|_{\Gamma_N}$, so
$\nu$ can be extended to $\nu \: \mathbb{R} \to
\mathrm{A}(\mathbb{R}^n)$ via $\nu(t) := \mu(t \cdot t_1)$. Further,
by Corollary \ref{Wang eind}, we have $\Gamma_N \cong \mathbb{Z}^n$.
\q
\N Hence we have seen, that the existence of a lattice in an almost
nilpotent Lie group implies that a certain Lie group automorphism
must preserve a lattice in the (nilpotent) nilradical. The next
theorem deals with such automorphisms.

\begin{Thm} \label{1ParamSL}
Let $N$ be a connected and simply-connected nilpotent Lie group with
Lie algebra $\mathfrak{n}$, $f_* \in \mathrm{A}(\mathfrak{n})$, and
$f := \exp^N \circ f_* \circ \log^N \in \mathrm{A}(N)$, i.e.\ $d_e f
= f_*$. Assume that $f$ preserves a lattice $\Gamma$ in $N$.

Then there exists a basis $\mathfrak{X}$ of $\mathfrak{n}$ such that
$M_{\mathfrak{X}}(f_*) \in \mathrm{GL}(n,\mathbb{Z})$, where
$M_{\mathfrak{X}}(f_*)$ denotes the matrix of $f_*$ with respect to
$\mathfrak{X}$.

Moreover, if there are a one-parameter group $\mu \: \mathbb{R} \to
\mathrm{A}(N)$ and $t_0 \ne 0$ such that $\mu(t_0) = f$, i.e.\
$d_e(\mu(t_0)) = f_*$, then $\det\big(d_e(\mu(\ldots))\big) \equiv
1$.
\end{Thm}

\textit{Proof.} By Theorem \ref{rat Strukturkonst} (ii),
$$\mathcal{L} := \langle \log^{N}(\Gamma) \rangle_{\mathbb{Z}} = \{
\sum_{i=1}^m k_i \, V_i \,|\, m \in \mathbb{N}_+, k_i \in
\mathbb{Z}, V_i \in \log^{N}(\Gamma) \}$$ is a lattice in
$\mathfrak{n}$. Therefore, there exists a basis $\mathfrak{X} = \{
X_1 , \ldots , X_n \}$ of $\mathfrak{n}$ such that $\mathcal{L} =
\langle \mathfrak{X} \rangle_{\mathbb{Z}}$.

Since  $f(\Gamma) \subset \Gamma$, we have $f_* \big( \log^N
(\Gamma) \big) \subset \log^N (\Gamma)$. This implies
$f_*(\mathcal{L}) \subset \mathcal{L}$ and hence,
$M_{\mathfrak{X}}(f_*) \in \mathrm{GL}(n,\mathbb{Z})$.

Further, if $\mu(t_0) = f$ with $\mu$, $t_0 \ne 0$ as in the
statement of the theorem, then the map $\Delta := \det \circ
d_e(\mu(\ldots)) \: (\mathbb{R},+) \to (\mathbb{R} \setminus \{0\} ,
\cdot)$ is a continuous group homomorphism with $\Delta(0) = 1$ and
$\Delta(t_0) = \pm 1,$ i.e.\ $\Delta \equiv 1.$ \q

\begin{Rem}
The basis $\mathfrak{X}$ in the last theorem has rational structure
constants.
\end{Rem}

Obviously, a one-parameter group $\mu$ in the automorphism group of
an abelian Lie group with $\mu(t_0)$ integer valued for $t_0 \ne 0$
defines a lattice in $\mathbb{R} \ltimes_{\mu} \mathbb{R}^n$. It is
easy to compute the first Betti number of the corresponding
solvmanifold, as the next proposition will show. Before stating it,
we mention that the situation becomes more complicated in the case
of a non-abelian and nilpotent group $N$.

Let a one-parameter group $\mu \: \mathbb{R} \to \mathrm{A}(N)$ be
given and $t_0 \ne 0$ such that $d_e(\mu(t_0))$ is an integer matrix
with respect to a basis $\mathfrak{X}$ of the Lie algebra
$\mathfrak{n}$ of $N$. In general, this does not enable us to define
a lattice in $\mathbb{R} \ltimes_{\mu} N$. But if $\Gamma_N :=
\exp^N(\langle \mathfrak{X} \rangle_{\mathbb{Z}})$ is a lattice in
$N$, i.e.\ $\Gamma_N$ is a lattice group, then this is possible.

\begin{Prop}  \label{b1 fast abelsch}
Let $\mu \: \mathbb{R} \to \mathrm{SL}(n,\mathbb{R})$ be a
one-parameter group such that $\mu(1) = (m_{ij})_{i,j} \in
\mathrm{SL}(n,\mathbb{Z})$.

Then $M := (\mathbb{R} \ltimes_{\mu} \mathbb{R}^n) / (\mathbb{Z}
\ltimes_{\mu} \mathbb{Z}^n)$ is a solvmanifold with
\begin{eqnarray*} \pi_1(M) = \langle e_0, e_1, \ldots , e_n ~ | &
\forall_{i \in \{1,\ldots,n\}} ~ e_0 e_i e_0^{-1} = e_1^{m_{1i}}
\cdots e_n^{m_{ni}} & \\ & \forall_{i,j \in \{1,\ldots,n\}} ~ [e_i ,
e_j] = 1 & \rangle
\end{eqnarray*}
and $b_1(M) = n+1 - \mathrm{rank}\, \big( \mu(1) - \mathrm{id}
\big).$
\end{Prop}

\textit{Proof.} The statement about the fundamental group is clear.
Therefore, we get
\begin{eqnarray*}
H_1(M , \mathbb{Z}) = \langle e_0, e_1, \ldots , e_n ~ | &
\forall_{i \in \{1,\ldots,n\}} ~ e_1^{m_{1i}} \cdots e_i^{m_{ii}-1}
\cdots e_n^{m_{ni}} = 1 \\ & \forall_{i,j \in \{0,\ldots,n\}} ~ [e_i
, e_j] = 1 & \rangle
\end{eqnarray*}
and this group is the abelianisation of $$\mathbb{Z} \oplus \langle
e_1, \ldots , e_n \, |\, \forall_{i \in \{1,\ldots,n\}} \,
e_1^{m_{1i}} \cdots e_i^{m_{ii}-1} \cdots e_n^{m_{ni}} = 1
\rangle.$$ Now, the proof of the theorem about finitely generated
abelian groups (see e.g.\ \cite{Bosch}) shows $H_1(M , \mathbb{Z}) =
\mathbb{Z}^{n-k+1} \oplus \bigoplus_{i=1}^k \mathbb{Z}_{d_i},$ where
$d_1, \ldots , d_k \in \mathbb{N}_+$ denote the elementary divisors
of $\mu(1) - \mathrm{id}$. The proposition follows. \q
\N We finally mention a result of Gorbatsevich. In view of Theorem
\ref{min chev eil} (iii), it enables us to compute the minimal model
of a wide class of solvmanifolds which are discrete quotients of
almost abelian Lie groups.

\begin{Thm}[{\cite[Theorem 4]{Gorb}}] \label{Gorb pi}
Let $\mu \: \mathbb{R} \to \mathrm{SL}(n,\mathbb{R})$ be a
one-parameter group such that $\mu(1) =
\exp^{SL(n,\mathbb{R})}(\dot{\mu}(0)) \in
\mathrm{SL}(n,\mathbb{Z})$. Denote by $\lambda_1, \ldots, \lambda_n$
the (possibly not pairwise different) roots of $\dot{\mu}(0)$. Then
$\Gamma := (\mathbb{Z} \ltimes_{\mu} \mathbb{Z}^n)$ is a lattice in
$G := (\mathbb{R} \ltimes_{\mu} \mathbb{R}^n)$.

The Zariski closures of $\mathrm{Ad}\,(\Gamma)$ and
$\mathrm{Ad}\,(G)$ coincide if and only if the number $\pi i$ is not
representable as a linear combination of the numbers $\lambda_k$
with rational coefficients. \q
\end{Thm}

In distinction from the nilpotent case, criteria for the existence
of a lattice in connected and simply-connected solvable Lie groups
have rather cumbersome formulations. The criterion that we present
is due to Auslander \cite{Aus} and makes use of the concept of
semisimple splitting.

Let $G$ be a connected and simply-connected Lie group. We call a
connected and simply-connected solvable Lie group $G_s = T
\ltimes_{\nu_s} N_s$ a \emph{semisimple splitting for}
$G$\index{Semisimple Splitting} if the following hold:
\begin{itemize}
\item[(i)] $N_s$ is the nilradical of $G_s$ -- the so called
\emph{nilshadow of }$G$\index{Nilshadow} -- and $T \cong
\mathbb{R}^k$ for $k = \dim G_s - \dim N_s$,
\item[(ii)] $T$ acts on $N_s$ via $\nu_s$ by semisimple
automorphisms,
\item[(iii)] $G$ is a closed normal subgroup of $G_s$ and $G_s =
T \ltimes_{\varpi} G$,
\item[(iv)] $N_s = Z_{N_s}(T) \cdot (N_s \cap G)$, where
$Z_{N_s}(T)$ denotes the centralizer of $T$ in $N_s$.
\end{itemize}
This definition then implies (see e.g. \cite[Lemma 5.2]{Dekimpe})
that $N_s$ is a connected and simply-connected nilpotent Lie group,
$N = N_s \cap G$, and $G/N \cong T$.

\begin{Thm}[{\cite[Theorems 5.3 and 5.4]{Dekimpe}}]
Let $G$ be a connected and simply-connected solvable Lie group. Then
$G$ admits a unique semisimple splitting.
\end{Thm}

We shall not give the whole proof of this theorem that can be found
in \cite{Dekimpe}. But we shortly describe the construction of the
semisimple splitting. In order to do this, we recall the Jordan
decomposition\index{Jordan Decomposition} of certain morphisms:

Let $\varphi$ be an endomorphism of a finite-dimensional vector
space over a field of characteristic zero. There is a unique
\emph{Jordan sum decomposition}
\begin{eqnarray*}
\varphi = \varphi_s + \varphi_n, && \varphi_s \circ \varphi_n =
\varphi_n \circ \varphi_s,
\end{eqnarray*}
with $\varphi_s$ semisimple and $\varphi_n$ nilpotent. They are
called respectively the \emph{semisimple part} and the
\emph{nilpotent part of} $\varphi$. If $\varphi$ is an automorphism,
it also has a unique \emph{Jordan product decomposition}
\begin{eqnarray*}
\varphi = \varphi_s \circ \varphi_u, && \varphi_s \circ \varphi_u =
\varphi_u \circ \varphi_s,
\end{eqnarray*}
with $\varphi_s$ semisimple and $\varphi_u$ unipotent; $\varphi_s$
is the same as in the sum decomposition and $\varphi_u = \mathrm{id}
+ (\varphi_s^{-1} \circ \varphi_n)$. The latter is called the
\emph{unipotent part of} $\varphi$.

Note, if $\varphi$ is a derivation resp.\ an automorphism of a Lie
algebra, then the semisimple and the nilpotent resp.\ unipotent part
of $\varphi$ are also derivations resp.\ automorphisms of the Lie
algebra.

Now, let $G$ be a connected and simply-connected Lie group and $f \:
G \to G$ a Lie group automorphism.

Then $f_* := d_ef$ is a Lie algebra automorphism which has a Jordan
product decomposition $f_* = (f_*)_s \circ (f_*)_u = (f_*)_u \circ
(f_*)_s$. The \emph{semisimple} and \emph{unipotent part of} $f$ are
by definition the unique Lie group automorphisms $f_s,f_u \: G \to
G$ with $d_ef_s = (f_*)_s$ and $d_ef_u = (f_*)_u$.

\textit{Construction of the semisimple
splitting.}\label{Konstruktsemisimple} Let $G$ be a connected and
simply-connected solvable Lie group. Denote by $N$ the nilradical of
$G$.

By \cite[Proposition 3.3]{Dekimpe}, there exists a connected and
simply-connected nilpotent Lie subgroup $H$ of $G$ such that $G = H
\cdot N$. Fix such an $H$ and consider the well-defined (!) action
$\widetilde{\varpi} \: H \to \mathrm{A}(G)$ given by
$\widetilde{\varpi}(a)(h \cdot n) := h \cdot (I_a|_N)_s (n)$, where
$(I_a|_N)_s$ is the semisimple part of the automorphism of $N$ which
is obtained by conjugating every element of $N$ by $a$.

Define $T := H / (N \cap H) \cong HN/N \cong G/N \cong
\mathbb{R}^k$. Note that there is an action $\varpi$ of $T$ on $G$
making the following diagram commutative:
$$
\begin{diagram}
H & && \rTo^{\widetilde{\varpi}} && &\mathrm{A}(G) \\
~~~~~\pi \! \! \! & \rdTo && &&  \ruTo & \! \! \! \varpi~~~~~ \\
& && T = H / (H \cap N) && & \\
\end{diagram}
$$
Set $G_s := T \ltimes_{\varpi} G$.

One calculates that $$N_s := \{\pi(h^{-1}) \cdot h \,|\, h \in H \}
\cdot N = \{ (\pi(h^{-1}), h \cdot n) \,|\, h \in H, n \in N \}
\subset T \ltimes G = G_s$$ is the nilradical of $G_s$. Furthermore,
we have $T \cdot N_s = G_s$ and $T \cap N_s = \{e\}$. For $t \in T$,
$h \in H$, $n \in N$ and every $h_t \in \pi^{-1}(\{t\})$ holds
\begin{eqnarray*}
\nu_s(t)\big( \pi(h^{-1}) \cdot (h \cdot n) \big) & := &  t \cdot
\pi(h^{-1}) \cdot (h \cdot n) \cdot t^{-1} \\ &\;=& \pi(h^{-1})
\cdot (h \cdot \widetilde{\varpi}(h_t) (n)),
\end{eqnarray*}
i.e.\ $\nu_s(t)$ is a semisimple automorphism and $G_s = T
\ltimes_{\nu_s} N_s$.
\begin{Rem}
As usual, we denote the Lie algebras of the above Lie groups by the
corresponding small German letters. In \cite[Chapter III]{DER} can
be found:

There exists a vector subspace $V$ of $|\mathfrak{g}|$ with
$|\mathfrak{g}| = V \oplus |\mathfrak{n}|$ as vector spaces and
$\forall_{A,B \in V}~ \mathrm{ad}(A)_s (B) = 0$, where
$\mathrm{ad}(A)_s$ denotes the semisimple part of $\mathrm{ad}(A)$.

Let $\mathfrak{v}$ be a copy of $V$, considered as abelian Lie
algebra. Then the Lie algebra of the semisimple splitting for $G$ is
$\mathfrak{g}_s = \mathfrak{v} \ltimes_{ad(\ldots)_s} \mathfrak{g}$,
i.e.\ $$\forall_{(A,X),(B,Y) \in \mathfrak{g}_s}~ [(A,X),(B,Y)] =
\big( 0, [X,Y] + \mathrm{ad}(A)_s (Y) - \mathrm{ad}(B)_s(X) \big),$$
with nilradical $\mathfrak{n}_s = \{ (-X_V,X) \,|\, X \in
\mathfrak{g} \}$, where $X_V$ denotes the component of $X$ in $V$.
\end{Rem}

Now we state the announced criterion for the existence of lattices
in solvable Lie groups.

\begin{Thm}[{\cite[p.\ 248]{Aus}}] \label{GitterAuslander}
Let $G$ be a connected and simply-connected solvable Lie group with
nilradical $N$ and semisimple splitting $G_s = T \ltimes_{\nu_s}
N_s$, where $N_s$ is the nilshadow of $G$.

Then $G/N$ is contained as a subgroup in $G_s / N = T \times
(N_s/N)$ and the projections $\pi_1 \: G/N \to T$, $\pi_2 \: G/N \to
N_s/N$ are isomorphisms of abelian Lie groups.

Moreover, $G$ admits a lattice if and only if the following
conditions are satisfied:
\begin{itemize}
\item[(i)] There exists a basis $\mathfrak{X} := \{X_1,\ldots,X_n,\ldots,X_m\}$
with rational structure constants of the Lie algebra
$\mathfrak{n}_s$ of $N_s$ such that $\{X_1,\ldots,X_n\}$ is a basis
of the Lie algebra $\mathfrak{n}$ of $N$.

We write $\mathfrak{n}_s(\mathbb{Q})$ for the rational Lie algebra
$\langle \mathfrak{X} \rangle_{\mathbb{Q}}$ and $N_s(\mathbb{Q})$
for its image under the exponential map.
\item[(ii)] There exists a lattice subgroup $\Gamma_T$ of $T$ with
$\Gamma_T \subset \pi_1 \circ \pi_2^{-1} \big( N_s(\mathbb{Q}) /
N  \big)$ such that the natural action $\Gamma_T \to
\mathrm{A}(\mathfrak{n}_s(\mathbb{Q}))$ is described by integer
matrices in an appropriate basis of $\mathfrak{n}_s(\mathbb{Q})$. \q
\end{itemize}
\end{Thm}
\section{Three-dimensional solvmanifolds}
The only one- and two-dimensional solvmanifolds are tori. Therefore,
we begin our studies of low-dimensional solvmanifolds in dimension
three.

\begin{Prop}[\cite{AGH}]  \label{3 hat Gitter}
Every $3$-dimensional connected and simply-connected solvable
non-nilpotent Lie group $G$ that possesses a lattice $\Gamma$ has a
$2$-dimensional nilradical. The Lie group can be written as $G =
\mathbb{R} \ltimes_{\mu} \mathbb{R}^2$ and the lattice as $\Gamma =
\mathbb{Z} \ltimes_{\mu} \mathbb{Z}^2$. 
\end{Prop}

\textit{Proof.} This is a direct consequence of Proposition \ref{dim Nil} and
Corollary \ref{Gitter in fastabelsch}. \q

\begin{Thm}
A three-dimensional solvmanifold $G / \Gamma$ is non-formal if and
only if $b_1(G / \Gamma) = 2.$ In this case, $G / \Gamma$ is
diffeomorphic to a  nilmanifold.
\end{Thm}

\textit{Proof.} By Theorem \ref{formale Nilmgf}, it suffices to
consider the case when $G$ is solvable and non-nilpotent. The last
proposition implies that there is a map $\nu \: \mathbb{Z} \to
\mathrm{SL}(2,\mathbb{Z})$ such that $\Gamma = \mathbb{Z}
\ltimes_{\nu} \mathbb{Z}^2$.

If none of the roots of $\nu (1)$ equals $1$, Proposition \ref{b1
fast abelsch} implies $b_1 = 1$, so $G / \Gamma$ is formal by
Theorem \ref{nicht sympl}.

Assume that $\nu(1)$ possesses the double root $1$. Then Proposition \ref{b1
fast abelsch} implies $b_1 = 3$ if $\nu(1)$ is diagonalisable and
$b_1 = 2$ if $\nu(1)$ is not diagonalisable.

Case A: $\nu(1)$ is diagonalisable

\nopagebreak \noindent Recall that a solvmanifold is uniquely
determined by its fundamental group. Therefore, we can assume $G =
\mathbb{R} \ltimes_{\mu} \mathbb{R}^2$ and $\Gamma = \mathbb{Z}
\ltimes_{\mu} \langle v_1 , v_2 \rangle_{\mathbb{Z}}$ with linearly
independent $v_1 , v_2 \in \mathbb{R}^2$ and $\mu(t) \equiv
\mathrm{id}$. In this case, $G / \Gamma$ is a torus which is formal.

Case B: $\nu(1)$ is not diagonalisable

\nopagebreak \noindent In this case, we can assume $G = \mathbb{R}
\ltimes_{\mu} \mathbb{R}^2$ as well as $\Gamma = \mathbb{Z}
\ltimes_{\mu} \langle v_1 , v_2 \rangle_{\mathbb{Z}}$ with linearly
independent $v_1 , v_2 \in \mathbb{R}^2$ and
$$\mu(t) = \left(
\begin{array}{cc}
1 & t \\ 0 & 1
\end{array}
\right) .$$ The Lie algebra $\mathfrak{g} = \langle T,X,Y \,|\,
[T,Y] = X \rangle$ of $G$ is nilpotent, so $G / \Gamma$ is a
nilmanifold with $b_1 = 2$. Therefore, it cannot be a torus and is
not formal by Theorem \ref{formale Nilmgf}. \q
\N In \cite[Chapter III §3]{AGH} the three-dimensional solvmanifolds
which have no nilmanifold structure are examined. This, together
with the last theorem, yields a ``cohomological'' classification of
three-dimensional solvmanifolds.

\begin{Thm}
Every $3$-dimensional solvmanifold $G / \Gamma$ is contained in
Table \ref{Solv3}.
\begin{table}[htb] \centering
\caption{3-dimensional solvmanifolds} \label{Solv3}
\begin{tabular}{|c|c|c|c|c|} \hline \hline
& $b_1(G / \Gamma)$ & $G / \Gamma$ formal & Nilmfd. \footnotemark[3] & c.s. \footnotemark[4] %
\\ \hline
\hline a) & $3$ & yes & Torus & yes \\
\hline b) & $2$ & no & yes & yes \\
\hline c) & $1$ & yes & no & yes \\
\hline d) & $1$ & yes & no & no \\
\hline
\end{tabular}
\end{table}
In particular, $G / \Gamma$ is non-formal if and only if it is a
non-toral nilmanifold.\q
\end{Thm}
\footnotetext[3]{possesses the structure of a solvmanifold as
quotient of a nilpotent Lie group}%
\footnotetext[4]{possesses the structure of a solvmanifold as
quotient of a \textbf{c}ompletely \textbf{s}olvable Lie group}%
\addtocounter{footnote}{2}%
\begin{Ex} \label{Ex3}
The torus $\mathbb{R}^3 / \mathbb{Z}^3$ is a solvmanifold with $b_1
= 3$, and examples of $3$-dimensional solvmanifolds with $b_1 = 2$
will be given in the next theorem.

For $i \in \{1,2\}$ consider the Lie groups $G_i = \mathbb{R}
\ltimes_{\mu_i} \mathbb{R}^2$, where $\mu_i$ is given by
$\mu_1(t)(x,y) = (e^{t} \, x , e^{-t} \, y)$, $\mu_2(t)(x,y) =
(\cos(t) \, x + \sin(t) \,y , -\sin(t) \, x + \cos(t) \, y)$.

$G_1$ is completely solvable and possesses the lattice
\begin{equation*}
\Gamma_1 := t_1 \, \mathbb{Z} \ltimes_{\mu_1} \langle \left(
\begin{array}{c} 1 \\ 1 \end{array} \right) , \left( \begin{array}{c}
\frac{18 + 8 \sqrt{5}}{7 + 3 \sqrt{5}} \\ \frac{2}{3 + \sqrt{5}}
\end{array} \right) \rangle_{\mathbb{Z}},
\end{equation*}
where $t_1 = \ln (\frac{3 + \sqrt{5}}{2})$. Note that the following
equation implies that $\Gamma_1$ really is a lattice
\begin{equation} \label{expt1}
\left( \begin{array}{cc} e^{t_1} & 0 \\
0 & e^{-t_1}
\end{array} \right) = \left( \begin{array}{cc} 1 & \frac{18 + 8
\sqrt{5}}{7 + 3 \sqrt{5}} \\ 1 & \frac{2}{3 + \sqrt{5}} \end{array}
\right) \left(
\begin{array}{cc} 0 & -1 \\ 1 & 3 \end{array} \right) \left(
\begin{array}{cc} 1 & \frac{18 + 8 \sqrt{5}}{7 + 3
\sqrt{5}} \\ 1 & \frac{2}{3 + \sqrt{5}}
\end{array} \right)^{-1}.
\end{equation}
$G_2$ is not completely solvable and contains the lattice
$$ \Gamma_2 = \pi \, \mathbb{Z} \ltimes_{\mu_2} \mathbb{Z}^2.$$

A short computation yields that the abelianisations of $\Gamma_i$
have rank one, i.e.\ we have constructed examples of type c) and d)
in Table \ref{Solv3}.
\end{Ex}

\begin{Thm} \label{Nil3}
Every lattice in the unique $3$-dimensional connected and
simply-connected non-abelian nilpotent Lie group
$$U_3(\mathbb{R}) := \{ \left(
\begin{array}{ccc}
1 & x & z \\ 0 & 1 & y \\ 0 & 0 & 1
\end{array}
\right) \,|\, x,y,z \in\mathbb{R} \}$$ is isomorphic to
$\Gamma_{3,n} := \Gamma_{3,n}(\mathbb{Z}) := \{ \left(
\begin{array}{ccc}
1 & x & \frac{z}{n} \\ 0 & 1 & y \\ 0 & 0 & 1
\end{array}
\right) \,|\, x,y,z \in\mathbb{Z} \}$ with $n \in \mathbb{N}_+$.

Therefore, any three-dimensional nilmanifold with $b_1 = 2$ is of
the form $U_3(\mathbb{R}) / \Gamma_{3,n}(\mathbb{Z})$.

$\Gamma_{3,n}(\mathbb{Z})$ is presented by $\langle e_1, e_2, e_3
\,|\, [e_1,e_2] =  e_3^n \mbox{ and }  e_3 \mbox{ central } \rangle
.$
\end{Thm}

\textit{Proof.} $U_3(\mathbb{R})$ is the only connected and
simply-connected non-abelian nilpotent Lie group of dimension three.
By \cite[Chapter III § 7]{AGH}, each lattice in it is isomorphic to
$\Gamma_{3,n}$. The other assertions follow trivially. \q
\N Sometimes, we shall write $(x,y,z)$ for the corresponding matrix
in $U_3(\mathbb{R})$.

For later applications, we are going to determine the Lie group
automorphisms and the one-parameter groups of $U_3(\mathbb{R})$. In
order to do this, we start with the following proposition. Note that
$Z(G)$ denotes the center of a group $G$.
\begin{Prop} \label{U3,Gamma3n} $\,$
\begin{itemize}
\item[(i)] $[U_3(\mathbb{R}),U_3(\mathbb{R})] = Z(U_3(\mathbb{R}))
= \{ (0,0,z) \,|\, z \in \mathbb{R} \}$, $U_3(\mathbb{R}) /
Z(U_3(\mathbb{R})) \cong \mathbb{R}^2$
\item[(ii)] Every Lie group homomorphism $f \: U_3(\mathbb{R}) \to
U_3(\mathbb{R})$ induces natural Lie group homomorphisms $$f_Z \:
Z(U_3(\mathbb{R})) \longrightarrow Z(U_3(\mathbb{R}))$$ and
\begin{eqnarray*}
\overline{f} \: ~~ U_3(\mathbb{R}) / Z(U_3(\mathbb{R})) ~~&
\longrightarrow & ~~~~~~~~ U_3(\mathbb{R}) / Z(U_3(\mathbb{R})). \\
{[}(x,y,0){]} = {[} (x,y,z) {]} & \longmapsto & {[} f \big( (x,y,z)
\big) {]} = {[}(f_1(x,y,0),f_2(x,y,0),0){]}
\end{eqnarray*}
$\overline{f}$ uniquely determines $f_Z$, and $\overline{f}$ is an
automorphism if and only if $f$ is such.
\item[(iii)] Let $\gamma_1 = (a_1,b_1,\frac{c_1}{n}), \gamma_2 =
(a_2,b_2,\frac{c_2}{n}) \in \Gamma_{3,n}$. Then there is a unique
homomorphism $g \: \Gamma_{3,n} \to \Gamma_{3,n}$ such that
$g\big((1,0,0)\big) = \gamma_1$ and $g\big((0,1,0)\big) = \gamma_2$.
Moreover, $g\big((0,0,\frac{1}{n})\big) = \big(0,0,\frac{1}{n}(a_1
b_2 - a_2 b_1)\big)$.

One has $\Gamma_{3,n} / Z(\Gamma_{3,n}) \cong \mathbb{Z}^2$, and $g$
is an isomorphism if and only if $$\overline{g} \: \Gamma_{3,n} /
Z(\Gamma_{3,n}) \longrightarrow \Gamma_{3,n} / Z(\Gamma_{3,n})$$ is
an isomorphism, i.e.\ $a_1 b_2 - a_2 b_1 = \pm 1$.
\end{itemize}
\end{Prop}

\textit{Proof.} (i) is trivial.

ad (ii): Let $f \: U_3(\mathbb{R}) \to U_3(\mathbb{R})$ be a Lie
group homomorphism. Then
\begin{equation} \label{f(0,0,z)}
f \big( (0,0,z) \big) = 
[ f \big( (z,0,0) \big), f \big( (0,1,0) \big) ] \in
Z(U_3(\mathbb{R})),
\end{equation}
i.e.\ $f \big( Z(U_3(\mathbb{R})) \big) \subset
Z(U_3(\mathbb{R})).$ Moreover, one has for 
$(a,b,c) := f \big( (x,y,0) \big)$
$$(a,b,0)^{-1} \cdot (a,b,c) = (-a,-b,0) \cdot (a,b,c) = (0,0,-ab +
c ) \in Z(U_3(\mathbb{R})),$$ and therefore $[(a,b,0)] =
\overline{f} ([(x,y,0)])$. Now, (\ref{f(0,0,z)}) implies that $f_Z$
is uniquely determined by $\overline{f}$.

Assume, $f$ is an isomorphism. Then (\ref{f(0,0,z)}) also holds for
$f^{-1}$ and we get $f \big( Z(U_3(\mathbb{R})) \big) =
Z(U_3(\mathbb{R})),$ i.e.\ $f_Z$ is an isomorphism of the additive
group $\mathbb{R}$. Since $f$ is continuous, there exists $m \in
\mathbb{R} \setminus \{0\}$ such that $f_Z \big( (0,0,z) \big) =
(0,0,mz)$. Denote by $(f_{ij})_{1 \le i,j \le 2}$ the matrix of
$\overline{f} \: \mathbb{R}^2 \to \mathbb{R}^2$ with respect to the
basis $\{ \left(
\begin{array}{c}
1 \\ 0
\end{array}
\right) , \left(
\begin{array}{c}
0 \\ 1
\end{array}
\right) \}$ of the vector space $\mathbb{R}^2$. One calculates
\begin{eqnarray*}
(0,0, \det (f_{ij}) ) & = & [(f_{11},f_{21},0), (f_{12},f_{22},0)] =
[f \big( (1,0,0) \big) , f \big( (0,1,0) \big) ] \\ &
\stackrel{(\ref{f(0,0,z)})}{=} & (0,0,m),
\end{eqnarray*}
so $\overline{f}$ is an automorphism, since $m \ne 0$.

Conversely, if $\overline{f}$ is an automorphism, then the
homomorphism $f_Z$ is given by $f_Z \big( (0,0,z) \big) = (0,0,
\det(\overline{f}) z )$ which is even an automorphism. Therefore,
the $5$-Lemma implies that $f$ is an automorphism.

ad (iii): Let $\gamma_1 , \gamma_2$ be as in (iii). Then
$[\gamma_1,\gamma_2] = 
\big(0,0,\frac{1}{n}(a_1 b_2 - a_2 b_1)\big)^n$ and this implies the
existence of the (unique) homomorphism $g$ with the mentioned
properties.

If $g$ is an isomorphism, then $g(Z(\Gamma_{3,n})) = Z(\Gamma_{3,n})
= \{(0,0,\frac{z}{n}) \,|\, z \in \mathbb{Z}\},$ and therefore $|a_1
b_2 - a_2 b_1| = 1$. Since the matrix of $\overline{g}$ has
determinant $a_1 b_2 - a_2 b_1$, $\overline{f}$ is an isomorphism.

Again, the converse is trivial. \q

\begin{Thm}
As a set, the group of Lie group automorphisms
$\mathrm{A}(U_3(\mathbb{R}))$ equals $\mathrm{GL}(2,\mathbb{R})
\times \mathbb{R}^2$, the group law is given by
\begin{eqnarray} \label{Formel1 Aut}
(A,a) \circ (B,b) &\longmapsto& \big(AB, \det(B) B^{-1} a + \det(A)
b) \big),
\end{eqnarray}
and for $f = ( A = \left(
\begin{array}{cc}
\alpha & \beta \\ \gamma & \delta
\end{array}
\right) , \left(
\begin{array}{cc}
u  \\ v
\end{array}
\right) ) \in \mathrm{A}(U_3(\mathbb{R}))$ and $(x,y,z) \in
U_3(\mathbb{R})$ we have

\parbox{12cm}{$$
\begin{array}{ccl} f\big((x,y,z)\big) &=& \big(\alpha x + \beta y, \gamma x + \delta
y, \\ && ~\det(A) z + \beta \gamma xy + \frac{\alpha \gamma}{2} x^2
+ \frac{\beta \delta}{2} y^2 + uy - vx \big).
\end{array}$$} \hfill \parbox{8mm}{\begin{eqnarray} \label{Formel
Aut} \end{eqnarray}}
\end{Thm}

\textit{Proof.} Let $f \in \mathrm{A}(U_3(\mathbb{R}))$ and $(x,y,z)
\in U_3(\mathbb{R})$ be given. We have to show that there is
$(\left(
\begin{array}{cc}
\alpha & \beta \\ \gamma & \delta
\end{array}
\right) , \left(
\begin{array}{cc}
u  \\ v
\end{array}
\right) ) \in \mathrm{GL}(2,\mathbb{R}) \times \mathbb{R}^2$ such
that $f\big( (x,y,z) \big)$ satisfies (\ref{Formel Aut}). Then a
short computation yields (\ref{Formel1 Aut}).

Let $\left(
\begin{array}{cc}
\alpha & \beta \\ \gamma & \delta
\end{array}
\right) \in \mathrm{GL}(2,\mathbb{R})$ be the matrix of
$\overline{f}$ with respect to the canonical basis of
$\mathbb{R}^2$. We showed in the last proof $f\big((0,0,z)\big) =
(0,0\det(\overline{f})z)$.

There exist smooth functions $f_1,f_2 \: \mathbb{R} \to \mathbb{R}$
with
\begin{eqnarray*}
f\big((x,0,0)\big) & = & (\alpha x , \gamma x, f_1(x)), \\
f\big((0,y,0)\big) & = & (\beta y , \delta y, f_2(y)).
\end{eqnarray*}
We set $u := f_2^{'}(0)$ and $v := - f_1^{'}(0)$. The homomorphism
property of $f$ implies
\begin{eqnarray*}
\frac{1}{h}(f_1(x + h) - f_1(x)) & = & \frac{f_1(x) - f_1(0)}{h} + \alpha \gamma x ,\\
\frac{1}{h}(f_2(y + h) - f_2(y)) & = & \frac{f_2(y) - f_2(0)}{h} +
\beta \delta y,
\end{eqnarray*}
and this yields
\begin{eqnarray*}
f_1(x) & = & -v x + \frac{\alpha \gamma}{2} x^2, \\
f_2(y) & = & u y + \frac{\beta \delta}{2} y^2.
\end{eqnarray*}
Using $(x,y,z) = (0,y,0) (x,0,0) (0,0,z)$, one computes (\ref{Formel
Aut}). \q

\begin{Cor}  \label{Nil3 1Param}
$f = ( A , \left(
\begin{array}{cc}
u  \\ v
\end{array}
\right) ) \in \mathrm{A}(U_3(\mathbb{R}))$ with $A = \left(
\begin{array}{cc}
\alpha & \beta \\ \gamma & \delta
\end{array}
\right)$ lies on a one-parameter group of
$\mathrm{A}(U_3(\mathbb{R}))$ if and only if $A$ lies one a
one-parameter group of $GL(2,\mathbb{R})$.

If $\nu_t = \left(
\begin{array}{cc}
\alpha_t & \beta_t \\ \gamma_t & \delta_t
\end{array}
\right)$ denotes a one-parameter group with $\nu_1 = A$, then the
map $\mu_t \: \mathbb{R} \to \mathrm{A}(U_3(\mathbb{R}))$ defined by
$$\begin{array}{ccl} \mu_t \big((x,y,z)\big) &=& \big(\alpha_t x + \beta_t y, \gamma_t x +
\delta_t y, \\ && ~\underbrace{(\alpha_t \delta_t - \beta_t
\delta_t)}_{\D =1} z + \beta_t \gamma_t xy + \frac{\alpha_t
\gamma_t}{2} x^2 + \frac{\beta_t \delta_t}{2} y^2 + tuy - tvx \big)
\end{array}$$
is a one-parameter group with $\mu_1 = f$.
\end{Cor}

\textit{Proof.} The only claim that is not obvious is the fact that
$\mu_t$ defines a one-parameter group. Using $\nu_{t+s} = \nu_t
\circ \nu_s$, this can be seen by a short calculation. \q

\section{Four-dimensional solvmanifolds}\label{KapSolvmgfen in Dim4}
As we have done in the three-dimensional case, we are going to give
a ``cohomological'' classification of four-dimensional
solvmanifolds. We shall consider all isomorphism classes of lattices
that can arise in a four-dimensional connected and simply-connected
solvable Lie group. The next proposition describes such lattices in
the case of a non-nilpotent group.

\begin{Prop} \label{nu forsetzbar}
Every $4$-dimensional connected and simply-connected solvable
non-nilpotent Lie group $G$ that possesses a lattice $\Gamma$ has a
$3$-dimensional nilradical $N$ which is either $\mathbb{R}^3$ or
$U_3(\mathbb{R})$. Therefore, $G / \Gamma$ fibers over $S^1$ (this
is the Mostow bundle) and the Lie group can be written as $G =
\mathbb{R} \ltimes_{\mu} N$. If $N$ is abelian, a basis
transformation yields $\Gamma = \mathbb{Z}
\ltimes_{\mu|_{\mathbb{Z}^3}} \mathbb{Z}^3$. Otherwise, $\Gamma$ is
isomorphic to $\mathbb{Z} \ltimes_{\nu} \Gamma_{3,n}$, where $\nu \:
\mathbb{Z} \to \mathrm{Aut}(\Gamma_{3,n})$ is a group homomorphism
with
$$\begin{array}{ccl} \nu(1)(x,y,\frac{z}{n}) ~= & \big( a_1 x + a_2
y,~~ b_1 x + b_2 y, & a_2 b_1 x y + a_1 b_1 \frac{x(x-1)}{2} + a_2
b_2 \frac{y(y-1)}{2} \\ && + \frac{1}{n}(c_1 x + c_2 y  + (a_1 b_2 -
a_2 b_1)z) \big),
\end{array} $$
where $c_1,c_2, \in \mathbb{Z}$, and $\left(
\begin{array}{cc}
a_1 & a_2 \\ b_1 & b_2
\end{array}
\right) \in \mathrm{GL}(2,\mathbb{Z})$ is the matrix of
$\overline{\nu(1)}$ with respect to the canonical basis of the
$\mathbb{Z}$-module $\mathbb{Z}^2 = \Gamma_{3,n} / Z(\Gamma_{3,n})$.
Moreover, $\overline{\nu}(1)$ lies on a one-parameter group
$\mathbb{R} \to \mathrm{A}(U_3(\mathbb{R})/Z(U_3(\mathbb{R}))) =
\mathrm{GL}(2,\mathbb{R})$, i.e.\ $\overline{\nu}(1) \in
\mathrm{SL}(2,\mathbb{R})$.
\end{Prop}

\textit{Proof.} From \cite[Theorem 3.1.10]{TO} follows $\dim N = 3$
and $G = \mathbb{R} \ltimes_{\mu} N$. If $N$ is abelian, Corollary
\ref{Gitter in fastabelsch} implies that we can assume $\Gamma =
\mathbb{Z} \ltimes_{\mu|_{ \mathbb{Z}^3}} \mathbb{Z}^3$.

Assume now that $N$ is not abelian, i.e.\ $N = U_3(\mathbb{R})$.
$\Gamma_N = \Gamma \cap N$ is a lattice in $N$ and by Theorem
\ref{Nil3}, we have $\Gamma_N = \Gamma_{3,n}$. By Corollary
\ref{Gitter in fastabelsch}, there is a homomorphism $\nu \:
\mathbb{Z} \to \mathrm{Aut}(\Gamma_{3,n})$ with $\Gamma \cong
\mathbb{Z} \ltimes_{\nu} \Gamma_{3,n}$. Proposition
\ref{U3,Gamma3n}(iii) implies that $\nu(1)$ is determined by
$(a_1,b_1, \frac{c_1}{n}) := \nu(1) \big( (1,0,0) \big)$ and
$(a_2,b_2, \frac{c_2}{n}) := \nu(1) \big( (0,1,0) \big) \in
\Gamma_{3,n}$. Since $(x,y,\frac{z}{n}) = (0,1,0)^y (1,0,0)^x
(0,0,\frac{1}{n})^z$, a short computation yields the claimed formula
for $\nu(1) \big( (x,y,\frac{z}{n}) \big)$.

Further, Corollary \ref{Nil3 1Param} implies that
$\overline{\nu}(1)$ lies on a one-parameter group. \q
\begin{Thm}  \label{Solvmgf Dim4}
Every $4$-dimensional solvmanifold $G / \Gamma$ is contained in
Table \ref{Solv4}. 
\begin{table}[h]
\centering \caption{\label{Solv4} $4$-dimensional solvmanifolds}
\begin{tabular}{|c|c|c|c|c|c|c|c|} \hline \hline
& $b_1(G/\Gamma)$ & $G / \Gamma$ formal & symplectic & complex &
K\"ahler & Nilmfd. \footnotemark[5] & c.s. \footnotemark[6] %
\\ \hline
\hline a) & $4$ & yes & yes & Torus & yes & Torus & yes \\
\hline b) & $3$ & no & yes & PKS \footnotemark[7] & no & yes & yes \\
\hline c) & $2$ & yes & yes & no & no & no & yes \\
\hline d) & $2$ & yes & yes & HS \footnotemark[8] & yes & no & no \\
\hline e) & $2$ & no & yes & no & no & yes & yes \\
\hline f) &$1$ & yes & no & no & no & no & yes \\
\hline g) & $1$ & yes & no & I$S^0$ \footnotemark[9] & no & no & no \\
\hline h) & $1$ & yes & no & I$S^+$ \footnotemark[10] & no & no & yes \\
\hline i) &$1$ & yes & no & SKS \footnotemark[11] & no & no & no \\
\hline
\end{tabular}
\end{table}
\footnotetext[5]{possesses the structure of a solvmanifold as
quotient of a nilpotent Lie group}%
\footnotetext[6]{possesses the structure of a solvmanifold as
quotient of a \textbf{c}ompletely \textbf{s}olvable Lie group}%
\footnotetext[7]{\textbf{P}rimary\textbf{K}odaira \textbf{S}urface}%
\footnotetext[8]{\textbf{H}yperelliptic \textbf{S}urface}%
\footnotetext[9]{\textbf{I}noue Surface of Type $\mathbf{S^0}$}%
\footnotetext[10]{\textbf{I}noue Surface of Type $\mathbf{S^+}$}%
\footnotetext[11]{\textbf{S}econdary \textbf{K}odaira
\textbf{S}urface}%
\addtocounter{footnote}{7}%
In particular, $G / \Gamma$ is non-formal if and only if it is a
non-toral nilmanifold.
\end{Thm}

\textit{Proof.} Apart from the column on formality the theorem
follows from works of Geiges \cite{Ge4} and Hasegawa \cite{Has4}.
(Attention: In \cite{Has4} a more general notion of solvmanifold is
used!)

A decomposable four-dimensional connected and simply-connected
nilpotent Lie group is abelian or has a two-dimensional center. The
only connected and simply-connected indecomposable nilpotent Lie
group of dimension four has a two-dimensional commutator. By
Propositions \ref{Zentrum} and \ref{lokal kompakt}, the
corresponding nilmanifolds have the structure of orientable
$T^2$-bundles over $T^2$. (The orientability follows from the total
spaces' orientability.)

From a result of Geiges \cite[Theorems 1 and 3]{Ge4} follows that
they are contained in Table \ref{Solv4}. (Recall that a nilmanifold
is formal if and only if it is a torus.) In particular, every
four-dimensional nilmanifold is symplectic.

Now, we regard a lattice $\Gamma = \mathbb{Z} \ltimes_{\nu}
\Gamma_N$, $\Gamma_N \in \{ \mathbb{Z}^3 , \Gamma_{3,n}(\mathbb{Z})
\}$, in a Lie group $G = \mathbb{R} \ltimes_{\mu} N$ as in the last
proposition.

We expand Hasegawa's argumentation in \cite{Has4} by the aspect of
formality and consider the ``roots'' of $\nu(1)$. Recall, Corollary
\ref{Wang eind} implies that a solvmanifold is determined by its
fundamental group. Below, we shall use this fact several times.
\\[2.5mm] \textbf{Case A.:} $\Gamma_N = \mathbb{Z}^3$

\nopagebreak \noindent By Proposition \ref{nu forsetzbar}, $\nu$
extends to a one-parameter group $\mathbb{R} \to
\mathrm{SL}(3,\mathbb{R})$. Denote by $\lambda_1 , \lambda_2,
\lambda_3 \in \mathbb{C}$ the roots of $\nu(1) \in
\mathrm{SL}(3,\mathbb{Z})$, i.e.\ $\lambda_1 \cdot \lambda_2 \cdot
\lambda_3 = 1$. We get from Theorem \ref{Bild exp} and Lemma
\ref{Z3doppelt} that the following subcases can occur:
\begin{description}
\item{A.1.)} $\lambda_1, \lambda_2, \lambda_3 \in \mathbb{R}_+$
\begin{description}
\item{A.1.1.)} $\exists_{i_0} ~ \lambda_{i_0} = 1$ (w.l.o.g. $\lambda_1 =
1$)
\begin{description}
\item{A.1.1.1.)} $\lambda_2 = \lambda_3 =1 $
\item{A.1.1.2.)} $\lambda_2 = \lambda_3^{-1} \in \mathbb{R} \setminus \{ 1 \}$
\end{description}
\item{A.1.2.)} $\forall_i~ \lambda_i \ne 1$
\begin{description}
\item{A.1.2.1.)} $\nu(1)$ is diagonalisable
\item{A.1.2.2.)} $\nu(1)$ is not diagonalisable
\end{description}
\end{description}
\item{A.2.)} $\lambda_1= 1$, $\lambda_2 = \lambda_3 = -1$ and
$\nu(1)$ is diagonalisable
\item{A.3.)} $\exists_{i_0}~ \lambda_{i_0} \in \mathbb{C} \setminus
\mathbb{R}$ (w.l.o.g. $\lambda_2 = \overline{\lambda_3} \in
\mathbb{C} \setminus \mathbb{R}$ and $\lambda_1 \in \mathbb{R}_+$)
\begin{description}
\item{A.3.1.)} $\lambda_1 = 1$
\item{A.3.2.)} $\lambda_1 \ne 1$
\end{description}
\end{description}

\textbf{Case A.1.1.1.:} $\lambda_1 = \lambda_2 = \lambda_3 = 1$

\nopagebreak \noindent If $\nu(1)$ is diagonalisable, then
$G/\Gamma$ is clearly a torus. This is case a). If $\nu(1)$ is not
diagonalisable we can assume $G = \mathbb{R} \ltimes_{\mu}
\mathbb{R}^3$ and $\Gamma = \mathbb{Z} \ltimes_{\mu} \langle v_1 ,
v_2 , v_3 \rangle_{\mathbb{Z}}$, where $\mu(t)$ is one of the
following one-parameter groups
\begin{eqnarray*}
\exp \Big(t \left(
\begin{array}{ccc} 0& 0 & 0 \\ 0 & 0 & 1 \\ 0 & 0 & 0
\end{array} \right) \Big) & = & \left(
\begin{array}{ccc} 1& 0 & 0 \\ 0 & 1 & t \\ 0 & 0 & 1
\end{array} \right), \\ \exp \Big(t \left(
\begin{array}{ccc} 0 & 1 & - \frac{1}{2} \\ 0 & 0 & 1 \\ 0 & 0 & 0
\end{array} \right) \Big) & = & \left(
\begin{array}{ccc} 1& t & \frac{1}{2} (t^2 -t) \\ 0 & 1 & t \\ 0 & 0 & 1
\end{array} \right).
\end{eqnarray*}
In both cases $G/\Gamma$ is a $4$-dimensional nilmanifold and
therefore symplectic. In the first case, we have a primary Kodaira
surface with $b_1 =3$, see \cite[Section 2.2.3)]{Has4}; in the
second case the nilmanifold has $b_1 = 2$ and no complex structure,
see \cite[Example 2]{HasCK}. Being non-toral nilmanifolds, both are
not formal and we get the cases b) and e).

\textbf{Cases A.1.1.2. and A.1.2.1.:} The $\lambda_i$ are positive
and pairwise different or two of them are equal but $\nu(1)$ is
diagonalisable. (The latter cannot happen by Lemma \ref{Z3doppelt}.)

\nopagebreak \noindent We can assume $G = \mathbb{R} \ltimes_{\mu}
\mathbb{R}^3$ and $\Gamma = \mathbb{Z} \ltimes_{\mu} \langle v_1 ,
v_2 , v_3 \rangle_{\mathbb{Z}}$ with linearly independent $v_1 , v_2
, v_3 \in \mathbb{R}^3,$ where $\mu(t) = \left(
\begin{array}{ccc} \exp(t \ln(\lambda_1)) & 0 & 0 \\ 0 & \exp(t
\ln(\lambda_2)) & 0 \\ 0 & 0 & \exp(t \ln(\lambda_3))
\end{array} \right) .$ By \cite[Example 2]{HasCK}, the solvmanifold
$G / \Gamma$ does not admit a complex structure.

One computes the Lie algebra of $G$ as $$\mathfrak{g} = \langle \,
T,X,Y,Z \,|\, [T,X] = \ln(\lambda_1) X,\, [T,Y] = \ln(\lambda_2)
Y,\, [T,Z] = \ln(\lambda_3) Z \, \rangle $$ which is completely
solvable and non-nilpotent. Therefore, the minimal model of the
Chevalley-Eilenberg complex is the minimal model of $G / \Gamma$.

If none of the roots $\lambda_i$ is one, we see by Proposition
\ref{b1 fast abelsch} that $b_1(G/\Gamma) = 1$. Since $G / \Gamma$
is parallelisable, this implies $b_2(G / \Gamma) = 0$, so the space
cannot be symplectic. Further it is formal by Theorem \ref{nicht
sympl}. This is case f) in Table \ref{Solv4}.

If one of the roots is one (w.l.o.g. $\lambda_1 = 1$), we have
$b_1(G/\Gamma) = 2$ and the Chevalley-Eilenberg complex is $$\big(
\bigwedge ( \tau , \alpha , \beta , \gamma ) \, ,\, d \tau = d
\alpha = 0 ,~ d \beta = -\ln(\lambda_2) \, \tau \wedge \beta ,~ d
\gamma = - \ln(\lambda_3) \, \tau \wedge \gamma \big).$$ $\tau
\wedge \alpha + \alpha \wedge \beta + \alpha \wedge \gamma - \beta
\wedge \gamma$ defines a symplectic form on $G / \Gamma$. Further,
one computes the first stage of the minimal model of the
Chevalley-Eilenberg complex as
$$\mathcal{M}^{\le 1} = \bigwedge (x_1 , x_2),~ d x_i = 0 .$$
Therefore, $G/\Gamma$ is $1$-formal and by Theorem \ref{formal = n-1
formal} formal. This is case c) in Table \ref{Solv4}.

\textbf{Case A.1.2.2.:} $\lambda_i \in \mathbb{R}_+ \setminus \{1\}$
and $\nu(1)$ is not diagonalisable

\nopagebreak \noindent In this case two roots must be equal
(w.l.o.g. $\lambda_2 = \lambda_3$) and the third is different from
the others, i.e.\ $\lambda_1 = \frac{1}{\lambda_2^2} \ne \lambda_2$.
Since the characteristic polynomial of $\nu(1)$ has integer
coefficients, Lemma \ref{Z3doppelt} implies $\lambda_2 = \pm 1$ and
this is a contradiction.

\textbf{Cases A.2. and A.3.1.:} $\lambda_1 = 1, \lambda_2 =
\overline{\lambda_3} = \exp(i \varphi) \in \mathbb{C} \setminus
\mathbb{R}, \varphi \in ]0,2\pi[$

\nopagebreak \noindent We can assume $G = \mathbb{R} \ltimes_{\mu}
\mathbb{R}^3$ and $\Gamma = \mathbb{Z} \ltimes_{\mu} \langle v_1 ,
v_2 , v_3 \rangle_{\mathbb{Z}}$ with linearly independent $v_1 , v_2
, v_3 \in \mathbb{R}^3,$ where $\mu(t) = \left(
\begin{array}{ccc} 1 & 0 & 0 \\ 0 & \cos(t
\varphi) & -\sin(t \varphi) \\ 0 & \sin(t \varphi) & \cos(t \varphi)
\end{array} \right) .$ Thus $G / \Gamma$ is a hyperelliptic surface
(see \cite[Section 3.3.]{Has4}) which is K\"ahlerian and has $b_1 =
2$. The Lie algebra of $G$ is not completely solvable and we are in
case d).

\textbf{Case A.3.2.:} $\lambda_1 \ne 1, \lambda_2 =
\overline{\lambda_3} = |\lambda_2| \exp(i \varphi) \in \mathbb{C}
\setminus \mathbb{R}, \varphi \in ]0,2\pi[ \setminus \{\pi\}$

\nopagebreak \noindent We can assume $G = \mathbb{R} \ltimes_{\mu}
\mathbb{R}^3$ and $\Gamma = \mathbb{Z} \ltimes_{\mu} \langle v_1 ,
v_2 , v_3 \rangle_{\mathbb{Z}}$ with linearly independent $v_1 , v_2
, v_3 \in \mathbb{R}^3$, where $\mu(t) = \left(
\begin{array}{ccc} \lambda_1^t & 0 & 0 \\ 0 & |\lambda_2|^t \cos(t
\varphi) & -|\lambda_2|^t \sin(t \varphi) \\ 0 & |\lambda_2|^t
\sin(t \varphi) & |\lambda_2|^t \cos(t \varphi)
\end{array} \right) .$ Thus $G / \Gamma$ is a Inoue surface
of type $S^0$ (see \cite[Section 3.6.]{Has4}), which is not
symplectic and has $b_1 = 1$ (by Proposition \ref{b1 fast abelsch},
since $1$ is no root of $\mu(1)$). By Theorem \ref{nicht sympl},
$G/\Gamma$ is formal. The Lie algebra of $G$ is not completely
solvable and this yields case g) of Table \ref{Solv4}.
\\[2.5mm] \textbf{Case B.:} $\Gamma_N = \Gamma_{3,n}(\mathbb{Z})$

\nopagebreak \noindent In this case we have a homomorphism $\nu \:
\mathbb{Z} \to \mathrm{Aut}(\Gamma_{3,n}(\mathbb{Z})).$ We shall
write $N$ for $U_3(\mathbb{R})$ and $\Gamma_N$ for
$\Gamma_{3,n}(\mathbb{Z}).$ The automorphism $\nu(1)$ induces an
automorphism $\overline{\nu}(1)$ of $\Gamma_N / Z(\Gamma_N) =
\mathbb{Z}^2$ which lies by Proposition \ref{nu forsetzbar} on a
one-parameter group $\mathbb{R} \to \mathrm{A}(U_3(\mathbb{R}) /
\Gamma_N) = \mathrm{GL}(2,\mathbb{R})$. Denote the roots of
$\overline{\nu}(1) \in \mathrm{GL}(2,\mathbb{Z})$ by
$\widetilde{\lambda}_1, \widetilde{\lambda}_2.$ Theorem \ref{Bild
exp} shows that the following cases are possible:
\begin{description}
\item{B.1.)} $\widetilde{\lambda}_1, \widetilde{\lambda}_2 \in \mathbb{R}_+$
\begin{description}
\item{B.1.1.)} $\widetilde{\lambda}_1 = \widetilde{\lambda}_2 = 1$
\item{B.1.2.)} $\widetilde{\lambda}_1  = \widetilde{\lambda}_2^{-1} \ne 1$
\end{description}
\item{B.2.)} $\widetilde{\lambda}_1 = \widetilde{\lambda}_2 = -1$ and
$\overline{\nu}(1)$ is diagonalisable
\item{B.3.)} $\widetilde{\lambda}_1 = \overline{\widetilde{\lambda}_2} \in \mathbb{C}
\setminus \mathbb{R}$
\end{description}
$\nu(1)$ also induces an automorphism $\nu(1)_Z$ of the center
$Z(\Gamma_N)$ of $\Gamma_N$ which equals $\det \big(
\overline{\nu}(1) \big) \cdot \mathrm{id} = \mathrm{id}$ by
Proposition \ref{U3,Gamma3n}(iii).

\textbf{Case B.1.1.:} $\widetilde{\lambda}_1 = \widetilde{\lambda}_2
= 1$

\nopagebreak \noindent By \cite[Lemma 1]{SF}, we can assume
$\overline{\nu}(1) = \left(
\begin{array}{cc}
1 &  k \\ 0 & 1
\end{array}
\right) \in \mathrm{SL}(2, \mathbb{Z})$ with $k \in \mathbb{N}$.
Then Proposition \ref{nu forsetzbar} yields $\nu(1) \big(
(x,y,\frac{z}{n}) \big) = (x + ky ,y, k \frac{y(y-1)}{2} + \frac{c_1
x + c_2 y +z}{n})$ and this implies $$\Gamma = \mathbb{Z}
\ltimes_{\nu} \Gamma_{3,n} = \langle e_0, \ldots , e_3 \,|\,
[e_0,e_1] = e_3^{c_1}, [e_2^{-1},e_0] = e_1^k e_3^{c_2} , [e_0,e_3]
= 1 , [e_1,e_2] = e_3^n \rangle.$$ This is a discrete torsion-free
nilpotent group, which can be embedded as a lattice in a connected
and simply-connected nilpotent Lie group by \cite[Theorem
2.18]{Rag}. Since a solvmanifold is uniquely determined by its
fundamental group, $G / \Gamma$ is diffeomorphic to a nilmanifold.

As at the beginning of the proof, we conclude that $G/\Gamma$ is the
total space of a $T^2$-bundle over $T^2$ and occurs in our list. The
quotient $G/\Gamma$ is of type b) if $k = 0$ and of type e) if $k
\ne 0$.

\textbf{Case B.1.2.:} $\widetilde{\lambda}_1 =
\widetilde{\lambda}_2^{-1} \in \mathbb{R}_+ \setminus \{1\}$

\nopagebreak \noindent We have $\overline{\nu}(1) = \left(
\begin{array}{cc}
a_1 &  a_2 \\ b_1 & b_2
\end{array} \right) \in \mathrm{SL}(2, \mathbb{Z}),$
and Proposition \ref{nu forsetzbar} implies
$$\begin{array}{ccl} \nu(1)(x,y,\frac{z}{n}) ~= & \big( a_1 x + a_2
y,~~ b_1 x + b_2 y, & a_2 b_1 x y + a_1 b_1 \frac{x(x-1)}{2} + a_2
b_2 \frac{y(y-1)}{2} \\ && + \frac{1}{n}(c_1 x + c_2 y  + (a_1 b_2 -
a_2 b_1)z) \big)
\end{array} $$
for certain $c_1,c_2 \in \mathbb{Z}$.

Choose eigenvectors $\left(
\begin{array}{c}
v_1 \\ v_ 2
\end{array}
\right), \left(
\begin{array}{c}
w_1 \\ w_ 2
\end{array}
\right) \in \mathbb{R}^2 \setminus \{0\}$ with respect to the
eigenvalues $\widetilde{\lambda}_1$ resp.\ $\widetilde{\lambda}_2$
of $^{\tau} \overline{\nu}(1)$ (where $^{\tau}$\! denotes the
transpose). There exist $u_1, u_2, u_3 \in \mathbb{R}$ such that for
$\gamma_i := (v_i,w_i,u_i)$, $i \in \{1,2\}$, and $\gamma_3 :=
(0,0,u_3) \in U_3(\mathbb{R})$ we have
\begin{eqnarray*}
& [\gamma_1,\gamma_2] = \gamma_3^n ,& \\ & \widetilde{\mu}(1)
(\gamma_1) = \gamma_1^{a_1} \gamma_2^{b_1} \gamma_3^{c_1}, ~~
\widetilde{\mu}(1) (\gamma_2) = \gamma_1^{a_2} \gamma_2^{b_2}
\gamma_3^{c_2}, &
\end{eqnarray*}
where $\widetilde{\mu}(t) \big( (x,y,z) \big) = \big( \exp (t
\ln(\widetilde{\lambda}_1)) \, x , \exp (t
\ln(\widetilde{\lambda}_2)) \, y , z \big).$

Then $G/\Gamma$ is diffeomorphic to the solvmanifold $\widetilde{G}
/ \widetilde{\Gamma}$, where $\widetilde{G} = \mathbb{R}
\ltimes_{\widetilde{\mu}} \nolinebreak U_3(\mathbb{R})$ and
$\widetilde{\Gamma} = \mathbb{Z} \ltimes_{\widetilde{\mu}} \langle
\gamma_1 ,\gamma_2 ,\gamma_3 \rangle$, i.e.\ $G / \Gamma$ is a Inoue
surface of type $S^+$, see \cite[Section 3.7]{Has4}. The Lie algebra
of $\widetilde{G}$, $$\widetilde{\mathfrak{g}} = \langle \, T,X,Y,Z
\,|\, [T,X] = X, \, [T,Y] = -Y, \, [X,Y] = Z \, \rangle,$$ is
completely solvable and not nilpotent. Further, the knowledge of
$\widetilde{\mathfrak{g}}$ implies $b_1(G/\Gamma) = 1$. By Theorem
\ref{nicht sympl}, $G/\Gamma$ is formal. Therefore, this is a
solvmanifold of type h) in Table \ref{Solv4}.

\textbf{Case B.2.:} $\widetilde{\lambda}_1 = \widetilde{\lambda}_2 =
-1$ and $\overline{\nu}(1)$ is diagonalisable

\nopagebreak \noindent \cite[Lemma 1]{SF} implies that we can assume
$\overline{\nu}(1) = \left(
\begin{array}{cc}
-1 &  0 \\ 0 & -1
\end{array} \right) \in \mathrm{SL}(2, \mathbb{Z})$.
So Proposition \ref{nu forsetzbar} implies $\nu(1) \big(
(x,y,\frac{z}{n}) \big) = (-x ,-y, \frac{c_1 x + c_2 y +z}{n})$ for
certain integers $c_1,c_2 \in \mathbb{Z}$. Moreover, $G/\Gamma$ is
diffeomorphic to $\widetilde{G}/\widetilde{\Gamma}$, where
$\widetilde{G} := \mathbb{R} \ltimes_{\widetilde{\mu}}
U_3(\mathbb{R})$,
$$\widetilde{\mu}(t) \big( ( x, y, z) \big) = \big(
\cos (t \pi) \,x - \sin (t \pi) \,y , \sin (t \pi) \,x + \cos (t
\pi) \,y, z + h_t(x,y) \big),$$ $h_t(x,y) = \frac{1}{2} \sin(t \pi)
\Big( cos(t \pi) \big( x^2 - y^2 \big) - 2 \sin(t \pi) x y \Big)$
and $\widetilde{\Gamma} = \mathbb{Z} \ltimes_{\widetilde{\mu}}
\langle \gamma_1 ,\gamma_2 ,\gamma_3 \rangle$ such that
$[\gamma_1,\gamma_2] = \gamma_3^n$, $\widetilde{\mu}(1) (\gamma_1) =
\gamma_1^{-1} \gamma_3^{c_1}$ and $\widetilde{\mu}(1) (\gamma_2) =
\gamma_2^{-1} \gamma_3^{c_2}$. (Using the addition theorems for
$\sin$ and $\cos$, one calculates that $\widetilde{\mu}$ is a
one-parameter group in $\mathrm{A}(U_3(\mathbb{R}))$.) By
\cite[Section 3.5]{Has4}, $\widetilde{G}/\widetilde{\Gamma}$ is a
secondary Kodaira surface.

Obviously, the Lie algebra of $\widetilde{G}$ is not completely
solvable and we cannot use its Chevalley-Eilenberg complex for
computing $b_1(G/\Gamma)$. But since
$$\Gamma = \langle e_0, \ldots, e_3 \,|\, e_0 e_1 e_0^{-1} = e_1^{-1} e_3^{c_1}
, e_0 e_2 e_0^{-1} = e_2^{-1} e_3^{c_1}, [e_0,e_3]=1,
[e_1,e_2]=e_3^n \rangle, $$ we see $b_1(G/\Gamma) = \mathrm{rank} \,
\Gamma_{ab} = 1$ and $G/\Gamma$ belongs to the last row in Table
\ref{Solv4}.

\textbf{Case B.3.:} $\widetilde{\lambda}_1 =
\overline{\widetilde{\lambda}_2} = \exp(i \varphi) \in \mathbb{C}
\setminus \mathbb{R}, \varphi \in ]0,2\pi[ \setminus \{\pi\}$

\nopagebreak \noindent This case is similar to the last one. We have
$| \mathrm{tr} \, \overline{\nu}(1) | \le |\widetilde{\lambda}_1| +
|\widetilde{\lambda}_2| = 2$ and \cite[Lemma 1]{SF} implies that we
can assume $\nu(1)$ to be $\left(
\begin{array}{cc}
0 & -1 \\ 1 & 0
\end{array}
\right)$ or $\pm \left(
\begin{array}{cc}
0 & -1 \\ 1 & 1
\end{array}
\right)$. In each case, one computes $b_1(G/\Gamma) = \mathrm{rank}
\, \Gamma_{ab} = 1$, as above. Moreover, one embeds a lattice
$\widetilde{\Gamma}$ isomorphic to $\Gamma$ in the Lie group
$\widetilde{G} := \mathbb{R} \ltimes_{\widetilde{\mu}}
U_3(\mathbb{R})$, where
$$\widetilde{\mu}(t) \big( ( x, y, z) \big) = \big(
\cos (t \varphi) \,x - \sin (t \varphi) \,y , \sin (t \varphi) \,x +
\cos (t \varphi) \,y, z + h_t(x,y) \big),$$ $h_t(x,y) = \frac{1}{2}
\sin(t \varphi) \Big( cos(t \varphi) \big( x^2 - y^2 \big) - 2
\sin(t \varphi) x y \Big)$. Again, $\widetilde{G} /
\widetilde{\Gamma}$ is a secondary Kodaira surface and $G/\Gamma$ is
an example for case i). For more details see \cite[Section
3.5]{Has4}. \q
\N Below, we give examples for each of the nine types of
four-\-di\-men\-sio\-nal solvmanifolds. The Lie algebras of the
connected and simply-connected four-\-di\-men\-sio\-nal solvable Lie
groups that admit lattices are listed in Table \ref{Sol4} in
Appendix \ref{Liste Solv}.

\begin{Ex}
The following manifolds belong to the corresponding row in Table
\ref{Solv4}.
\begin{itemize}
\item[a)] $\mathbb{R}^4 / \mathbb{Z}^4$
\item[b)] $(\mathbb{R} \ltimes_{\mu_b} \mathbb{R}^3 )/ (\mathbb{Z}
\ltimes_{\mu_b} \mathbb{Z}^3),~ \mu_b(t) = \left(
\begin{array}{ccc} 1 & 0 & 0 \\ 0 & 1 & t \\ 0 & 0 & 1
\end{array} \right)$
\item[c)] $(\mathbb{R} \ltimes_{\mu_c} \mathbb{R}^3) / \Gamma_c$ with
\begin{equation*}
\Gamma_c = \mathbb{Z} \ltimes_{\mu_c} \langle \left(
\begin{array}{c} 1 \\ 0 \\ 0
\end{array} \right), \left(
\begin{array}{c} 0 \\ 1 \\ 1 \end{array} \right) , \left( \begin{array}{c}
0 \\ \frac{18 + 8 \sqrt{5}}{7 + 3 \sqrt{5}} \\ \frac{2}{3 +
\sqrt{5}} \end{array} \right) \rangle_{\mathbb{Z}},
\end{equation*}
$t_1 = \ln(\frac{3 + \sqrt{5}}{2})$ and $\mu_c(t) = \left(
\begin{array}{ccc} 1 & 0 & 0 \\ 0 & e^{t \, t_1} & 0  \\ 0 & 0 &
e^{-t \, t_1} \end{array} \right);$ the proof that this is really a
solvmanifold is analogous to that in the example on page
\pageref{Ex3}.
\item[d)] $(\mathbb{R} \ltimes_{\mu_d} \mathbb{R}^3) / (\pi \, \mathbb{Z}
\ltimes_{\mu_d} \mathbb{Z}^3),~ \mu_d(t) = \left(
\begin{array}{ccc} 1 & 0 & 0 \\ 0 & \cos(t) & -\sin(t) \\ 0 & \sin(t) &
\cos(t) \end{array} \right)$
\item[e)] $(\mathbb{R} \ltimes_{\mu_e} \mathbb{R}^3 )/ (\mathbb{Z}
\ltimes_{\mu_e} \mathbb{Z}^3),~ \mu_e(t) = \left(
\begin{array}{ccc} 1 & t & \frac{1}{2} (t^2 -t) \\ 0 & 1 & t \\ 0 & 0 & 1
\end{array} \right)$
\item[f)] Consider $A := \left(
\begin{array}{ccc} 0 & 0 & 1 \\ 1 & 0 & -11 \\ 0 & 1 &
8 \end{array} \right) \in \mathrm{SL}(3,\mathbb{Z})$. $A$ has $X^3 - 8
X^2 + 11 X -1$ as characteristic polynomial which possesses three
pairwise different real roots $t_1 \approx 6,271$, $t_2 \approx
1,631$ and $t_3 \approx 0,098$. Therefore, $A$ is conjugate to
$\mu_f(1)$, where $\mu_f(t) = \left(
\begin{array}{ccc} e^{t \, \ln(t_1)} & 0 & 0 \\ 0 &  e^{t \, \ln(t_2)} & 0
\\ 0 & 0 &  e^{t \, \ln(t_3)} \end{array} \right)$, and this implies
the existence of a lattice $\Gamma_f$ in the completely solvable Lie
group $\mathbb{R} \ltimes_{\mu_f} \mathbb{R}^3$.
\item[g)] Let $A := \left(
\begin{array}{ccc} 0 & 0 & 1 \\ 1 & 0 & -8 \\ 0 & 1 &
4 \end{array} \right) \in \mathrm{SL}(3,\mathbb{Z})$. The
characteristic polynomial of $A$ is $X^3 - 4 X^2 + 8 X -1$ which has
three pairwise different roots $t_1 \approx 0,134$ and $t_{2,3} =
(1/\sqrt{t_1}) \, ( \cos(\varphi) \pm i \sin(\varphi) ) \approx
1,933 \pm 1,935 \, i$. So $A$ is conjugate to $\mu_g(1)$, where
$\mu_g(t) = \left(
\begin{array}{ccc} e^{t \, \ln(t_1)} & 0 & 0 \\ 0 &
e^{t \, \ln(|t_2|)} \cos(t \, \varphi) & -e^{t \, \ln(|t_2|)} \sin(t
\, \varphi)
\\ 0 & e^{t \, \ln(|t_2|)} \sin(t
\, \varphi) &  e^{t \, \ln(|t_2|)} \cos(t \, \varphi) \end{array}
\right)$, and this implies the existence of a lattice $\Gamma_g$ in
the Lie group $\mathbb{R} \ltimes_{\mu_g} \mathbb{R}^3$.
\item[h)] Using Theorem \ref{rat Strukturkonst}, one shows that
\begin{eqnarray*}
\gamma_1   &:=& (1,1,-\frac{1+\sqrt{5}}{3+\sqrt{5}}),
\\ \gamma_2 &:=& (- \frac{2(2+\sqrt{5})}{3+\sqrt{5}},
\frac{1+\sqrt{5}}{3+\sqrt{5}} ,- \frac{11+5\sqrt{5}}{7+3\sqrt{5}}),
\\  \gamma_3 &:=& (0,0,\sqrt{5})
\end{eqnarray*}
generate a lattice $\Gamma$ in $U_3(\mathbb{R})$ with
$[\gamma_1,\gamma_2] = \gamma_3$ and $\gamma_3$ central.

Define the one-parameter group $\mu_h \: \mathbb{R} \to
\mathrm{A}(U_3(\mathbb{R}))$ by $$\mu_h (t) \big( (x,y,z) \big) =
(e^{-t \, t_1} x, e^{t \, t_1} y, z), $$ where $t_1 := \ln(
\frac{3+\sqrt{5}}{2} )$. Then $\mu_h(1)$ preserves the lattice
$\Gamma$ with $$\mu_h(1)(\gamma_1) = \gamma_1^2 \, \gamma_2, ~
\mu_h(1) (\gamma_2) = \gamma_1 \, \gamma_2, ~ \mu_h(1)(\gamma_3) =
\gamma_3$$ and therefore, $\mathbb{Z} \ltimes_{\mu_h} \Gamma$ is a
lattice in $\mathbb{R} \ltimes_{\mu_h} U_3(\mathbb{R})$.
\item[i)] Consider the Lie group $\widetilde{G}$ and the one-parameter
group $\widetilde{\mu}$ of Case B.2 from the proof of the last
theorem. Setting $\gamma_1 = (1,0,0), ~ \gamma_2 = (0,1,0)$ as well
as $\gamma_3 = (0,0,1), ~ n=1$ and $c_1=c_2=0$, one explicitly gets
an example.
\end{itemize}
\end{Ex}

The manifolds of type c) show that formal spaces with the same
minimal model as a K\"ahler manifold need not be K\"ahlerian. This
was proved by Fern\'andez and Gray.

\begin{Thm}[\cite{FG90}]
Let $M$ be one of the symplectic solvmanifolds of type c) in the
last theorem, i.e.\ $M$ is formal and possesses no complex
structure. $M$ has the same minimal model as the K\"ahler manifold
$T^2 \times S^2$. \q
\end{Thm}

\section{Five-dimensional solvmanifolds}
We study the five-dimensional solvmanifolds by regarding lattices in
the corresponding connected and simply-connected Lie groups. By
Proposition \ref{unimod}, their Lie algebras have to be unimodular.
These are listed in Appendix \ref{Liste Solv}.

\subsection{Nilpotent and decomposable solvable Lie algebras}
There are nine classes of nilpotent Lie algebras in dimension five,
see Table \ref{Nil5}. Each of them has a basis with rational
structure constants. By Theorem \ref{rat Strukturkonst}, the
corresponding connected and simply-connected Lie groups admit
lattices and accordingly to Theorem \ref{formale Nilmgf}, the
associated nilmanifolds are formal if and only if they are tori. For
$i \in \{4,5,6\}$ the connected and simply-connected nilpotent Lie
group with Lie algebra $\mathfrak{g}_{5.i}$ possesses the
left-invariant contact form $x_1$ (where $x_1$ is dual to the basis
element $X_1 \in \mathfrak{g_i}$ as in Table \ref{Nil5}). Therefore,
the corresponding nilmanifolds are contact.

The eight classes of decomposable unimodular non-nilpotent solvable
Lie algebras are listed in Table \ref{UniDecSol5}. Except for
$\mathfrak{g}_{4.2} \oplus \mathfrak{g}_1$, the corresponding
connected and simply-connected Lie groups admit lattices since both
of their factors admit lattices.

\begin{SThm}
The connected and simply-connected Lie group $G_{4.2} \times
\mathbb{R}$ with Lie algebra $\mathfrak{g}_{4.2} \oplus
\mathfrak{g}_1$ possesses no lattice.
\end{SThm}

\textit{Proof.} Write $G$ for $G_{4.2} \times \mathbb{R}$ and
$$\mathfrak{g} = \langle X_1, \ldots, X_5 \,|\, [X_1,X_4] = -2 X_1,\,
[X_2,X_4] = X_2,\, [X_3,X_4] = X_2 + X_3 \rangle$$ for its Lie
algebra which has $\mathfrak{n} = \mathbb{R}^4_{X_1,X_2,X_3,X_5}$ as
nilradical. Therefore, $G$ can be written as almost abelian Lie
group $\mathbb{R} \ltimes_{\mu} \mathbb{R}^4$ with
$$ \mu(t) = \exp^{GL(n,\mathbb{R})} ( t\, \mathrm{ad}(X_4)) = \left(
\begin{array}{cccc}
e^{2t} & 0 & 0 & 0 \\ 0 & e^{-t} & -t e^{-t} & 0 \\ 0 & 0 & e^{-t} &
0 \\ 0 & 0 & 0 & 1
\end{array}
\right).$$ By Corollary \ref{Gitter in fastabelsch}, the existence
of a lattice in $G$ would imply that there is $t_1 \in \mathbb{R}
\setminus \{0\}$ such that $\mu(t_1)$ is conjugate to an element of
$\mathrm{SL}(4,\mathbb{Z})$. Clearly, the characteristic polynomial
of $\mu(t_1)$ is $\, P(X) = (X-1) \, \widetilde{P}(X) \,$, where the
polynomial $\widetilde{P}(X) = X^3 - k X^2 + l X - 1 \in
\mathbb{Z}[X]$ has the double root $e^{-t_1}$. Lemma \ref{Z3doppelt}
then implies $e^{-t_1} = 1$, i.e.\ $t_1 = 0$ which is a
contradiction. \q

\begin{SProp}
If $\Gamma$ is a lattice in a five-dimensional completely solvable
non-nilpotent connected and simply-connected decomposable Lie group
$G$, then $G/\Gamma$ is formal.
\end{SProp}

\textit{Proof.} Let $G$, $\Gamma$ be as in the proposition. As
usual, we denote by $\mathfrak{g}$ the Lie algebra of $G$. We have
$\mathfrak{g} = \mathfrak{h} \, \oplus \, k \mathfrak{g}_1$ with $k
\in \{1,2\}$ and a certain $(5-k)$-dimensional completely solvable
non-nilpotent Lie algebra $\mathfrak{h}$, see Tables
\ref{UniDecSol5} and \ref{Sol4}. By completely solvability and
Theorem \ref{min chev eil} (ii), $G/\Gamma$ and the
Chevalley-Eilenberg complex of $\mathfrak{h} \, \oplus \,k
\mathfrak{g}_1$ share their minimal model $\mathcal{M}$. The lower
dimensional discussion above shows that for all $\mathfrak{h}$ which
can arise in the decomposition of $\mathfrak{g}$ the algebras
$\mathcal{M}_{(\bigwedge \mathfrak{h}^*, \delta_{\mathfrak{h}})}$
and $\mathcal{M}_{(\bigwedge k \mathfrak{g}_1^*, \delta = 0)} =
(\bigwedge k \mathfrak{g}_1^*, \delta = 0)$ are formal. This implies
the formality of $\mathcal{M}= \mathcal{M}_{(\bigwedge
\mathfrak{h}^*, \delta_{\mathfrak{h}})} \otimes
\mathcal{M}_{(\bigwedge k \mathfrak{g}_1^*, \delta = 0)}$. \q
\subsection{Indecomposable non-nilpotent Lie algebras}
There are $19$ classes of indecomposable non-nilpotent Lie algebras
in dimension five which are unimodular. These are listed in Tables
\ref{UniAlab5} -- \ref{Uni2ab5}. Instead of the small German letters
for the Lie algebras in the mentioned tables, we use capital Latin
letters (with the same subscripts) for the corresponding connected
and simply-connected Lie groups.

We want to examine which of them admit lattices and where
appropriate, whether the quotients are formal. The non-existence
proofs of lattices in certain almost abelian Lie groups below are
taken from Harshavardhan's thesis \cite{Har}. Some of the existence
proofs of lattices in certain almost abelian Lie groups are sketched
in \cite[pp.\ 29 and 30]{Har}.

\subsubsection*{Almost abelian algebras}
We now consider the almost abelian Lie groups $G_{5.i} = \mathbb{R}
\ltimes_{\mu_i} \mathbb{R}^4$. We write $\mu(t) = \mu_i(t) =
\exp^{GL(4,\mathbb{R})}( t \, \mathrm{ad}(X_5) )$, where $X_5 \in
\mathfrak{g}_{5.i}$ is as in Table \ref{UniAlab5} ($X_5$ depends on
$i$). We know by Corollary \ref{Gitter in fastabelsch}, Theorem
\ref{1ParamSL} and Proposition \ref{b1 fast abelsch} that there is a
lattice $\Gamma$ in $G_{5.i}$ if and only if there exists $t_1 \ne
0$ such that $\mu(t_1)$ is conjugate to $\widetilde{\mu}(1) \in
\mathrm{SL}(4,\mathbb{Z})$ and $\Gamma = \mathbb{Z}
\ltimes_{\widetilde{\mu}} \mathbb{Z}^4$. This will be used in the
proof of the following propositions.

Methods to obtain integer matrices with given characteristic
polynomial and necessary conditions for their existence are given in
Appendix \ref{AppZ}.

\begin{SProp} \label{G5.7}
Let $p,q,r \in \mathbb{R}$ with $-1 \le r \le q \le p \le 1$, $pqr
\ne 0$ and $p + q + r = -1$. If the completely solvable Lie group
$G_{5.7}^{p,q,r}$ admits a lattice and $M$ denotes the corresponding
solvmanifold, then $M$ is formal, $b_1(M) = 1$ and one of the
following conditions holds:
\begin{itemize}
\item[(i)] $b_2(M) = 0$,
\item[(ii)] $b_2(M) = 2$, i.e.\ $r = -1,~ p = -q \in \, ]0,1[$ or
\item[(iii)] $b_2(M) = 4$, i.e.\ $r = q = -1,~ p=1$.
\end{itemize}
Moreover, there exist $p,q,r$ as above satisfying (i), (ii) resp.\
(iii) such that $G_{5.7}^{p,q,r}$ admits a lattice.
\end{SProp}

\textit{Proof.} We suppress the sub- and superscripts of $G$ and
$\mathfrak{g}$.

a) Assume, there is a lattice in $G$ and denote the corresponding
solvmanifold by $M$. Since $\mathfrak{g}$ is completely solvable,
the inclusion of the Chevallier-Eilenberg complex $\big(
\bigwedge(x_1, \ldots, x_5),\delta \big)$ into the forms on $M$
induces an isomorphism on cohomology. Moreover, the minimal model of
$\big( \bigwedge(x_1, \ldots, x_5),\delta \big)$ is isomorphic to
the minimal model of $M$.

$\delta$ is given by $$\delta x_1 = - x_{15} ,\, \delta x_2 = -p\,
x_{25} ,\, \delta x_3 = -q \, x_{35} ,\, \delta x_4 = -r \, x_{45}
,\, \delta x_5 = 0.$$ (Here we write $x_{ij}$ for $x_i x_j$.) This
implies $b_1(M) = 1$.

One computes the differential of the non-exact generators of degree
two in the Chevalley-Eilenberg complex as
$$\begin{array}{l@{~~~~}l@{~~~~}l}
\delta x_{12} = (1+p)\, x_{125}, & \delta x_{13} = (1+q)\, x_{135},
& \delta x_{14} = (1+r)\, x_{145}, \\ \delta x_{23} = (p+q)\,
x_{235}, & \delta x_{24} = (p+r)\, x_{245}, & \delta x_{34} =
(q+r)\, x_{345}.
\end{array}$$
$-1 \le r \le q \le p \le 1$, $pqr \ne 0$ and $p + q + r = -1$
implies $p \ne -1$ and $q \ne -r$ and a short computation yields
that either (i), (ii) or (iii) holds.

In each case, we determine the $2$-minimal model, i.e.\ the minimal
model up to generators of degree two and will see, that these
generators are closed. By Definition \ref{s-form}, the minimal model
then is $2$-formal and Theorem \ref{formal = n-1 formal} implies the
formality of $M$.

If we are in case (i), the minimal model has  one closed generator
of degree one, and no generator of degree two.

If we are in case (ii), we have $r = -1,~ p = -q \in \, ]0,1[$,
$$\begin{array}{l@{~~~~}l@{~~~~}l}
\delta x_{12} = (1+p) x_{125} \ne 0, &  \delta x_{13} = (1-p)
x_{135} \ne 0, & \delta x_{14} = 0, \\  \delta x_{23} = 0, & \delta
x_{24} = (p-1)\, x_{245} \ne 0, & \delta x_{34} = (-1-p)\, x_{345}
\ne 0,
\end{array}$$ $$
\begin{array}{ccc}
H^1(M) & \cong & \langle [x_5] \rangle, \\
H^2(M) & \cong & \langle [x_{14}],[x_{23}] \rangle,
\end{array}$$
and the $2$-minimal model $\rho \: (\bigwedge V^{\le 2} , d) \to
\big( \bigwedge(x_1, \ldots, x_5),\delta \big)$ is given by
$$\begin{array}{c@{~~~~}c@{~~~~}c}
\rho(y) = x_5, & |y|= 1, & dy = 0; \\
\rho(z_1) = x_{14}, & |z_1|= 2, & dz_1 = 0; \\
\rho(z_2) = x_{23}, & |z_2|= 2, & dz_2 = 0.
\end{array}$$
Note, further generators of degree $\le 2$ do not occur, since $y^2
= 0$ (by graded commutativity) and $\rho(y z_i)$ is closed and
non-exact. Here we use the construction of the minimal model that we
have given in the proof of Theorem \ref{Konstruktion d m M}.

Case (iii) is similar to case (ii).

b) Now, we show that there are examples for each of the three cases.
In case (i), we follow \cite{Har} and consider the matrix $\left(
\begin{array}{cccc}
1 & 0 & 0 & -2 \\ 1 & 2 & 0 & -3 \\ 0 & 1 & 3 & 5 \\
0 & 0 & 1 & 2
\end{array}
\right).$ It suffices to show that there are $t_1 \in \mathbb{R}
\setminus \{0\}$, $-1 < r < q < p < 1$ with $pqr \ne 0$, $p \ne -q$,
$p \ne -r$, $q \ne -r$ and $p + q + r = -1$ such that $$\mu(t_1) =
\exp^{GL(4,\mathbb{R})}( t_1 \, \mathrm{ad}(X_5) ) = \left(
\begin{array}{cccc}
e^{-t_1} & 0 & 0 & 0 \\ 0 & e^{-pt_1} & 0 & 0
\\ 0 & 0 & e^{-qt_1} & 0 \\ 0 & 0 & 0 & e^{-rt_1}
\end{array}
\right)$$ is conjugate to the matrix above, which has $P(X) = X^4 -
8 X^3 +18 X^2 - 10 X +1$ as characteristic polynomial. $P$ has four
distinct roots $\lambda_1, \ldots, \lambda_4$ with $\lambda_1
\approx 0,12$, $\lambda_2 \approx 0,62$, $\lambda_3 \approx 2,79$
and $\lambda_4 \approx 4,44$. Define $t_1 := - \ln(\lambda_1)$ and
$p,q,r$ by $e^{-pt_1} = \lambda_2$, $e^{-qt_1} = \lambda_3$ and
$e^{-rt_1} = \lambda_4$. Then $t_1,p,q,r$ have the desired
properties.

In case (ii), regard the matrix $\left(
\begin{array}{cccc}
0 & 0 & 0 & -1 \\ 1 & 0 & 0 & 10 \\ 0 & 1 & 0 & -23 \\
0 & 0 & 1 & 10
\end{array}
\right)$ which is conjugate to $\mu(t_1) = \left(
\begin{array}{cccc}
e^{-t_1} & 0 & 0 & 0 \\ 0 & e^{-pt_1} & 0 & 0
\\ 0 & 0 & e^{pt_1} & 0 \\ 0 & 0 & 0 & e^{t_1}
\end{array}
\right)$ for $t_1 = 2 \ln(\frac{3 + \sqrt{5}}{2})$ and $p =
\frac{1}{2}$ since both matrices have the same characteristic
polynomial which has four distinct real roots.

In case (iii), regard the matrix $\left(
\begin{array}{cccc}
3 & 0 & -1 & 0 \\ 0 & 3 & 0 & -1 \\ 1 & 0 & 0 & 0 \\
0 & 1 & 0 & 0
\end{array}
\right)$ which is conjugate to $\mu(t_1) = \left(
\begin{array}{cccc}
e^{-t_1} & 0 & 0 & 0 \\ 0 & e^{-t_1} & 0 & 0
\\ 0 & 0 & e^{t_1} & 0 \\ 0 & 0 & 0 & e^{t_1}
\end{array}
\right)$ for $t_1 = \ln(\frac{3 + \sqrt{5}}{2})$ since both matrices
have the same minimal polynomial by Proposition \ref{Matrix} (ii).
\q
\N We have seen that a non-formal solvmanifold is a non-toral
nilmanifold in dimensions three and four. In higher dimensions this
is no longer true as the following proposition shows:

\begin{SProp}\label{G58}
The completely solvable Lie group $G_{5.8}^{-1}$ admits a lattice.

Moreover, for each lattice $\Gamma$ the corresponding solvmanifold
$M = G_{5.8}^{-1} / \Gamma$ has $b_1(M) = 2$ and is not formal.
\end{SProp}

\textit{Proof.} Again, we suppress the sub- and superscripts. $G$
admits a lattice since $\mu(t) = \exp^{GL(4,\mathbb{R})}( t \,
\mathrm{ad}(X_5) ) = \left(
\begin{array}{cccc}
1 & -t & 0 & 0 \\ 0 & 1 & 0 & 0
\\ 0 & 0 & e^{-t} & 0 \\ 0 & 0 & 0 & e^{t}
\end{array}
\right)$ and $\left(
\begin{array}{cccc}
0 & 0 & 0 & -1 \\ 1 & 0 & 0 & 5
\\ 0 & 1 & 0 & -8 \\ 0 & 0 & 1 & 5
\end{array}
\right)$ are conjugated for $t_1 = \ln(\frac{3 + \sqrt{5}}{2})$.
Note that the transformation matrix $T \in
\mathrm{GL}(4,\mathbb{R})$ with $T A T^{-1} = \mu(t_1)$ is
\begin{align*}
T = \left(
\begin{array}{cccc}
 1 & 0 & -1 & -2 \\
 \frac{1}{\ln(\frac{3 + \sqrt{5}}{2})} & \frac{1}{\ln(\frac{3 +
\sqrt{5}}{2})}
& \frac{1}{\ln(\frac{3 + \sqrt{5}}{2})} & \frac{1}{\ln(\frac{3 + \sqrt{5}}{2})} \\
 -\frac{5 + 3\sqrt{5}}{10} & -\frac{1}{\sqrt{5}} & \frac{5 -
3\sqrt{5}}{10} & \frac{3}{2} - \frac{7}{2 \sqrt{5}} \\
 \frac{-5 + 3\sqrt{5}}{10} & \frac{1}{\sqrt{5}} & \frac{5 +
3\sqrt{5}}{10} & \frac{3}{2} + \frac{7}{2 \sqrt{5}}
\end{array} \right).
\end{align*}

Now, let $\Gamma$ be an arbitrary lattice in $G$. By completely
solvability and Theorem \ref{min chev eil} (ii), we get the minimal
model of $M = G / \Gamma$ as the minimal model $\mathcal{M}$ of the
Chevalley-Eilenberg complex $(\bigwedge \mathfrak{g}^* , \delta )$.
The latter is given by
$$\delta x_1 = - x_{25}, \, \delta x_2 = 0, \, \delta x_3 = - x_{35},\,
\delta x_4 = x_{45} ,\, \delta x_5 = 0,$$ which implies $b_1(M) =
2$. Further, the minimal model $\rho \: (\bigwedge V, d) \to
(\bigwedge \mathfrak{g}^* , \delta )$ must contain two closed
generators $y_1,y_2$ which map to $x_2$ and $x_5$. Then we have
$\rho(y_1 y_2) = x_{25} = - \delta x_1$ and the minimal model's
construction in the proof of Theorem \ref{Konstruktion d m M}
implies that there is another generator $u$ of degree one such that
$\rho(u) = - x_1$ and $d u = y_1 y_2$. Since $\rho(u y_1) = -
x_{12}$ and $\rho(u y_2) = -x_{15}$ are closed and non-exact, there
are no further generators of degree one in $V$. But this implies 
that $(u+c) \, y_1$ is closed and non-exact in $\mathcal{M}$ for 
each closed element $c$ of degree one. Using the
notation of Theorem \ref{nichtformal}, we have $u \in N^1, y_1 \in
V^1$ and $\mathcal{M}$ is not formal. \q

\begin{SProp}
The completely solvable Lie group $G_{5.9}^{p,-2-p}$, $p \ge -1$,
does not admit a lattice.
\end{SProp}

\textit{Proof.} The first half of the proof is taken from
\cite{Har}. Assume there is a lattice. $\mu(t) = \left(
\begin{array}{cccc}
e^{-t} & -t e^{-t} & 0 & 0 \\
0 & e^{-t} & 0 & 0 \\
0 & 0 & e^{-tp} & 0 \\
0 & 0 & 0 & e^{t(2+p)}
\end{array} \right)$ is conjugate to an element of $\mathrm{SL}(4,\mathbb{Z})$
for $t = t_1 \ne 0$ and has roots $e^{-t_1}, e^{-t_1}, e^{-t_1p}$
and $e^{t_1(2+p)}$. By Proposition \ref{Z4doppelt}, this can occur
if and only if $p = -1$. Therefore, for the remainder of the proof
we assume $p = -1$.

The Jordan form of $\mu(t_1)$ is $\left(
\begin{array}{cccc}
e^{-t} & 1 & 0 & 0 \\
0 & e^{-t} & 0 & 0 \\
0 & 0 & e^{t} & 0 \\
0 & 0 & 0 & e^{t}
\end{array} \right)$, i.e.\ the characteristic and the minimal
polynomial of $\mu(t_1)$ are
\begin{eqnarray*}
P(X) & = & (X - e^{-t_1})^2 (X - e^{t_1})^2 \\ & = & X^4 - 2
(e^{-t_1} + e^{t_1}) X^3 + (e^{-2t_1} + e^{2t_1} + 4) X^2 - 2
(e^{-t_1} + e^{t_1}) X + 1, \\ m(X) & = & (X - e^{-t_1})^2 (X -
e^{t_1}) \\ & = & X^3 - (2 e^{-t_1} + e^{t_1}) X^2 + (e^{-2t_1} +2)
X - e^{-t_1}.
\end{eqnarray*}
Since $\mu(t_1)$ is conjugate to an integer matrix, we have $P(X) ,
m(X) \in \mathbb{Z}[X]$ by Theorem \ref{minpol}. This is impossible
for $t_1 \ne 0$. \q

\begin{SProp}[\cite{Har}]
The completely solvable Lie group $G_{5.11}^{-3}$ does not admit a
lattice.
\end{SProp}

\textit{Proof.} If the group admits a lattice, there exists $t_1 \in
\mathbb{R} \setminus \{0\}$ such that the characteristic polynomial
of $\mu(t_1) = \left(
\begin{array}{cccc}
e^{-t_1} & -t_1 e^{-t_1} & \frac{t_1^2}{2} e^{-t_1} & 0 \\
0 & e^{-t_1} & -t_1 e^{-t_1} & 0 \\
0 & 0 & e^{-t_1} & 0 \\
0 & 0 & 0 & e^{3t_1}
\end{array} \right)$ is a monic integer polynomial with a three-fold root
$e^{-t_1}$ and a simple root $e^{3t_1}$. By Proposition
\ref{Z4doppelt}, this is impossible for $t_1 \ne 0$. \q

\begin{SProp}\label{G513a}
There are $q,r \in \mathbb{R}$ with $-1 \le q < 0,\, q \ne
-\frac{1}{2}, r \ne 0$ such that $G_{5.13}^{-1-2q,q,r}$ admits a
lattice.
\end{SProp}

\textit{Proof.} We have $\mu_{q,r}(t) = \left(
\begin{array}{cccc}
e^{-t} & 0 & 0 & 0 \\
0 & e^{t+2qt} & 0 & 0 \\
0 & 0 & e^{-qt} \cos(r t) & -e^{-qt} \sin(r t) \\
0 & 0 & e^{-qt} \sin(r t) & e^{-qt} \cos(r t)
\end{array} \right)$ and as claimed in \cite{Har}, there exist
$t_1 \ne 0$, $q_0, r_0$ such that $\mu_{q_0,r_0}(t_1)$ is conjugate
to $A := \left(
\begin{array}{cccc}
1 & 0 & 0 & 1 \\
1 & 2 & 0 & 2 \\
0 & 1 & 3 & 0 \\
0 & 0 & 1 & 0
\end{array} \right)$ which implies the existence of a lattice
$\Gamma_A$ in $G_{5.13}^{-1-2q_0,q_0,r_0}$.

If $\lambda_1 \approx 0,15 < \lambda_2 \approx 3,47$ denote the real
roots and $\lambda_{3,4} \approx 1,17 \pm i \, 0,67$ the non-real
roots of $P_A(X) = X^4 - 6 X^3 + 11 X^2 - 8 X + 1$, then $t_1 =
-\ln(\lambda_1) \approx 1,86$, $q_0 =
\frac{1}{2}(\frac{\ln(\lambda_2)}{t_1} - 1) \approx -0,16$ and $r_0
= \frac{1}{t_1}\arccos\big(\mathrm{Re}(\lambda_3) e^{q_0t_1}\big)
\approx 0,27$. \q

\begin{Rem}
If the real number $\frac{\pi}{t_1 r_0}$ is not rational, then
Theorems \ref{min chev eil} (iii) and \ref{Gorb pi} enable us to
show that the manifold $G_{5.13}^{-1-2q_0,q_0,r_0} / \Gamma_A$ has
$b_1 = 1$ and is formal.
\end{Rem}

\begin{SProp}\label{G513b}
There exists $r \in \mathbb{R} \setminus \{0\}$ such that
$G_{5.13}^{-1,0,r}$ admits a lattice.
\end{SProp}

\textit{Proof.} Let $t_1 = \ln(\frac{3+\sqrt{5}}{2})$, $r = \pi/t_1$
and $A = \left(
\begin{array}{cccc}
3 & 1 & 0 & 0 \\
-1 & 0 & 0 & 0 \\
0 & 0 & -1 & 0 \\
0 & 0 & 0 & -1
\end{array} \right)$. Then $A$ is conjugate to $\mu_{0,r}(t_1) = \left(
\begin{array}{cccc}
e^{-t_1} & 0 & 0 & 0 \\
0 & e^{t_1} & 0 & 0 \\
0 & 0 & \cos(r t_1) & \sin(r t_1) \\
0 & 0 & \sin(r t_1) & \cos(r t_1)
\end{array} \right)$ and this implies the existence of a lattice.

Note that we have $T A T^{-1} = \mu_{0,r}(t_1)$, where $T = \left(
\begin{array}{cccc}
1 & \frac{18 + 8 \sqrt{5}}{7 + 3 \sqrt{5}} & 0 & 0 \\
1 & \frac{2}{3 + \sqrt{5}} & 0 & 0 \\
0 & 0 & 1 & 0 \\
0 & 0 & 0 & 1
\end{array} \right)$. \q

\begin{Rem}
Since the abelianisation of the lattice in the last proof is
isomorphic to $\mathbb{Z} \oplus \mathbb{Z}_2{}^2$, the constructed
solvmanifold has $b_1 = 1$.
\end{Rem}

\begin{SProp}\label{G514}
$G_{5.14}^{0}$ admits a lattice.
\end{SProp}

\textit{Proof.} We have $\mu(t) = \left(
\begin{array}{cccc}
1 & -t & 0 & 0 \\
0 & 1 & 0 & 0 \\
0 & 0 & \cos(t) & -\sin(t) \\
0 & 0 & \sin(t) & \cos(t)
\end{array} \right)$. Let $t_1 = \frac{\pi}{3}$, then
$\mu(t_1)$ is conjugate to $\left(
\begin{array}{cccc}
1 & 0 & 0 & 0 \\
1 & 1 & 0 & 0 \\
0 & 1 & 1 & -1 \\
0 & 0 & 1 & 0
\end{array} \right)$, so there is a lattice. Note that the matrix
$T = \left(
\begin{array}{cccc}
-1 & 1 & 0 & 0 \\
-\frac{\pi}{3} & 0 & 0 & 0 \\
-\frac{1}{\sqrt{3}} & -\frac{1}{\sqrt{3}} & \frac{2}{\sqrt{3}} &
-\frac{1}{\sqrt{3}} \\ 1 & -1 & 0 & 1
\end{array} \right) \in \mathrm{GL}(4,\mathbb{R})$ satisfies $T A T^{-1} = \mu(t_1)$. \q

\begin{Rem}
The abelianisation of the lattice in the last proof is isomorphic to
$\mathbb{Z}^2$, i.e.\ the corresponding solvmanifold has $b_1 =2$.
\end{Rem}

\begin{SProp}
If there is a lattice $\Gamma$ in the Lie group $G := G_{5.14}^0$
such that $b_1(G/\Gamma) = 2$, then $G/\Gamma$ is not formal.
\end{SProp}

\textit{Proof.} By Theorem \ref{min chev eil}(i), the natural
inclusion of the Chevalley-Eilenberg complex $(\bigwedge
\mathfrak{g}^* , \delta) \to (\Omega(G/\Gamma), d)$ induces an
injection on cohomology. $(\bigwedge \mathfrak{g}^* , \delta)$ is
given by
$$\delta x_1 = - x_{25},\, \delta x_2 = 0,\, \delta x_3 = - x_{45},\,
\delta x_4 = x_{35},\, \delta x_5 = 0.$$ This implies $b_1(\bigwedge
\mathfrak{g}^* , \delta) = 2$, hence $H^1(G/\Gamma,d) = \langle
[x_2], [x_5] \rangle.$ Therefore $$[x_2] \cdot H^1(G/\Gamma) +
H^1(G/\Gamma) \cdot [x_5] = \langle [x_{25}] \rangle = \langle
[\delta x_1] \rangle = 0,$$ and in the Massey product $\langle
[x_2], [x_2], [x_5] \rangle = [-x_{15}]$ is no indeterminacy. Since
$x_{15}$ is closed and not exact, $G/\Gamma$ cannot be formal. \q

\begin{SProp}  \label{G5.15}
The completely solvable Lie group $G_{5.15}^{-1}$ admits a lattice.
For each lattice the corresponding solvmanifold satisfies $b_1 = 1$
and is non-formal.
\end{SProp}

\textit{Proof.} As we have done above, we suppress the sub- and
superscripts. First, we follow \cite{Har} and consider the matrix $A
:= \left(
\begin{array}{cccc}
2 & 0 & 0 & -1 \\
1 & 2 & 0 & 2 \\
0 & 1 & 1 & 2 \\
0 & 0 & 1 & 1
\end{array} \right)$ which is conjugate to $\mu(t_1) = \left(
\begin{array}{cccc}
e^{-t_1} & -t_1 e^{-t_1} & 0 & 0 \\
0 & e^{-t_1} & 0 & 0 \\
0 & 0 & e^{t_1} & -t_1 e^{t_1} \\
0 & 0 & 0 & e^{t_1}
\end{array} \right)$ for $t_1 = \ln(\frac{3 + \sqrt{5}}{2})$.
This implies the existence of a lattice. The transformation matrix
$T \in \mathrm{GL}(4,\mathbb{R})$ with $T A T^{-1} = \nolinebreak
\mu(t_1)$ is
$$T = \left(
\begin{array}{cccc}
\frac{2}{5 \sqrt{5}} & -\frac{1}{5 \sqrt{5}} & -\frac{2}{5
\sqrt{5}} & \frac{1}{2} + \frac{3}{50} \sqrt{5} \\
 -\frac{3 + \sqrt{5}}{10  \ln(\frac{3 + \sqrt{5}}{2})} & \frac{2
+ \sqrt{5}}{5  \ln(\frac{3 + \sqrt{5}}{2})} & -\frac{(3 +
\sqrt{5})^2}{20 \ln(\frac{3 + \sqrt{5}}{2})} &
\frac{2 + \sqrt{5}}{5  \ln(\frac{3 + \sqrt{5}}{2})} \\
 -\frac{2}{5 \sqrt{5}} & \frac{1}{5 \sqrt{5}} & \frac{2}{5
\sqrt{5}}
& \frac{1}{2} - \frac{3}{50} \sqrt{5} \\
 -\frac{2}{5(3 + \sqrt{5}) \ln(\frac{3 + \sqrt{5}}{2})} &
-\frac{-1+\sqrt{5}}{5(3 + \sqrt{5}) \ln(\frac{3 + \sqrt{5}}{2})} &
 \frac{-3 + \sqrt{5}}{5 (3 + \sqrt{5}) \ln(\frac{3 +
\sqrt{5}}{2})} & -\frac{-1+\sqrt{5}}{5(3 + \sqrt{5}) \ln(\frac{3 +
\sqrt{5}}{2})}
\end{array} \right). $$

Let $\Gamma$ be an arbitrary lattice in $G$. By completely
solvability and Theorem \ref{min chev eil} (ii), we get the minimal
model of $M = G / \Gamma$ as the minimal model $\mathcal{M}$ of the
Chevalley-Eilenberg complex $(\bigwedge \mathfrak{g}^* , \delta )$.
The latter is given by
$$\delta x_1 = - x_{15} - x_{25},~ \delta x_2 = - x_{25},
~ \delta x_3 = x_{35} - x_{45},~ \delta x_4 = x_{45} ,~ \delta x_5 =
0,$$ which implies $b_1(M) = 1$.

One computes the differential of the non-exact generators of degree
two in the Chevalley-Eilenberg complex as
\begin{eqnarray*}
\delta x_{12} = 2 x_{125}, & \delta x_{13} = x_{145} + x_{235}, &
\delta x_{14} = x_{245}, \\ \delta x_{23} = x_{245}, & \delta x_{24}
= 0, & \delta x_{34} = -2 x_{345}.
\end{eqnarray*}
The minimal model $\rho \: (\bigwedge V, d) \to (\bigwedge
\mathfrak{g}^* , \delta )$ must contain three closed generators
$y,z_1,z_2$ which map to $x_5,x_{14}-x_{23}$ and $x_{24}$. We see
$\rho(y z_1) = x_{145}-x_{235}$ is closed and non-exact, $\rho(y
z_2) = x_{245} = \delta x_{23}$ and the minimal model's construction
in the proof of Theorem \ref{Konstruktion d m M} implies that there
is another generator $u$ of degree two such that $\rho(u) = x_{23}$
and $d u = y z_2$. Since $\rho(uy) = x_{235}$ is closed and
non-exact, there are no further generators of degree less than or
equal to two in $V$. But this implies that $(u+c) \, y$ is closed and
non-exact in $\mathcal{M}$ for each closed element $c$ of degree two.
Using the notation of Theorem \ref{nichtformal}, we have $u \in N^2, 
y \in V^1$ and $\mathcal{M}$ is not formal. \q

\begin{SProp}[\cite{Har}]
$G_{5.16}^{-1,q}$, $q \ne 0$, does not admit a lattice.
\end{SProp}

\textit{Proof.} If the group admits a lattice, there exists $t_1 \in
\mathbb{R} \setminus \{0\}$ such that the characteristic polynomial
of $\mu(t_1) = \left(
\begin{array}{cccc}
e^{-t_1} & -t_1 e^{-t_1} & 0 & 0 \\
0 & e^{-t_1} & 0 & 0 \\
0 & 0 & e^{t_1} \cos(t_1 q) & -e^{t_1} \sin(t_1 q) \\
0 & 0 & e^{t_1} \sin(t_1 q) & e^{t_1} \cos(t_1 q)
\end{array} \right)$ is a monic integer polynomial with simple roots
$e^{t_1} (\cos(t_1 q) \pm i \sin(t_1 q))$ and a double root
$e^{-t_1}$. By Proposition \ref{Z4doppelt}, this is impossible for
$t_1 \ne 0$. \q

\begin{SProp}\label{G517a}
There are $p,r \in \mathbb{R}$, $p \ne 0$, $r \notin \{0,\pm1\}$,
such that $G_{5.17}^{p,-p,r}$ admits a lattice.
\end{SProp}

\textit{Proof.} We follow \cite{Har}.

$A:= \left(
\begin{array}{cccc}
2 & 0 & 0 & -11 \\
1 & 2 & 0 & -9 \\
0 & 1 & 1 & -1 \\
0 & 0 & 1 & 1
\end{array} \right)$ is conjugate to $$\mu(t_1) = \left(
\begin{array}{cccc}
e^{-t_1 p} \cos(t_1) & -e^{-t_1 p} \sin(t_1) & 0 & 0 \\
e^{-t_1 p} \sin(t_1) & e^{-t_1 p} \cos(t_1) & 0 & 0 \\
0 & 0 & e^{t_1 p} \cos(t_1 r) & -e^{t_1 p} \sin(t_1 r) \\
0 & 0 & e^{t_1 p} \sin(t_1 r) & e^{t_1 p} \cos(t_1 r)
\end{array} \right)$$ for certain $t_1,p,r \ne 0$, i.e.\ there is a lattice.

If $\lambda_{1,2} \approx 0,306 \pm i \, 0,025$ and $\lambda_{3,4}
\approx 2,694 \pm i \, 1,83 $ denote the roots of $P_A(X) = X^4 - 6
X^3 + 14 X^2 - 7 X + 1$, one has $t_1 p = -\ln(|\lambda_1|) \approx
1,181$, hence \linebreak $t_1 = \arccos\big(\mathrm{Re}(\lambda_1)
e^{t_1p}\big) \approx 0,062$, $p \approx 14,361$. $t_1 r =
\arccos\big(\mathrm{Re}(\lambda_3) e^{t_1 p}\big) \approx 0,597$
implies $r \approx 7,259$. \q

\begin{Rem}
Since the abelianisation of the lattice in the last proof is
isomorphic to $\mathbb{Z} \oplus \mathbb{Z}_3$, the corresponding
solvmanifold has $b_1 =1$.
\end{Rem}

\begin{SProp}\label{G517b}
There exists $p \in \mathbb{R} \setminus \{0\}$ such that
$G_{5.17}^{p,-p,\pm 1}$ admits a lattice.
\end{SProp}

\textit{Proof.} Let $p := \frac{1}{\pi} \ln(\frac{3+\sqrt{5}}{2})$,
$t_1 := \pi$ and $A := \left(
\begin{array}{cccc}
0 & -1 & 0 & 0 \\
1 & -3 & 0 & 0 \\
0 & 0 & 0 & -1 \\
0 & 0 & 1 & -3
\end{array} \right)$. Then $\mu(t_1) = \left(
\begin{array}{cccc}
e^{-t_1 p} \cos(t_1) & -e^{-t_1 p} \sin(t_1) & 0 & 0 \\
e^{-t_1 p} \sin(t_1) & e^{-t_1 p} \cos(t_1) & 0 & 0 \\
0 & 0 & e^{t_1 p} \cos(\pm t_1) & -e^{t_1 p} \sin(\pm t_1) \\
0 & 0 & e^{t_1 p} \sin(\pm t_1) & e^{t_1 p} \cos(\pm t_1)
\end{array} \right)$ is conjugate to $A$ and this implies the
existence of a lattice. Note that we have $T A T^{-1} = \mu(t_1)$
with $T := \left(
\begin{array}{cccc}
\frac{1}{\sqrt{5}} & \frac{5-3\sqrt{5}}{10} & 0 & 0 \\
0 & 0 & \frac{1}{\sqrt{5}} & \frac{5-3\sqrt{5}}{10} \\
-\frac{1}{\sqrt{5}} & \frac{5+3\sqrt{5}}{10} & 0 & 0 \\
0 & 0 & -\frac{1}{\sqrt{5}} & \frac{5+3\sqrt{5}}{10}
\end{array} \right)$. \q

\begin{Rem}
The abelianisation of the lattice in the last proof is $\mathbb{Z}
\oplus \mathbb{Z}_3{}^2$, i.e.\ the corresponding solvmanifold has
$b_1 = 1$.
\end{Rem}

\begin{SProp}\label{G517c}
There exists $r \in \mathbb{R} \setminus \{ 0, \pm1 \}$ such that
$G_{5.17}^{0,0,r}$ admits a lattice.
\end{SProp}

\textit{Proof.} Let $r \in \{2,3\}$. Then $\mu(t) = \left(
\begin{array}{cccc}
\cos(t) & \sin(t) & 0 & 0 \\
\sin(t) & \cos(t) & 0 & 0 \\
0 & 0 & \cos(t r) & \sin(t r) \\
0 & 0 & \sin(t r) & \cos(t r)
\end{array} \right)$ is an integer matrix for $t = \pi$.
This implies the existence of a lattice. \q

\begin{Rem}
If we chose in the last proof $r = 2$, then the corresponding
solvmanifold has $b_1 = 3$. For $r = 3$ we obtain a solvmanifold
with $b_1 = 1$.
\end{Rem}

\begin{SProp}\label{G517d}
$G_{5.17}^{0,0,\pm 1}$ admits a lattice.
\end{SProp}

\textit{Proof.} $\mu(t) = \left(
\begin{array}{cccc}
\cos(t) & \sin(t) & 0 & 0 \\
\sin(t) & \cos(t) & 0 & 0 \\
0 & 0 & \cos(\pm t) & \sin(\pm t) \\
0 & 0 & \sin(\pm t) & \cos(\pm t)
\end{array} \right)$ is an integer matrix for $t = \pi$.
This implies the existence of a lattice.\q

\begin{Rem}
The first Betti number of the solvmanifold induced by the lattice of
the last proof equals one.
\end{Rem}

\begin{SProp}\label{G518}
$G_{5.18}^0$ admits a lattice.
\end{SProp}

\textit{Proof.} Again, we follow \cite{Har}. The matrix $\left(
\begin{array}{cccc}
2 & 0 & 0 & -9 \\
1 & 0 & 0 & -4 \\
0 & 1 & 0 & -3 \\
0 & 0 & 1 & 0
\end{array} \right)$ is conjugate to $\mu(t_1) = \left(
\begin{array}{cccc}
\cos(t_1) & - \sin(t_1) & -t_1 \cos(t_1) & t_1 \sin(t_1) \\
\sin(t_1) & \cos(t_1) & -t_1 \sin(t_1) & -t_1 \cos(t_1) \\
0 & 0 & \cos(t_1) & - \sin(t_1) \\
0 & 0 & \sin(t_1) & \cos(t_1)
\end{array} \right)$ for $t_1 =
\frac{\pi}{3}$. This implies the existence of a lattice.

Note, $T = \left(
\begin{array}{cccc}
\frac{4}{3 \sqrt{3}} & - \frac{2}{\sqrt{3}} & 0 & - \frac{1}{\sqrt{3}} \\
0 & 0 & 0 & 1 \\
\frac{\sqrt{3}}{\pi} & - \frac{2\sqrt{3}}{\pi} &  -
\frac{\sqrt{3}}{\pi} &
\frac{\sqrt{3}}{\pi} \\
\frac{1}{\pi} & 0 &  - \frac{3}{\pi} & - \frac{3}{\pi}
\end{array} \right) \in \mathrm{GL}(4,\mathbb{R})$ is the transformation matrix with $T A T^{-1} =
\mu(t_1)$. \q

\begin{Rem}
The abelianisation of the lattice in the last proof is isomorphic to
$\mathbb{Z}$, i.e.\ the corresponding solvmanifold has $b_1 =1$.
\end{Rem}

\subsubsection*{Algebras with nilradical $\mathfrak{n} :=
\mathfrak{g}_{3.1} \oplus \mathfrak{g}_1 = \langle X_1,\ldots,X_4
\,|\, [X_2,X_3] = X_1 \rangle$}
We now regard the unimodular almost-nilpotent Lie groups $G_{5.i}$
with nilradical $N := U_3(\mathbb{R}) \times \mathbb{R}$, i.e.\ $i
\in \{19,20,23,25,26,28\}$. We can identify $N$ with $\mathbb{R}^4$
as a manifold and the group law given by
$$(a,b,c,r) \cdot (x,y,z,w) = (a+x+bz\, , \, b+y \, , \, c+z \, , \,r+w).$$

The Lie algebras of the unimodular Lie groups $G_{5.i} = \mathbb{R}
\ltimes_{\mu_i} N$ with nilradical $N$ are listed in Table
\ref{UniAlnil15}. We have $\mu_i(t) = \exp^N \circ
\exp^{A(\mathfrak{n})}(t \, \mathrm{ad}(X_5)) \circ \log^N$, where
$X_5$ depends on $i$.

Assume there is a lattice $\Gamma$ in $G_{5.i}$. By Corollary
\ref{Gitter in fastabelsch}, there are $t_1 \ne 0$ and an inner
automorphism $I_{n_1}$ of $N$ such that $\nu_i := \mu_i(t_1) \circ
I_{n_1}, \nu_i^{-1} \in \mathrm{A}(N)$ preserve the lattice
$\Gamma_N := \Gamma \cap N$ in $N$. For $n_1 = (a,b,c,r)$ one
calculates
\begin{equation}\label{In1} I_{n_1}(x,y,z,w) = (x + bz - yc\,,\,y\,,\,z\,,\,w).
\end{equation} $\Gamma_{N{'}} := \Gamma_N \cap N{'} \cong \mathbb{Z}$
is a lattice in $N{'} := [N,N] = \{(x,0,0,0) \,|\, x \in
\mathbb{R}\} \cong \mathbb{R}$ by Theorem \ref{Zentrum} and since
$\nu_i(\Gamma_{N{'}}), \nu_i^{-1}(\Gamma_{N{'}}) \subset
\Gamma_{N{'}}$, we have $\nu_i|_{\Gamma_{N{'}}} \in
\mathrm{Aut}(\mathbb{Z})$. This implies $\nu_i|_{\Gamma_{N{'}}} =
\pm \mathrm{id}$ and hence $\mu_i(t_1)|_{[N,N]} = \pm \mathrm{id}$
(a cause of (\ref{In1}) and the shape of $[N,N]$). Moreover, we have
$[\mathfrak{n},\mathfrak{n}] = \langle X_1 \rangle$ and since
$\exp^{\mathbb{R}}$ is the identity,
$$\pm \mathrm{id} = \mu_i(t_1)|_{[N,N]} = \exp^{A(\mathfrak{n})}(t_1 \,
\mathrm{ad}(X_5)|_{\langle X_1 \rangle})|_{[N,N]}.$$ (Note that
$\exp^N([\mathfrak{n},\mathfrak{n}]) = [N,N]$ by \cite[Theorem
3.6.2]{V}.) Therefore, $t_1 [X_5,X_1]$ has no component in $\langle
X_1 \rangle$ and since $t_1 \ne 0$, this means that $[X_1,X_5]$ has
no component in $X_1$-direction. The list of Lie algebras in Table
\ref{UniAlnil15} implies:

\begin{SProp} \label{nur526}
The only connected and simply-connected solvable Lie groups with
nilradical $U_3(\mathbb{R}) \times \mathbb{R}$ that can contain a
lattice are $G_{5.20}^{-1}$ and $G_{5.26}^{0,\pm 1}$. \q
\end{SProp}

\begin{Rem}
In a previous version of this article, the group $G_{5.20}^{-1}$ is
absent. It was added, after the author had read \cite{DF} in April 2009. 
\end{Rem}

\begin{SProp}\label{G520}
The completely solvable Lie group $G_{5.20}^{-1}$ admits a lattice. 
For each lattice the corresponding solvmanifold admits a contact form, 
is formal and has $b_1 = 2$.
\end{SProp}

\textit{Proof.} Using Theorem \ref{rat Strukturkonst}, one shows that
\begin{eqnarray*}
\gamma_1   &:=& (\frac{20 + 9 \sqrt{5}}{9 + 4 \sqrt{5}} ,0,0,0),
\\ \gamma_2 &:=& ( \frac{181 + 81 \sqrt{5}}{47+ 21 \sqrt{5}},
\frac{18+8 \sqrt{5}}{7+3\sqrt{5}} ,\frac{2}{3+\sqrt{5}},0),
\\ \gamma_3 &:=& (\frac{181 + 81 \sqrt{5}}{47+ 21 \sqrt{5}}, 1, 1, 0),
\\ \gamma_4 &:=& (0,0,0,-\frac{20 + 9 \sqrt{5}}{(9 + 4 \sqrt{5})
\ln(\frac{3 + \sqrt{5}}{2})})
\end{eqnarray*}
generate a lattice $\Gamma_N$ in $N$ with
$[\gamma_2,\gamma_3] = \gamma_1$ and $\gamma_1, \gamma_4$ central.

A short calculation yields that $\mu(t) \big((x,y,z,w)\big) = 
(x-tw, e^{-t}y,e^tz,w)$ defines a one-parameter group in 
$\mathrm{A}(N)$. Moreover, for $t_1 = \ln(\frac{3+\sqrt{5}}{2})$ holds
$\mu(t_1)(\gamma_1) = \gamma_1$, $\mu(t_1)(\gamma_2) = \gamma_3$, 
$\mu(t_1)(\gamma_3) = \gamma_2^{-1} \gamma_3^3$
and $\mu(t_1)(\gamma_4) = \gamma_1 \gamma_4$.

This implies the existence of a lattice in $G := G_{5.20}^{-1} =
\mathbb{R} \ltimes_{\mu} N$. 
 
Let $\Gamma$ be an arbitrary lattice in $G$. By completely
solvability and Theorem \ref{min chev eil} (ii), we get the minimal
model of $M = G / \Gamma$ as the minimal model $\mathcal{M}$ of the
Chevalley-Eilenberg complex $(\bigwedge \mathfrak{g}^* , \delta )$.
The latter is given by
$$\delta x_1 = - x_{23} - x_{45},~ \delta x_2 = - x_{25},
~ \delta x_3 = x_{35} ,~ \delta x_4 = \delta x_5 =
0,$$ which implies $b_1(M) = 2$. Moreover, $x_1$ defines a 
left-invariant contact form on $G/\Gamma$.

One computes the differential of the non-exact generators of degree
two in the Chevalley-Eilenberg complex as
\begin{eqnarray*}
\delta x_{12} = x_{125} - x_{245}, & \delta x_{13} = - x_{135} - x_{345},&
\delta x_{14} = - x_{234}, \\ \delta x_{15} = - x_{235}, & \delta x_{23}
= 0, & \delta x_{24} = x_{245}, \\ \delta x_{34} = - x_{345}, & 
\delta x_{45} =0 . &
\end{eqnarray*}
The minimal model $\rho \: (\bigwedge V, d) \to (\bigwedge
\mathfrak{g}^* , \delta )$ must contain two closed generators
$y_1,y_2$ which map to $x_4$ and $x_{5}$. We see
$\rho(y_1 y_2) = x_{45}$ is closed and non-exact. Since $b_2(G/\Gamma)=1$,
the minimal model's construction
in the proof of Theorem \ref{Konstruktion d m M} implies that there
are no further generators of degree less than or
equal to two in $V$. This implies that $G/\Gamma$ is formal. \q

\begin{SProp}\label{G526}
$G_{5.26}^{0,\varepsilon}$ admits a lattice for $\varepsilon = \pm
1$. For each lattice the corresponding solvmanifold is contact and
has $b_1 \ge 2$.
\end{SProp}

\textit{Proof.} One calculates that $\mu \: \mathbb{R} \to
\mathrm{A}(N)$ defined by \pagebreak
\begin{eqnarray*}
&& \mu(t) \big( ( x, y, z, w) \big) \\ && ~~  =  \big( x + h_t(y,z)
- \varepsilon t w, \cos (t \pi) \,y - \sin (t \pi) \,z , \sin (t
\pi) \,y + \cos (t \pi) \,z , w \big),
\end{eqnarray*}
where $h_t(y,z) = \frac{1}{2} \sin(t \pi) \Big( cos(t \pi) \big( y^2
- z^2 \big) - 2 \sin(t \pi) y z \Big)$, is a one-parameter group.

Then we have $G := G_{5.26}^{0,\varepsilon} = \mathbb{R}
\ltimes_{\mu} N$ and $\mathbb{Z} \ltimes_{\mu} \{(x,y,z,w) \in N
\,|\, x,y,z,w \in \mathbb{Z}\}$ is a lattice in $G$ since $\mu (1)
\big( (x,y,z,w) \big) = ( x - \varepsilon w , - y , -z , w ).$

Using $d_e\big(\mu(t)\big) = \log^N \circ \mu(t) \circ \exp^N$, we
obtain the Lie algebra $\mathfrak{g}$ of $G$ as
$$\langle X_1, \dots X_5 \,|\, [X_2,X_3] = X_1,\,
[X_2,X_5] = X_3,\, [X_3,X_5] = - X_2,\, [X_4,X_5] = \varepsilon X_1
\rangle.$$ Denote $\{x_1,\ldots,x_5\}$ the basis of $\mathfrak{g}^*$
which is dual to $\{X_1,\ldots,X_5\}$, i.e.\ the $x_i$ are
left-invariant $1$-forms on $G$. One calculates that $x_1$ is a
left-invariant contact form on $G$, so it descends to a contact form
on the corresponding solvmanifold.

The statement about the first Betti number follows from Theorem
\ref{min chev eil}(i). \q

\begin{Rem}
Since the abelianisation of the lattice in the last proof is
isomorphic to $\mathbb{Z}^2 \oplus \mathbb{Z}_2^2$, the
corresponding solvmanifold has $b_1 = 2$.
\end{Rem}

\subsubsection*{Algebras with nilradical $\mathfrak{g}_{4.1} = \langle
 X_1,\ldots,X_4 \,|\, [X_2,X_4] = X_1,\, [X_3,X_4] = X_2 \rangle$}
\begin{SProp}\label{garnixG41}
No connected and simply-connected solvable Lie group $G_{5.i}$ with
nilradical $N := G_{4.1}$ admits a lattice.
\end{SProp}

\textit{Proof.} There is only one unimodular connected and
simply-connected solvable Lie group with nilradical $G_{4.1}$,
namely the completely solvable group $G := G_{5.30}^{-\frac{4}{3}}$.
We show that it admits no lattice.

The group $N$ is $\mathbb{R}^4$ as a manifold  with multiplication
given by
$$(a,b,c,r) \cdot (x,y,z,w) = (a+ x + wb + \frac{1}{2} w^2 c\,,\, b
+y + wc \,,\, c+z \,,\, r+ w),$$ and one calculates for $n_1 =
(a,b,c,r)$
$$ I_{n_1}(x,y,z,w) = (x + wb + \frac{1}{2} w^2 c - ry - rwc
+ \frac{1}{2} r^2 z\,,\, y + wc - rz\,,\,z\,,\,w).$$

Let $G = \mathbb{R} \ltimes_{\mu} N$, where $ \mu(t) = \exp^N \circ
\exp^{A(\mathfrak{n})}(t \, \mathrm{ad}(X_5)) \circ \log^N$ and
assume there is a lattice $\Gamma$ in $G$. By Corollary \ref{Gitter
in fastabelsch}, there are $t_1 \ne 0$ and $n_1 \in N$ such that
$\nu := \mu(t_1) \circ I_{n_1} \in \mathrm{A}(N)$ preserves the
lattice $\Gamma_N := \Gamma \cap N$ in $N$.

$\Gamma_{N'} := N' \cap \Gamma_N$ is a lattice in $N{'} := [N,N] =
\{(x,y,0,0) \in N \,|\, x,y \in \mathbb{R}\} \cong \mathbb{R}^2$ by
Theorem \ref{Zentrum}, and since $\nu (N{'}) \subset N{'}$, this
lattice is preserved by $\nu|_{N'}$. This and $\exp^{\mathbb{R}^2} =
\mathrm{id}$ imply
$$\pm 1 = \det (\nu|_{N'}) = \det \big( \exp^{A(\mathfrak{n})}(t_1 \,
\mathrm{ad}(X_5)|_{[\mathfrak{n},\mathfrak{n}]})|_{[N,N]} \big)
\cdot \underbrace{\det (I_{n_1}|_{N'})}_{=1},$$ i.e.\
$\mathrm{ad}(X_5)|_{[\mathfrak{n},\mathfrak{n}]}$ has trace equal to
zero. This contradicts $\mathfrak{g}_{5.30}^{-\frac{4}{3}}$, see
Table \ref{UniAlnil25}. \q

\subsubsection*{Non-almost nilpotent algebras}
Now, there remain two unimodular connected and simply-connected
solvable Lie groups in dimension five, namely $G_{5.33}^{-1,-1}$ and
$G_{5.35}^{-2,0}$. Unfortunately, we do not know whether the former
group admits a lattice or not. Note, Harshavardhan's argumentation
in \cite[p.\ 33]{Har} is not sufficient.

\begin{Rem}
If the completely solvable Lie group $G_{5.33}^{-1,-1}$ admits a
lattice, one easily proves that the corresponding solvmanifold
admits a contact form (since $G_{5.33}^{-1,-1}$ possesses the
left-invariant contact form $x_1+x_2+x_3$ with $x_i$ dual to $X_i
\in \mathfrak{g}_{5.33}^{-1,-1}$ as in Table \ref{Uni2ab5}), is
formal and has $b_1=2$.
\end{Rem}

\begin{Rem}
In April 2009, A.\ Diatta and B.\ Foreman \cite{DF} proved that 
$G_{5.33}^{-1,-1}$ possesses a lattice.
\end{Rem}

\begin{SProp}\label{G535}
$G_{5.35}^{-2,0}$ contains a lattice. For each lattice the
corresponding solvmanifold is contact and has $b_1 \ge 2$.
\end{SProp}

\textit{Proof.} A lattice and a contact form were constructed by
Geiges in \cite{Ge5}. One has the left-invariant contact form $x_1 +
x_2$ on the Lie group, where $x_1,x_2$ are dual to the
left-invariant vector fields as in Table \ref{Uni2ab5}. Hence the
form descends to each compact quotient by a discrete subgroup.

The statement about the first Betti number follows from Theorem
\ref{min chev eil}(i). \q

\subsubsection*{Conclusion}
We have seen that each connected and simply-connected
$5$-dimensional solvable Lie group admits a lattice if it is
nilpotent or decomposable with the exception of $G_{4.2} \times
\mathbb{R}$. If an indecomposable non-nilpotent group $G_{5.i}$
gives rise to a solvmanifold it is contained in Table \ref{Solv5}.
Recall, by Theorem \ref{min chev eil}, we always have a lower bound
for the solvmanifold's Betti numbers and in some cases the exact
value. These can be read of in the second and the third column. The
last column refers to the examples that we have constructed above.
``yes'' means that we have such for certain parameters that satisfy
the conditions of the column ``Comment''. Except for $i=33$ we have
examples for all possible values of $i$.
\begin{table}[htb]
\centering \caption{\label{Solv5} $5$-dimensional indecomposable
non-nilmanifolds}
\begin{tabular}{|c|c|c|c|c|c|} \hline \hline
& $b_1$ & $b_2$ & formal & Comment & Example
\\ \hline
\hline $G_{5.7}^{p,q,r}$ & $1$ & $0$ & yes & $-1 < r < p < q <1,$&
\ref{G5.7} (i)
\\ &&&& $pqr \ne 0$, & \\ &&&& $p+q+r=-1$ &\\
\hline $G_{5.7}^{p,q,-1}$ & $1$ & $2$ & yes & $p=-q \in ]0,1[$ & \ref{G5.7} (ii)\\
\hline $G_{5.7}^{1,-1,-1}$ & $1$ & $4$ & yes & & \ref{G5.7} (iii)\\
\hline $G_{5.8}^{-1}$ & $2$ & $3$ & no & & \ref{G58}\\
\hline $G_{5.13}^{-1-2q,q,r}$ & $\ge 1$ & $\ge 0$ & ? & $q\in[-1,0]
\setminus\{\frac{1}{2}\},$ & \ref{G513a} \\ &&&& $r \ne 0$ &\\
\hline $G_{5.13}^{-1,0,r}$ & $\ge 1$ & $\ge 2$ & ? & $r \ne 0$ & \ref{G513b} \\
\hline $G_{5.14}^{0}$ & $\ge 2$ & $\ge 3$ & ?  & & \ref{G514}\\
\hline $G_{5.15}^{-1}$ & $1$ & $2$ & no & & \ref{G5.15}\\
\hline $G_{5.17}^{p,-p,r}$ & $\ge 1$ & $\ge 0$ & ? & $p\ne 0,~
r \notin \{0,\pm1\}$& \ref{G517a} \\
\hline $G_{5.17}^{p,-p,\pm1}$ & $\ge 1$ & $\ge 2$ & ? & $p\ne 0$ & \ref{G517b}\\
\hline $G_{5.17}^{0,0,r}$ & $\ge 1$ & $\ge 2$ & ? & $r \notin
\{0,\pm1\}$& \ref{G517c} \\
\hline $G_{5.17}^{0,0,\pm1}$ & $\ge 1$ & $\ge 4$ & ? & & \ref{G517d}\\
\hline $G_{5.18}^{0}$ & $\ge 1$ & $\ge 2$ & ? & & \ref{G518}\\
\hline $G_{5.20}^{-1}$ & $2$ & $1$ & yes & & \ref{G520}\\
\hline $G_{5.26}^{0,\pm 1}$ & $\ge 2$ & $\ge 1$ & ? & & \ref{G526}\\
\hline $G_{5.33}^{-1,-1}$ & $2$ & $1$ & yes & & no\\
\hline $G_{5.35}^{-2,0}$ & $\ge 2$ & $\ge 1$ & ? & & \ref{G535} \\
\hline
\end{tabular}
\end{table}

Assuming that there is a lattice in one the non-completely solvable
Lie groups $G_{5.i}$, i.e.\ $i \in \nolinebreak
\{13,14,17,18,26,35\}$, such that the inequalities in the above
table are equalities, then one can calculate that such quotients are
formal for $i \in \{ 13, 17, 26, 35\}$ and not formal for $i \in
\{14,18\}$. The assumptions about the Betti numbers are needed to
ensure that the Lie algebra cohomology is isomorphic to the
solvmanifold's cohomology.
\subsection{Contact structures}
Some of the connected and simply-connected five-dimensional solvable
Lie groups $G_{5.i}$ which admit a lattice $\Gamma$ possess a
left-invariant contact form. Obviously, it also defines a contact
form on the corresponding solvmanifold. By this way, we showed that
the manifolds $G_{5.i}/\Gamma$ for $i \in \{4,5,6\}$ and quotients
of almost nilpotent groups with non-abelian nilradical (i.e.\ $i \ge
19$) by lattices are contact.

But $\mathbb{R}^5$, $U_3(\mathbb{R}) \times \mathbb{R}^2$, $G_{4.1}
\times \mathbb{R}$ and $G_{5.i}$ do not have a left-invariant
contact form for $i \in \{1,2,3,7,\ldots,18\}$, see e.g.\
\cite{Diatta}. For some of the nilmanifolds, we can provide a
contact structure by another approach.

\begin{SThm}  \label{5 contact}
Let $G \in \{ \mathbb{R}^5, U_3(\mathbb{R}) \times \mathbb{R}^2,
G_{4.1} \times \mathbb{R}, G_{5.1}, G_{5.3} \}$ and $\Gamma$ a
lattice $G$. Then $G / \Gamma$ admits a contact structure.
\end{SThm}

\textit{Proof.} For $G$ chosen as in the theorem, the dimension of
the center is greater than or equal to two. Therefore, we can find a
two-dimensional closed normal subgroup that lies in the center such
that its intersection with $\Gamma$ is a lattice in it. By Theorem
\ref{SolvHauptfaserbdl}, $G/\Gamma$ has the structure of a principal
$T^2$-bundle over a three dimensional closed orientable manifold.
Then the following result of Lutz implies the claim. \q

\begin{SThm}[\cite{Lutz}]
The total space of a principal $T^2$-bundle over a closed orientable
$3$-manifold admits a contact form. \q
\end{SThm}

Unfortunately, we did not find a contact structure on the manifold
of Proposition \ref{G5.15}. If such exists, this yields a
five-dimensional non-formal contact solvmanifold with $b_1 = 1$.
\section{Six-dimensional solvmanifolds}
There are $164$ types of connected and simply-connected
indecomposable solvable Lie groups in dimension six, most of them
depending on parameters. For classifying six-dimensional
solvmanifolds, we restrict ourselves to the following types:
\begin{itemize}
\item[(1)] nilmanifolds
\item[(2)] symplectic solvmanifolds that are quotients of
indecomposable non-nilpotent groups
\item[(3)] products of lower-dimensional solvmanifolds
\end{itemize}

Although we have to make some restrictions to get a manageable number
of cases, one certainly has to consider types (1) and (3). Concerning the
third type, the reader can even ask the legitimate question why we do not
consider arbitrary lattices in products of lower dimensional
solvable Lie groups $G_1, G_2$, instead of direct products $\Gamma_1
\times \Gamma_2$ of lattices $\Gamma_i$ in the factors $G_i$. The
reason is that we have no tool to construct arbitrary lattices or
disprove their existence, unless we can ensure that they contain the
semidirect factor $\mathbb{Z}$. (When we wanted to investigate
$G_{5.33}^{-1,-1}$, we already had this problem.)

The further restriction in (2) is justified by the large number of
indecomposable non-nilpotent solvable Lie algebras in dimension six:
There are $140$ types of it. The author has decided to consider the
most interesting among them. Since we are not able to refute a
symplectic form's existence in the non-completely solvable case, we
shall partly make even more restrictions.
\subsection{Nilmanifolds}\label{nil6} There are $34$ isomorphism classes of
nilpotent Lie algebras in dimension six. Each of them possesses a
basis with rational structure constants and therefore determines a
nilmanifold. They are listed on page \pageref{6Nilmgf} in Table
\ref{6Nilmgf} which is taken from \cite{Salamon}. The corresponding
Lie algebras are listed in Appendix \ref{Liste Solv}. Among the $34$
classes of nilmanifolds, there are $26$ which admit a symplectic
form.

\begin{table}[p!]
\centering \caption{\label{6Nilmgf} $6$-dimensional nilmanifolds}
\begin{tabular}{|c|c|c|c|} \hline \hline
$b_1(G/\Gamma)$ & $b_2(G/\Gamma)$ & Comment & $\mathfrak{g}$ %
\\ \hline
\hline $6$ & $15$ & Torus, symplectic & $6 \mathfrak{g}_1$ \\
\hline \hline $5$ & $11$ & symplectic & $\mathfrak{g}_{3.1} \oplus
3 \mathfrak{g}_1$ \\
\hline \hline $5$ & $9$ & not symplectic & $\mathfrak{g}_{5.4}
\oplus
\mathfrak{g}_1$ \\
\hline \hline $4$ & $9$ & symplectic & $\mathfrak{g}_{5.1} \oplus
\mathfrak{g}_1$ \\
\hline \hline $4$ & $8$ & symplectic & $2 \mathfrak{g}_{3.1}$ \\
\hline $4$ & $8$ & symplectic & $\mathfrak{g}_{6.N4}$ \\
\hline $4$ & $8$ & symplectic & $\mathfrak{g}_{6.N5}$ \\
\hline \hline $4$ & $7$ & symplectic & $\mathfrak{g}_{5.5} \oplus
\mathfrak{g}_1$ \\
\hline $4$ & $7$ & symplectic & $\mathfrak{g}_{4.1} \oplus
2 \mathfrak{g}_1$ \\
\hline \hline $4$ & $6$ & not symplectic & $\mathfrak{g}_{6.N12}$ \\
\hline \hline $3$ & $8$ & symplectic & $\mathfrak{g}_{6.N3}$ \\
\hline \hline $3$ & $6$ & symplectic & $\mathfrak{g}_{6.N1}$ \\
\hline $3$ & $6$ & symplectic & $\mathfrak{g}_{6.N6}$ \\
\hline $3$ & $6$ & symplectic & $\mathfrak{g}_{6.N7}$ \\
\hline \hline $3$ & $5$ & symplectic & $\mathfrak{g}_{5.2} \oplus
\mathfrak{g}_1$ \\
\hline $3$ & $5$ & not symplectic & $\mathfrak{g}_{5.3} \oplus
\mathfrak{g}_1$ \\
\hline $3$ & $5$ & symplectic & $\mathfrak{g}_{5.6} \oplus
\mathfrak{g}_1$ \\
\hline $3$ & $5$ & symplectic & $\mathfrak{g}_{6.N8}$ \\
\hline $3$ & $5$ & symplectic & $\mathfrak{g}_{6.N9}$ \\
\hline $3$ & $5$ & symplectic & $\mathfrak{g}_{6.N10}$ \\
\hline $3$ & $5$ & not symplectic & $\mathfrak{g}_{6.N13}$ \\
\hline $3$ & $5$ & not symplectic & $\mathfrak{g}_{6.N14}^1$ \\
\hline $3$ & $5$ & not symplectic & $\mathfrak{g}_{6.N14}^{-1}$ \\
\hline $3$ & $5$ & symplectic & $\mathfrak{g}_{6.N15}$ \\
\hline $3$ & $5$ & symplectic & $\mathfrak{g}_{6.N17}$ \\
\hline \hline $3$ & $4$ & symplectic & $\mathfrak{g}_{6.N16}$ \\
\hline \hline $2$ & $4$ & symplectic & $\mathfrak{g}_{6.N11}$ \\
\hline $2$ & $4$ & symplectic & $\mathfrak{g}_{6.N18}^1$ \\
\hline $2$ & $4$ & symplectic & $\mathfrak{g}_{6.N18}^{-1}$ \\
\hline \hline $2$ & $3$ & symplectic & $\mathfrak{g}_{6.N2}$ \\
\hline $2$ & $3$ & symplectic & $\mathfrak{g}_{6.N19}$ \\
\hline $2$ & $3$ & symplectic & $\mathfrak{g}_{6.N20}$ \\
\hline \hline $2$ & $2$ & not symplectic & $\mathfrak{g}_{6.N21}$ \\
\hline $2$ & $2$ & not symplectic & $\mathfrak{g}_{6.N22}$  \\
\hline
\end{tabular}
\end{table}
Recall that a nilmanifold is formal or K\"ahlerian if and only if
the corresponding Lie algebra is abelian.
\subsection{Candidates for the existence of lattices}
\label{Kandidaten} Among the $61$ types of indecomposable unimodular
almost nilpotent Lie algebras in dimension six that are listed in
Tables \ref{UniAlab6} -- \ref{letzteAlgebra}, there are some that
cannot be the Lie algebra of a connected and simply-connected Lie
group which admits a lattice.

Instead of the small German letters for the Lie algebras in the
mentioned tables, we use again capital Latin letters with the same
subscripts for the corresponding connected and simply-connected Lie
groups. If any, we chose the same designation for the parameters
$a,b,c,h,s,\varepsilon$ of $G_{6.i}$ as for their Lie algebras.

\begin{SProp} Let $i \in \{13,\ldots,38\}$, i.e.\ $\mathrm{Nil}(G_{6.i}) =
U_3(\mathbb{R}) \times \mathbb{R}^2$. Then it is necessary for
$G_{6.i}$ to contain a lattice that one of the following conditions
holds:
$$\begin{array}{l@{~~~~}l@{~~~~}l@{~~~~}l}
i = 15, & i=18 \wedge a =-1, & i=21 \wedge a=0,& i=23 \wedge a=0, \\
i=25 \wedge b=0, & i=26, & i=29 \wedge b=0, & i=32 \wedge a=0, \\
i=33 \wedge a=0, & i=34 \wedge a=0, & i=35 \wedge a=-b, & i=36
\wedge a=0, \\ i=37 \wedge a=0, & i=38. & &
\end{array}$$
\end{SProp}

\textit{Proof.} This can be seen analogous as in the proof of
Proposition \ref{nur526}. Denote $\{X_1,\ldots,X_6\}$ the basis used
for the description of the Lie algebra in Appendix \ref{Liste Solv}.
Then the existence of a lattice implies that $[X_6,X_1]$ has no
component in $X_1$-direction and this yields the claim. \q

\begin{SProp}
Let $i \in \{39,\ldots,47\}$, i.e.\ the nilradical of $G_{6.i}$ is
$G_{4.1} \times \mathbb{R}$. If $G_{6.i}$ admits a lattice, then
holds $i = 39 \wedge h = -3$ or $i = 40$.
\end{SProp}

\textit{Proof.} Use the designation $X_1,\ldots,X_6$ as above. Then
$\langle X_1,X_2 \rangle$ is the commutator of the nilradical of
$\mathfrak{g}_{6.i}$. Analogous as in the proof of Proposition
\ref{garnixG41}, one shows that $\mathrm{ad}(X_6)|_{\langle X_1,X_2
\rangle}$ has trace equal to zero. This is only satisfied for $i =
39 \wedge h = -3$ or $i = 40$. \q

\begin{SProp} $\,$
\begin{itemize}
\item[(i)] Let $i \in \{54,\ldots,70\}$, i.e.\ the nilradical of
$G_{6.i}$ is $G_{5.1}$. If $G_{6.i}$ admits a lattice, then holds $i
= 54 \wedge l = -1$, $i = 63$, $i=65 \wedge l=0$ or $i=70 \wedge
p=0$.
\item[(ii)] No connected and simply-connected almost nilpotent Lie
group with nilradical $G_{5.2}$ or $G_{5.5}$ admits a lattice. \q
\end{itemize}
\end{SProp}

\textit{Proof.} This follows in the same manner as the last
proposition. The trace of $\mathrm{ad}(X_6)$ restricted to the
commutator of the nilradical must be zero. \q

\subsection[Symplectic solvmanifolds with $b_1 = 1$]{Symplectic
solvmanifolds whose first Betti number equals one}%
If we are looking for solvmanifolds with $b_1 = 1$, it is necessary
that the corresponding Lie algebra is unimodular, almost nilpotent
and has $b_1 = 1$ itself. Note that the latter forces the algebra to
be indecomposable. In Tables \ref{b16} -- \ref{b16cc} on pages
\pageref{b16} -- \pageref{b16cc} we have listed all possible values
that can arise as $b_1$ for the classes of unimodular non-nilpotent
solvable indecomposable Lie algebras in dimension six.

Since we are mainly interested in symplectic 6-manifolds, we now
investigate which Lie algebras contained in Tables \ref{UniAlab6c}
-- \ref{letzteAlgebra} that satisfy $b_1 =1$ are
\emph{cohomologically symplectic}\index{Cohomologically Symplectic
Lie Algebra}, i.e.\ there is a closed element $\omega \in
\bigwedge^2 \mathfrak{g}^*$ such that $\omega^3$ is not exact.

Note, if a unimodular Lie algebra is cohomologically symplectic,
then each compact quotient of the corresponding Lie group by a
lattice is symplectic. If the Lie algebra is completely solvable,
this is even necessary for the quotient to be symplectic.

\begin{SProp} \label{Sympl61}
Let $\mathfrak{g}_{6.i}$ be a unimodular almost-nilpotent Lie
algebra with $b_1(\mathfrak{g}_{6.i}) = 1$. Then we have:

$\mathfrak{g}_{6.i}$ is cohomologically symplectic if and only if $i
\in \{15,38,78\}$.
\end{SProp}

\textit{Proof.} For $i \in \{15,38,78\}$ one computes all symplectic
forms up to exact summands as
\begin{description}
\item[\,$i=15:$] $\omega = (\lambda + \mu) \, x_{16} + \lambda \, x_{25} - \mu
\, x_{34},~~ \lambda,\mu \in \mathbb{R} \setminus \{0\}, \lambda \ne
- \mu,$
\item[\,$i=38:$] $\omega = \lambda \, x_{16} + \mu \, x_{24} +
\frac{\lambda}{2} \, x_{25} - \frac{\lambda}{2} \, x_{34} + \mu \,
x_{35},~~ \lambda, \mu \in \mathbb{R}, \lambda \ne 0,
-\frac{3}{2}\lambda^3 \ne 2 \lambda \mu^2,$
\item[\,$i=78:$] $\omega = \lambda \, x_{14} + \lambda \, x_{26} + \lambda \,
 x_{35},~~ \lambda \in \mathbb{R} \setminus \{0\}.$
\end{description}
If $i \notin \{15,38,78\}$, then the conditions on the parameters of
$\mathfrak{g}_{6.i}$ to ensure its unimodularity and
$b_1(\mathfrak{g}_{6.i}) = 1$ imply that there are no closed
elements of $\bigwedge^2 \mathfrak{g}_{6.i}^*$ without exact
summands which contain one of the elements $x_{16}, x_{26}, x_{36},
x_{46}$ or $x_{56}$. Therefore, $\mathfrak{g}_{6.i}$ cannot be
cohomologically symplectic. \q

\begin{Rem}
We give an explicit example of the argumentation in the last proof
for $i=2$:

$\mathfrak{g}_{6.2}$ depends on three parameters $a,c,d \in
\mathbb{R}$ with $0 < |d| \le |c| \le 1$ and the brackets are given
in Table \ref{UniAlab6c} as
$$ \begin{array}{l@{~~~~}l@{~~~~}l}
[X_1,X_6] = a \, X_1, &  [X_2,X_6] = X_1 + a \, X_2, & [X_3,X_6] = X_3, \\
{[}X_4,X_6] = c \, X_4, & [X_5,X_6] = d \, X_5. &
\end{array}$$
The condition of unimodularity implies $2a + c +d = -1$. Moreover,
if first the Betti number equals one, we see in Tabular \ref{b16}
that $a \ne 0$.

The Chevalley-Eilenberg complex is given by
$$ \begin{array}{l@{~~~~}l@{~~~~}l}  \delta
x_1 = - a \, x_{16} - x_{26}, & \delta x_2 = -a \, x_{26}, & \delta
x_3 = - x_{36}, \\  \delta x_4 = -c \, x_{46} , & \delta x_5 = -d \,
x_{56}, & \delta x_6 = 0
\end{array} $$
and since $a,c,d \ne 0$, $x_{26}, x_{36}, x_{46}, x_{56}$ are exact.
Moreover, $x_{16} = \delta (- \frac{1}{a} x_1 + \frac{1}{a^2} x_2)$
is exact, too. This implies the claim.
\end{Rem}

We now examine the three Lie groups that have cohomologically
symplectic Lie algebras.

The next theorem was announced in \cite{ipse}. It provides an
example of a symplectic non-formal $6$-manifold with $b_1 = 1$.
Since it is a solvmanifold, this manifold is symplectically
aspherical. Hence, we found an example for which K\c{e}dra, Rudyak
and Tralle looked in \cite[Remark 6.5]{KRT}.

\begin{SThm} $\,$ \label{G615}
\begin{itemize}
\item[(i)] The completely solvable Lie group $G_{6.15}^{-1}$ contains
a lattice.
\item[(ii)] If $\Gamma$ is any lattice in $G:= G_{6.15}^{-1}$, then $M :=
G/\Gamma$ is a symplectic and non-formal manifold with $b_1(M) = 1$
and $b_2(M) = 2$.
\end{itemize}
\end{SThm}

\textit{Proof.} ad (i): Let $N = U_3(\mathbb{R}) \times
\mathbb{R}^2$ denote the nilradical of $G$. We can identify $N$ with
$\mathbb{R}^5$ as a manifold and the multiplication given by
\begin{eqnarray*}
(a,b,c,r,s) \cdot (x,y,z,v,w) & = & (a+x+bz,b+y,c+z,r+v,s+w),
\end{eqnarray*}
i.e.\ $[N,N] = \{(x,0,0,0,0) \,|\, x \in \mathbb{R}\} \cong
\mathbb{R}$ and $\overline{N} := N/[N,N] \cong \mathbb{R}^4$.

By definition of $G$, we have $G = \mathbb{R} \ltimes_{\mu} N$,
where
\begin{eqnarray} \label{expG615}
\forall_{t \in \mathbb{R}} & \mu(t) = \exp^N \circ
\exp^{A(\mathfrak{n})} (t \left( \begin{array}{ccccc}
0 & 0 & 0 & 0 & 0 \\
0 & -1 & 0 & 0 & 0 \\
0 & 0 & 1 & 0 & 0 \\
0 & -1 & 0 & -1 & 0 \\
0 & 0 & -1 & 0 & 1
\end{array} \right) ) \circ
\log^N, &
\end{eqnarray}
and since $\exp^{\mathbb{R}^4} = \mathrm{id}$, the induced maps
$\overline{\mu}(t) \: \overline{N} \to \overline{N}$ are given by
\begin{eqnarray*}
\overline{\mu}(t) \big( (y,z,v,w) \big) &=& \exp^{GL(4,\mathbb{R})}
(t \left(
\begin{array}{cccc}
-1 & 0 & 0 & 0 \\
0 & 1 & 0 & 0 \\
-1 & 0 & -1 & 0 \\
0 & -1 & 0 & 1
\end{array} \right) )
\left( \begin{array}{c} y \\ z \\ v \\ w
\end{array} \right)
\\ &=& \left( \begin{array}{cccc}
e^{-t} & 0 & 0 & 0 \\
0 & e^{t} & 0 & 0 \\
-t e^{-t} & 0 & e^{-t} & 0 \\
0 & -t e^{t} & 0 & e^{t}
\end{array} \right)
\left( \begin{array}{c} y \\ z \\ v \\ w
\end{array} \right).
\end{eqnarray*}
One calculates that $\widetilde{\mu} \: \mathbb{R} \to
\mathrm{A}(N)$ given by
\begin{eqnarray} \label{expG615=}
\forall_{t \in \mathbb{R}} & \forall_{(x,y,z,v,w) \in N} &
\widetilde{\mu}(t) \big((x,y,z,v,w)\big) = \big( x, \overline{\mu}
(t) \big((y,z,v,w) \big) \big)
\end{eqnarray}
is a one-parameter group, and since the derivations of
(\ref{expG615}) and (\ref{expG615=}) in zero are equal, we have $\mu
\equiv \widetilde{\mu}$.

Let $t_1 = \ln(\frac{3 + \sqrt{5}}{2})$, then $\overline{\mu}(t_1)$
is conjugate to $A:= \left(
\begin{array}{cccc}
2 & 1 & 0 & 0 \\
1 & 1 & 0 & 0 \\
2 & 1 & 2 & 1 \\
1 & 1 & 1 & 1
\end{array} \right)$. The transformation matrix $T \in
\mathrm{GL}(4,\mathbb{R})$ with $T A T^{-1} = \overline{\mu}(t_1)$
is
$$T = \left(
\begin{array}{cccc}
1 & -\frac{2 (2 + \sqrt{5})}{3 + \sqrt{5}} & 0 & 0 \\
1 & \frac{1 + \sqrt{5}}{3 + \sqrt{5}} & 0 & 0 \\
0 & 0 & \ln(\frac{2}{3 + \sqrt{5}}) &
\frac{2 ( 2+ \sqrt{5}) \ln(\frac{3 + \sqrt{5}}{2})}{3 + \sqrt{5}} \\
0 & 0 & \ln(\frac{2}{3 + \sqrt{5}}) & - \frac{(1 + \sqrt{5}) \ln(
\frac{3 + \sqrt{5}}{2})}{3 + \sqrt{5}}
\end{array} \right). $$
Denote by $\{b_1,\ldots,b_4\}$ the basis of $\mathbb{R}^4$ for which
$\overline{\mu}(t_1)$ is represented by $A$, i.e.\ $b_i$ is the
$i$-th column of $T$. One calculates
\begin{eqnarray*}
b_{11} b_{22} - b_{12} b_{21} &=&  \sqrt{5}, \\
 b_{i1} b_{j2} - b_{i2} b_{j1} &=&  0 ~~\mbox{ for } i<j, (i,j) \ne
(1,2).
\end{eqnarray*}
This implies that we have for $\gamma_0 :=
(\sqrt{5},0_{\mathbb{R}^4})$, $\gamma_i := (b_{i0}, b_i) \in N$ with
arbitrary $b_{i0} \in \mathbb{R}$, $i=1,\ldots,4$,
$$[\gamma_1,\gamma_2] = \gamma_0,~ [\gamma_1,\gamma_3] =
[\gamma_1,\gamma_4] = [\gamma_2,\gamma_3] = [\gamma_2,\gamma_4] =
[\gamma_3,\gamma_4] = e_N.$$

We can choose the $b_{i0}$ such that the following equations hold:

\parbox{12cm}{$$\begin{array}{ccccccc}
\mu(t_1) (\gamma_0) &=& \gamma_0, &&&&\\
\mu(t_1) (\gamma_1) &=& & \gamma_1^2 & \gamma_2 & \gamma_3^2 & \gamma_4,\\
\mu(t_1) (\gamma_2) &=& & \gamma_1 & \gamma_2 & \gamma_3 & \gamma_4,\\
\mu(t_1) (\gamma_3) &=& &&& \gamma_3^2 & \gamma_4, \\
\mu(t_1) (\gamma_4) &=& &&& \gamma_3 & \, \gamma_4.
\end{array} $$} \hfill \parbox{8mm}{\begin{eqnarray} \label{Mu615}
\end{eqnarray}}
Note that (\ref{Mu615}) leads to the equation
$(\mathrm{id} - \, ^{\tau} \! A) \left( \begin{array}{c} b_{10} \\ b_{20} \\
b_{30} \\ b_{40} \end{array} \right) = \left( \begin{array}{c} 1
+ \frac{2(1 + \sqrt{5})}{3+\sqrt{5}} \\ \frac{1+ \sqrt{5}}{3+\sqrt{5}} \\
0 \\ 0 \end{array} \right)$ which has the (unique) solution $b_{10}
= - \frac{1 + \sqrt{5}}{3+\sqrt{5}}$, $b_{20} = - \frac{11 + 5
\sqrt{5}}{7+ 3 \sqrt{5}}$ and $b_{30} = b_{40} = 0$.

We claim that $t_1 \mathbb{Z} \ltimes_{\mu} \langle \exp^N \big(
\mathrm{Span}_{\mathbb{Z}} \log^N ( \{\gamma_0,\ldots, \gamma_4 \} )
\big) \rangle$ defines a lattice in $G$:

It suffices to show that $\langle \exp^N \big(
\mathrm{Span}_{\mathbb{Z}} \log^N ( \{\gamma_0,\ldots, \gamma_4 \} )
\big) \rangle$ defines a lattice in $N$, so let us prove this
assertion.
There exist uniquely $Y_0, \ldots , Y_4 \in \mathfrak{n}$ with
$\exp^N(Y_i) = \gamma_i$ for $i \in \{0,\ldots,4\}$. If we prove
that $\mathfrak{Y} := \{Y_0,\ldots,Y_4\}$ is a basis of
$\mathfrak{n}$ with rational structure constants, then Theorem
\ref{rat Strukturkonst} (i) implies that $\langle \exp^N (
\mathrm{Span}_{\mathbb{Z}} \mathfrak{Y}) \rangle$ is a lattice in
$N$.

We identify $\mathfrak{n}$ with $\mathbb{R}^5$ and brackets given by
the Campbell-Hausdorff formula, see e.g.\ \cite[Chapter 2.15]{V}.
Since $\mathfrak{n}$ is $2$-step nilpotent (and $\exp^N$ is a
diffeomorphism), the formula yields for all $V,W \in \mathfrak{n}$
$$\log^N \big( \exp^N(V) \exp^N(W) \big) = V + W + \frac{1}{2}
[V,W].$$ Since $U_3(\mathbb{R})$ can be considered as a group of
matrices, one can easily calculate its exponential map. Then, its
knowledge implies that the exponential map resp.\ the logarithm of
$N$ is given by
\begin{eqnarray*}
& \exp^N \big( (x, \, y, \, z, \, v, \, w ) \big) =  ( x +
\frac{1}{2} y z, \, y, \, z, \, v , \, w), & \\ & \log^N \big( ( x,
\, y, \, z, \, v , \, w ) \big) = ( x - \frac{1}{2} y z, \, y, \, z,
\, v, \, w), &
\end{eqnarray*}
and we obtain $ Y_0 = (\sqrt{5}, 0_{\mathbb{R}^4})$, $Y_1 = (b_{10}
- \frac{1}{2}, b_1)$, $Y_2 = (b_{20} +
\frac{(2+\sqrt{5})(1+\sqrt{5})}{(3 + \sqrt{5})^2}, b_2)$, $Y_3 = (0,
b_3)$, $Y_4 = (0, b_4)$, $[Y_1,Y_2] = Y_0$. The other brackets
vanish.

ad (ii): Let $\Gamma$ be an arbitrary lattice in $G$. By completely
solvability and Theorem \ref{min chev eil} (ii), we get the minimal
model of $M = G / \Gamma$ as the minimal model $\mathcal{M}$ of the
Chevalley-Eilenberg complex $(\bigwedge \mathfrak{g}^* , \delta )$.
The latter has the closed generator $x_6$ and the non-closed
generators satisfy
$$\delta x_1 = - x_{23}, \, \delta x_2 = - x_{26},
\, \delta x_3 = x_{36}, \, \delta x_4 = -x_{26} - x_{46} , \, \delta
x_5 = -x_{36} + x_{56},$$ which implies $b_1(M) = 1$.

One computes the differential of the non-exact generators of degree
two in the Chevalley-Eilenberg complex as
$$\begin{array}{l@{~~~~}l@{~~~~}l} \delta x_{12} = x_{126}, & \delta x_{13} =
-x_{136}, & \delta x_{14} = x_{126} + x_{146} - x_{234},  \\ \delta
x_{15} = x_{136} - x_{156} - x_{235}, & \delta x_{16} = x_{236}, &
\delta x_{24} = 2 x_{246}, \\ \delta x_{25} = x_{236}, & \delta
x_{34} = - x_{236}, & \delta x_{35} = -2 x_{356}, \\ \delta x_{45} =
x_{256} - x_{346}, &&
\end{array}$$
i.e.\ $b_2(M) =2$.

The minimal model $\rho \: (\bigwedge V, d) \to (\bigwedge
\mathfrak{g}^* , \delta )$ must contain three closed generators
$y,z_1,z_2$ which map to $x_6,x_{16}+x_{25}$ and  $x_{16}-x_{34}$.
$\rho(y z_1) = x_{256}$ and $\rho(y z_2) = -x_{346}$ are closed and
not exact. But in the generation of $y,z_1$ and $z_2$ is one (and up
to a scalar only one) element that maps onto an exact form, namely
$\rho(y(z_1 +z_2)) = \delta x_{45}$. The minimal model's
construction in the proof of Theorem \ref{Konstruktion d m M}
implies that there is another generator $u$ of degree two such that
$\rho(u) = x_{45}$ and $d u = y (z_1 + z_2)$. Since $\rho(yu) =
x_{456}$ is closed and non-exact, there are no further generators of
degree less than or equal to two in $V$. But this implies for each 
closed element $c$ of degree two that $y\,(u+c)$
is closed and non-exact in $\mathcal{M}$. Using the notation of
Theorem \ref{nichtformal}, we have $u \in N^2, y \in V^1$ and
$\mathcal{M}$ is not formal.

Finally, the existence of a symplectic form on $G/\Gamma$ follows
from Proposition \ref{Sympl61}.\q

\begin{SProp}$\,$
\begin{itemize}
\item[(i)] Each quotient of the Lie group $G^0_{6.38}$ by a
lattice is symplectic. $G^0_{6.38}$ contains a lattice $\Gamma$ with
$b_1(G^0_{6.38}/\Gamma)=1$.
\item[(ii)] If the Lie group $G^0_{6.38}$ contains a lattice $\Gamma$
such that $M := G^0_{6.38}/\Gamma$ satisfies $b_1(M) = 1$ and
$b_2(M) = 2$, then $M$ is a symplectic and non-formal manifold.
\end{itemize}
\end{SProp}

\textit{Proof.} The proof is similar to that of the last theorem.
Therefore, we just give a sketch of the proof.

ad (i): The existence of a symplectic form on each quotient of $G :=
G^0_{6.38}$ by a lattice follows from Proposition \ref{Sympl61}.

The nilradical $N$ of $G$ is the same as in the proof of Theorem
\ref{G615}, so we have $[N,N] = \mathbb{R}$ and $\overline{N} =
N/[N,N] = \mathbb{R}^4$. If $\overline{\mu}(t) \: \overline{N} \to
\overline{N}$ is defined by
\begin{eqnarray*}
\overline{\mu}(t) \big( (y,z,v,w) \big) &=& \exp^{GL(4,\mathbb{R})}
(t \left(
\begin{array}{cccc}
0 & 1 & 0 & 0 \\
-1 & 0 & 0 & 0 \\
-1 & 0 & 0 & 1 \\
0 & -1 & -1 & 0
\end{array} \right) )
\left( \begin{array}{c} y \\ z \\ v \\ w
\end{array} \right)
\\ &=& \left( \begin{array}{cccc}
\cos(t) & \sin(t) & 0 & 0 \\
- \sin(t) & cos(t) & 0 & 0 \\
-t \cos(t) & -t \sin(t) & \cos(t) & sin(t) \\
t \sin(t) & -t \cos(t) & - \sin(t) & \cos(t)
\end{array} \right)
\left( \begin{array}{c} y \\ z \\ v \\ w
\end{array} \right),
\end{eqnarray*} one calculates that
$\mu \: \mathbb{R} \to \mathrm{A}(N)$ given by
\begin{eqnarray*}
\mu(t) \big((x,y,z,v,w)\big) & = & \big( x - \sin^2(t) y z +
\frac{\sin(t) \cos(t)}{2} (z^2 - y^2) + t \frac{\sqrt{3}}{8} (y - z)
, \\ && ~~ \overline{\mu} (t) \big((y,z,v,w) \big) \,\big)
\end{eqnarray*}
is a one-parameter group with $d_e \big(\mu(t)\big) =
\exp^{A(\mathfrak{n})} (t \underbrace{\left(
\begin{array}{ccccc}
0 & 0 & 0 & 0 & 0 \\
0 & -1 & 0 & 0 & 0 \\
0 & 0 & 1 & 0 & 0 \\
0 & -1 & 0 & -1 & 0 \\
0 & 0 & -1 & 0 & 1
\end{array} \right)}_{\D = \mathrm{ad}(X_6)} )$, i.e.\ $G = N \ltimes_{\mu} \mathbb{R}$.
(Here $X_6$ is chosen as in the last line of Table
\ref{UniAlnil16cc} on page \pageref{UniAlnil16cc}.) For $t_1 :=
\frac{\pi}{3}$ we have
\begin{eqnarray*}
&\mu(t_1) \big((x,y,z,v,w)\big) = \big( x - \frac{3}{4} y z +
\frac{\sqrt{3}}{8} (z^2 - y^2) + \frac{\pi}{8 \sqrt{3}} (y - z), \,
\overline{\mu} (t) \big((y,z,v,w) \big) \big), &
\end{eqnarray*}
and in order to construct a lattice in $G$, it is enough to
construct a lattice in $N$ that is preserved by $\mu(t_1)$.
$\overline{\mu}(t_1)$ is conjugate to $A:= \left(
\begin{array}{cccc}
-1 & -3 & 0 & 0 \\
1 & 2 & 0 & 0 \\
-2 & -3 & -1 & -3 \\
1 & 1 & 1 & 2
\end{array} \right)$ and the transformation matrix $T \in
\mathrm{GL}(4,\mathbb{R})$ with $T A T^{-1} = \overline{\mu}(t_1)$
is
$$T = \left(
\begin{array}{cccc}
\frac{\sqrt{3}}{\pi} & 0 & 0 & 0 \\
- \frac{3}{\pi} & - \frac{6}{\pi} & 0 & 0 \\
0 & 0 & - \frac{2}{\sqrt{3}} & - \sqrt{3} \\
0 & 0 & 0 & 1
\end{array} \right). $$
Denote by $\{b_1,\ldots,b_4\}$ the basis of $\mathbb{R}^4$ for which
$\overline{\mu}(t_1)$ is represented by $A$, i.e.\ $b_i$ is the
$i$-th column of $T$. One calculates
\begin{eqnarray*}
b_{11} b_{22} - b_{12} b_{21} &=&  \frac{-6 \sqrt{3}}{\pi^2}, \\
 b_{i1} b_{j2} - b_{i2} b_{j1} &=&  0 ~~\mbox{ for } i<j, (i,j) \ne
(1,2).
\end{eqnarray*}
This implies that we have for $\gamma_0 := ( b_{11} b_{22} - b_{12}
b_{21} ,0_{\mathbb{R}^4})$, $\gamma_i := (b_{i0}, b_i) \in N$ with
arbitrary $b_{i0} \in \mathbb{R}$, $i=1,\ldots,4$,
$$[\gamma_1,\gamma_2] = \gamma_0,~ [\gamma_1,\gamma_3] =
[\gamma_1,\gamma_4] = [\gamma_2,\gamma_3] = [\gamma_2,\gamma_4] =
[\gamma_3,\gamma_4] = e_N.$$

If we set $b_{10} = \frac{1488 \sqrt{3} + 72 \sqrt{3} \pi - 19
\sqrt{3} \pi^2 + 4 \pi^3}{128 \pi^2}$, $b_{20} = \frac{2736 \sqrt{3}
+ 216 \sqrt{3} \pi - 25 \sqrt{3} \pi^2 + 12 \pi^3}{128 \pi^2}$ and
$b_{30} = b_{40} = 0$, we obtain
$$\begin{array}{ccccccc}
\mu(t_1) (\gamma_0) &=& \gamma_0, &&&&\\
\mu(t_1) (\gamma_1) &=& & \gamma_1^{-1} & \gamma_2 & \gamma_3^{-2} & \gamma_4,\\
\mu(t_1) (\gamma_2) &=& & \gamma_1^{-3} & \gamma_2^2 & \gamma_3^{-3} & \gamma_4,\\
\mu(t_1) (\gamma_3) &=& &&& \gamma_3^{-1} & \gamma_4, \\
\mu(t_1) (\gamma_4) &=& &&& \gamma_3^{-3} & \, \gamma_4^{2}.
\end{array} $$
Then $\langle \exp^N \big( \mathrm{Span}_{\mathbb{Z}} \log^N (
\{\gamma_0,\ldots, \gamma_4 \} ) \big) \rangle$ is a lattice in $N$.
This can be seen by a similar computation as in the proof of the
last theorem. Finally, one checks that the abelianisation of this
lattice is isomorphic to $\mathbb{Z}$, hence the corresponding
solvmanifold has $b_1 = 1$.

ad (ii): Let $\Gamma$ be a lattice in $G$ such that
$b_1(G/\Gamma)=1$ and $b_2(G/\Gamma) = 2$.

The Chevalley-Eilenberg complex $(\bigwedge \mathfrak{g}^* , \delta
)$ has the closed generator $x_6$ and $\delta$ is given on the
non-closed generators by
$$\delta x_1 = - x_{23}, \, \delta x_2 = x_{36},
\, \delta x_3 = -x_{26}, \, \delta x_4 = -x_{26} + x_{56} , \,
\delta x_5 = -x_{36} - x_{46},$$ which implies $H^1(\bigwedge
\mathfrak{g}^* , \delta) = \langle [x_6] \rangle$.

One computes the differential of the non-exact generators of degree
two in the Chevalley-Eilenberg complex as
$$\begin{array}{l@{~~~~}l} \delta x_{12} = -x_{136}, & \delta x_{13} =
x_{126}, \\ \delta x_{14} = x_{126} - x_{156} - x_{234},  & \delta
x_{15} = x_{136} + x_{146} - x_{235}, \\ \delta x_{16} = -x_{236}, &
\delta x_{24} = -x_{256} - x_{346}, \\ \delta x_{25} = x_{236} +
x_{246} -x_{356}, & \delta x_{34} = - x_{236} + x_{246} - x_{356},
\\ \delta x_{35} = x_{256} + x_{346}, &
\delta x_{45} = x_{256} - x_{346},
\end{array}$$
i.e.\ $H^2(\bigwedge \mathfrak{g}^* , \delta) = \langle [x_{16} +
\frac{1}{2} x_{25} - \frac{1}{2} x_{34}], [x_{24} + x_{35}]
\rangle$.

This implies that $G/\Gamma$ and $(\bigwedge \mathfrak{g}^* ,
\delta)$ have the same Betti numbers and therefore, by Theorem
\ref{min chev eil}, they share their minimal model.

The minimal model $\rho \: (\bigwedge V, d) \to (\bigwedge
\mathfrak{g}^* , \delta )$ must contain three closed generators
$y,z_1,z_2$ which map to $x_6,x_{16} + \frac{1}{2} x_{25} -
\frac{1}{2} x_{34}$ and  $x_{24} + x_{35}$. $\rho(y z_2) = x_{246} +
x_{356}$ is closed and not exact. But $\rho(y z_1) = \frac{1}{2}
(x_{256} - x_{346}) = \frac{1}{2} \delta x_{45}$ is exact. Hence the
minimal model's construction in the proof of Theorem
\ref{Konstruktion d m M} implies that there is another generator $u$
of degree two such that $\rho(u) = \frac{1}{2} x_{45}$ and $d u = y
z_1$. Since $\rho(yu) = \frac{1}{2} x_{456}$ is closed and
non-exact, there are no further generators of degree less than or
equal to two in $V$. Using the notation of Theorem \ref{nichtformal},
we have $u \in N^2, y \in V^1$, $(u+c) \, y$ is closed and not exact
for each $c \in C^2$ and $(\bigwedge V, d)$ is not formal. 
\q
\begin{SThm}\label{G678} $\,$
\begin{itemize}
\item[(i)] The completely solvable Lie group $G:=G_{6.78}$ possesses
a lattice.
\item[(ii)] For each lattice the corresponding quotient is a
symplectic and formal manifold with $b_1 = b_2 = 1$.
\end{itemize}
\end{SThm}

\textit{Proof.} ad (i): By definition, we have $G = \mathbb{R}
\ltimes_{\mu} N$ with $N = G_{5.3}$ and $\mu(t) = \exp^N \circ
\exp^{A(\mathfrak{n})}(t \, \mathrm{ad}(X_6)) \circ \log^N,$ where
$\{X_1,\ldots,X_6\}$ denotes a basis of $\mathfrak{g}$ as in the
second row of Table \ref{UniAlnil56}. Note that $\{X_1, \ldots
,X_5\}$ is a basis for the nilradical $\mathfrak{n}$. One computes
\begin{equation}\label{mu78}
\mu(t)_* := d_e \big( \mu(t) \big) = \exp^{A(\mathfrak{n})}(t \,
\mathrm{ad}(X_6)) = \left( \begin{array}{ccccc}
e^{t} & 0 & 0 & 0 & 0 \\
0 & 1 & 0 & 0 & 0 \\
0 & 0 & e^{-t} & -t e^{-t} & 0 \\
0 & 0 & 0 & e^{-t} & 0 \\
0 & 0 & 0 & 0 & e^{t}
\end{array} \right).
\end{equation}

Using $$\mathfrak{n} = \langle X_5 \rangle \ltimes_{ad(X_5)} \big(
\langle X_1 \rangle \oplus \langle X_2,X_3,X_4 \, | \, [X_2,X_4] =
X_3 \rangle \big)$$ with $\mathrm{ad}(X_5)(X_2)=-X_1$,
$\mathrm{ad}(X_5)(X_4)=-X_2$,
$\mathrm{ad}(X_5)(X_1)=\mathrm{ad}(X_5)(X_3)=0$ and $$\langle
X_2,X_3,X_4 \, | \, [X_2,X_4] = X_3 \rangle \cong
\mathfrak{g}_{3.1},$$ we can determine the Lie group $N$.

As a smooth manifold $N$ equals $\mathbb{R}^5$, and the
multiplication is given by
\begin{eqnarray*}
\lefteqn{(a, b, c, r, s) \cdot (x, y, z, v, w)} \\ &=& ( a + x + b w
+ \frac{r w^2}{2} \,,\, b + y + r w \,,\, c + z + b v + \frac{r^2
w}{2} + r v w \,,\, r + v \,,\, s + w ) .
\end{eqnarray*}

Now, Theorem \ref{expsemidirekt} enables us to compute the
exponential map of $N$ as
\begin{eqnarray*}
\lefteqn{\exp^N(x X_1 + y X_2 + z X_3 + v X_4 + w X_5)} \\ &=& \big(
x + \frac{y w}{2} + \frac{v w^2}{6} \,,\, y + \frac{v w}{2} \,,\, z
+  \frac{y v}{2} + \frac{v^2 w}{3} \,,\, v \,,\, w \big),
\end{eqnarray*}
and therefore, we also obtain the logarithm of $N$
\begin{eqnarray*}
\lefteqn{\log^N \big((x, y, z, v , w) \big)} \\ &=&  (x - \frac{y
w}{2} + \frac{v w^2}{12}) X_1 + (y - \frac{v w}{2}) X_2 + (z -
\frac{y v}{2} - \frac{v^2 w}{12}) X_3 + v X_4 + w X_5.
\end{eqnarray*}

Finally, a short computation shows that (\ref{mu78}) implies
$$\mu(t) \big( (x,y,z,v,w) \big) = (e^t x, y, e^{-t} (z- tw), e^{-t}
v, e^t w). $$

Let $t_1 := \ln(\frac{3+ \sqrt{5}}{2})$, $b_0 := - \frac{2 t_1}{1+
\sqrt{5}}$ and consider for $t \in \mathbb{R}$ the automorphisms
$I(t) \: N \to N$ given by
\begin{eqnarray*}
\lefteqn{I(t) \big( (x,y,z,v,w) \big)} \\ &=& (0,t b_0,0,0,0)
(x,y,z,v,w) (0,t b_0,0,0,0)^{-1} = (x + t b_0 w, y, z +  t b_0 v, v,
w),
\end{eqnarray*}
and $\nu(t) := \mu(t) \circ I(t) \: N \to N$. It is easy to see that
$\nu \: \mathbb{R} \to \mathrm{A}(N)$ is a one-parameter group in
$N$.

We shall show that there exists a lattice $\Gamma_N$ in $N$
preserved by $\nu(t_1)$, and this then implies the existence of a
lattice in $G_{6.78}$, namely $t_1 \mathbb{Z} \ltimes_{\nu}
\Gamma_N$.

For the remainder of the proof, we identify $\mathfrak{n} \equiv
\mathbb{R}^5$ with respect to the basis $\{X_1,\ldots,X_5\}$ of
$\mathfrak{n}$. Under this identification, consider the basis
$\{Y_1,\ldots,Y_5\}$ of $\mathfrak{n}$, $Y_i$ being the $i$-th
column of $T = (T_{ij}) \in \mathrm{GL}(5,\mathbb{R})$, where $T$
has the following entries:
\begin{align*}
T_{11} &=  \frac{10 (161 + 72 \sqrt{5})
\ln(\frac{3+\sqrt{5}}{2})^2}{1165 + 521 \sqrt{5}}, \\
T_{12} &= 0, \\
T_{13} &= \frac{5 (2 + \sqrt{5}) (161 + 72 \sqrt{5})
\ln(\frac{3+\sqrt{5}}{2})^2}{1525
+ 682 \sqrt{5}}, \\
T_{14} &=  \frac{328380 + 146856 \sqrt{5} - (159975 + 71543
\sqrt{5}) \ln(\frac{3 + \sqrt{5}}{2})^2}{202950 +  90762 \sqrt{5}}, \\
T_{15} &= 1,\\
T_{21} &= 0,\\
T_{22} &= -\frac{(5 + 3 \sqrt{5}) \ln(\frac{3 + \sqrt{5}}{2})}{3 +
\sqrt{5}},\\
T_{23} &= 0,\\
T_{24} &= -\frac{(158114965 + 70711162 \sqrt{5}) \ln(\frac{3 +
\sqrt{5}}{2})}{141422324 + 63245986 \sqrt{5}},\\
T_{25} &= \frac{5 (3940598 + 1762585 \sqrt{5}) \ln(\frac{3 +
\sqrt{5}}{2})}{17622890 + 7881196 \sqrt{5}},\\
T_{31} &= \frac{1}{2} (5 + \sqrt{5}) \ln(\frac{3 + \sqrt{5}}{2}),\\
T_{32} &= 0,\\
T_{33} &= T_{22},\\
T_{34} &= 1,\\
T_{35} &= -\frac{597 + 267 \sqrt{5} + (3808 + 1703 \sqrt{5})
\ln(\frac{3 + \sqrt{5}}{2})}{369 + 165 \sqrt{5}},\\
T_{41} &= 0,\\
T_{42} &= 0,\\
T_{43} &= 0,\\
T_{44} &= 1,\\
T_{45} &= -\frac{2 (2 + \sqrt{5})}{3 + \sqrt{5}},\\
T_{51} &= 0,\\
T_{52} &= 0,\\
T_{53} &= 0,\\
T_{54} &= \ln(\frac{2}{3+\sqrt{5}}),\\
T_{55} &= -\frac{2 \ln(\frac{3 + \sqrt{5}}{2})}{1 + \sqrt{5}}.
\end{align*}
Let $\gamma_i := \exp^N(Y_i)$ for $i \in \{1, \ldots 5\}$ and
\begin{align*}
S_1 &= \frac{92880525355200 +  41537433696024
\sqrt{5}}{57403321562460 + 25671545829588 \sqrt{5}} \\
&~~~~ - \frac{(3591421616495 + 1606132574069 \sqrt{5}) \ln(\frac{3 +
\sqrt{5}}{2})^2}{57403321562460 + 25671545829588 \sqrt{5}},\\
S_2 &= -\frac{(228826127 + 102334155 \sqrt{5}) \ln(\frac{3 +
\sqrt{5}}{2})}{141422324 + 63245986 \sqrt{5}}, \\
S_3 &= 1 - \frac{(757189543 + 338625458 \sqrt{5}) \ln(\frac{3 +
\sqrt{5}}{2})}{848533944 + 379475916 \sqrt{5}}, \\
S_4 &= \frac{724734510 + 324111126 \sqrt{5} - (325041375 + 145362922
\sqrt{5}) \ln(\frac{3 + \sqrt{5}}{2})^2}{724734510 + 324111126 \sqrt{5}},& \\
S_5 &= \frac{(120789085 + 54018521 \sqrt{5}) \ln(\frac{3 +
\sqrt{5}}{2})}{74651760 + 33385282 \sqrt{5}}, \\
S_6 &= -\frac{466724522940 + 208725552012 \sqrt{5}}{24
(12018817440 + 5374978561 \sqrt{5})} \\
&~~~~ + \frac{(3393446021605 + 1517595196457 \sqrt{5}) \ln(\frac{3 +
\sqrt{5}}{2})}{24 (12018817440 + 5374978561 \sqrt{5})}.
\end{align*}
One computes $\gamma_1 = (T_{11},0,T_{31},0,0)$, $\gamma_2 =
(0,T_{22},0,0,0)$, $\gamma_3 = (T_{13},0,T_{33},0,0)$, $\gamma_4 =
(S_{1},S_{2},S_{3},T_{44},T_{54})$ and $\gamma_5 =
(S_{4},S_{5},S_{6},T_{45},T_{55})$.

Moreover, if $A$ denotes the matrix $\left(
\begin{array}{ccccc}
1 & 0 & 1 & \frac{13}{6} & \frac{11}{6} \\
0 & 1 & 0 & 0 & -\frac{1}{2} \\
1 & 0 & 2 & -\frac{5}{6} & -\frac{1}{3} \\
0 & 0 & 0 & 2 & 1 \\
0 & 0 & 0 & 1 & 1
\end{array} \right)$, we can calculate  $T A T^{-1} = \nu(t_1)_* :=
d_e(\nu(t))$. Since $\nu(t_1) = \exp^N \circ \nu(t_1)_* \circ
\log^N$, this yields  $\nu(t_1)(\gamma_1) = \gamma_1 \, \gamma_3$,
$\nu(t_1)(\gamma_2) = \gamma_2$, $\nu(t_1)(\gamma_3) = \gamma_1 \,
\gamma_3^2$, $\nu(t_1)(\gamma_4) = \gamma_1^2 \, \gamma_2^{-2} \,
\gamma_4^2 \, \gamma_5$ and $\nu(t_1)(\gamma_5) = \gamma_1^2 \,
\gamma_2^{-1} \, \gamma_{4} \, \gamma_5$.

Therefore, we have shown that $\nu(t_1)$ preserves the subgroup
$\Gamma_N$ of $N$ which is generated by $\gamma_1, \ldots ,
\gamma_5$. In order to complete the proof of (i), it suffices to
show that $\Gamma_N$ is a lattice in $N$.

Since $\mathfrak{n}$ is $3$-step nilpotent, the
Baker-Campbell-Hausdorff formula (see e.g.\ \cite[Chapter 2.15]{V})
yields for all $V,W \in \mathfrak{n}$
$$\log^N \big( \exp^N(V) \exp^N(W)
\big) = V + W + \frac{1}{2} [V,W] + \frac{1}{12} ([[V,W],W] -
[[V,W],V]). $$ Therefore, we obtain by a short calculation
$[Y_2,Y_4] = Y_3$, $[Y_2 , Y_5] = Y_1$ and $[Y_4,Y_5] = \frac{1}{2}
Y_1 + Y_2 + \frac{1}{2} Y_3$, i.e.\ the basis $\{Y_1,\ldots,Y_5\}$
has rational structure constants. Theorem \ref{rat Strukturkonst}
then implies that $\Gamma_N$ is a lattice in $N$.

ad (ii): Let $\Gamma$ be a lattice in $G := G_{6.78}$. By completely
solvability and Theorem \ref{min chev eil} (ii), the minimal model
of $M = G / \Gamma$ is the same as the minimal model $\mathcal{M}$
of the Chevalley-Eilenberg complex $(\bigwedge \mathfrak{g}^* ,
\delta )$. In view of Theorem \ref{formal = n-1 formal}, it suffices
to prove that the latter is $2$-formal. On the non-closed generators
of $(\bigwedge \mathfrak{g}^* , \delta )$ the differential is  given
by
$$\delta x_1 = x_{16} - x_{25}, \, \delta x_2 = -x_{45},
\, \delta x_3 = -x_{24} -x_{36} -x_{46}, \, \delta x_4 = -x_{46}, \,
\delta x_5 = x_{56},$$ i.e.\ $H^1(\bigwedge \mathfrak{g}^* , \delta
) = \langle [x_6] \rangle$. Further, one calculates $H^2(\bigwedge
\mathfrak{g}^* , \delta ) = \langle x_{14} + x_{26} + \nolinebreak
x_{35} \rangle$. The minimal model $\rho \: (\bigwedge V, d) \to
(\bigwedge \mathfrak{g}^* , \delta )$ then must contain two closed
generators $y,z$ which map to $x_6$ and  $x_{14} + x_{26} + x_{35}$.
Since $\rho(y z) = x_{146} + x_{356}$ is closed and non-exact, there
are no other generators of degree two in $(\bigwedge V, d)$, hence
up to degree two, all generators are closed. This implies the
minimal model's $2$-formality.

Moreover, $x_{14} + x_{26} + x_{35}$ defines a symplectic form. \q

\begin{Rem}
In order to determine a lattice in $G_{6.78}$, the author also found
a lattice of the completely solvable Lie group $G_{6.76}^{-1}$. One
can show that the corresponding solvmanifold is formal and has first
Betti number equal to one. Unfortunately, it is not symplectic by
Proposition \ref{Sympl61}.
\end{Rem}
\begin{Rem}
Besides the mentioned groups above, the following non-completely
solvable Lie groups $G_{6.i}$ could give rise to a symplectic
solvmanifold $G_{6.i} / \Gamma$ with $b_1(G_{6.i} / \Gamma) = 1$:
\begin{equation}\label{nvavs} \begin{array}{l@{~~~}l@{~~~}l} i=8; & i=9,~ b\ne0;
& i = 10 ,~ a\ne0; \\ i=11;  & i=12; & i = 32,~ a = \varepsilon = 0 < h; \\
i= 37,~ a = 0; & i= 88,~ \mu_0 \nu_0 \ne 0; & i=89,~ \nu_0 s \ne 0; \\
i=90,~ \nu_0 \ne 0; & i=92,~ \mu_0 \nu_0 \ne  0; & i=92^*; \\
i=93,~ |\nu_0|>\frac{1}{2}. &&
\end{array}\end{equation}
But then the cohomology class of the symplectic form cannot lie in
the image of the inclusion $H^* (\bigwedge \mathfrak{g}^*_{6.i},
\delta)\hookrightarrow H^*(G_{6.i}/\Gamma,d)$ by Proposition \ref{Sympl61}.
\end{Rem}
\subsection[Symplectic solvmanifolds with $b_1 > 1$]{Symplectic
solvmanifolds whose first Betti number is greater than
one}\label{6Sgr1}
In this section, we examine which Lie groups $G$ can give rise to a
six-dimensional solvmanifold $G/\Gamma$ with $b_1(G/\Gamma) >1$.
Again, we just consider indecomposable connected and
simply-connected solvable Lie groups. The nilradical of such a group
has not dimension equal to three, see e.g.\ \cite{Mub63b}.
Proposition \ref{dim Nil} then tells us that indecomposable solvable
Lie groups have nilradicals of dimension greater than three.
Moreover, the nilpotent ones were considered in Section \ref{nil6},
hence we can assume that $G$ is non-nilpotent, i.e.\ $\dim
\mathrm{Nil}(G) \in \{4,5\}$. The corresponding Lie algebras are
listed in Tables \ref{UniAlab6} -- \ref{Nilrad42c}.

In Section \ref{Kandidaten}, we have excluded some groups $G$ since
they cannot admit lattices. Clearly, we omit them in the
considerations below.

By Theorem \ref{min chev eil}(ii), we have in the completely
solvable case an isomorphism from Lie algebra cohomology to the
solvmanifold's cohomology, i.e.\ the Lie algebra $\mathfrak{g}$ must
satisfy $b_1(\mathfrak{g}) > 1$, too. In the last section, we saw
that $\mathfrak{g}_{6.38}^0$ is the only non-completely solvable but
cohomologically symplectic Lie algebra with $b_1(\mathfrak{g}) = 1$.
Therefore, for each lattice $\Gamma$ in $G_{6.38}^0$ with
$b_1(G_{6.38}^0 / \Gamma) >1$, the quotient is symplectic. We now
turn to Lie algebras with $b_1(\mathfrak{g}) > 1$. The possible
values of $b_1$ can be read of in Tables \ref{b16} -- \ref{b16cc}.

The remaining algebras to examine are $\mathfrak{g}_{6.i}$ with
\begin{equation}\label{uebrig}
\begin{array}{l@{~~~}l@{~~~}l} i=2,~ a
=0; & i=3,~ d=-1; & i = 6 ,~ a=-\frac{1}{2},~ b=0; \\
i= 9,~ b=0;  & i=10,~ a=0; & i = 21,~ a=0; \\ i= 23,~ a = 0; & i=
25,~ b=0; & i=26; \\ i=29,~ b=0; & i=33,~ a=0; & i = 34,~ a=0; \\
i=36,~ a =0; & i=54,~ l=-1; & i =63; \\ i=65,~ l=0; & i= 70,~ p=0; &
i = 83,~l=0; \\ i=84; & i=88, ~ \mu_0 = \nu_0 = 0; & i=89,~ \nu_0 s
= 0; \\ i=90,~ \nu_0 =0; & i=92,~ \nu_0 \mu_0 = 0; & i=93,~ \nu_0
=0; \\ i= 102; & i=105; & i=107; \\ i=113; & i=114; & i =115; \\
i=116;& i=118; & i=120; \\ i=125; & i=129; & i =135.
\end{array}
\end{equation}
As above, we just consider such Lie algebras that are
cohomologically symplectic, although this condition is only in the
completely solvable case necessary for the existence of a symplectic
form on $G / \Gamma$.

\begin{SProp}\label{Symplgg1}
Let $\mathfrak{g}_{6.i}$ be one of the Lie algebras listed in
(\ref{uebrig}).

Then $\mathfrak{g}_{6.i}$ is cohomologically symplectic if and only
if it is contained in the following list:

$$\begin{array}{l@{~~~}l@{~~~}l@{~~~}l} b_1=2: & i=3 ~ \wedge ~ d
=-1, & i=10 ~ \wedge ~ a=0, & i = 21 ~ \wedge ~ a=0, \\ & i= 36 ~
\wedge ~ a =0,  & i=54 ~ \wedge ~ l=-1, \\ & i = 70  ~ \wedge ~ p=0,
& i= 118 ~ \wedge ~ b=\pm1. & \\ b_1=3: & i= 23  ~ \wedge ~ a = 0~
\wedge~ \varepsilon \ne 0, & i=29 ~ \wedge ~ b=0. &
\end{array}$$
\end{SProp}

\textit{Proof.} This is done by a case by case analysis as described
in the proof of Proposition \ref{Sympl61}. We list the symplectic
forms for the Lie algebras that are cohomologically symplectic. In
the cases with $b_1=2$, the symplectic forms are given by
$$\begin{array}{l@{:~~}ll}
i = 3,~d=-1 & \omega = \lambda \, x_{16} + \mu \, x_{23} +
\nu \, x_{45},& \lambda \mu \nu \ne 0, \\
i = 10,~a=0 & \omega = \lambda \, x_{16} + \mu \, x_{23} +
\nu \, x_{45},& \lambda \mu \nu \ne 0, \\
i = 21,~a=0 & \omega = \lambda \, x_{12} + \mu \, x_{36} +
\nu \, x_{45},& \lambda \mu \nu \ne 0, \\
i = 36,~a=0 & \omega = \lambda \, x_{12} + \mu \, x_{36} +
\nu \, x_{45},& \lambda \mu \nu \ne 0, \\
i = 54,~l=-1 & \omega = \lambda \, (x_{12} + x_{23}) + \mu
\, x_{34} + \nu \, x_{56},& \lambda \nu \ne 0, \\
i = 70,~p=0 & \omega = \lambda \, (x_{13} + x_{24}) + \mu
\, x_{34} + \nu \, x_{56},& \lambda \nu \ne 0, \\
i = 118,~b=\pm1 & \omega = \lambda \, (x_{13} \pm x_{24}) + \mu
\, (x_{14} - x_{23}) + \nu \, x_{56},& (\lambda^2 + \mu^2) \nu \ne 0.\\
\end{array}$$
In the cases with $b_1 =3$, we have the symplectic forms
$$ \omega = \lambda \, (x_{12} +
\varepsilon \, x_{35})  + \mu \, (x_{16} + x_{24}) + \nu \, (x_{23}
- \varepsilon \, x_{56}) + \rho \, x_{25} + \sigma \, x_{46}$$ with
$\lambda \mu \nu \ne 0$ for $i = 23,~a=0,~ \varepsilon \ne 0$,
$$ \omega = \lambda \, (x_{13} +
\varepsilon \, x_{45})  + \mu \, (x_{16} + x_{24}) + \nu \, (x_{23}
- \varepsilon \, x_{56}) + \rho \, x_{26} + \sigma \, x_{34}$$ with
$\lambda \ne 0, \rho \ne \frac{(\lambda + \varepsilon) \mu
\nu}{\lambda}$ for $i = 29,~b=0,~ \varepsilon \ne 0$ and
$$ \omega = \lambda \, x_{12} + \mu \, x_{13} + \nu \, (x_{16} +
x_{24}) + \rho \, x_{26} + \sigma \, x_{34} + \tau \, x_{56}$$ with
$\nu(\nu \sigma + \mu \tau) \ne 0$ for $i = 29,~b=0,~ \varepsilon
\ne 0$. \q
\N Provided there is a lattice in one of the ten Lie groups
$G_{6.i}$ in the last proposition whose Lie algebras are
cohomologically symplectic, we can ensure that the corresponding
solvmanifold is symplectic. In the completely solvable case, i.e.\
$i \in \{3,21,23,29,54\}$, we can determine cohomological properties
of the potential solvmanifolds.

\begin{SProp}$\,$ \label{G63}
\begin{itemize}
\item[(i)] There is a lattice in the completely solvable Lie group
$G_{6.3}^{0,-1}$.
\item[(ii)] For each lattice the corresponding solvmanifold is symplectic, not
formal and satisfies $b_1= 2$ as well as $b_2 = 3$.
\end{itemize}
\end{SProp}

\textit{Proof.} ad (i) : We have $G := G_{6.3}^{0,-1} = \mathbb{R}
\ltimes_{\mu} \mathbb{R}^4$ with $\mu(t) = \exp^{GL(4,\mathbb{R})}(t
\, \mathrm{ad}(X_6))$, where $X_6 \in \mathfrak{g}_{6.3}^{0,-1}$ is
chosen as in Table \ref{UniAlab6c}, i.e.\ $$\mu(t) = \left(
\begin{array}{ccccc}
1 & -t & \frac{t^2}{2} & 0 & 0 \\
0 & 1 & -t & 0 & 0 \\
0 & 0 & 1 & 0 & 0 \\
0 & 0 & 0 & e^{-t} & 0 \\
0 & 0 & 0 & 0 & e^{t}
\end{array} \right).$$
Set $t_1 := \ln(\frac{3+\sqrt{5}}{2})$. Then $\mu(t_1)$ is conjugate
to $\left(
\begin{array}{ccccc}
1 & 1 & 0 & 0 & 0 \\
0 & 1 & 1 & 0 & 0 \\
0 & 0 & 1 & 0 & 0 \\
0 & 0 & 0 & 0 & -1 \\
0 & 0 & 0 & 1 & 3
\end{array} \right)$. This follows from (\ref{expt1}) and the fact
that the Jordan form of the upper left block of $\mu(t_1)$ is
$\left(
\begin{array}{ccc}
1 & 1 & 0  \\
0 & 1 & 1  \\
0 & 0 & 1
\end{array} \right)$. Hence $G$ admits a lattice.

ad (ii): By completely solvability and Theorem \ref{min chev eil}
(ii), the solvmanifold's minimal model is the same as the minimal
model of the Chevalley-Eilenberg complex $(\bigwedge
(x_1,\ldots,x_6) , \delta )$. In view of Theorem \ref{formal = n-1
formal}, it suffices to prove that the latter is not $2$-formal.

Using the knowledge of the Chevalley-Eilenberg complex, one can
compute $\rho \: (\bigwedge (y_1,\ldots,y_4,z) , d ) \to
(\bigwedge(x_1,\ldots,x_6),\delta)$ as the minimal model up to
generators of degree two, where
\begin{eqnarray*}
&\rho(y_1) = x_3,~ \rho(y_2) = x_6,~ \rho (y_3) = -x_2,~ \rho(y_4) =
-x_1,~ \rho(z) = x_4 x_5, & \\ &d y_1 = d y_2 = 0 ,~ d y_3 = y_1
y_2,~ d y_4 = y_2 y_3 ,~ d z = 0.&
\end{eqnarray*}
This obviously implies the statement about the Betti numbers.
Moreover, using the notation of Theorem \ref{nichtformal}, we have
$C^1= \langle y_1, y_2 \rangle$, $N^1 = \langle y_3, y_4 \rangle$,
and $y_1 \, (y_3+c)$ is closed but not exact for each $c \in C^1$.
Hence the minimal model is not $1$-formal. \q

\begin{SProp}$\,$ \label{G621}
\begin{itemize}
\item[(i)] There is a lattice in the completely solvable Lie group
$G_{6.21}^{0}$.
\item[(ii)] For each lattice the corresponding solvmanifold is symplectic,
not formal and satisfies $b_1= 2$ as well as $b_2 = 3$.
\end{itemize}
\end{SProp}

\textit{Proof.} The proof of (ii) is completely analogous to that of
(ii) in the last proposition. But this time, the minimal model is
given by
\begin{eqnarray*}
&\rho(y_1) = x_2,~ \rho(y_2) = x_6,~ \rho (y_3) = -x_3,~ \rho(y_4) =
x_1,~ \rho(z) = x_4 x_5, & \\ &d y_1 = d y_2 = 0 ,~ d y_3 = y_1
y_2,~ d y_4 = y_1 y_3 ,~ d z = 0,&
\end{eqnarray*}
and $y_1 \, (y_4+c)$ is closed but non-exact for each closed $c$ 
of degree one.

ad (i): In order to prove the existence of a lattice, we use the
same argumentation as in the proof of Theorem \ref{G615} (i). (Note
that $G_{6.15}^{-1}$ and $G := G_{6.21}^{0}$ share their nilradical
$N$.) But of course, we now have a different initial data: $G =
\mathbb{R} \ltimes_{\mu} N$ with $\mu(t) = \exp^N \circ
\exp^{A(\mathfrak{n})} (t \left(
\begin{array}{ccccc}
0 & 0 & 0 & 0 & 0 \\
0 & 0 & 0 & 0 & 0 \\
0 & -1 & 0 & 0 & 0 \\
0 & 0 & 0 & -1 & 0 \\
0 & 0 & 0 & 0 & 1
\end{array} \right) ) \circ
\log^N$ and \begin{eqnarray*} \overline{\mu}(t) \big( (y,z,v,w)
\big) &=& \exp^{GL(4,\mathbb{R})} (t \left(
\begin{array}{cccc}
0 & 0 & 0 & 0 \\
-t & 0 & 0 & 0 \\
0 & 0 & -t & 0 \\
0 & 0 & 0 & t
\end{array} \right) )
\left( \begin{array}{c} y \\ z \\ v \\ w
\end{array} \right)
\\ &=& \left( \begin{array}{cccc}
1 & 0 & 0 & 0 \\
-t & 1 & 0 & 0 \\
0 & 0 & e^{-t} & 0 \\
0 & 0 & 0 & e^{t}
\end{array} \right)
\left( \begin{array}{c} y \\ z \\ v \\ w
\end{array} \right).
\end{eqnarray*}
Arguing analogous as in (\ref{expG615=}), one obtains $$ \mu(t)
\big( (x,y,z,v,w) \big) = \big (x-\frac{t}{2}y^2 , \overline{\mu}(t)
\big( (x,y,z,v,w) \big) \big).$$ Let $t_1 =
\ln(\frac{3+\sqrt{5}}{2})$, $A := \left(
\begin{array}{cccc}
1 & 1 & 0 & 0 \\
0 & 1 & 0 & 0 \\
0 & 0 & 0 & -1 \\
0 & 0 & 1 & 3
\end{array} \right)$ and $T = \left(
\begin{array}{cccc}
0 & -\frac{1}{t_1} & 0 & 0 \\
1 & 0 & 0 & 0 \\
0 & 0 & \frac{18 + 8\sqrt{5}}{7 + 3\sqrt{5}} & 1 \\
0 & 0 & \frac{2}{3 + \sqrt{5}} & 1
\end{array} \right)$. Then we have $T A T^{-1} = \overline{\mu}(t_1)$.
Denote the $i$-th column of T by $b_i$. Analogous calculations as in
loc.\ cit.\ imply the existence of a lattice generated by $\gamma_0
:= \nolinebreak (\frac{1}{t_1},0_{\mathbb{R}^4})$ and $\gamma_i
\nolinebreak := \nolinebreak (b_{i0}, b_i)$, $i \in \{1,\ldots,4\}$,
where $b_{20} \in \mathbb{R}$ arbitrary and $b_{10}= - \frac{1}{2
t_1}$ as well as $b_{30}=b_{40} = 0$. \q

\begin{SProp} $\,$
\begin{itemize}
\item[(i)] Let $\varepsilon = \pm 1$. There is a lattice in the completely
solvable Lie group $G_{6.23}^{0,0,\varepsilon}$.
\item[(ii)] If there is a lattice in $G_{6.23}^{0,0,\varepsilon}, \, \varepsilon \ne 0$,
then the corresponding solvmanifold is symplectic, non-formal and
satisfies $b_1= 3$ as well as $b_2 = 5$.
\end{itemize}
\end{SProp}

\textit{Proof.} ad (i): $G_{6.23}^{0,0,\varepsilon}$ has the same
nilradical $N$ as $G_{6.15}^{-1}$ and the latter is described at the
beginning of the proof of Theorem \ref{G615}.

By definition, $G_{6.23}^{0,0,\varepsilon} = \mathbb{R}
\ltimes_{\mu} N$ with
\begin{eqnarray*}
\mu(t) & = & \exp^N \circ \exp^{A(\mathfrak{n})} (t \left(
\begin{array}{ccccc}
0 & 0 & 0 & 0 & -\varepsilon \\
0 & 0 & 0 & 0 & 0 \\
0 & -1 & 0 & 0 & 0 \\
0 & 0 & -1 & 0 & 0 \\
0 & 0 & 0 & 0 & 0
\end{array} \right) ) \circ
\log^N.
\end{eqnarray*}
The functions $\exp^N, \log^N$ also can be found in the proof of
Theorem \ref{G615}. Using their knowledge, we calculate
$$ \mu(t)\big( (x,y,z,v,w) \big) = \big( x - \frac{t}{2} y^2 - t
\varepsilon \, , \,  y \, , \, z - ty \, , \, \frac{t^2}{2} y - t z
+ v \, , \, w \big).
$$ If $\varepsilon = \pm 1$, then the map $\mu(2)$ preserves the
lattice
$$\{ (x,y,z,v,w) \in N \, | \, x,y,z,v,w \in \mathbb{Z} \} \subset N.$$
Therefore, $G_{6.23}^{0,0,\varepsilon}$ admits a lattice.

ad (ii): By completely solvability, the Betti numbers of the
Chevalley-Eilenberg complex coincide with the solvmanifold's Betti
numbers. A short calculation yields the first Betti numbers of the
former as $b_1 = 3$ and $b_2 =5$.

As above, the knowledge of the Chevalley-Eilenberg complex
$(\bigwedge(x_1,\ldots,x_6),\delta)$ enables us to compute the first
stage of the minimal model. It is given by $\rho \: (\bigwedge
(y_1,\ldots,y_6) , d ) \to (\bigwedge(x_1,\ldots,x_6),\delta)$ with
\begin{eqnarray*}
&\rho(y_1) = x_2,~ \rho(y_2) = x_5,~ \rho (y_3) = x_6,~ \rho(y_4) =
-x_3,~ \rho(x_5) = x_1,~ \rho(y_6) = -x_4, & \\ &d y_1 = d y_2 =
dy_3 = 0 ,~ d y_4 = y_1 y_3,~ d y_5 = y_1 y_4 - \varepsilon \, y_2
y_3 ,~ d y_6 = y_3 y_4.&
\end{eqnarray*} Since $y_3\,(y_6+c)$ is closed and non-exact for each closed $c$ 
of degree one, the minimal model is not $1$-formal. \q

\begin{SProp} $\,$
\begin{itemize}
\item[(i)] Let $\varepsilon \in \{0,\pm 1\}$. There is a lattice in the
completely solvable Lie group $G_{6.29}^{0,0,\varepsilon}$.
\item[(ii)] If there is a lattice in $G_{6.29}^{0,0,\varepsilon}, \,
\varepsilon \in \mathbb{R}$, then the corresponding solvmanifold is
symplectic, non-formal and has $b_1= 3$ as well as $b_2 = \left\{
\begin{array}{c} 5,~ \mbox{if } \varepsilon \ne  0 \\ 6,~ \mbox{if }
\varepsilon = 0 \end{array} \right\}$.
\end{itemize}
\end{SProp}

\textit{Proof.} The argumentation is analogous to the last proof,
but this time we have
$$ \mu(t)\big( (x,y,z,v,w) \big) = \big( x - \frac{\varepsilon}{6} t^3 z
+ \frac{\varepsilon}{2} t^2 v - \varepsilon t w \, , \,  y \, , \, z
\, , \, - t z + v \, , \, \frac{1}{2} t^2 z -tv +w  \big).
$$
(Note that there is no misprint. The maps $\exp^N \circ
\exp^{A(\mathfrak{n})}(t \, \mathrm{ad}(X_6)) \circ \log^N$ and
$\exp^{A(\mathfrak{n})}(t \, \mathrm{ad}(X_6))$ have the same form.)
For $\varepsilon \in \{0,\pm1\}$, $\mu(6)$ preserves the integer
lattice mentioned in the last proof. This implies (i).

In order to prove (ii), we consider the minimal model. Up to
generators of degree one, it is given by
\begin{eqnarray*}
&\rho(y_1) = x_2,~ \rho(y_2) = x_3,~ \rho (y_3) = x_6,~ \rho(y_4) =
-x_4,~ \rho(x_5) = -x_5,~ \rho(y_6) = -x_1, & \\ &d y_1 = d y_2 =
dy_3 = 0 ,~ d y_4 = y_2 y_3,~ d y_5 = y_3 y_4 ,~ d y_6 = y_1 y_2 +
\varepsilon \, y_3 y_5, &
\end{eqnarray*}
if $\varepsilon \ne 0$, and
\begin{eqnarray*}
&\rho(y_1) = x_2,~ \rho(y_2) = x_3,~ \rho (y_3) = x_6,~ \rho(y_4) =
-x_1,~ \rho(x_5) = -x_4,~ \rho(y_6) = -x_5, & \\ &d y_1 = d y_2 =
dy_3 = 0 ,~ d y_4 = y_1 y_2,~ d y_5 = y_2 y_3 ,~ d y_6 = y_3 y_5,&
\end{eqnarray*}if $\varepsilon = 0$. In both cases $y_2 \,(y_4+c)$ is closed
and non-exact for all closed $c$ of degree one. \q
\N The following result is due to Fern\'andez, de L\'eon and
Saralegui.  Its proof can be found in \cite[Section 3]{FLS}. Note
that the cohomological results are independent of the choice of the
lattice, since the Lie group in the proposition is completely
solvable.

\begin{SProp} \label{G654}
The completely solvable Lie group $G_{6.54}^{0,-1}$ admits a
lattice. For each such, the corresponding solvmanifold is
symplectic, non-formal and satisfies $b_1=2$ as well as $b_2=5$. \q
\end{SProp}

Summing up the results concerning completely solvable Lie groups
that admit symplectic quotients, we obtain:

\begin{SThm}
All six-dimensional symplectic solvmanifolds that can be written as
quotient of a non-nilpotent completely solvable indecomposable Lie
group are contained in one of the last five propositions, Theorem
\ref{G615} or Theorem \ref{G678}. \q
\end{SThm}

To end this section, we consider the four cohomologically symplectic
Lie algebras $\mathfrak{g}_{6.i}$ of Proposition \ref{Symplgg1} that
are not completely solvable, this means $i=10 \wedge a=0$, $i=36
\wedge a=0$, $i=70 \wedge p=0$ or $i=118 \wedge b= \pm 1$. Clearly,
the existence of a lattice implies that the corresponding
solvmanifold is symplectic. But in order to make a statement about
cohomological properties, one needs an assumption about the first
two Betti numbers to ensure the knowledge of the cohomology algebra.

\begin{SProp}$\,$
\begin{itemize}
\item[(i)] Each quotient of the Lie group $G :=
G_{6.10}^{0,0}$ by a lattice is symplectic and $G$ admits a lattice
$\Gamma$ with $b_1(G/\Gamma) = 2$.
\item[(ii)] If there is a lattice in $G$ such that the corresponding
solvmanifold satisfies $b_1= 2$ and $b_2= 3$, then it is symplectic
and not formal.
\end{itemize}
\end{SProp}

\textit{Proof.} We have $G = \mathbb{R} \ltimes_{\mu} \mathbb{R}^4$
with $\mu(t) = \exp^{GL(4,\mathbb{R})}(t \, \mathrm{ad}(X_6))$ and
$X_6 \in \mathfrak{g}_{6.10}^{0,0}$ chosen as in Table
\ref{UniAlab6c}, i.e.\
$$\mu(t) = \left(
\begin{array}{ccccc}
1 & -t & \frac{t^2}{2} & 0 & 0 \\
0 & 1 & -t & 0 & 0 \\
0 & 0 & 1 & 0 & 0 \\
0 & 0 & 0 & \cos(t) & - \sin(t) \\
0 & 0 & 0 & \sin(t) & \cos(t)
\end{array} \right).$$
$\mu(\pi)$ is conjugate to $\left(
\begin{array}{ccccc}
1 & 1 & 0 & 0 & 0 \\
0 & 1 & 1 & 0 & 0 \\
0 & 0 & 1 & 0 & 0 \\
0 & 0 & 0 & -1 & 0 \\
0 & 0 & 0 & 0 & -1
\end{array} \right)$. This follows from the fact that the Jordan form
of the upper left block of $\mu(\pi)$ is $\left(
\begin{array}{ccc}
1 & 1 & 0  \\
0 & 1 & 1  \\
0 & 0 & 1
\end{array} \right)$. Hence $G$ admits a lattice $\Gamma$.

A short calculation yields that the abelianisation of this lattice
is isomorphic to $\mathbb{Z}^2 \oplus \mathbb{Z}_2{}^2$, i.e.\
$b_1(G/\Gamma) = 2$.

Using the assumptions of (ii), one calculates the minimal model up
to generators of degree one as
\begin{eqnarray*}
&\rho(y_1) = x_3,~ \rho(y_2) = x_6,~ \rho (y_3) = -x_2,~ \rho(y_4) =
-x_1,& \\ &d y_1 = d y_2 = 0 ,~ d y_3 = y_1 y_2,~ d y_4 = y_2 y_3, &
\end{eqnarray*}
and $y_1 \, (y_3+c)$ is closed but not exact for each closed $c$ of 
degree one. \q

\begin{SProp} $\,$
\begin{itemize}
\item[(i)] Each quotient of the Lie group $G :=
G_{6.36}^{0,0}$ by a lattice is symplectic and $G$ admits a lattice
$\Gamma$ with $b_1(G/\Gamma) = 2$.
\item[(ii)] If there is a lattice in the Lie group $G$ such that
the corresponding solvmanifold satisfies $b_1=2$ and $b_2=3$, then
it is symplectic and not formal.
\end{itemize}
\end{SProp}

\textit{Proof.} The proof of (ii) is analogous to the last one. Up
to generators of degree one, the minimal model is given by
\begin{eqnarray*}
&\rho(y_1) = x_2,~ \rho(y_2) = x_6,~ \rho (y_3) = -x_3,~ \rho(y_4) =
x_1,& \\ &d y_1 = d y_2 = 0 ,~ d y_3 = y_1 y_2,~ d y_4 = y_1 y_3, &
\end{eqnarray*}
and $y_1 \, (y_4+c)$ is closed but not exact for each closed $c$ of 
degree one.

ad (i): Using another initial data, we argue as in the proof of
Proposition \ref{G621}. We now have $ \mu(t) \big( (x,y,z,v,w) \big)
= \big (x-\frac{t}{2}y^2 , \overline{\mu}(t) \big( (x,y,z,v,w) \big)
\big)$ with \begin{eqnarray*} \overline{\mu}(t) \big( (y,z,v,w)
\big) &=& \exp^{GL(4,\mathbb{R})} (t \left(
\begin{array}{cccc}
0 & 0 & 0 & 0 \\
-t & 0 & 0 & 0 \\
0 & 0 & 0 & t \\
0 & 0 & -t & 0
\end{array} \right) )
\left( \begin{array}{c} y \\ z \\ v \\ w
\end{array} \right)
\\ &=& \left( \begin{array}{cccc}
1 & 0 & 0 & 0 \\
-t & 1 & 0 & 0 \\
0 & 0 & \cos(t) & \sin(t) \\
0 & 0 & - \sin(t) & \cos(t)
\end{array} \right)
\left( \begin{array}{c} y \\ z \\ v \\ w
\end{array} \right).
\end{eqnarray*}
Let $t_1 = \pi$, $A := \left(
\begin{array}{cccc}
1 & 1 & 0 & 0 \\
0 & 1 & 0 & 0 \\
0 & 0 & -1 & 0 \\
0 & 0 & 0 & -1
\end{array} \right)$ and $T := \left(
\begin{array}{cccc}
0 & -\frac{1}{t_1} & 0 & 0 \\
1 & 0 & 0 & 0 \\
0 & 0 & 1 & 0 \\
0 & 0 & 0 & 1
\end{array} \right)$. Then we have $T A T^{-1} = \overline{\mu}(t_1)$.
Denote the $i$-th column of T by $b_i$. Analogous calculations as in
loc.\ cit.\ lead to a lattice generated by $\gamma_0 :=
(\frac{1}{t_1},0_{\mathbb{R}^4})$ and $\gamma_i := (b_{i0}, b_i)$
for $i \in \{1,\ldots,4\}$, where $b_{20} \in \mathbb{R}$ arbitrary
and $b_{10}= - \frac{1}{2 t_1}, b_{30}=b_{40} = 0$.

Obviously, this lattice is represented by $$\langle \tau , \gamma_0
, \ldots, \gamma_4 \,|\, [\tau, \gamma_1] = 1, \, [\tau, \gamma_2]
=\gamma_1, \, [\tau, \gamma_3] = \gamma_3^{-2}, \, [\tau, \gamma_4]
= \gamma_4^{-2}, \, [\gamma_1, \gamma_2] = \gamma_0 \rangle$$ and
its abelianisation is $\mathbb{Z}^2 \oplus \mathbb{Z}_2{}^2$, i.e.\
the solvmanifold's first Betti number equals two. \q

\begin{SProp} $\,$
\begin{itemize}
\item[(i)] Each quotient of the Lie group $G :=
G_{6.70}^{0,0}$ by a lattice is symplectic and $G$ admits a lattice
$\Gamma$ with $b_1(G/\Gamma) = 2$.
\item[(ii)] If there is a lattice $\Gamma$ in $G$ such that
$b_1(G/\Gamma) = 2$ and $b_2(G/\Gamma) = 3$, then $G/\Gamma$ is
formal.
\end{itemize}
\end{SProp}

\textit{Proof.} ad (i): By definition, we have $G = \mathbb{R}
\ltimes_{\mu} N$ with $N = G_{5.1}$ and $\mu(t) = \exp^N \circ
\exp^{A(\mathfrak{n})}(t \, \mathrm{ad}(X_6)) \circ \log^N,$ where
$\{X_1,\ldots,X_6\}$ denotes a basis of $\mathfrak{g}$ as in the
second row of Table \ref{UniAlnil36c}. Note that $\{X_1, \ldots
,X_5\}$ is a basis of the nilradical $\mathfrak{n}$. One computes
\begin{eqnarray*}
\mu(t)_* &:=& d_e \big( \mu(t) \big) = \exp^{A(\mathfrak{n})}(t \,
\mathrm{ad}(X_6)) \\ &=& \left( \begin{array}{ccccc}
\cos(t) & \sin(t) & 0 & 0 & 0 \\
- \sin(t) & \cos(t) & 0 & 0 & 0 \\
0 & 0 & \cos(t) & \sin(t) & 0 \\
0 & 0 & - \sin(t) & \cos(t) & 0 \\
0 & 0 & 0 & 0 & 1
\end{array} \right).
\end{eqnarray*}

Using $\mathfrak{n} = \langle X_5 \rangle \ltimes_{ad} \langle X_1,
\ldots ,X_4 \, | \,  \rangle$, $\mathrm{ad}(X_5)(X_3)=-X_1$,
$\mathrm{ad}(X_5)(X_4)=-X_2$ and
$\mathrm{ad}(X_5)(X_1)=\mathrm{ad}(X_5)(X_2)=0$, we can determine
the Lie group $N$.

As a smooth manifold $N$ equals $\mathbb{R}^5$, and the
multiplication is given by
\begin{eqnarray*}
(a, b, c, r, s) \cdot (x, y, z, v, w) &=& \big( a + x + c w \,,\, b
+ y + r w \,,\, c + z \,,\, r + v \,,\, s + w \big) .
\end{eqnarray*}

By Theorem \ref{expsemidirekt}, we can obtain the exponential map of
$N$ as
\begin{eqnarray*}
\exp^N(x X_1 + y X_2 + z X_3 + v X_4 + w X_5) &=& ( x + \frac{z
w}{2} \,,\, y + \frac{v w}{2} \,,\, z \,,\, v \,,\, w ),
\end{eqnarray*}
and obviously, this implies
\begin{eqnarray*}
\log^N \big((x, y, z, v , w) \big) &=& (x - \frac{w z}{2}) X_1 + (y
- \frac{v w}{2}) X_2 + z  X_3 + v X_4 + w X_5.
\end{eqnarray*}

From $\mu(t) = \exp^N \circ \mu(t)_* \circ \log^N$ we get
\begin{eqnarray*}
\mu(t) \big( (x,y,z,v,w) \big) &=& (\cos(t) \, x + \sin(t) \, y
\,,\, -\sin(t) \, x + \cos(t) \, y \,,\\ && ~ \cos(t) \, z + \sin(t)
\, v \,,\, -\sin(t) \, z+ \cos(t) \, v \,,\, w),
\end{eqnarray*}
and $\mu(\pi)$ preserves the lattice $\{(x,y,z,v,w) \in N \,|\,
x,y,z,v,w \in \mathbb{Z}\}$.

The corresponding solvmanifold has $b_1 = 2$ since the
abelianisation of this lattice is isomorphic to $\mathbb{Z}^2 \oplus
\mathbb{Z}_2{}^4$.

ad (ii): Up to generators of degree two, the minimal model is given
by
\begin{eqnarray*}
&\rho(y_1) = x_5,~ \rho(y_2) = x_6,~ \rho (z_1) = x_{13} + x_{24},~
\rho(z_2) = x_{34},& \\ &d y_1 = d y_2 = 0 ,~ d z_1 = d z_2 = 0,&
\end{eqnarray*}
hence the minimal is $2$-formal. By Theorem \ref{formal = n-1
formal}, the solvmanifold is formal. \q

\begin{SProp} $\,$
\begin{itemize}
\item[(i)] $G := G_{6.118}^{0,\pm1,-1}$ admits a lattice such that
the first Betti number of the corresponding solvmanifold equals two
and the second Betti number equals five.
\item[(ii)] If there is a lattice $\Gamma$ in $G$
such that $b_1(G/\Gamma) = 2$ and $b_2(G/\Gamma) = 3$, then
$G/\Gamma$ is symplectic and formal.
\end{itemize}
\end{SProp}

\textit{Proof.} The construction of the lattices mentioned in (i)
can be found in \cite{Yam}. In loc.\ cit.\ $G_{6.118}^{0,1,-1}$ is
denoted by $G_3$ and $G_{6.118}^{0,-1,-1}$ by $G_1$, respectively.
The Betti numbers of the quotient of $G_{118}^{0,-1,-1}$ are
determined explicitly. In the case of $G_{118}^{0,1,-1}$, one can
make an analogous computation.

Assume there is a lattice that satisfies the condition of (ii). Up
to generators of degree two, the solvmanifold's minimal model is
given by
\begin{eqnarray*}
&\rho(y_1) = x_5,~ \rho(y_2) = x_6,~ \rho (z_1) = x_{13} \pm
x_{24},~ \rho(z_2) = x_{14} \mp x_{23},& \\ &d y_1 = d y_2 = 0 ,~ d
z_1 = d z_2 = 0,&
\end{eqnarray*}
hence it is $2$-formal. Theorem \ref{formal = n-1 formal} then
implies formality. \q

\begin{Rem}
$G_{6.118}^{0,-1,-1}$ is the underlying real Lie group of the unique
connected and simply-connected complex three-dimensional Lie group
that is solvable and not nilpotent. Its compact quotients by
lattices are classified in \cite[Theorem 1]{Nak}. They always
satisfy $b_1 = 2$ and moreover, for the Hodge number $h^{0,1}$ holds
either $h^{0,1}=1$ or $h^{0,1}=3$.
\end{Rem}

\begin{Rem}
Besides the groups mentioned in this section, the following solvable
but not completely solvable Lie groups $G_{6.i}$ could give rise to
a symplectic solvmanifold with $b_1 > 1$. But then the cohomology
class of the symplectic form cannot lie in the image of the
inclusion $H^* (\bigwedge \mathfrak{g}^*_{6.i},
\delta)\hookrightarrow H^*(G_{6.i}/\Gamma,d)$.

There are the sixteen classes of groups in (\ref{nvavs}) and
\begin{equation*} \begin{array}{l@{~~~}l@{~~~}l} i=9,~ b=0;
& i = 33 ,~ a=0; & i=34,~ a=0; \\ i=35,~ a=-b; &
i= 89,~ s=0,~ \nu_0 \ne 0; & i=92,~ \nu_0 \mu_0 = 0; \\
i \in \{107,113,\ldots,116,125,135\}. &&
\end{array}\end{equation*}
\end{Rem}

\subsection{Decomposable solvmanifolds}
The six-dimensional decomposable solvmanifolds $G / \Gamma = H_1 /
\Gamma_1 \times H_2 / \Gamma_2$ being not a nilmanifold are
contained in Table \ref{6decSolvmgf} on page \pageref{6decSolvmgf}.
Using Theorem \ref{min chev eil}, one can deduce the statement about
the Betti numbers. The results on the existence of a symplectic form
were mostly made by Campoamor-Stursberg in \cite{CS}. He examined
whether the Lie algebra admits a symplectic form. Note that in
\cite{CS} the symplectic forms $$\lambda \, x_{12} + \mu \, x_{15} +
\nu \, x_{26} + \rho \, x_{34} + \sigma \, x_{56},~~ \rho \ne 0,~
\lambda \sigma \ne \mu \nu,$$ on $\mathfrak{g}_{5.14}^0 \oplus
\mathfrak{g}_1$ are absent.

Since there is a monomorphism from the Lie algebra cohomology to the
solvmanifold's cohomology, the existence of a symplectic form with
non-exact cubic on the Lie algebra implies the existence of such an
on the solvmanifold. Recall that the Lie algebra is generated by the
left-invariant one-forms on the Lie group. If the Lie algebra
cohomology is isomorphic to the solvmanifold's
cohomology\footnote{E.g.\ this happens if the Lie algebra is
completely solvable or if the above monomorphism must be an
isomorphism by dimension reasons.}, one knows whether the
solvmanifold is symplectic or not. Up to exact summands the
symplectic forms are listed in Table \ref{6decSolvmgfsympl} with
respect to the dual of the Lie algebra's bases given in Appendix
\ref{Liste Solv}. In the column ``isom.'', we mark whether there is
an isomorphism of the cohomology algebras.

%
We do not claim that Table \ref{6decSolvmgf} contains all connected
and simply-connected decomposable solvable and non-nilpotent Lie
groups which admit a lattice -- just  those Lie groups admitting a
lattice such that the corresponding solvmanifold is a product of
lower-dimensional ones.

\begin{table}[p!]
\centering \caption{\label{6decSolvmgf} Decomposable
non-nil-solvmanifolds $G/\Gamma = H_1/\Gamma_1 \times H_2/\Gamma_2$}
\begin{tabular}{|c|c|c|c|c|c|} \hline \hline
$G$ & $b_1(G/\Gamma)$ & $b_2(G/\Gamma)$ & formal & sympl. & Comment %
\\ \hline
\hline $G_{5.7}^{p,q,r} \times \mathbb{R}$ & $2$ & $1$ & yes & no
& $-1 < r < q < p < 1,$ \\ &&&&& $pqr \ne 0,$ \\ &&&&&  $p+q+r = -1$ \\
\hline $G_{5.7}^{p,-p,-1} \times \mathbb{R}$ & $2$ & $3$ & yes & yes
& $p \in ]0,1[$\\
\hline $G_{5.7}^{1,-1,-1} \times \mathbb{R}$ & $2$ & $5$ & yes & yes & \\
\hline $G_{5.8}^{-1} \times \mathbb{R}$ & $3$ & $5$ & no & yes & \\
\hline $G_{5.13}^{-1-2q,q,r} \times \mathbb{R}$ & $\ge 2$ & $\ge 1$
& ?
& ? & $q \in [-1,0[,$ \\ &&&&& $q \ne -\frac{1}{2}, r \ne 0 $\\
\hline $G_{5.13}^{-1,0,r} \times \mathbb{R}$ & $\ge 2$ & $\ge 3$ &
? & yes & $r \ne 0$\\
\hline $G_{5.14}^{0} \times \mathbb{R}$ & $\ge 3$ & $\ge 5$ &
? & yes &\\
\hline $G_{5.15}^{-1} \times \mathbb{R}$ & $2$ & $3$ &
no &yes&\\
\hline $G_{5.17}^{p,-p,r} \times \mathbb{R}$ & $\ge 2$ & $\ge 1$ & ?
& ? & $p \ne 0,~ r\not\in\{0,\pm 1\},$ \\
\hline $G_{5.17}^{p,-p,r} \times \mathbb{R}$ & $\ge 2$ & $\ge 3$ &?&
yes & $(p\ne 0, r=\pm1)$\\ &&&&& or $(p=0, r\not\in\{0,\pm1\})$\\
\hline $G_{5.17}^{0,0,\pm 1} \times \mathbb{R}$ & $\ge 2$ & $\ge 5$
&?& yes & \\
\hline $G_{5.18}^{0} \times \mathbb{R}$ & $\ge 2$ & $\ge 3$ & ? &
yes & \\
\hline $G_{5.20}^{-1} \times \mathbb{R}$ & $3$ & $3$ & yes & no & \\
\hline $G_{5.26}^{0,\pm1} \times \mathbb{R}$ & $\ge 3$ & $\ge 3$ & ?
& ? & \\
\hline $G_{5.33}^{-1,-1} \times \mathbb{R}$ & $3$ & $3$ & yes &
no & 
\\
\hline $G_{5.35}^{-2,0} \times \mathbb{R}$ & $\ge 3$ & $\ge 3$ & ? &
? & \\
\hline \hline $G_{4.5}^{p,-p-1} \times \mathbb{R}^2$ & $3$ & $3$ &
yes &
no & $p \in [-\frac{1}{2} ,0[$\\
\hline $G_{4.6}^{-2p,p} \times \mathbb{R}^2$ & $3$ & $3$ & yes &
no & $p > 0$\\
\hline $G_{4.8}^{-1} \times \mathbb{R}^2$ & $3$ & $3$ & yes &
no & \\
\hline $G_{4.9}^{0} \times \mathbb{R}^2$ & $3$ & $3$ & yes &
no &\\
\hline \hline $G_{3.4}^{-1} \times \mathbb{R}^3$ & $4$ & $7$ & yes &
yes & \\
\hline $G_{3.5}^{0} \times \mathbb{R}^3$ & $4$ & $7$ & yes &
yes &\\
\hline $G_{3.1} \times G_{3.4}^{-1}$ & $3$ & $5$ & no & yes &  \\
\hline $G_{3.1} \times G_{3.5}^0$ & $3$ & $5$ & no & yes &  \\
\hline $G_{3.4}^{-1} \times G_{3.4}^{-1}$ & $2$ & $3$ & yes
& yes &  \\
\hline $G_{3.4}^{-1} \times G_{3.5}^0$ & $2$ & $3$ & yes & yes &  \\
\hline $G_{3.5}^0 \times G_{3.5}^0$ & $2$ & $3$ & yes & yes &  \\
\hline
\end{tabular}
\end{table}
%
\begin{table}[p!]
\centering \caption{\label{6decSolvmgfsympl} Symplectic forms on
$G/\Gamma = H_1/\Gamma_1 \times H_2/\Gamma_2$}
\begin{tabular}{|c|l|c|} \hline \hline
$\mathfrak{g}$ & symplectic forms & isom. %
\\ \hline
\hline $\mathfrak{g}_{5.7}^{p,-p,-1} \oplus \mathfrak{g}_1$ & $a \,
x_{14} + b \, x_{23} + c \, x_{56},~ abc \ne0$ & yes\\
\hline $\mathfrak{g}_{5.7}^{1,-1,-1} \oplus \mathfrak{g}_1$ & $a \,
x_{13} + b \, x_{14} + c \, x_{23} + d \, x_{24} + e \, x_{56},~
e(bc - ad) \ne 0$ & yes  \\
\hline $\mathfrak{g}_{5.8}^{-1} \oplus \mathfrak{g}_1$ & $a \,
x_{12} + b \, x_{15} + c \, x_{26} + d \, x_{34} + e \, x_{56},~
d(ae - bc) \ne 0$ & yes \\
\hline $\mathfrak{g}_{5.13}^{-1,0,r} \oplus \mathfrak{g}_1$ & $a \,
x_{12} + b \, x_{34} + c \, x_{56},~ abc \ne 0$ & ? \\
\hline $\mathfrak{g}_{5.14}^{0} \oplus \mathfrak{g}_1$ & $a \,
x_{12} + b \, x_{15} + c \, x_{26} + d \, x_{34} + e \, x_{56},~
d(ae-bc) \ne 0$ & ? \\
\hline $\mathfrak{g}_{5.15}^{-1} \oplus \mathfrak{g}_1$ & $a \,
(x_{14} - x_{23}) + b \, x_{24} + c \, x_{56},~
abc \ne 0$ & yes\\
\hline $\mathfrak{g}_{5.17}^{p,-p,\pm1} \oplus \mathfrak{g}_1$ & $a
\, (x_{13} \pm x_{24}) + b \, (x_{14} \mp x_{23}) + c \, x_{56},~
abc \ne 0$ & ? \\
$p \ne 0$ && \\
\hline $\mathfrak{g}_{5.17}^{0,0,r} \oplus \mathfrak{g}_1$ & $a \,
x_{12} + b \, x_{34} + c \, x_{56},~ abc \ne 0$ & ? \\ $r \ne \pm 1$
&& \\
\hline $\mathfrak{g}_{5.17}^{0,0,\pm 1} \oplus \mathfrak{g}_1$ & $a
\, x_{12} + b \, (x_{13} \pm x_{24}) + c \, (x_{14} \mp x_{23}) + d
\, x_{34} + e \, x_{56},$ & ? \\ &$e(ad \mp (b^2+c^2)) \ne 0$& \\
\hline $\mathfrak{g}_{5.18}^{0} \oplus \mathfrak{g}_1$ & $a \,
(x_{13} + x_{24}) + b \, x_{24} + c \, x_{56},~ ac \ne0$ & ? \\
\hline \hline $\mathfrak{g}_{3.4}^{-1} \oplus 3 \mathfrak{g}_1$ & $a
\, x_{12} + b \, x_{34} + c \, x_{35} + d \, x_{36} + e \, x_{45} +
f \, x_{46} + g \, x_{56},$ & yes \\ & $ a(de - cf + bg) \ne 0$ & \\
\hline $\mathfrak{g}_{3.5}^{0} \oplus 3 \mathfrak{g}_1$ & $a \,
x_{12} + b \, x_{34} + c \, x_{35} + d \, x_{36} + e \, x_{45} + f
\, x_{46} + g \, x_{56},$ & yes \\ & $ a(de - cf + bg) \ne 0$ & \\
\hline $\mathfrak{g}_{3.1} \oplus \mathfrak{g}_{3.4}^{-1}$ & $a \,
x_{12} + b \, x_{13} + c \, x_{26} + d \, x_{36} + e \, x_{45},~
e(ad - bc) \ne 0$ & yes \\
\hline $\mathfrak{g}_{3.1} \oplus \mathfrak{g}_{3.5}^0$ & $a \,
x_{12} + b \, x_{13} + c \, x_{26} + d \, x_{36} + e \, x_{45},~
e(ad - bc) \ne 0$ & yes \\
\hline $\mathfrak{g}_{3.4}^{-1} \oplus \mathfrak{g}_{3.4}^{-1}$ & $a
\, x_{12} + b \, x_{36} + c \, x_{45},~
abc \ne 0$ & yes \\
\hline $\mathfrak{g}_{3.4}^{-1} \oplus \mathfrak{g}_{3.5}^0$ & $a \,
x_{12} + b \, x_{36} + c \, x_{45},~
abc \ne 0$ & yes \\
\hline $\mathfrak{g}_{3.5}^0 \oplus \mathfrak{g}_{3.5}^0$ & $a \,
x_{12} + b \, x_{36} + c \, x_{45},~ abc \ne 0$ & yes \\ \hline
\end{tabular}
\end{table}
%
%

\section{Relations with the Lefschetz property}
We have seen in Section \ref{fskl} that a compact K\"ahler manifold
is formal, Hard Lefschetz and its odd-degree Betti numbers are even.
Even if a manifold has a complex structure, these conditions are not
sufficient as the following theorem which is mentioned in
\cite{HasCK} shows. Recall, we have seen above that
$G_{5.7}^{1,-1,-1}$ admits a lattice.

\begin{Thm}  \label{alles nur nicht Kaehler}
Let $\Gamma$ be an arbitrary lattice in $G_{5.7}^{1,-1,-1}$. Then
the solvmanifold $M := G_{5.7}^{1,-1,-1} / \Gamma \times S^1$ is
formal, Hard Lefschetz and has even odd-degree Betti numbers.
Moreover, $M$ possesses a complex structure but it cannot be
K\"ahlerian.
\end{Thm}

\textit{Proof.} From Proposition \ref{G5.7} follows that the Lie
group $G := G_{5.7}^{1,-1,-1} \times \mathbb{R}$ possesses a lattice
$\Gamma$. The Chevalley-Eilenberg complex of its Lie algebra
$$\langle \, X_1, \ldots, X_6 \,|\, [X_1,X_5] = X_1,\,
[X_2,X_5] = X_2,\, [X_3,X_5] = - X_3,\, [X_4,X_5] = - X_4 \,
\rangle$$ is given by
$$ \delta x_1 = -x_{15},~ \delta x_2 = -x_{25},~ \delta x_3 =
x_{35},~ \delta x_4 = x_{45},~ \delta x_5 = \delta x_6 = 0,$$ where
$\{x_1,\ldots,x_6\}$ is a basis of the left-invariant one-forms on
$G$. Since $G$ is completely solvable, Theorem \ref{min chev eil}
(ii) enables us to compute the cohomology of $M$ as
\begin{eqnarray}
\nonumber H^1(M,\mathbb{R}) &=& \langle [x_5], [x_6] \rangle, \\
\nonumber H^2(M,\mathbb{R}) &=& \langle [x_{13}], [x_{14}],
[x_{23}], [x_{24}], [x_{56}] \rangle,
\\ \label{H357S} H^3(M,\mathbb{R}) &=& \langle [x_{135}], [x_{136}], [x_{145}],
[x_{146}], [x_{235}], [x_{236}], [x_{245}], [x_{246}] \rangle,  \\
\nonumber H^4(M,\mathbb{R}) &=& \langle [x_{1234}], [x_{1356}],
[x_{1456}], [x_{2356}], [x_{2456}] \rangle, \\ \nonumber
H^5(M,\mathbb{R}) &=& \langle [x_{12345}], [x_{12346}] \rangle.
\end{eqnarray}
Let $[\omega] \in H^2(M,\mathbb{R})$ represent a symplectic form on
$M$. A short calculation shows that there are $a,b,c,d,e \in
\mathbb{R}$ with $e(bc - ad) \ne 0$ and $$[\omega] = a [x_{13}] + b
[x_{14}] + c [x_{23}] + d [x_{24}] + e [x_{56}]. $$ Since $[x_5]
\cup [\omega]^2 = 2 (bc-de) [x_{12345}] \ne 0$ and $[x_6] \cup
[\omega]^2 = 2 (bc-de) [x_{12346}] \ne 0 $, the homomorphism $L^2 \:
H^1(M,\mathbb{R}) \to H^5(M,\mathbb{R})$ is an isomorphism.

In the basis (\ref{H357S}), the homomorphism $L^1 \:
H^2(M,\mathbb{R}) \to H^4(M,\mathbb{R})$ is represented by the
matrix $\left(
\begin{array}{ccccc}
-d & c & -b & -a & 0 \\
e & 0 & 0 & 0 & a \\
0 & e & 0 & 0 & b \\
0 & 0 & e & 0 & c \\
0 & 0 & 0 & e & d
\end{array} \right)$ which has $2 e^3 (ad-bc)\ne 0$ as determinant,
hence $M$ is Hard Lefschetz.

We define an almost complex structure $J$ on $G$ by
$$JX_1 = X_2,~ JX_2 = -X_1,~ JX_3 = X_4,~ JX_4 = - X_3,~ JX_5 =
X_6,~ JX_6 = -X_5,$$ which induces an almost complex structure on
$M$. It is easy to see that the Nijenhuis tensor vanishes, hence $M$
is a complex manifold.

$M$ is a non-toral solvmanifold which is a quotient of a completely
solvable Lie group. Therefore, $M$ cannot be K\"ahlerian by Theorem
\ref{Kaehler vollst aufl}. \q
\N The authors of \cite{IRTU99} considered the relations between the
above three properties for closed symplectic manifolds. We want to
try to complete \cite[Theorem 3.1 Table 1]{IRTU99} in the case of
symplectic solvmanifolds. Actually, the mentioned table deals with
symplectically aspherical closed manifolds, but note that symplectic
solvmanifolds are symplectically aspherical.

We start our investigations by the examination of the Lefschetz
property in dimension four.

\begin{Thm}
A four-dimensional symplectic solvmanifold is not (Hard) Lefschetz
if and only if it is a non-toral nilmanifold. Especially, the (Hard)
Lefschetz property is independent of the choice of the symplectic
form.
\end{Thm}

\textit{Proof.} By Theorem \ref{Solvmgf Dim4}, there are five
classes of four-dimensional symplectic solvmanifolds. Three of them
are nilmanifolds and satisfy the claim by Corollary
\ref{NilLefschetzgdw}.

There remain two non-nilmanifolds to consider. We start with a
quotient of the Lie group which has $\mathfrak{g}_{3.1}^{-1} \oplus
\mathfrak{g}_1$ as Lie algebra, see Table \ref{Sol4}.  The Lie group
is completely solvable, hence the Lie algebra cohomology is
isomorphic to the solvmanifold's cohomology. If $x_1,\ldots,x_4$
denote the left-invariant one-forms which are dual to the basis
given in Table \ref{Sol4}, one computes
\begin{eqnarray} \label{KohomSol4}
H^1 = \langle [x_3], [x_4] \rangle, & H^2 = \langle [x_{12}],
[x_{34}] \rangle, & H^3 = \langle [x_{123}], [x_{124}] \rangle.
\end{eqnarray}
The class representing a symplectic form must be of the form $[a \,
x_{12} + b \, x_{34}]$ with $a,b \ne 0$ and obviously, the Lefschetz
map with respect to this class is an isomorphism.

Now, consider a solvmanifold $G/\Gamma$ such that the Lie algebra of
$G$ is $\mathfrak{g}_{3.5}^0 \oplus \mathfrak{g}_1$ and
$b_1(G/\Gamma) = 2$. A short computation yields that the Lie algebra
cohomology of $\mathfrak{g}_{3.5} \oplus \mathfrak{g}_1$ is the same
as in (\ref{KohomSol4}). Since $G / \Gamma$ is compact and
parallelisable, we see further $b_i(G/\Gamma) = 2$ for $i \in
\{1,2,3\}$, and Theorem \ref{min chev eil} (i) implies that
(\ref{KohomSol4}) also gives the cohomology of $G/\Gamma$. We have
yet seen that a symplectic four-manifold with this cohomology is
Hard Lefschetz. \q
\N Denote $KT$ ``the'' four-dimensional symplectic nilmanifold with
$b_1(KT) = 3$. We have seen that $KT$ is not formal and not
Lefschetz. Its square has the following properties:

\begin{Thm}[\cite{IRTU99}]
$KT \times KT$ is not formal, not Lefschetz and has even odd-degree
Betti numbers. \q
\end{Thm}

Next, we are looking for an example of a formal manifold that is not
Lefschetz and has even odd degree Betti numbers resp.\ an odd odd
degree Betti number.

\begin{Thm}
The Lie group $G_{6.78}$ admits a lattice $\Gamma$, see above. $M :=
G_{6.78} / \Gamma$ is a formal solvmanifold  with $b_1(M) = 1$ that
admits a symplectic form $\omega$ such that $(M,\omega)$ is not Hard
Lefschetz. Moreover, $(M \times M,\omega \times \omega)$ is a formal
symplectic manifold with even odd degree Betti numbers that is not
Hard Lefschetz.
\end{Thm}

\textit{Proof.} By Theorem \ref{G678}, $M$ is a formal symplectic
manifold with Betti numbers $b_1(M) = b_2(M) =1$. Note that this
implies that $M \times M$ is symplectic and formal (the latter
property by Proposition \ref{prod}).

Corollary \ref{b1 Lefschetz} forces $M$ to be not Lefschetz and
since \cite[Proposition 4.2]{FMDon} says that a product is Lefschetz
if and only if both factors are Lefschetz, $M \times M$ is not
Lefschetz.

$M$ is a six-dimensional solvmanifold and so it is parallelisable.
Hence the fact $b_0(M) = b_1(M) = b_2(M) =1$ implies $b_3(M) = 2$.
This and Poincar\'e Duality imply $b_1(M \times M) = b_{11}(M \times
M) = 2$, $b_3(M \times M) = b_9(M \times M) = 6$ and $b_5(M \times
M) = b_7(M \times M) = 4$. \q
\N In 1990, Benson and Gordon \cite[Example 3]{BG88} constructed an
eight-dimensional non-exact symplectic and completely solvable Lie
algebra that does not satisfy the Hard Lefschetz property, but they
did not know whether the corresponding connected and
simply-connected Lie group $G^{BG}$ admits a lattice.

Fern\'andez, de Le\'on and Saralegui computed in \cite[Proposition
3.2]{FLS} the minimal model of the complex of the left-invariant
differential forms on $G^{BG}$. It is formal and its cohomology of
odd degree is even-dimensional. If $G^{BG}$ admits a lattice, by
completely solvability, the corresponding solvmanifold would be a
symplectic and formal manifold with even odd degree Betti numbers
that violates the Hard Lefschetz property.

In 2000, Tralle \cite{Tralle} claimed that a lattice does not exist
but Sawai and Yamada noted 2005 Tralle's proof would contain
calculatory errors and constructed a lattice \cite[Theorem 1]{SY}.
This proves the next theorem.

\begin{Thm}
There exists an eight-dimensional symplectic and formal solvmanifold
$M^{BG}$ with even odd degree Betti numbers that is not Hard
Lefschetz. \q
\end{Thm}

We sum up the above results in Table \ref{relokp} on page \pageref{relokp}.
It is an enlargement of \cite[Theorem 3.1 Table 1]{IRTU99}.

\begin{table}[ht]
\centering \caption{Relations of the K\"ahler properties}
\label{relokp}
\begin{tabular}{|c|c|c|c|} \hline \hline
Formality & Hard Lefschetz & $b_{2i+1} \equiv 0(2)$ & Example %
\\ \hline
\hline yes & yes & yes & K\"ahler \\
\hline yes & yes & no & impossible\\
\hline yes & no & yes & $M^{BG}$, $G_{6.78}/\Gamma \times G_{6.78}/\Gamma$ \\
\hline yes & no & no & $G_{6.78}/\Gamma$ \\
\hline no & yes & yes & ? \\
\hline no & yes & no & impossible \\
\hline no & no & yes & $KT \times KT$ \\
\hline no & no & no & $KT$ \\ \hline
\end{tabular}
\end{table}

Unfortunately, the missing example does not arise among
the six-dimensional solvmanifolds. In order to see
this, one has to examine which of them satisfy the (Hard) Lefschetz
property. We briefly mention the results.

By Corollary \ref{b1 Lefschetz}, a manifold with odd first Betti
number cannot be Lefschetz. We now examine such indecomposable
solvmanifolds whose first Betti number is even; in the completely
solvable case, these are quotients of $G_{6.3}^{-1}$, $G_{6.21}^0$
and $G_{6.54}^{-1}$.

The proof of the next two propositions is done analogous as that of
Theorem \ref{alles nur nicht Kaehler}. By complete solvability, we
know the solvmanifolds' cohomology and all possible symplectic forms
were determined in the proof of Proposition \ref{Symplgg1}.
Therefore, one can compute the image of the Lefschetz maps.

\begin{Prop}
Let a lattice in $G_{6.3}^{-1}$ or $G_{6.21}^0$ be given. (Wee have
seen above that such exists.) Then the corresponding (non-formal)
symplectic solvmanifold (with $b_1 = 2, \, b_2 = 3$) is not
Lefschetz, independent of the choice of the symplectic form. \q
\end{Prop}

\begin{Prop} \label{G654Lef}
Let a lattice in $G_{6.54}^{-1}$ be given. (Such exists by
Proposition \ref{G654}.) The corresponding (non-formal) symplectic
solvmanifold (with $b_1=2, \linebreak b_2=5$) is Lefschetz but not
Hard Lefschetz, independent of the choice of the symplectic form. \q
\end{Prop}

\begin{Rem}
The existence of a lattice in $G_{6.54}^{-1}$ was proven by
Fern\'andez, de L\'eon and Saralegui in \cite{FLS}. They also
computed the Betti numbers of the corresponding solvmanifold, showed
that it is not formal and does not satisfy the Hard Lefschetz
property with respect to a certain symplectic form. Moreover,
Fern\'andez and Mu\~noz proved in \cite[Example 3]{FMDon} that the
manifold is Lefschetz. (Analogous calculations work for other
symplectic forms.)
\end{Rem}

In the non-completely solvable case, the situation becomes a little
more complicated. If we are willing to make a statement about the
Lefschetz property, we have to know the cohomology and need
therefore assumptions on the Betti numbers.

\begin{Prop}
If there is a lattice in one of the non-completely solvable groups
$G_{6.i}^{0,0}$, $i \in \{10,36,70\}$ resp.\ $G_{6.118}^{0,\pm1,-1}$
such that the cohomology of the corresponding solvmanifold $M_i$ is
isomorphic to the Lie algebra cohomology of $\mathfrak{g}_{6.i}$
(i.e.\ the cohomology is as small as possible), then one computes
that the following hold, independent of the choice of the symplectic
forms provided by Proposition \ref{Symplgg1}:
\begin{itemize}
\item $M_{10}$ and $M_{36}$ are not formal and not Lefschetz.
\item $M_{70}$ is formal and Lefschetz but not Hard Lefschetz.
\item $M_{118}$ is formal and Hard Lefschetz.
\end{itemize}
(The statements on formality follow from the propositions at the end
of Section \ref{6Sgr1}.) \q
\end{Prop}

Finally, we consider the decomposable symplectic solvmanifolds
listed in Table \ref{6decSolvmgf}.

\begin{Prop}
Let $G/\Gamma = H_1/\Gamma_1 \times H_2/\Gamma_2$ be one of the
symplectic solvmanifolds listed in Table \ref{6decSolvmgf} such that
in the corresponding row of the table arises no $\ge$-sign.

Then $G/\Gamma$ is formal if and only if it is Hard Lefschetz
(independent of the special choice of the symplectic form).
\end{Prop}

\textit{Sketch of the proof.} One has an isomorphism from the Lie
algebra cohomology to the solvmanifold's cohomology for each
manifold as in the theorem. Then an explicit calculation as in the
proof of Theorem \ref{alles nur nicht Kaehler} yields that the Hard
Lefschetz manifolds among the considered are exactly the formal
ones.

Note, if $b_1$ is not even, the claim follows directly from Theorem
\ref{b2i+1 Lefschetz}. \q
\begin{Rem}
$G_{5.15}^{-1}/\Gamma_1 \times S^1$ is Lefschetz. The other
manifolds in the last proposition are even not Lefschetz if they are
not Hard Lefschetz.

$G_{5.15}^{-1}/\Gamma_1 \times S^1$ is a Lefschetz manifold that is
not formal and has even odd degree Betti numbers.
\end{Rem}

A similar result as the last proposition holds for the manifolds in
Table \ref{6decSolvmgf} such that in the corresponding row of the
table arises a $\ge$-sign. But we again must make an assumption that
enables us to compute the whole cohomology.

\begin{Prop} $\,$
\begin{itemize}
\item[(i)] Let $M = G_{5.i}/\Gamma \times \mathbb{R} / \mathbb{Z}$ be a
symplectic manifold such that one of the following conditions holds:
\begin{itemize}
\item[a)] $i=13$ with $q=0$ and $b_1(M)=2$ as well as $b_2(M) = 3$,
\item[b)] $i=17$ with $p \ne 0, \; r = \pm 1$ and $b_1(M)=2$ as well as $b_2(M)=3$,
\item[c)] $i=17$ with $p = 0, \; r \not\in \{0,\pm 1\}$ and $b_1(M)=2$ as well as $b_2(M)=3$,
\item[d)] $i=17$ with $p = 0, \; r = \pm 1$ and $b_1(M)=2$ as well
as $b_2(M)=5$.
\end{itemize}

Then $M$ is formal and Hard Lefschetz (independent of the special
choice of the symplectic form).
\item[(ii)] Let $\Gamma$ be a lattice in $G_{5.14}^0$ such that
$M = G_{5.14}^0/\Gamma \times \mathbb{R} / \mathbb{Z}$ satisfies
$b_1(M)=3$ and $b_2(M)=5$.

Then $M$ is not formal and not Lefschetz (independent of the special
choice of the symplectic form).

\item[(iii)] Let $\Gamma$ be a lattice in $G_{5.18}^0$ such that
$M = G_{5.18}^0/\Gamma \times \mathbb{R} / \mathbb{Z}$ satisfies
$b_1(M)=2$ and $b_2(M)=3$.

Then $M$ is not formal and Lefschetz but not Hard Lefschetz
(independent of the special choice of the symplectic form). \q
\end{itemize}
\end{Prop}

\begin{appendix}

\hyphenation{mani-fold mani-folds nil-mani-fold nil-mani-folds
solv-mani-fold solv-mani-folds}

\section{Lists of Lie Algebras} \label{Liste Solv}

\begin{table}[h]
In Table \ref{Sol4}, we give the isomorphism classes of Lie algebras
of the simply-connected solvable Lie groups up to dimension four
that possesses lattices. The designation $\mathfrak{g}_{i,j}$ means
the $j$-th indecomposable solvable Lie algebra of dimension $i$. The
choice of the integer $j$ bases on the notation of \cite{Mub63}. The
superscripts, if any, give the values of the continuous parameters
on which the algebra depends.
\end{table}

\begin{table}[h]
\centering \caption{\label{Sol4} Solvmanifolds up to dimension four}
\begin{tabular}{|c|c|c|} \hline \hline
& $[X_i,X_j]$ & cpl. solv.
\\ \hline \hline $\mathfrak{g}_1$ &  & abelian
\\ \hline \hline $2 \mathfrak{g}_1$&
& abelian
\\ \hline \hline $3 \mathfrak{g}_1$ & &
abelian
\\ \hline $\mathfrak{g}_{3.1}$ & $ [X_2 , X_3 ] = X_1$
& nilpotent
\\ \hline $\mathfrak{g}_{3.4}^{-1}$ & $[X_1,X_3] = X_1,~
[X_2,X_3] = - X_2$ & yes
\\ \hline $\mathfrak{g}_{3.5}^{0}$ & $[X_1,X_3] = - X_2,~
[X_2,X_3] =  X_1$ & no
\\ \hline \hline $4 \mathfrak{g}_1$ & &
abelian
\\ \hline $\mathfrak{g}_{3.1} \oplus \mathfrak{g}_1$
& $ [X_2 , X_3 ] = X_1$ & nilpotent
\\ \hline $\mathfrak{g}_{3.4}^{-1} \oplus \mathfrak{g}_1$ & $[X_1,X_3] = X_1,~
[X_2,X_3] = - X_2$ & yes
\\ \hline $\mathfrak{g}_{3.5}^0 \oplus \mathfrak{g}_1$  & $[X_1,X_3] = - X_2,~
[X_2,X_3] =  X_1$ & no
\\ \hline $\mathfrak{g}_{4.1}$ & $[X_2,X_4] = X_1,~ [X_3,X_4] =
X_2$ & nilpotent
\\ \hline $\mathfrak{g}_{4.5}^{p,-p-1}$ & $[X_1,X_4] = X_1,~ [X_2,X_4] =
p X_2,$ & yes \\ &  $[X_3,X_4] = (-p-1) X_3,~ 
- \frac{1}{2} \le p < 0$ &
\\ \hline $\mathfrak{g}_{4.6}^{-2p,p}$ & $[X_1,X_4] = -2p X_1,~ [X_2,X_4] = p
X_2 - X_3,$ & no \\ & $[X_3,X_4] = X_2 + p X_3,~ p>0$ &
\\ \hline $\mathfrak{g}_{4.8}^{-1}$ & $[X_2,X_3] = X_1,~ [X_2,X_4] =
X_2,~ [X_3,X_4] = -X_3 $ & yes
\\ \hline $\mathfrak{g}_{4.9}^0$ & $[X_2,X_3] = X_1,~ [X_2,X_4] =
-X_3,~ [X_3,X_4] = X_2 $& no
\\ \hline
\end{tabular}
\end{table}

\begin{table}[b!]
We do not claim that the corresponding Lie groups admit a lattice
for all parameters. We just know that there exist such for certain
parameters! \\[15mm] The set of non-isomorphic five dimensional nilpotent
Lie algebras is exhausted by three types of decomposable algebras
and six indecomposables which are listed in the next table. The
designation is taken from \cite{Mub63a}.
\end{table}

\begin{table}[t]
\centering \caption{\label{Nil5} $5$-dimensional nilpotent algebras}
\begin{tabular}{|c|c|} \hline \hline
& $[X_i,X_j]$
\\ \hline \hline $5 \mathfrak{g}_1$ & abelian
\\ \hline $\mathfrak{g}_{3.1} \oplus 2 \mathfrak{g}_1$ &
$[X_2 , X_3 ] = X_1$
\\ \hline $\mathfrak{g}_{4.1} \oplus \mathfrak{g}_1$ & $[X_2,X_4] = X_1,~
[X_3,X_4] = X_2$
\\ \hline $\mathfrak{g}_{5.1}$ & $[X_3,X_5] = X_1,~ [X_4,X_5] = X_2$
\\ \hline $\mathfrak{g}_{5.2}$ & $[X_2,X_5] = X_1,~ [X_3,X_5] = X_2, ~
[X_4,X_5] = X_3$
\\ \hline $\mathfrak{g}_{5.3}$ & $[X_2,X_4] = X_3,~ [X_2,X_5] = X_1, ~
[X_4,X_5] = X_2$
\\ \hline $\mathfrak{g}_{5.4}$ & $[X_2,X_4] = X_1,~ [X_3,X_5] = X_1$
\\ \hline $\mathfrak{g}_{5.5}$ & $[X_3,X_4] = X_1,~ [X_2,X_5] = X_1, ~
[X_3,X_5] = X_2$
\\ \hline $\mathfrak{g}_{5.6}$ & $[X_3,X_4] = X_1,~ [X_2,X_5] = X_1, ~
[X_3,X_5] = X_2,~ [X_4,X_5] = X_3$
\\ \hline
\end{tabular}
\end{table}

\begin{table}[hb]
There are $24$ classes of solvable and non-nilpotent decomposable
Lie algebras in dimension five. The unimodular among them are the
ones in Table \ref{UniDecSol5}.
\end{table}

\begin{table}[ht]
\centering \caption{\label{UniDecSol5} $5$-dimensional decomposable
unimodular non-nilpotent algebras}
\begin{tabular}{|c|c|c|} \hline \hline
& $[X_i,X_j]$ & cpl. solv.
\\ \hline \hline $\mathfrak{g}_{3.4}^{-1} \oplus 2 \mathfrak{g}_1$ & $[X_1,X_3] = X_1,~
[X_2,X_3] = - X_2$ & yes
\\ \hline $\mathfrak{g}_{3.5}^{0} \oplus 2 \mathfrak{g}_1$  & $[X_1,X_3] = - X_2,~
[X_2,X_3] =  X_1$ & no
\\ \hline $\mathfrak{g}_{4.2}^{-2} \oplus \mathfrak{g}_1$ &
$[X_1,X_4] = -2 X_1,~ [X_2,X_4] = X_2,$ & yes \\ & $[X_3,X_4] = X_2
+ X_3$ &
\\ \hline $\mathfrak{g}_{4.5}^{p,-p-1} \oplus \mathfrak{g}_1$ & $[X_1,X_4] = X_1,~ [X_2,X_4] =
p X_2,$ & yes \\ &  $[X_3,X_4] = (-p-1) X_3,~ - \frac{1}{2} \le p <
0$ &
\\ \hline $\mathfrak{g}_{4.6}^{-2p,p} \oplus \mathfrak{g}_1$ & $[X_1,X_4] = -2p X_1,~ [X_2,X_4] = p
X_2 - X_3,$ & no \\ & $[X_3,X_4] = X_2 + p X_3,~ p>0$ &
\\ \hline $\mathfrak{g}_{4.8}^{-1} \oplus \mathfrak{g}_1$ & $[X_2,X_3] = X_1,~ [X_2,X_4] =
X_2,~ [X_3,X_4] = -X_3 $ & yes
\\ \hline $\mathfrak{g}_{4.9}^0 \oplus \mathfrak{g}_1$ & $[X_2,X_3] = X_1,~ [X_2,X_4] =
-X_3,~ [X_3,X_4] = X_2 $& no
\\ \hline
\end{tabular}
\end{table}

\begin{table}[h]
Except for $\mathfrak{g}_{4.2} \oplus \mathfrak{g}_1$, to each class
of algebras there is a connected and simply-connected solvable Lie
group admitting a lattice and has a Lie algebra belonging to the
class.
\end{table}

\begin{table}[h]
Mubarakzjanov's list in \cite{Mub63a} contains $33$ classes of
five-dimensional indecomposable non-nilpotent solvable Lie algebras,
namely $\mathfrak{g}_{5.7}, \ldots ,\mathfrak{g}_{5.39}$. We list
the unimodular among them in Tables \ref{UniAlab5} to \ref{Uni2ab5}.

~~ Note that there is a minor misprint in \cite{Mub63a} which has
been corrected in the list below.
\end{table}

\begin{table}[ht]
\centering \caption{\label{UniAlab5} $5$-dimensional indecomposable
unimodular almost abelian algebras }
\begin{tabular}{|c|c|c|} \hline \hline
& $[X_i,X_j]$ & cpl. solv.
\\ \hline \hline $\mathfrak{g}_{5.7}^{p,q,r}$ & $[X_1,X_5] = X_1,~
[X_2,X_5] = p X_2,$  & yes \\ & $[X_3,X_5] = q X_3,~ [X_4,X_5] = r
X_4,$ & \\ & $ -1 \le r \le q \le p \le 1,~ pqr \ne 0,~ p+q+r = -1$&
\\ \hline $\mathfrak{g}_{5.8}^{-1}$ & $[X_2,X_5] = X_1,~ [X_3,X_5] =
X_3,~ [X_4,X_5] = - X_4,$ & yes
\\ \hline $\mathfrak{g}_{5.9}^{p,-2-p}$ & $[X_1,X_5] = X_1,~ [X_2,X_5] =
X_1 + X_2,~ [X_3,X_5] = p X_3,$ & yes \\ & $[X_4,X_5] = (-2-p) X_4,~
p \ge -1$&
\\ \hline $\mathfrak{g}_{5.11}^{-3}$ & $[X_1,X_5] = X_1,~ [X_2,X_5] =
X_1 + X_2,$ & yes \\ & $[X_3,X_5] = X_2 + X_3,~ [X_4,X_5] = -3
X_4,$&
\\ \hline $\mathfrak{g}_{5.13}^{-1-2q,q,r}$ & $[X_1,X_5] = X_1,~ [X_2,X_5] =
(-1-2q) X_2, $ & no \\ & $[X_3,X_5] = q X_3 - r X_4,~ [X_4,X_5] = r
X_3 + q X_4,$ & \\ & $ -1 \le q \le 0,~ q \ne - \frac{1}{2},~ r \ne
0 $&
\\ \hline $\mathfrak{g}_{5.14}^{0}$ & $[X_2,X_5] = X_1,~ [X_3,X_5] =
- X_4,~ [X_4,X_5] = X_3$ & no
\\ \hline $\mathfrak{g}_{5.15}^{-1}$ & $[X_1,X_5] = X_1,~ [X_2,X_5] =
X_1 + X_2,$ & yes \\ & $[X_3,X_5] = - X_3,~ [X_4,X_5] = X_3 - X_4$ &
\\ \hline $\mathfrak{g}_{5.16}^{-1,q}$ & $[X_1,X_5] = X_1,~ [X_2,X_5] =
X_1 + X_2,$ & no \\ & $[X_3,X_5] = - X_3 - q X_4 ,~ [X_4,X_5] = q
X_3 - X_4$ & \\ & $q \ne 0$ &
\\ \hline $\mathfrak{g}_{5.17}^{p,-p,r}$ & $[X_1,X_5] = p X_1 - X_2 ,~ [X_2,X_5] =
X_1 + p X_2,$ & no \\ & $[X_3,X_5] = - p X_3 - r X_4 ,~ [X_4,X_5] =
r X_3 - p X_4$ & \\ & $r \ne 0$ &
\\ \hline $\mathfrak{g}_{5.18}^0$ & $[X_1,X_5] = - X_2 ,~ [X_2,X_5] =
X_1,$ & no \\ & $[X_3,X_5] = X_1 - X_4 ,~ [X_4,X_5] = X_2 + X_3$ &
\\ \hline
\end{tabular}
\end{table}

\begin{table}[ht]
\centering \caption{\label{UniAlnil15} $5$-dimensional
indecomposable unimodular algebras with nilradical
$\mathfrak{g}_{3.1} \oplus \mathfrak{g}_1$}
\begin{tabular}{|c|c|c|} \hline \hline
& $[X_i,X_j]$ & cpl. solv.
\\ \hline \hline $\mathfrak{g}_{5.19}^{p,-2p-2}$ & $[X_2,X_3] = X_1,~
[X_1,X_5] = (1+p) X_1,~ [X_2,X_5] = X_2,$  & yes \\ & $[X_3,X_5] = p
X_3,~ [X_4,X_5] = (-2p-2) X_4,~ p \ne -1$&
\\ \hline $\mathfrak{g}_{5.20}^{-1}$ & $[X_2,X_3] = X_1,~
[X_2,X_5] = X_2,~ [X_3,X_5] = - X_3,$  & yes \\ & $[X_4,X_5] = X_1$&
\\ \hline $\mathfrak{g}_{5.23}^{-4}$ & $[X_2,X_3] = X_1,~
[X_1,X_5] = 2 X_1,~ [X_2,X_5] = X_2 + X_3,$  & yes \\ & $[X_3,X_5] =
X_3,~ [X_4,X_5] = -4 X_4$&
\\ \hline $\mathfrak{g}_{5.25}^{p,4p}$ & $[X_2,X_3] = X_1,~
[X_1,X_5] = 2p X_1,~ [X_2,X_5] = p X_2 + X_3,$  & no \\ & $[X_3,X_5]
= - X_2 + p X_3,~ [X_4,X_5] = -4p X_4,~ p \ne 0$&
\\ \hline $\mathfrak{g}_{5.26}^{0,\varepsilon}$ & $[X_2,X_3] = X_1,~
[X_2,X_5] = X_3,~ [X_3,X_5] = - X_2 ,$  & no \\ & $ [X_4,X_5] =
\varepsilon X_1,~ \varepsilon = \pm 1$&
\\ \hline $\mathfrak{g}_{5.28}^{-\frac{3}{2}}$ & $[X_2,X_3] = X_1,~
[X_1,X_5] = - \frac{1}{2} X_1,~ [X_2,X_5] = -\frac{3}{2} X_2,$ & yes
\\ & $[X_3,X_5] = X_3 + X_4,~ [X_4,X_5] = X_4$ &
\\ \hline
\end{tabular}
\end{table}

\begin{table}[ht]
\centering \caption{\label{UniAlnil25} $5$-dimensional
indecomposable unimodular algebras with nilradical
$\mathfrak{g}_{4.1}$}
\begin{tabular}{|c|c|c|} \hline \hline
& $[X_i,X_j]$ & cpl. solv.
\\ \hline \hline $\mathfrak{g}_{5.30}^{-\frac{4}{3}}$ & $[X_2,X_4] =
X_1,~ [X_3,X_4] = X_2,~ [X_1,X_5] = \frac{2}{3} X_1,$  & yes \\
& $[X_2,X_5] = - \frac{1}{3} X_2,~ [X_3,X_5] = -\frac{4}{3} X_3,~
[X_4,X_5] = X_4,$&
\\ \hline
\end{tabular}
\end{table}

\begin{table}[ht]
\centering \caption{\label{Uni2ab5} $5$-dimensional indecomposable
unimodular algebras with nilradical $3 \mathfrak{g}_1$}
\begin{tabular}{|c|c|c|} \hline \hline
& $[X_i,X_j]$ & cpl. solv.
\\ \hline \hline $\mathfrak{g}_{5.33}^{-1,-1}$ & $[X_1,X_4] = X_1,~
[X_3,X_4] = - X_3,$  & yes \\ & $[X_2,X_5] = X_2,~ [X_3,X_5] = -
X_3$&
\\ \hline $\mathfrak{g}_{5.35}^{-2,0}$ & $[X_1,X_4] = -2 X_1,~
[X_2,X_4] = X_2,~[X_3,X_4] = X_3$ & no \\ & $[X_2,X_5] = -X_3,~
[X_3,X_5] = X_2$ &
\\ \hline
\end{tabular}
\end{table}

\begin{table}[ht]
There are ten classes of decomposable nilpotent Lie algebras in
dimension six: $6 \mathfrak{g}_1$, $\mathfrak{g}_{3.1} \oplus 3
\mathfrak{g}_1$, $2 \mathfrak{g}_{3.1}$, $\mathfrak{g}_{4.1} \oplus
2 \mathfrak{g}_1$ and $\mathfrak{g}_{5.i} \oplus \mathfrak{g}_1$ for
$i \in \{1,\ldots 6\}$.
\end{table}

\begin{table}[ht]
Tables \ref{Nil6} and \ref{Nil62} contain the six-dimensional
indecomposable nilpotent real Lie algebras. They base on Morozov's
classification in \cite{Mor}, where nilpotent algebras over a field
of characteristic zero are determined. Note that over $\mathbb{R}$,
there is only one isomorphism class of Morozov's indecomposable type
$5$ resp.\ type $10$ and type $14$ resp.\ $18$ splits into two
non-isomorphic ones.
\end{table}

\begin{table}[ht]
\centering \caption{\label{Nil6} $6$-dimensional indecomposable
nilpotent algebras}
\begin{tabular}{|c|c|} \hline \hline
& $[X_i,X_j]$
\\ \hline \hline $\mathfrak{g}_{6.N1}$ & $[X_1,X_2] = X_3,~ [X_1,X_3] =
X_4,~ [X_1,X_5] = X_6$
\\ \hline $\mathfrak{g}_{6.N2}$ & $[X_1,X_2] = X_3,~ [X_1,X_3] = X_4,~
[X_1,X_4] = X_5,~ [X_1,X_5] = X_6$
\\ \hline $\mathfrak{g}_{6.N3}$ & $[X_1,X_2] = X_6,~
[X_1,X_3] = X_4,~ [X_2,X_3] = X_5$
\\ \hline $\mathfrak{g}_{6.N4}$ & $[X_1,X_2] = X_5,~
[X_1,X_3] = X_6,~ [X_2,X_4] = X_6$
\\ \hline $\mathfrak{g}_{6.N5}$ & $[X_1,X_3] = X_5,~
[X_1,X_4] = X_6,~ [X_2,X_3] = - X_6,~ [X_2,X_4] = X_5$
\\ \hline $\mathfrak{g}_{6.N6}$ & $[X_1,X_2] = X_6,~
[X_1,X_3] = X_4,~ [X_1,X_4] = X_5,~ [X_2,X_3] = X_5$
\\ \hline $\mathfrak{g}_{6.N7}$ & $[X_1,X_3] = X_4,~
[X_1,X_4] = X_5,~ [X_2,X_3] = X_6$
\\ \hline $\mathfrak{g}_{6.N8}$ & $[X_1,X_2] = X_3 + X_5,~
[X_1,X_3] = X_4,~ [X_2,X_5] = X_6$
\\ \hline $\mathfrak{g}_{6.N9}$ & $[X_1,X_2] = X_3,~ [X_1,X_3] = X_4,~ [X_1,X_5] =
X_6,~ [X_2,X_3] = X_5$
\\ \hline $\mathfrak{g}_{6.N10}$ & $[X_1,X_2] = X_3,~
[X_1,X_3] = X_5,~ [X_1,X_4] = X_6,$ \\ & $[X_2,X_3] = - X_6,~
[X_2,X_4] = X_5$
\\ \hline $\mathfrak{g}_{6.N11}$ & $[X_1,X_2] = X_3,~ [X_1,X_3] = X_4, ~
[X_1,X_4] = X_5,~ [X_2,X_3] = X_6$
\\ \hline $\mathfrak{g}_{6.N12}$ & $[X_1,X_3] = X_4,~ [X_1,X_4] = X_6, ~
[X_2,X_5] = X_6$
\\ \hline $\mathfrak{g}_{6.N13}$ & $[X_1,X_2] = X_5,~ [X_1,X_3] = X_4, ~
[X_1,X_4] = X_6,~ [X_2,X_5] = X_6$
\\ \hline $\mathfrak{g}_{6.N14}^{\pm 1}$ & $[X_1,X_3] = X_4,~ [X_1,X_4] = X_6, ~
[X_2,X_3] = X_5,~ [X_2,X_5] = \pm X_6$
\\ \hline $\mathfrak{g}_{6.N15}$ & $[X_1,X_2] = X_3 + X_5,~ [X_1,X_3] = X_4, ~
[X_1,X_4] = X_6,~ [X_2,X_5] = X_6$
\\ \hline
\end{tabular}
\end{table}
\clearpage
\begin{table}[ht]
\centering \caption{\label{Nil62} $6$-dimensional indecomposable
nilpotent algebras (continued)}
\begin{tabular}{|c|c|} \hline \hline
& $[X_i,X_j]$
\\ \hline \hline
$\mathfrak{g}_{6.N16}$ & $[X_1,X_3] = X_4,~ [X_1,X_4] = X_5, ~
[X_1,X_5] = X_6,$ \\ & $[X_2,X_3] = X_5,~ [X_2,X_4] = X_6$
\\ \hline $\mathfrak{g}_{6.N17}$ & $[X_1,X_2] = X_3,~ [X_1,X_3] = X_4, ~
[X_1,X_4] = X_6,~ [X_2,X_5] = X_6$
\\ \hline $\mathfrak{g}_{6.N18}^{\pm 1}$ & $[X_1,X_2] = X_3,~ [X_1,X_3] = X_4,
~ [X_1,X_4] = X_6,$ \\ & $[X_2,X_3] = X_5,~ [X_2,X_5] = \pm X_6$
\\ \hline $\mathfrak{g}_{6.N19}$ & $[X_1,X_2] = X_3,~ [X_1,X_3] = X_4,
~ [X_1,X_4] = X_5,$ \\ & $[X_1,X_5] = X_6,~ [X_2,X_3] = X_6$
\\ \hline $\mathfrak{g}_{6.N20}$ & $[X_1,X_2] = X_3,~ [X_1,X_3] = X_4,
~ [X_1,X_4] = X_5,$ \\ & $[X_1,X_5] = X_6,~ [X_2,X_3] = X_5,~
[X_2,X_4] = X_6$
\\ \hline $\mathfrak{g}_{6.N21}$ & $[X_1,X_2] = X_3,~ [X_1,X_5] = X_6,
~ [X_2,X_3] = X_4,$ \\ & $[X_2,X_4] = X_5,~ [X_3,X_4] = X_6$
\\ \hline $\mathfrak{g}_{6.N22}$ & $[X_1,X_2] = X_3,~ [X_1,X_3] = X_5,
~ [X_1,X_5] = X_6,$ \\ & $[X_2,X_3] = X_4,~ [X_2,X_4] = X_5,~
[X_3,X_4] = X_6$
\\ \hline
\end{tabular}
\end{table}
\vspace{10mm}

\begin{table}[ht]
Mubarakzjanov's list in \cite{Mub63b} contains $99$ classes of
six-dimensional indecomposable almost nilpotent Lie algebras, namely
$\mathfrak{g}_{6.1}, \ldots ,\mathfrak{g}_{6.99}$.

~~ As first remarked by Turkowski, there is one algebra missing. The
complete (and partly corrected) list can be found in the article
\cite{CS2} of Campoamor-Stursberg\footnotemark[13], where the
missing algebra is denoted by $\mathfrak{g}_{6.92}^*$.

~~ We list the unimodular among this $100$ algebras in Tables
\ref{UniAlab6} to \ref{letzteAlgebra} (where some minor misprints
have been corrected). Note that there is no table with Lie algebras
with nilradical $\mathfrak{g}_{5.6}$ since the only such algebra is
not unimodular. \vspace{6mm}

\centering \caption{\label{UniAlab6} $6$-dimensional indecomposable
unimodular almost abelian algebras}
\begin{tabular}{|c|c|c|} \hline \hline
& $[X_i,X_j]$ & cpl. solv.
\\ \hline \hline $\mathfrak{g}_{6.1}^{a,b,c,d}$ & $[X_1,X_6] = X_1,~
[X_2,X_6] = a X_2,~ [X_3,X_6] = b X_3,$  & yes \\ & $[X_4,X_6] = c
X_4,~ [X_5,X_6] = d X_5$ & \\ & $ 0 < |d| \le |c| \le |b| \le |a|
\le 1,~ a+b+c+d = -1$&
\\ \hline
\end{tabular}
\end{table}
\footnotetext[13]{The author wishes to express his gratitude to R.
Campoamor-Stursberg for providing him with copies of \cite{CS2} and
\cite{Mub63b}.}
\clearpage
\begin{table}[t]
\centering \caption{\label{UniAlab6c} $6$-dimensional indecomposable
unimodular almost abelian algebras (continued)}
\begin{tabular}{|c|c|c|} \hline \hline
& $[X_i,X_j]$ & c. s.
\\ \hline \hline $\mathfrak{g}_{6.2}^{a,c,d}$ & $[X_1,X_6] = a X_1,~
[X_2,X_6] = X_1 + a X_2,~ [X_3,X_6] = X_3,$& yes \\ &$[X_4,X_6] = c
X_4,~ [X_5,X_6] = d X_5$&\\&$ 0 < |d| \le |c| \le 1,~ 2a+c+d = -1$&
\\ \hline $\mathfrak{g}_{6.3}^{-\frac{d+1}{3},d}$ & $[X_1,X_6] =
-\frac{d+1}{3} X_1,~ [X_2,X_6] = X_1 -\frac{d+1}{3} X_2,$ & yes \\
&$[X_3,X_6] = X_2 -\frac{d+1}{3} X_3,~ [X_4,X_6] = X_4,~ [X_5,X_6] = d X_5$& \\
&$ 0 < |d| \le 1,$&
\\ \hline $\mathfrak{g}_{6.4}^{-\frac{1}{4}}$ & $[X_1,X_6] = -\frac{1}{4}
X_1,~ [X_2,X_6] = X_1 - \frac{1}{4} X_2,$& yes \\ &$[X_3,X_6] = X_2
- \frac{1}{4} X_3,~ [X_4,X_6] = X_3 - \frac{1}{4} X_4,~ [X_5,X_6] =
X_5$&
\\ \hline $\mathfrak{g}_{6.6}^{a,b}$ & $[X_1,X_6] = X_1,~
[X_2,X_6] = a X_2,~ [X_3,X_6] = X_2 + a X_3,$  & yes \\ & $[X_4,X_6]
= b X_4,~ [X_5,X_6] = X_4 + b X_5,~ a \le b,~ a+b = -\frac{1}{2}$ &
\\ \hline $\mathfrak{g}_{6.7}^{a,-\frac{2}{3}a}$ & $[X_1,X_6] = a X_1,~
[X_2,X_6] = X_1 + a X_2,~ [X_3,X_6] = X_2 +a X_3,$  & yes \\ &
$[X_4,X_6] = -\frac{3}{2} a X_4,~ [X_5,X_6] = X_4 - \frac{3}{2} a
X_5,~ a \ne 0$ &
\\ \hline $\mathfrak{g}_{6.8}^{a,b,c,p}$ & $[X_1,X_6] = a X_1,~
[X_2,X_6] = b X_2,~ [X_3,X_6] = c X_3,$  & no \\ & $[X_4,X_6] = p
X_4 - X_5,~ [X_5,X_6] = X_4 + p X_5$ & \\ & $ 0 < |c| \le |b| \le
|a|,~ a+b+c+2p = 0$&
\\ \hline $\mathfrak{g}_{6.9}^{a,b,p}$ & $[X_1,X_6] = a X_1,~
[X_2,X_6] = b X_2,~ [X_3,X_6] = X_2 + b X_3,$  & no \\ & $[X_4,X_6]
= p X_4 - X_5,~ [X_5,X_6] = X_4 + p X_5,$ & \\ & $a \ne 0,~ a+2b+2p
= 0$&
\\ \hline $\mathfrak{g}_{6.10}^{a,-\frac{3}{2}a}$ & $[X_1,X_6] = a X_1,~
[X_2,X_6] = X_1 + a X_2,~ [X_3,X_6] = X_2 + a X_3,$  & no \\ &
$[X_4,X_6] = -\frac{3}{2} a X_4 - X_5,~ [X_5,X_6] = X_4 -
\frac{3}{2} a X_5$ &
\\ \hline $\mathfrak{g}_{6.11}^{a,p,q,s}$ & $[X_1,X_6] = a X_1,~ [X_2,X_6] =
p X_2 - X_3,~ [X_3,X_6] = X_2 + p X_3,$ & no \\ & $[X_4,X_6] = q X_4
- s X_5,~ [X_5,X_6] = s X_4 + q X_5,$&\\& $as \ne 0, a+2p+2q = 0$&
\\ \hline $\mathfrak{g}_{6.12}^{-4p,p}$ & $[X_1,X_6] = -4p X_1,~ [X_2,X_6] =
p X_2 - X_3,$ & no \\ & $[X_3,X_6] = X_2 + p X_3,~ [X_4,X_6] = X_2 +
p X_4 - X_5,$ & \\ &$[X_5,X_6] = X_3 + X_4 + p X_5,~ p \ne 0$&
\\ \hline
\end{tabular}
\end{table}
\begin{table}[b]
\centering \caption{\label{UniAlnil16} $6$-dimensional
indecomposable unimodular algebras with nilradical
$\mathfrak{g}_{3.1} \oplus 2 \mathfrak{g}_1$}
\begin{tabular}{|c|c|c|} \hline \hline
& $[X_i,X_j]$ & cpl. solv.
\\ \hline \hline $\mathfrak{g}_{6.13}^{a,b,h}$ & $[X_2,X_3] = X_1,~
[X_1,X_6] = (a+b) X_1,~ [X_2,X_6] = a X_2,$  & yes \\ & $[X_3,X_6] =
b X_3,~ [X_4,X_6] = X_4,~ [X_5,X_6] = h X_5,$& \\ & $a \ne 0,~ 2a +
2b +h = -1$&
\\ \hline $\mathfrak{g}_{6.14}^{a,b}$ & $[X_2,X_3] = X_1,~
[X_1,X_6] = (a+b) X_1,~ [X_2,X_6] = a X_2,$  & yes \\ & $[X_3,X_6] =
b X_3,~ [X_4,X_6] = X_4,~ [X_5,X_6] = X_1 + (a+b) X_5,$& \\ & $a \ne
0,~ a + b = -\frac{1}{3}$&
\\ \hline
\end{tabular}
\end{table}
\clearpage
\begin{table}[ht]
\centering \caption{\label{UniAlnil16c} $6$-dimensional
indecomposable unimodular algebras with nilradical
$\mathfrak{g}_{3.1} \oplus 2 \mathfrak{g}_1$ (continued)}
\begin{tabular}{|c|c|c|} \hline \hline
& $[X_i,X_j]$ & c. s.
\\ \hline \hline $\mathfrak{g}_{6.15}^{-1}$ & $[X_2,X_3] = X_1,~
[X_2,X_6] = X_2 + X_4,$& yes \\ & $[X_3,X_6] = - X_3 + X_5,~
[X_4,X_6] = X_4,~ [X_5,X_6] = - X_5,$&
\\ \hline $\mathfrak{g}_{6.17}^{-\frac{1}{2},0}$ & $[X_2,X_3] = X_1,~
[X_1,X_6] = -\frac{1}{2} X_1,~ [X_2,X_6] = -\frac{1}{2} X_2,$& yes \\
& $[X_3,X_6] = X_4,~ [X_5,X_6] = X_5,$&
\\ \hline $\mathfrak{g}_{6.18}^{a,-2a-3}$ & $[X_2,X_3] = X_1,~
[X_1,X_6] = (1+a) X_1,~ [X_2,X_6] = a X_2,$& yes \\ & $[X_3,X_6] =
X_3 + X_4,~ [X_4,X_6] = X_4,$& \\ & $[X_5,X_6] = -(2a+3) X_5, ~ a
\ne -\frac{3}{2}$&
\\ \hline $\mathfrak{g}_{6.19}^{-\frac{4}{3}}$ & $[X_2,X_3] = X_1,~
[X_1,X_6] = -\frac{1}{3} X_1,~ [X_2,X_6] = -\frac{4}{3} X_2,$& yes \\
& $[X_3,X_6] = X_3 + X_4,~ [X_4,X_6] = X_4,~ [X_5,X_6] = X_1
-\frac{1}{3} X_5$&
\\ \hline $\mathfrak{g}_{6.20}^{-3}$ & $[X_2,X_3] = X_1,~
[X_1,X_6] = X_1,~ [X_3,X_6] = X_3 + X_4,$& yes \\
& $[X_4,X_6] = X_1 + X_4,~ [X_5,X_6] = -3 X_5$&
\\ \hline $\mathfrak{g}_{6.21}^{a}$ & $[X_2,X_3] = X_1,~
[X_1,X_6] = 2 a X_1,~ [X_2,X_6] = a X_2 + X_3,$& yes \\
& $[X_3,X_6] = a X_3 ,~ [X_4,X_6] = X_4,~ [X_5,X_6] = -(4a+1) X_5$&
\\ &$ a \ne -\frac{1}{4}$&
\\ \hline $\mathfrak{g}_{6.22}^{-\frac{1}{6}}$ & $[X_2,X_3] = X_1,~
[X_1,X_6] = -\frac{1}{3} X_1,~ [X_2,X_6] = -\frac{1}{6} X_2 + X_3,$& yes \\
& $[X_3,X_6] = -\frac{1}{6} X_3 ,~ [X_4,X_6] = X_4,~ [X_5,X_6] = X_1
- \frac{1}{3} X_5$&
\\ \hline $\mathfrak{g}_{6.23}^{a,-7a,\varepsilon}$ & $[X_2,X_3] = X_1,~
[X_1,X_6] = 2 a X_1,~ [X_2,X_6] = a X_2 + X_3,$& yes \\
& $[X_3,X_6] = a X_3 + X_4 ,~ [X_4,X_6] = a X_4,$&\\&$[X_5,X_6] =
\varepsilon X_1 - 5a X_5,~ \varepsilon a = 0$&
\\ \hline $\mathfrak{g}_{6.25}^{b,-1-b}$ & $[X_2,X_3] = X_1,~
[X_1,X_6] = -b X_1,$& yes \\
& $[X_2,X_6] = X_2,~ [X_3,X_6] = -(1+b) X_3,$&\\&$[X_4,X_6] = b X_4
+ X_5,~ [X_5,X_6] = b X_5$&
\\ \hline $\mathfrak{g}_{6.26}^{-1}$ & $[X_2,X_3] = X_1,~
[X_2,X_6] = X_2,~ [X_3,X_6] = - X_3$& yes \\
& $[X_4,X_6] = X_5,~ [X_5,X_6] = X_1$&
\\ \hline $\mathfrak{g}_{6.27}^{-2b,b}$ & $[X_2,X_3] = X_1,~
[X_1,X_6] = -b X_1,~ [X_2,X_6] = -2b X_2,$& yes \\
& $[X_3,X_6] = b X_3 + X_4,~ [X_4,X_6] = b X_4 + X_5,$&\\&
$[X_5,X_6] = b X_5,~ b\ne0$&
\\ \hline $\mathfrak{g}_{6.28}^{-2}$ & $[X_2,X_3] = X_1,~
[X_1,X_6] = 2 X_1,~ [X_2,X_6] = X_2 + X_3,$& yes \\
& $[X_3,X_6] = X_3,~ [X_4,X_6] = -2 X_4 + X_5,~ [X_5,X_6] = -2 X_5$&
\\ \hline $\mathfrak{g}_{6.29}^{-2b,b,\varepsilon}$ & $[X_2,X_3] = X_1,~
[X_1,X_6] = -b X_1,~ [X_2,X_6] = -2b X_2,$& yes \\
& $[X_3,X_6] = b X_3 + X_4,~ [X_4,X_6] = b X_4 + X_5,$&\\&$
[X_5,X_6] = \varepsilon X_1 + b X_5,~ \varepsilon b = 0~ (?)$&
\\ \hline $\mathfrak{g}_{6.32}^{a,-6a-h,h,\varepsilon}$ & $[X_2,X_3] = X_1,~
[X_1,X_6] = 2a X_1,~ [X_2,X_6] = a X_2 + X_3,$& no \\
& $[X_3,X_6] = - X_2 + a X_3,~ [X_4,X_6] = \varepsilon X_1 + (2a +
h) X_4,$&\\&$ [X_5,X_6] = -(6a+h) X_5,~ a > -\frac{1}{4} h,~
\varepsilon h = 0$&
\\ \hline
\end{tabular}
\end{table}
\clearpage
\begin{table}[t]
\centering \caption{\label{UniAlnil16cc} $6$-dimensional
indecomposable unimodular algebras with nilradical
$\mathfrak{g}_{3.1} \oplus 2 \mathfrak{g}_1$ (continued)}
\begin{tabular}{|c|c|c|} \hline \hline
& $[X_i,X_j]$ & c. s.
\\ \hline \hline $\mathfrak{g}_{6.33}^{a,-6a}$ & $[X_2,X_3] = X_1,~
[X_1,X_6] = 2a X_1,~ [X_2,X_6] = a X_2 + X_3,$& no \\
& $[X_3,X_6] = - X_2 + a X_3,~ [X_4,X_6] = -6 a X_4,$&\\&$[X_5,X_6]
= X_1 + 2a X_5,~ a \ge 0$&
\\ \hline $\mathfrak{g}_{6.34}^{a,-4a,\varepsilon}$ & $[X_2,X_3] = X_1,~
[X_1,X_6] = 2a X_1,~ [X_2,X_6] = a X_2 + X_3,$ & no \\
& $[X_3,X_6] = -X_2 + a X_3,~ [X_4,X_6] = -2a X_4,$&\\&$ [X_5, X_6]=
\varepsilon X_1 - 2 a X_5,~ \varepsilon a =0$&
\\ \hline $\mathfrak{g}_{6.35}^{a,b,c}$ & $[X_2,X_3] = X_1,~
[X_1,X_6] = (a+b) X_1,~ [X_2,X_6] = a X_2,$ & no \\
& $[X_3,X_6] = b X_3,~ [X_4,X_6] = c X_4 + X_5,$&\\&$ [X_5, X_6]=
-X_4 + c X_5,~a+b+c=0,~a^2+b^2\ne 0$&
\\ \hline $\mathfrak{g}_{6.36}^{a,-2a}$ & $[X_2,X_3] = X_1,~
[X_1,X_6] = 2a X_1,~ [X_2,X_6] = a X_2 + X_3,$ & no \\
& $[X_3,X_6] = a X_3,~ [X_4,X_6] = -2a X_4 + X_5,$&\\&$ [X_5, X_6]=
-X_4 - 2a X_5$&
\\ \hline $\mathfrak{g}_{6.37}^{-a,-2a,s}$ & $[X_2,X_3] = X_1,~
[X_1,X_6] = 2a X_1,~ [X_2,X_6] = a X_2 + X_3,$ & no \\
& $[X_3,X_6] = -X_2 + a X_3,~ [X_4,X_6] = -2a X_4 + s X_5,$&\\&$
[X_5, X_6]= -s X_4 - 2a X_5,~ s \ne 0$&
\\ \hline $\mathfrak{g}_{6.38}^{0}$ & $[X_2,X_3] = X_1,~
[X_2,X_6] = X_3 + X_4 ,$ & no \\
& $[X_3,X_6] = -X_2 + X_5,~ [X_4,X_6] =  X_5,~ [X_5, X_6]= - X_4$&
\\ \hline
\end{tabular}
\end{table}

\begin{table}[b]
\centering \caption{\label{UniAlnil26} $6$-dimensional
indecomposable unimodular algebras with nilradical
$\mathfrak{g}_{4.1} \oplus \mathfrak{g}_1$}
\begin{tabular}{|c|c|c|} \hline \hline
& $[X_i,X_j]$ & c. s.
\\ \hline \hline $\mathfrak{g}_{6.39}^{-4-3h,h}$ & $[X_1,X_5] = X_2,~ [X_4,X_5] =
X_1,~ [X_1,X_6] = (1+h) X_1,$ & yes \\ & $[X_2,X_6] = (2+h) X_2,~
[X_3,X_6] = -(4+3h) X_3,$ & \\ &$[X_4,X_6] = h X_4,~ [X_5,X_6] =
X_5,~ h \ne -\frac{4}{3}$&
\\ \hline $\mathfrak{g}_{6.40}^{-\frac{3}{2}}$ & $[X_1,X_5] = X_2,~ [X_4,X_5] = X_1
$& yes \\ & $[X_1,X_6] = -\frac{1}{2} X_1,~ [X_2,X_6] = \frac{1}{2} X_2,$ & \\
&$[X_3,X_6] = X_2 + \frac{1}{2} X_3,~ [X_4,X_6] = -\frac{3}{2} X_4,~
[X_5,X_6] = X_5$&
\\ \hline $\mathfrak{g}_{6.41}^{-1}$ & $[X_1,X_5] = X_2,~ [X_4,X_5] = X_1
$& yes \\ & $[X_2,X_6] = X_2,~ [X_3,X_6] = - X_3,$ & \\
&$[X_4,X_6] = X_3 - X_4,~ [X_5,X_6] = X_5$&
\\ \hline
\end{tabular}
\end{table}

\begin{table}[ht]
\centering \caption{\label{UniAlnil26c} $6$-dimensional
indecomposable unimodular algebras with nilradical
$\mathfrak{g}_{4.1} \oplus \mathfrak{g}_1$ (continued)}
\begin{tabular}{|c|c|c|} \hline \hline
& $[X_i,X_j]$ & c. s.
\\ \hline \hline $\mathfrak{g}_{6.42}^{-\frac{5}{3}}$
& $[X_1,X_5] = X_2,~ [X_4,X_5] = X_1,$ & yes
\\ & $[X_1,X_6] = -\frac{2}{3} X_1,~[X_2,X_6] = \frac{1}{3} X_2,~
[X_3,X_6] = X_3$ & \\ &$[X_4,X_6] = -\frac{5}{3} X_4,~ [X_5,X_6] =
X_3 + X_5$&
\\ \hline $\mathfrak{g}_{6.44}^{-7}$ & $[X_1,X_5] = X_2,~ [X_4,X_5] =
X_1,$ & yes
\\ & $[X_1,X_6] = 2 X_1,~[X_2,X_6] =3 X_2,~
[X_3,X_6] = -7 X_3$ & \\ &$[X_4,X_6] = X_4,~ [X_5,X_6] = X_4 + X_5$&
\\ \hline $\mathfrak{g}_{6.47}^{-3,\varepsilon}$ & $[X_1,X_5] = X_2,
~ [X_4,X_5] = X_1,$ & yes \\ & $[X_1,X_6] = X_1,~[X_2,X_6] = X_2,~
[X_3,X_6] = -3 X_3$ & \\ &$[X_4,X_6] = \varepsilon X_2 + X_4,~
\varepsilon \in \{0, \pm1\}$&
\\ \hline
\end{tabular}
\end{table}
\begin{table}[ht]
\centering \caption{\label{UniAlnil36} $6$-dimensional
indecomposable unimodular algebras with nilradical
$\mathfrak{g}_{5.1}$}
\begin{tabular}{|c|c|c|} \hline \hline
& $[X_i,X_j]$ & c. s.
\\ \hline \hline $\mathfrak{g}_{6.54}^{2(1+l),l}$
& $[X_3,X_5] = X_1,~ [X_4,X_5] = X_2,$ & yes
\\ & $[X_1,X_6] = X_1,~ [X_2,X_6] = l X_2,~
[X_3,X_6] = (-1-2l) X_3$ & \\ &$[X_4,X_6] = (-2-l) X_4,~ [X_5,X_6] =
2(1+l) X_5$&
\\ \hline $\mathfrak{g}_{6.55}^{-4}$ & $[X_3,X_5] = X_1,~ [X_4,X_5] =
X_2,$ & yes
\\ & $[X_1,X_6] = X_1,~ [X_2,X_6] = -3 X_2,~
[X_3,X_6] = 4 X_3$ & \\ &$[X_4,X_6] = X_1 + X_4,~ [X_5,X_6] = -3
X_5$&
\\ \hline $\mathfrak{g}_{6.56}^{\frac{4}{3}}$ & $[X_3,X_5] = X_1,~ [X_4,X_5] =
X_2,$ & yes
\\ & $[X_1,X_6] = X_1,~ [X_2,X_6] = -\frac{1}{3} X_2 ,~
[X_3,X_6] = X_2 - \frac{1}{3} X_3$ & \\ &$[X_4,X_6] = - \frac{5}{3}
X_4,~ [X_5,X_6] = \frac{4}{3} X_5$&
\\ \hline $\mathfrak{g}_{6.57}^{-\frac{2}{3}}$ & $[X_3,X_5] = X_1,~ [X_4,X_5] =
X_2,$ & yes
\\ & $[X_1,X_6] = X_1,~ [X_2,X_6] = -\frac{4}{3} X_2 ,~
[X_3,X_6] = \frac{5}{3} X_3$ & \\ &$[X_4,X_6] = - \frac{2}{3} X_4,~
[X_5,X_6] = X_4 - \frac{2}{3} X_5$&
\\ \hline $\mathfrak{g}_{6.61}^{-\frac{3}{4}}$ & $[X_3,X_5] = X_1,~ [X_4,X_5] =
X_2,$ & yes
\\ & $[X_1,X_6] = 2 X_1,~ [X_2,X_6] = -\frac{3}{2} X_2 ,~
[X_3,X_6] = X_3$ & \\ &$[X_4,X_6] = - \frac{5}{2} X_4,~ [X_5,X_6] =
X_3 + X_5$&
\\ \hline $\mathfrak{g}_{6.63}^{-1}$ & $[X_3,X_5] = X_1,~ [X_4,X_5] =
X_2,$ & yes
\\ & $[X_1,X_6] = X_1,~ [X_2,X_6] = - X_2 ,~
[X_3,X_6] = X_3$ & \\ &$[X_4,X_6] = X_2 - X_4$&
\\ \hline $\mathfrak{g}_{6.65}^{4l,l}$ & $[X_3,X_5] = X_1,~ [X_4,X_5] =
X_2,$ & yes
\\ & $[X_1,X_6] = l X_1 + X_2,~ [X_2,X_6] = l X_2 ,~
[X_3,X_6] = -3l X_3 + X_4$ & \\ &$[X_4,X_6] = -3l X_4,~ [X_5,X_6] =
4l X_5$&
\\ \hline
\end{tabular}
\end{table}
\begin{table}[ht]
\centering \caption{\label{UniAlnil36c} $6$-dimensional
indecomposable unimodular algebras with nilradical
$\mathfrak{g}_{5.1}$ (continued)}
\begin{tabular}{|c|c|c|} \hline \hline
& $[X_i,X_j]$ & cpl. solv.
\\ \hline \hline $\mathfrak{g}_{6.70}^{4p,p}$
& $[X_3,X_5] = X_1,~ [X_4,X_5] = X_2,$ & no
\\ & $[X_1,X_6] = p X_1 + X_2,~ [X_2,X_6] = -X_1 + p X_2,$& \\
&$[X_3,X_6] = -3p X_3 + X_4,~ [X_4,X_6] = -X_3 - 3p X_4,$&\\&$
[X_5,X_6] = 4p X_5$&
\\ \hline
\end{tabular}
\end{table}
\begin{table}[ht]
\centering \caption{\label{UniAlnil46} $6$-dimensional
indecomposable unimodular algebras with nilradical
$\mathfrak{g}_{5.2}$}
\begin{tabular}{|c|c|c|} \hline \hline
& $[X_i,X_j]$ & cpl. solv.
\\ \hline \hline $\mathfrak{g}_{6.71}^{-\frac{7}{4}}$ &
$[X_2,X_5]=X_1,~ [X_3,X_5]=X_2,~ [X_4,X_5]=X_3,$& yes \\& $[X_1,X_6]
= \frac{5}{4} X_1,~ [X_2,X_6] = \frac{1}{4} X_2,~ [X_3,X_6] = -
\frac{3}{4} X_3$ & \\ &$[X_4,X_6] = -\frac{7}{4} X_4,~ [X_5,X_6] =
X_5$&
\\ \hline
\end{tabular}
\end{table}

\begin{table}[ht]
\centering \caption{\label{UniAlnil56} $6$-dimensional
indecomposable unimodular algebras with nilradical
$\mathfrak{g}_{5.3}$}
\begin{tabular}{|c|c|c|} \hline \hline
& $[X_i,X_j]$ & cpl. solv.
\\ \hline \hline $\mathfrak{g}_{6.76}^{-1}$ & $[X_2,X_4] = X_3,~ [X_2,X_5] = X_1, ~
[X_4,X_5] = X_2$ &  yes \\ &$[X_1,X_6] = - X_1,~ [X_3,X_6] = X_3, $& \\
&$[X_4,X_6] = X_4,~ [X_5,X_6] = - X_5$&
\\ \hline $\mathfrak{g}_{6.78}$ & $[X_2,X_4] = X_3,~ [X_2,X_5] = X_1, ~
[X_4,X_5] = X_2$ &  yes \\ &$[X_1,X_6] = - X_1,~ [X_3,X_6] = X_3, $& \\
&$[X_4,X_6] = X_3 + X_4,~ [X_5,X_6] = - X_5$&
\\ \hline
\end{tabular}
\end{table}
\begin{table}[ht]
\centering \caption{\label{UniAlnil66} $6$-dimensional
indecomposable unimodular algebras with nilradical
$\mathfrak{g}_{5.4}$}
\begin{tabular}{|c|c|c|} \hline \hline
& $[X_i,X_j]$ & cpl. solv.
\\ \hline \hline $\mathfrak{g}_{6.83}^{0,l}$ &$[X_2,X_4] = X_1,~
[X_3,X_5] = X_1,$ & yes \\ &$[X_2,X_6] = l X_2,~ [X_3,X_6] = l X_3, $& \\
&$[X_4,X_6] = -l X_4,~ [X_5,X_6] = -X_4 -  lX_5$&
\\ \hline $\mathfrak{g}_{6.84}$ &$[X_2,X_4] = X_1,~
[X_3,X_5] = X_1,$ & yes \\ &$[X_2,X_6] = X_2,~ [X_4,X_6] = -X_4,~
[X_5,X_6] = X_3$&
\\ \hline
\end{tabular}
\end{table}
\begin{table}[ht]
\centering \caption{\label{UniAlnil66c} $6$-dimensional
indecomposable unimodular algebras with nilradical
$\mathfrak{g}_{5.4}$ (continued)}
\begin{tabular}{|c|c|c|} \hline \hline
& $[X_i,X_j]$ & cpl. solv.
\\ \hline \hline $\mathfrak{g}_{6.88}^{0,\mu_0,\nu_0}$ &$[X_2,X_4] = X_1,~
[X_3,X_5] = X_1,$ & cpl. solv.  \\ &$[X_2,X_6] = \mu_0 X_2 + \nu_0
X_3,~
[X_3,X_6] =-\nu_0 X_2 + \mu_0 X_3, $& $\Updownarrow$ \\
&$[X_4,X_6] = -\mu_0 X_4 + \nu_0 X_5,~ [X_5,X_6] = -\nu_0 X_4 -
\mu_0 X_5$& $\nu_0 = 0$
\\ \hline $\mathfrak{g}_{6.89}^{0,\nu_0,s}$ &$[X_2,X_4] = X_1,~
[X_3,X_5] = X_1,$ & cpl. solv. \\ &$[X_2,X_6] = s X_2,~ [X_3,X_6] =\nu_0 X_5, $& $\Updownarrow$ \\
&$[X_4,X_6] = -s X_4,~ [X_5,X_6] = -\nu_0 X_3$& $\nu_0 =0$
\\ \hline $\mathfrak{g}_{6.90}^{0,\nu_0}$ &$[X_2,X_4] = X_1,~
[X_3,X_5] = X_1,$ & cpl. solv. \\ &$[X_2,X_6] = X_4,~ [X_3,X_6] = \nu_0 X_5, $& $\Updownarrow$ \\
&$[X_4,X_6] = X_2,~ [X_5,X_6] = -\nu_0 X_3,~ \nu_0 \ne 1$& $\nu_0 =
0$
\\ \hline $\mathfrak{g}_{6.91}$ &$[X_2,X_4] = X_1,~
[X_3,X_5] = X_1,$ & no \\ &$[X_2,X_6] = X_4,~ [X_3,X_6] = X_5, $& \\
&$[X_4,X_6] = X_2,~ [X_5,X_6] = - X_3$&
\\ \hline $\mathfrak{g}_{6.92}^{0,\mu_0,\nu_0}$ &$[X_2,X_4] = X_1,~
[X_3,X_5] = X_1,$ & no \\ &$[X_2,X_6] = \nu_0 X_3,~ [X_3,X_6] = -\mu_0 X_2, $& \\
&$[X_4,X_6] = \mu_0 X_5,~ [X_5,X_6] = -\nu_0 X_4$&
\\ \hline $\mathfrak{g}_{6.92^*}^{0}$ &$[X_2,X_4] = X_1,~
[X_3,X_5] = X_1,$ & no \\ &$[X_2,X_6] = X_4,~ [X_3,X_6] =  X_5, $& \\
&$[X_4,X_6] = - X_2,~ [X_5,X_6] = - X_3$&
\\ \hline $\mathfrak{g}_{6.93}^{0,\nu_0}$ &$[X_2,X_4] = X_1,~
[X_3,X_5] = X_1,$ & cpl. solv. \\ &$[X_2,X_6] = X_4 + \nu_0 X_5,~
[X_3,X_6] = \nu_0 X_4, $& $\Updownarrow$\\
&$[X_4,X_6] = X_2 - \nu_0 X_3,~ [X_5,X_6] = - \nu_0 X_2$& $|\nu_0|
\le \frac{1}{2}$
\\ \hline
\end{tabular}
\end{table}
\begin{table}[b]
\centering \caption{\label{letzteAlgebra} $6$-dimensional
indecomposable unimodular algebras with nilradical
$\mathfrak{g}_{5.5}$}
\begin{tabular}{|c|c|c|} \hline \hline
& $[X_i,X_j]$ & cpl. solv.
\\ \hline \hline $\mathfrak{g}_{6.94}^{-2}$ &$[X_3,X_4] = X_1,~ [X_2,X_5] = X_1, ~
[X_3,X_5] = X_2$& yes \\ &$[X_2,X_6] = - X_2,~ [X_3,X_6] = -2 X_3, $& \\
&$[X_4,X_6] = 2 X_4,~ [X_5,X_6] =  X_5$&
\\ \hline
\end{tabular}
\end{table}

\begin{table}[hb]
The six-dimensional solvable Lie algebras with four-dimensional
nilradical were classified by Turkowski in \cite{T}. We list the
unimodular among them in Tables \ref{Nilrad41} -- \ref{Nilrad42c}.
Note that there is no table with Lie algebras with nilradical
$\mathfrak{g}_{4.1}$ since the only such algebra is not unimodular.

The equations for the twenty-fifth algebra in Turkowoski's list
contain a minor misprint that we have corrected here.
\end{table}

\begin{table}[ht]
\centering \caption{\label{Nilrad41} $6$-dimensional indecomposable
unimodular algebras with nilradical $4 \mathfrak{g}_{1}$}
\begin{tabular}{|c|c|c|} \hline \hline
& $[X_i,X_j]$ & c. s.
\\ \hline \hline $\mathfrak{g}_{6.101}^{a,b,c,d}$ &$[X_5,X_1] = a X_1,~
[X_5,X_2] = c X_2,~ [X_5,X_4] = X_4,$ & yes \\ &$[X_6,X_1] = b X_1,~
[X_6,X_2] = d X_2,~ [X_6,X_3] =  X_3,$&
\\ &$a+c = -1,~ b+d = -1,~ ab \ne 0,~ c^2 + d^2 \ne 0$&
\\ \hline $\mathfrak{g}_{6.102}^{-1,b,-2-b}$ &$[X_5,X_1] = - X_1,~
[X_5,X_2] = X_2,~ [X_5,X_3] = X_4,$ & yes \\ &$[X_6,X_1] = b X_1,~
[X_6,X_2] = (-2-b) X_2,$ & \\ &$[X_6,X_3] =  X_3,~ [X_6,X_4] = X_4$&
\\ \hline $\mathfrak{g}_{6.105}^{-2,-1}$ &$[X_5,X_1] = -2 X_1,~
[X_5,X_3] = X_3 + X_4,$ & yes \\ &$[X_5,X_4] = X_4,~ [X_6,X_1] = -
X_1,~ [X_6,X_2] =  X_2$ &
\\ \hline $\mathfrak{g}_{6.107}^{-1,b,0}$ &$[X_5,X_1] = -X_1,~
[X_5,X_2] = -X_2,~ [X_5,X_3] = X_3 + X_4,$ & no \\ &$[X_5,X_4] =
X_4,~ [X_6,X_1] = X_2,~ [X_6,X_2] = - X_1$ &
\\ \hline $\mathfrak{g}_{6.113}^{a,b,-a,d}$ &$[X_5,X_1] = a X_1,~
[X_5,X_2] = -a X_2,~ [X_5,X_3] = X_4,$ & no \\ &$[X_6,X_1] = b X_1,~
[X_6,X_2] = d X_2,~ [X_6,X_3] = X_3$ & \\ &$[X_6,X_4] = X_4,~ a^2 +
b^2 \ne 0,~ a^2 + d^2 \ne 0,~ b+d=-2$&
\\ \hline $\mathfrak{g}_{6.114}^{a,-1,-\frac{a}{2}}$ &$[X_5,X_1] = a X_1,~
[X_5,X_3] = -\frac{a}{2} X_3 + X_4,$ & no \\ &$[X_5,X_4] = - X_3 +
\frac{a}{2} X_4,~ [X_6,X_1] = - X_1,$ & \\
&$[X_6,X_2] = X_2,~ a \ne 0$&
\\ \hline $\mathfrak{g}_{6.115}^{-1,b,c,-c}$ &$[X_5,X_1] = X_1,~
[X_5,X_2] = X_2,$ & no \\ &$[X_5,X_3] = - X_3 +
b X_4,~ [X_5,X_4] = -b X_3 - X_4,$ & \\
&$[X_6,X_1] = c X_1 + X_2,~ [X_6,X_2] = -X_1 + c X_2$& \\ &$[X_6,
X_3] = -c X_3,~ [X_6,X_4] = -c X_4,~ b \ne 0$&
\\ \hline $\mathfrak{g}_{6.116}^{0,-1}$ &$[X_5,X_1] = X_2,~
[X_5,X_3] = X_4,~ [X_5,X_4] = -X_3$ & no \\ &$[X_6,X_1] = X_1,~
[X_6,X_2] = X_2,$ & \\
&$[X_6,X_3] = - X_3,~ [X_6,X_4] = -X_4$&
\\ \hline $\mathfrak{g}_{6.118}^{0,b,-1}$ &$[X_5,X_1] = X_2,~
[X_5,X_2] = -X_1,~ [X_5,X_3] = b X_4,$ & no \\&$[X_5,X_4] = -b X_3,~
[X_6, X_1] = X_1,~ [X_6,X_2] = X_2$ & \\
&$[X_6,X_3] = - X_3,~ [X_6,X_4] = -X_4,~ b\ne0$&
\\ \hline $\mathfrak{g}_{6.120}^{-1,-1}$ & $[X_5,X_2] = - X_2,~
[X_5,X_4] = X_4,~ [X_5,X_6] = X_1,$ & yes \\ &$[X_6,X_2] = - X_1,~
[X_6,X_3] = X_3$&
\\ \hline $\mathfrak{g}_{6.125}^{0,-2}$ & $[X_5,X_3] = X_4,~
[X_5,X_4] = - X_3,~ [X_5,X_6] = X_1,$ & no \\ &$[X_6,X_2] = -2 X_2,~
[X_6,X_3] = X_3,~ [X_6,X_4] =X_4 $&
\\ \hline
\end{tabular}
\end{table}

\begin{table}[hb]
\centering \caption{\label{Nilrad42} $6$-dimensional indecomposable
unimodular algebras with nilradical $\mathfrak{g}_{3.1} \oplus
\mathfrak{g}_1$}
\begin{tabular}{|c|c|c|} \hline \hline
& $[X_i,X_j]$ & c. s.
\\ \hline \hline $\mathfrak{g}_{6.129}^{-2,-2}$ &$[X_2,X_3] = X_1,~
[X_5,X_1] = X_1,~ [X_5,X_2] = X_2,$ & yes \\ &$[X_5,X_4] = -2 X_4,~
[X_6,X_1] = X_1,$&
\\ &$[X_6,X_3] =  X_3,~ [X_6,X_4] = - 2 X_4$&
\\ \hline
\end{tabular}
\end{table}

\begin{table}[hb]
\centering \caption{\label{Nilrad42c} $6$-dimensional indecomposable
unimodular algebras with nilradical $\mathfrak{g}_{3.1} \oplus
\mathfrak{g}_1$ (continued)}
\begin{tabular}{|c|c|c|} \hline \hline
& $[X_i,X_j]$ & c. s.
\\ \hline \hline $\mathfrak{g}_{6.135}^{0,-4}$ &$[X_2,X_3] = X_1,~
[X_5,X_2] = X_3,~ [X_5,X_3] = -X_1,$ & no \\ &$[X_6,X_1] = 2 X_1,~
[X_6,X_2] = X_2,$&
\\ &$[X_6,X_3] =  X_3,~ [X_6,X_4] = -4 X_4$&
\\ \hline
\end{tabular}
\end{table}

\begin{table}[h]
In the introduction of \cite{Mub63b}, Mubarakzjanov quotes his own
result that a six-dimensional solvable Lie algebra with
three-dimensional nilradical is decomposable. Therefore, by
Proposition \ref{dim Nil}, we have listed all unimodular
indecomposable solvable Lie algebras of dimension six.

The first Betti numbers of the six-dimensional unimodular
indecomposable Lie algebras are listed in Tables \ref{b16} --
\ref{b16cc}. The word ``always'' means that the certain value arises
independent of the parameters on which the Lie algebra depends, but
we suppose that the parameters are chosen such that Lie algebra is
unimodular. The word ``otherwise'' in the tables means that this value
arises for all parameters such that the Lie algebra is unimodular
and the parameters are not mentioned in another column of the Lie
algebra's row.
\end{table}

\begin{table}[hb]
\centering \caption{\label{b16}$b_1(\mathfrak{g}_{6.i})$ for
$\mathfrak{g}_{6.i}$ unimodular}
\begin{tabular}{|c|c|c|c|} \hline \hline
$i$ & $b_1 = 1$ & $b_1 = 2$ & $b_1 = 3$ %
\\ \hline
\hline $1$ & always & - & - \\
\hline $2$ & $a \ne 0$ & $a = 0$ & -\\
\hline $3$ & $d \ne -1$ & $d = -1$ & - \\
\hline $4$ & always & - & -\\
\hline $6$ & $a,b \ne 0$ & $a=-\frac{1}{2} \wedge b=0$ & - \\
\hline $7$ & always & - & - \\
\hline $8$ & always & - & - \\
\hline $9$ & $b \ne 0$ & $b=0$ & - \\
\hline $10$ & $a \ne 0$ & $a=0$ & - \\
\hline $11$ & always & - & - \\
\hline $12$ & always & - & - \\
\hline $13$ & $b\ne0 \wedge h\ne0$ & otherwise & $a=-\frac{1}{2} \wedge b=h=0$ \\
\hline $14$ & otherwise & $a=-\frac{1}{3} \wedge b=0$ & - \\
\hline $15$ & always & - & - \\
\hline $17$ & - & always & - \\
\hline $18$ & $a \ne 0$ & $a=0$ & - \\
\hline $19$ & always & - & - \\
\hline $20$ & - & always & - \\
\hline
\end{tabular}
\end{table}
\begin{table}[hb]
\centering \caption{\label{b16c}$b_1(\mathfrak{g}_{6.i})$ for
$\mathfrak{g}_{6.i}$ unimodular (continued)}
\begin{tabular}{|c|c|c|c|c|c|} \hline \hline
$i$ & $b_1 = 1$ & $b_1 = 2$ & $b_1 = 3$ & $b_1 = 4$ & $b_1 =5$ %
\\ \hline
\hline $21$ & $a\ne 0$ & $a=0$ & - & - & -\\
\hline $22$ & always & - & - & - & -\\
\hline $23$ & $a \ne 0$ & - & $a=0$ & - & -\\
\hline $25$ & $b \notin \{-1,0\}$ & $b \in \{-1,0\}$ & - & - & -\\
\hline $26$ & - & always & - & - & -\\
\hline $27$ & always & - & - & - & -\\
\hline $28$ & always & - & - & - & -\\
\hline $29$ & $b \ne 0$ & - & $b=0$ & - & -\\
\hline $32$ & $h \notin \{-2a,-6a\}$ & otherwise & - & - & -\\
\hline $33$ & $a \ne 0$ & - & $a=0$ & - & -\\
\hline $34$ & $a \ne 0$ & - & $a=0$ & - & -\\
\hline $35$ & $a,b \ne 0$ & otherwise & - & - & -\\
\hline $36$ & $a \ne 0$ & $a=0$ & - & - & -\\
\hline $37$ & always & - & - & - & -\\
\hline $38$ & always & - & - & - & -\\
\hline $39$ & $h \ne 0$ & $h = 0$ & - & - & -\\
\hline $40$ & always & - & - & - & -\\
\hline $41$ & always & - & - & - & -\\
\hline $42$ & always & - & - & - & -\\
\hline $44$ & always & - & - & - & -\\
\hline $47$ & - & always & - & - & -\\
\hline $54$ & $l \notin \{-2,-1,-\frac{1}{2}\}$
&$l \in \{-2,-1,-\frac{1}{2}\}$& - & - & -\\
\hline $55$ & always & - & - & - & -\\
\hline $56$ & always & - & - & - & -\\
\hline $57$ & always & - & - & - & -\\
\hline $61$ & always & - & - & - & -\\
\hline $63$ & - & always & - & - & -\\
\hline $65$ & $l \ne 0$ & - & $l=0$ & - & -\\
\hline $70$ & $p \ne 0$ & $p=0$ & - & - & -\\
\hline $71$ & always & - & - & - & -\\
\hline $76$ & always & - & - & - & -\\
\hline $78$ & always & - & - & - & -\\
\hline $83$ & $l\ne0$& - & - & $l=0$& - \\
\hline $84$ & - & always & - & -& -\\
\hline $88$ & $\mu_0 \ne 0 \vee \nu_0 \ne 0$ & - & - & -
&  $\mu_0 = \nu_0 = 0$ \\
\hline $89$ &$\nu_0 \ne 0 \wedge s\ne0$ & - & otherwise & - & $\nu_0 = s = 0$\\
\hline $90$ & $\nu_0 \ne 0$ & - & $\nu_0 = 0$ & - & -\\

\hline
\end{tabular}
\end{table}

\begin{table}[t]
\centering \caption{\label{b16cc}$b_1(\mathfrak{g}_{6.i})$ for
$\mathfrak{g}_{6.i}$ unimodular (continued)}
\begin{tabular}{|c|c|c|c|c|c|} \hline \hline
$i$ & $b_1 = 1$ & $b_1 = 2$ & $b_1 = 3$ & $b_1=4$ & $b_1=5$%
\\ \hline \hline $91$ & always & - & - & - & -\\
\hline $92$ & $\mu_0 \ne 0 \wedge \nu_0  \ne0$ & - & otherwise & - & $\mu_0 = \nu_0 = 0$\\
\hline $92^*$ & always & - & - & - & -\\
\hline $93$ & $\nu_0 \ne 0$ & - & $\nu_0 = 0$ & - & -\\
\hline $94$ & always & - & - & - & -\\
\hline $101$ & - & always & - & - & -\\
\hline $102$ & - & always & - & - & -\\
\hline $105$ & - & always & - & - & -\\
\hline $107$ & - & always & - & - & -\\
\hline $113$ & - & always & - & - & -\\
\hline $114$ & - & always & - & - & -\\
\hline $115$ & - & always & - & - & -\\
\hline $116$ & - & always & - & - & -\\
\hline $118$ & - & always & - & - & -\\
\hline $120$ & - & always & - & - & -\\
\hline $125$ & - & always & - & - & -\\
\hline $129$ & - & always & - & - & -\\
\hline $135$ & - & always & - & - & -\\
\hline
\end{tabular}
\end{table}

\section{Integer Polynomials} \label{AppZ}
In this article, we often try to use necessary conditions for a
matrix to be conjugated to an integer matrix. We state briefly the
used results. Vice versa, we sometimes want to find integer matrices
with given minimal polynomial. We also present a few constructions.
\N Let be $n \in \mathbb{N}_+$, $\mathbb{K}$ a field and $A \in
\mathrm{M}(n,n;\mathbb{K})$. The \emph{characteristic
polynomial}\index{Characteristic Polynomial} of
$A$\index{Characteristic Polynomial} is the monic polynomial
$$P_A (X) := \det (X \, \mathrm{id} - A) \in \mathbb{K}[X],$$ and the
\emph{minimal polynomial}\index{Minimal Polynomial} $m_A(X)$ is the
unique monic divisor of lowest degree of $P_A(X)$ in $\mathbb{K}[X]$
such that $m_A(A) = 0$. (Note, by the theorem of Cayley-Hamilton,
one has $P_A(A) = 0$.)

If two matrices are conjugated, then they have the same
characteristic resp. minimal polynomials.

$\lambda \in \overline{\mathbb{K}}$ is called
\emph{root}\index{Root} of $A$ if $\lambda$ is a root of the
characteristic polynomial, considered as polynomial in
$\overline{\mathbb{K}}[X]$, where $\overline{\mathbb{K}}$ denotes
the algebraic closure of $\mathbb{K}$.

The next proposition follows directly from \cite[Corollaries
XIV.2.2, XIV.2.3]{Lang}.

\begin{Prop}
Let $n \in \mathbb{N}_+$. If $A \in \mathrm{M}(n,n;\mathbb{C})$ and
$B \in \mathrm{M}(n,n;\mathbb{Q})$ are conjugated via an element of
$\mathrm{GL}(n,\mathbb{C})$, then holds $P_A(X) = P_B(X) \in
\mathbb{Q}[X]$, $m_A(X) = m_B(X) \in \mathbb{Q}[X]$ and $m_A(X)$
divides $P_A(X)$ in $\mathbb{Q}[X]$.
\end{Prop}

\begin{Prop}
If $P(X) \in \mathbb{Z}[X]$, $m(X) \in \mathbb{Q}[X]$ are monic
polynomials and $m(X)$ divides $P(X)$ in $\mathbb{Q}[X]$, then holds
$m(X) \in \mathbb{Z}[X]$.
\end{Prop}

\textit{Proof.} Let $P(X), m(X)$ be as in the proposition and $f(X)
\in \mathbb{Q}[X]$ non-constant with $P(X) = f(X) \, m(X)$. There
exist $k,l \in \mathbb{Z} \setminus \{0\}$ such that $$k \, f(X) =
\sum_i a_i X^i, ~ l \, m(X) = \sum_j b_j X^j \in \mathbb{Z}[X]$$ are
primitive. (An integer polynomial is called \emph{primitive} if its
coefficients are relatively prime.) We have $$kl \, P(X) = (\sum_i
a_i X^i)(\sum_j b_j X^j)$$ and claim $kl = \pm 1$.

Otherwise, there is a prime $p \in \mathbb{N}$ that divides $kl$.
Since the coefficients of $k \, f(X)$ resp.\ $l \, m(X)$ are
relatively prime, there are minimal $i_0,j_0 \in \mathbb{N}$ such
that $p$ does not divide $a_{i_0}$ resp.\ $b_{j_0}$.

The coefficient of $X^{i_0 + j_0}$ of $k l \, f(X) \, m(X)$ is $$
a_{i_0} b_{j_0} + a_{i_0-1} b_{j_0+1} + a_{i_0+1} b_{j_0-1} + \ldots
$$ and $p$ divides each summand except the first. But since $p \,|\,
kl$, $p$ divides the whole sum. This is a contradiction. \q

\begin{Thm}  \label{minpol}
Let $n \in \mathbb{N}_+$ and $A \in \mathrm{M}(n,n;\mathbb{C})$ be
conjugated to an integer matrix. Then holds $P_A(X), m_A(X) \in
\mathbb{Z}[X]$.
\end{Thm}

\textit{Proof.} This follows from the preceding two propositions. \q

\begin{Lemma}[{\cite[Lemma 2.2]{Has4}}]  \label{Z3doppelt}
Let $\, P(X) = X^3 - k X^2 + l X - 1 \in \mathbb{Z}[X]$.

Then $P$ has a double root $\, X_0 \in \mathbb{R} \,$ if and only if
$\, X_0 = 1 \,$ or $\, X_0 = -1 \,$ for which $\, P(X) = X^3 - 3 X^2
+ 3 X - 1 \,$ or $\, P(X) = X^3 + X^2 - X - 1 \,$ respectively. \q
\end{Lemma}

\begin{Prop}[{\cite[Proposition 5]{Har}}] \label{Z3l1l2}
Let $\lambda_i \in \mathbb{R}_+$ with $\lambda_i +
\frac{1}{\lambda_i} = m_i \in \mathbb{N}_+$ and $m_i > 2$ for $i \in
\{1,2\}$.

Then there exists no element in $\mathrm{SL}(3,\mathbb{Z})$ with
roots $\lambda_1,\lambda_2,\frac{1}{\lambda_1 \lambda_2}$.
\end{Prop}

\begin{Prop}  \label{Z4doppelt}
Let $\, P(X) = X^4 - m X^3 + p X^2 - n X + 1 \in \mathbb{Z}[X]$.

Then $P$ has a root with multiplicity $> 1$ if and only if the zero
set of $P$ equals $\{1,1,a,a^{-1}\}$, $\{-1,-1,a,a^{-1}\}$,
$\{a,a^{-1},a,a^{-1}\}$ or $\{a,-a^{-1},a,-a^{-1}\}$ for fixed $a
\in \mathbb{C}$.
\end{Prop}

\textit{Proof.} The most part of the proof was done by Harshavardhan
in the proof of \cite[Propositon 2]{Har}.

We set $S := m^2 + n^2$ and $T := mn$ and get the discriminant $D$
of $P(X)$ as
\begin{eqnarray} \label{D=}
D & = & 16p^4 - 4Sp^3 + (T^2 -80 T -128)p^2 + 18S(T+8)p \\ &&
\nonumber + 256 - 192 T + 48T^2 - 4 T^3 -27 S^2.
\end{eqnarray}
Note that $P(X)$ has a root of multiplicity $>1$ if and only if
$D=0$. Solving $D=0$ for $S$, we see
\begin{eqnarray*}
S & = & -\frac{2}{27} p^3 + \frac{1}{3}pT + \frac{8}{3}p \pm
\frac{2}{27} \sqrt{(p^2 - 3T +12)^3},
\end{eqnarray*}
and since $S$ and $T$ are integers, there is $q \in \mathbb{Z}$ with
\begin{eqnarray*}
p^2 - 3T + 12 & = & q^2,
\end{eqnarray*}
which implies
\begin{eqnarray} \label{X= Y=}
S & = & 4p + \frac{1}{27}(p^3 - 3 pq^2 \pm 2 q^3)  \\ \nonumber T &
= & \frac{1}{3}(p^2-q^2+12) .
\end{eqnarray}
We first consider the plus sign in equation (\ref{X= Y=}). Then one
has
\begin{eqnarray*}
(m+n)^2 & = & S + 2T = \frac{1}{27}(p + 2q +6)(p-q+6)^2, \\
(m-n)^2 & = & S - 2T = \frac{1}{27}(p + 2q -6)(p-q-6)^2,
\end{eqnarray*}
and this implies the existence of $k_i,l_i \in \mathbb{N}$, $i =
1,2$, such that
\begin{eqnarray*}
3 k_1^2 & = & (p+2q+6) k_2^2, \\
3 l_1^2 & = & (p+2q-6) l_2^2.
\end{eqnarray*}
We shall show: $|m| = |n|$

[If $l_2 = 0$, the claim is proved. Therefore, we can assume $l_2
\ne 0$.

Case $1$: $k_2 = 0$

Then holds $k_1 = 0$ and this means $S+2T =0$, i.e.\ $(m+n)^2 = 0$,
so we have $m = -n$.

Case $2$: $k_2 \ne 0$

We write $k := \frac{k_1}{k_2}$ and $l := \frac{l_1}{l_2}$. Then
holds
\begin{eqnarray*}
3 k^2 & = &  p+2q+6 \in \mathbb{Z}, \\
3 l^2 & = &  p+2q-6 \in \mathbb{Z},
\end{eqnarray*}
and $3(k^2-l^2) = 12$. Therefore, we have $k^2 - l^2 = 4$, so $k^2 =
4$, $l^2=0$, i.e.\ $l_1 = 0$, $S-2T = 0$  and $m = n$. ]

Now, consider the minus sign in equation (\ref{X= Y=}). Then one has
\begin{eqnarray*}
(m+n)^2 & = & S + 2T = \frac{1}{27}(p - 2q +6)(p+q+6)^2, \\
(m-n)^2 & = & S - 2T = \frac{1}{27}(p - 2q -6)(p+q-6)^2,
\end{eqnarray*}
and shows analogously as above $|m| = |n|$.

We have shown: If $P(X)$ has a multiple root, then holds $m = \pm
n$.

If $m=n$, then one calculates the solutions of $D=0$ in (\ref{D=})
as the following
\begin{itemize}
\item[(i)] $p = -2 +2m$,
\item[(ii)] $p = -2-2m$,
\item[(iii)] $p = 2 + \frac{m^2}{4}$,
\end{itemize}
and if $m=-n$, then the real solution of $D=0$ in (\ref{D=}) is
\begin{itemize}
\item[(iv)] $p = - 2 + \frac{m^2}{4}$.
\end{itemize}

Moreover, a short computation yields the zero set of $P(X)$ in the
cases (i) -- (iv) as $\{1,1,a,a^{-1}\}$, $\{-1,-1,a,a^{-1}\}$,
$\{a,a^{-1},a,a^{-1}\}$,$\{a,-a^{-1},a,-a^{-1}\}$, respectively. \q

\begin{Prop}[{\cite[Proposition 4.4.14]{AdkWein}}]
Let $\mathbb{K}$ be a field and $$m(X) = X^n + a_{n-1} X^{n-1} +
\ldots + a_1 X^1 + a_0 \in \mathbb{K}[X]$$ a monic polynomial. Then
$\left(
\begin{array}{cccccc}
0 & 0 & \ldots & 0 & 0 & -a_0 \\
1 & 0 & \ldots & 0 & 0 & -a_1 \\
0 & 1 & \ldots & 0 & 0 & -a_2 \\
\vdots & \vdots & \ddots & \vdots & \vdots & \vdots \\
0 & 0 & \ldots & 1 & 0 & -a_{n-2} \\
0 & 0 & \ldots & 0 & 1 & -a_{n-1}
\end{array} \right)$ has minimal polynomial $m(X)$. \q
\end{Prop}

If one is willing to construct an integer matrix with given
characteristic and minimal polynomial, one always can chose any
matrix $M$ which has the desired polynomials and try to find an
invertible matrix $T$ such that $T^{-1} M T$ has integer entries. Of
course, this can be difficult. In the case of $4\times4$ - matrices
we have the following easy construction.
\begin{Prop}[{\cite[Section 2.3.1]{Har}}] $\,$ \label{Matrix}
\begin{itemize}
\item[(i)] Let integers $m,n,p \in \mathbb{Z}$ be given.

Choose $m_1, \ldots, m_4 \in \mathbb{Z}$ such that $\sum_{i=1}^4 m_i
= m$ and set
\begin{eqnarray*}
a & := & -m_1^2 p + m_1^3 m_2 + m_1^3 m_3 + m_1^3 m_4 + m_1 n - 1, \\
b & := & (-m_2 - m_1) p + m_1 m_2^2 + m_1 m_2 m_3 + m_1 m_2 m_4 +
m_2^2 m_3 + m_2^2 m_4 \\ && + m_1^2 m_2 + m_1^2 m_3 + m_1^2 m_4 + n,
\\ c & := & m_1 m_2 + m_1 m_3 + m_1 m_4 + m_2 m_3 + m_2 m_4 + m_3
m_4 - p.
\end{eqnarray*}
Then the matrix $\left(
\begin{array}{cccc}
m_1 & 0 & 0 & a \\
1 & m_2 & 0 & b \\
0 & 1 & m_3 & c \\
0 & 0 & 1 & m_4
\end{array}
\right)$ has $X^4 - m X^3 + p X^2 - n X +1$ as characteristic
polynomial.
\item[(ii)] Let $m \in 2\mathbb{Z}$ be an even integer. Then the matrix
$\left(
\begin{array}{cccc}
\frac{m}{2} & 0 & -1 & 0 \\
0 & \frac{m}{2} & 0 & -1 \\
1 & 0 & 0 & 0 \\
0 & 1 & 0 & 0
\end{array}
\right)$ has the characteristic polynomial $(X^2 - \frac{m}{2}X
+1)^2$, and $(X^2 - \frac{m}{2}X +1)$ as minimal polynomial. \q
\end{itemize}
\end{Prop}

\end{appendix}
$\,$\N
\begin{Ac}
The results presented in this paper were obtained in my dissertation
under the supervision of Prof.\ H.\ Geiges. I wish to express my
sincerest gratitude for his support during the last four years.
\end{Ac}

\textsc{Christoph Bock, Universit\"at zu K\"oln, Mathematisches Institut,
Weyertal 86--90, 50931 K\"oln, Germany}

\textit{e-mail:} \verb"bock@math.uni-koeln.de"

\begin{thebibliography}{99} \addcontentsline{toc}{chapter}{Bibliography}
\bibitem{AdkWein}W.\ A.\ Adkins, St.\ H.\ Weintraub: \textit{Algebra},
Springer (1992).
\bibitem{Aus}L.\ Auslander: \textit{An exposition of the structure of
solvmanifolds -- part I}, Bull.\ Amer.\ Math.\ Soc.\ \textbf{79}
(1973), no.\ 2, 227--261.
\bibitem{AGH}L.\ Auslander, L.\ Green, F.\ Hahn: \textit{Flows on Homogeneous
Spaces}, Princeton University Press (1963).
\bibitem{BG88}Ch.\ Benson, C.\ S.\ Gordon: \textit{K\"ahler and
symplectic structures on nilmanifolds}, Topology \textbf{27} (1988),
no.\ 4, 513--518.
\bibitem{BG90}Ch.\ Benson, C.\ S.\ Gordon: \textit{K\"ahler structures on
compact solvmanifolds}, Proc.\ Amer.\ Math.\ Soc.\ \textbf{108}
(1990), no.\ 4, 971--980.
\bibitem{ipse}Ch.\ Bock: \textit{Geography of non-formal symplectic
and contact manifolds}, arXiv:0812.1447.
\bibitem{Bosch}S.\ Bosch: \textit{Lineare Algebra}, Springer (2006).
\bibitem{Brown}K.\ S.\ Brown: \textit{Cohomology of Groups}, Springer
(1982).
\bibitem{CS}R.\ Campoamor-Stursberg: \textit{Symplectic forms on six
dimensional real solvable Lie algebras}, arXiv:math.DG/0507499.
\bibitem{CS2}R.\ Campoamor-Stursberg: \textit{Some remarks concerning
the invariants of rank one solvable real Lie algebras}, Algebra
Colloq. \textbf{12} (2005), no.\ 3, 497--518.
\bibitem{Chu}B.-Y.\ Chu: \textit{Symplectic homogeneous spaces},
Trans.\ Amer.\ Math.\ Soc.\ \textbf{197} (1974), 145--159.
\bibitem{CG}L.\ Corwin, F.\ P.\ Greenleaf: \textit{Representations of
Nilpotent Lie Groups and Their Applications}, Cambridge University
Press (1990).
\bibitem{Dekimpe}K.\ Dekimpe: \textit{Semi-simple splittings for solvable Lie groups and
polynomial structures}, Forum Math.\ \textbf{12} (2000), no.\ 1,
77--96.
\bibitem{DGMS}P.\ Deligne, P.\ Griffiths, J.\ Morgan, D.\ Sullivan:
\textit{Real homotopy theory of K\"ahler manifolds}, Invent.\ Math.\
\textbf{29} (1975), no.\ 3, 245--274.
\bibitem{Diatta}A.\ Diatta: \textit{Left invariant contact structures on Lie
groups}, Differential Geom. Appl. \textbf{26} (2008), no.\ 5,
544--552.
\bibitem{DF}A.\ Diatta, B.\ Foreman: \textit{Lattices in contact Lie groups and
5-dimensional contact solvmanifolds}, arXiv:0904.3113.
\bibitem{Dix}J.\ Dixmier: \textit{Sur les alg\`ebres d\'eriv\'ees
d'une alg\`ebre de Lie}, Proc.\ Cambride Philos.\ Soc.\ \textbf{51}
(1955), 541--544.
\bibitem{Dix2}J.\ Dixmier: \textit{L'application exponentielle dans
les groupes de Lie r\'esolubles}, Bull.\ Soc.\ Math.\ France
\textbf{85} (1957), 113--121.
\bibitem{DER}N.\ Dungey, A.\ F.\ M.\ ter Elst, D.\ W.\ Robinson:
\textit{Analysis on Lie Groups with Polynomial Growth}, Birkh\"auser
(2003).
\bibitem{FG90}M.\ Fern\'andez, A.\ Gray: \textit{Compact symplectic
solvmanifolds not admitting complex structures}, Geom.\ Dedicata
\textbf{34} (1990), no.\ 3, 295--299.
\bibitem{FLS}M.\ Fern\'andez, M.\ de Le\'on, M.\ Saralegui:
\textit{A  six-dimensional compact symplectic solvmanifold without
K\"ahler structures}, Osaka J.\ Math.\ \textbf{33} (1996), no.\ 1,
19--35.
\bibitem{FMDon}M.\ Fern\'andez, V.\ Mu\~noz: \textit{Formality of Donaldson
submanifolds},  Math.\ Z.\ \textbf{250} (2005), no.\ 1, 149--175.
\bibitem{FMGeo}M.\ Fern\'andez, V.\ Mu\~noz: \textit{The geography of non-formal
manifolds} in O.\ Kowalski, E.\ Musso, D.\ Perrone: \textit{Complex,
Contact and Symmetric Manifolds}, 121--129, Prog.\ Math.\ vol.\ 234,
Birkh\"auser (2005).
\bibitem{SF}S.\ Fukuhara, K.\ Sakamoto: \textit{Classification of
$T^2$-bundles over $T^2$}, Tokyo J.\ Math.\ \textbf{6} (1983), no.\
2, 311--327.
\bibitem{Ge4}H.\ Geiges: \textit{Symplectic structures on
$T^2$-bundles over $T^2$}, Duke Math.\ J.\ \textbf{67} (1992), no.\
3, 539--555.
\bibitem{Ge5}H.\ Geiges: \textit{Symplectic manifolds with
disconnected boundary of contact type}, Internat.\ Math.\ Res.\
Notices (1994), no.\ 1, 23--30.
\bibitem{Gorb}V.\ V.\ Gorbatsevich: \textit{Symplectic structures and
cohomologies on some solvmanifolds}, Sibiran Math.\ J.\ \textbf{44}
(2003), no.\ 2, 260--274.
\bibitem{GH}P.\ Griffiths, J.\ Harris: \textit{Principles of Algebraic
Geometry}, Wiley (1994).
\bibitem{Har}R.\ Harshavardhan: \textit{Geometric structures of Lie type
on 5-manifolds}, Ph.D.\ Thesis, Cambridge University (1996).
\bibitem{HasNil}K.\ Hasegawa: \textit{Minimal models of
nilmanifolds}, Proc.\ Amer.\ Math.\ Soc.\ \textbf{106} (1989), no.\
1, 65--71.
\bibitem{Has4}K.\ Hasegawa: \textit{Four-dimensional compact solvmanifolds with
and without complex analytic structures}, arXiv:math.CV/0401413.
\bibitem{HasCK}K.\ Hasegawa: \textit{Complex and K\"ahler structures on
compact solvmanifolds}, J.\ Symplectic Geom.\ \textbf{3} (2005),
no.\ 4, 749--767.
\bibitem{HasK}K.\ Hasegawa: \textit{A note on compact solvmanifolds with K\"ahler
structure}, Osaka J.\ Math.\ \textbf{43} (2006), no.\ 1, 131--135.
\bibitem{Hat}A.\ Hattori: \textit{Spectral sequence in the de Rham
cohomology of fibre bundles}, J.\ Fac.\ Sci.\ Univ.\ Tokyo Sect.\ I
\textbf{8} (1960), 289--331.
\bibitem{IRTU99}R.\ Ib\'a\~nez, Y.\ Rudiak, A.\ Tralle, L.\ Ugarte:
\textit{On certain geometric and homotopy properties of closed
symplectic manifolds}, Topology Appl.\ \textbf{127} (2003), no.\
1--2, 33--45.
\bibitem{KRT}J.\ K\c{e}dra, Y.\ Rudyak, A.\ Tralle:
\textit{Symplectically aspherical manifolds}, J.\ Fixed Point Theory
Appl.\ \textbf{3} (2008), no.\ 1, 1--21.
\bibitem{Lang}S.\ Lang: \textit{Algebra}, 3.\ ed., Addison-Wesley
(1993).
\bibitem{Lutz}R.\ Lutz: \textit{Sur la g\'eom\'etrie des structures
de contact invariantes}, Ann.\ Inst.\ Fourier (Grenoble),
\textbf{29} (1979), no.\ 1, xvii, 283--306.
\bibitem{Mil}J.\ Milnor: \textit{Curvature of left invariant metrics
on Lie groups}, Advances in Math.\ \textbf{21} (1976), no.\ 3,
293--329.
\bibitem{Mor}V.\ V.\ Morozov: \textit{Classifikaziya nilpotentnych
algebrach Lie shestovo poryadka}, Izv.\ Vys\u{s}.\ U\u{c}ebn.\
Zaved.\ Mathematika \textbf{5} (1958), no.\ 4, 161--171.
\bibitem{Mos}G.\ D.\ Mostow: \textit{Factor spaces of solvable
groups}, Ann.\ of Math.\ (2) \textbf{60} (1954), 1--27.
\bibitem{Mub63}G.\ M.\ Mubarakzjanov: \textit{O rasreshimych
algebrach Lie}, Izv.\ Vys\u{s}.\ U\u{c}ebn.\ Zaved.\ Mathematika
\textbf{32} (1963), no.\ 1, 114--123.
\bibitem{Mub63a}G.\ M.\ Mubarakzjanov: \textit{Classifkaziya
vehestvennych structur algebrach Lie pyatovo poryadko}, Izv.\
Vys\u{s}.\ U\u{c}ebn.\ Zaved.\ Mathematika \textbf{34} (1963), no.\
3, 99--106.
\bibitem{Mub63b}G.\ M.\ Mubarakzjanov: \textit{Classifkaziya
rasreshimych structur algebrach Lie shestovo poryadka a odnim
nenilpoentnym bazisnym elementom}, Izv.\ Vys\u{s}.\ U\u{c}ebn.\
Zaved.\ Mathematika \textbf{35} (1963), no.\ 4, 104--116.
\bibitem{Nak}I.\ Nakamura: \textit{Complex parallelisable manifolds
and their small deformations}, J.\ Differential Geometry \textbf{10}
(1975), 85--112.
\bibitem{NN}A.\ Newlander, L.\ Nirenberg: \textit{Complex analytic
coordinates in almost complex manifolds}, Ann.\ of Math.\ (2)
\textbf{65} (1957), 391--404.
\bibitem{N}K.\ Nomizu: \textit{On the cohomology of homogeneous
spaces of nilpotent Lie Groups}, Ann.\ of Math.\ (2) \textbf{59}
(1954), 531--538.
\bibitem{Nono}T.\ N\^ono: \textit{On the singularity of general
linear groups}, J.\ Sci.\ Hiroshima Univ.\ Ser.\ A \textbf{20}
(1956/1957), 115--123.
\bibitem{TO}J.\ Oprea, A.\ Tralle: \textit{Symplectic Manifolds with no K\"ahler
Structure}, Lecture Notes in Math.\ \textbf{1661}, Springer (1997).
\bibitem{OV}A.\ L.\ Onishchik, E.\ B.\ Vinberg: \textit{Lie Groups and Lie
Algebras III}, Springer (1994).
\bibitem{Rag}M.\ S.\ Raghunathan: \textit{Discrete Subgroups of Lie
Groups}, Springer (1972).
\bibitem{Saito}M.\ Sait\^o: \textit{Sur certaines groupes de Lie
r\'esolubles -- parties I et II}, Sci.\ Papers Coll.\ Gen.\ Ed.\
Univ.\ Tokyo \textbf{7} (1957), 1--11, 157--168.
\bibitem{Salamon}S.\ M.\ Salamon: \textit{Complex structures on
nilpotent Lie algebras}, J.\ Pure Appl.\ Algebra \textbf{157}
(2001), no.\ 2-3, 311--333.
\bibitem{SY}H. Sawai, T. Yamada: \textit{Lattices on Benson-Gordon
type solvable Lie groups}, Topology Appl.\ \textbf{149} (2005), no.\
1-3, 85--95.
\bibitem{Steen}N.\ Steenrod: \textit{The Topology of Fibre Bundles},
Princeton University Press (1951).
\bibitem{Sullivan}D.\ Sullivan: \textit{Infinitesimal computations in
topology}, Inst.\ Hautes \'Etudes Sci.\ Publ.\ Math.\ \textbf{47}
(1977), 269--331.
\bibitem{Thurston}W.\ P.\ Thurston: \textit{Some simple examples of
compact symplectic manifolds}, Proc.\ Amer.\ Math.\ Soc.\
\textbf{55} (1976), no.\ 2, 467--468.
\bibitem{Tralle}A.\ Tralle: \textit{On solvable Lie groups without lattices}
in M.\ Fern\'andez, J.\ A.\ Wolf: \textit{Global Differential
Geometry}, 437--441, Contemp.\ Math.\ \textbf{288}, Amer.\ Math.\
Soc.\ (2001).
\bibitem{T}P.\ Turkowski: \textit{Solvable Lie algebras of dimension
six}, J.\ Math.\ Phys.\ \textbf{31} (1990), no.\ 6, 1344--1350.
\bibitem{V}V.\ S.\ Vararadarajan: \textit{Lie Groups, Lie Algebras, and Their
Representations}, Springer (1984).
\bibitem{War}F.\ W.\ Warner: \textit{Foundations of Differentiable
Manifolds and Lie Groups}, Glenview (1971).
\bibitem{Weibel}Ch.\ A.\ Weibel: \textit{An Introduction to
Homological Algebra}, Cambridge University Press (1997).
\bibitem{Yam}T.\ Yamada: \textit{A pseudo-K\"ahler structure on a
nontoral compact complex parallelizable solvmanifold}, Geom.\
Dedicata \textbf{112} (2005), 115-122.
\end{thebibliography}
\end{document}